\documentclass[12pt,a4paper,leqno]{article}
\usepackage[toc]{appendix}
\usepackage{setspace} \usepackage{etoolbox} \usepackage{comment}
\usepackage{multicol}
\usepackage{array}
\usepackage[amsthm,thmmarks]{ntheorem}
\usepackage{amsmath,amsfonts,amssymb,mathtools}
\usepackage[bbgreekl]{mathbbol} \usepackage[a4paper,includeheadfoot,margin=2.54cm]{geometry}
\usepackage{times}
\usepackage{xypic}
\usepackage{pdflscape}
\usepackage[export]{adjustbox} 
\usepackage[cal=boondoxo,calscaled=.96,scr=rsfs]{mathalpha} \DeclareMathAlphabet{\mathscrbf}{OMS}{mdugm}{b}{n} \usepackage{hyphenat}
\usepackage[usenames,dvipsnames]{color}
\usepackage{rotating}
\usepackage{tikz-cd}
\usetikzlibrary{cd,decorations.pathmorphing}
\usepackage{longtable}
\usepackage{caption}
\usepackage[most]{tcolorbox}

\usepackage{epigraph}
\usepackage{enumitem}
\setlist[enumerate]{listparindent=0.5in}
\usepackage{mathrsfs}
\newcommand{\be}{\begin{equation}}
\newcommand{\ee}{\end{equation}}
\newcommand{\bes}{\begin{equation*}}
\newcommand{\ees}{\end{equation*}}
\newcommand{\bea}{\begin{eqnarray}}
\newcommand{\eea}{\end{eqnarray}}
\newcommand{\beas}{\begin{eqnarray}}
\newcommand{\eeas}{\end{eqnarray}}
\newcommand{\ben}{\begin{note}}
\newcommand{\een}{\end{note}}
\newcommand{\bexl}{\vskip0.1em\noindent\hrulefill\vskip1em\begin{ExerciseList}}
\newcommand{\eexl}{\end{ExerciseList}\hrulefill}

\newcommand{\bthm}{\begin{theorem}}
\newcommand{\ethm}{\end{theorem}}
\newcommand{\bpro}{\begin{prop}}
\newcommand{\epro}{\end{prop}}
\newcommand{\bcor}{\begin{corollary}}
\newcommand{\ecor}{\end{corollary}}
\newcommand{\bcon}{\begin{conjecture}}
\newcommand{\econ}{\end{conjecture}}
\newcommand{\bp}{\begin{proof}}
\newcommand{\ep}{\end{proof}}
\newcommand{\blem}{\begin{lemma}}
\newcommand{\elem}{\end{lemma}}
\newcommand{\bn}{\begin{note}}
\newcommand{\en}{\end{note}}
\newcommand{\benum}{\begin{enumerate}}
\newcommand{\eenum}{\end{enumerate}}
\newcommand{\bed}{\begin{defn}}
\newcommand{\eed}{\end{defn}}
\newcommand{\brem}{\begin{remark}}
\newcommand{\erem}{\end{remark}}

\newcommand{\btik}{\begin{tikzpicture}\begin{axis}[scale=0.5,axis y line=center, axis x line=middle]}
\newcommand{\etik}{\end{axis}\end{tikzpicture}}

\let\into=\hookrightarrow
\let\mapsto=\longmapsto
\let\cong=\equiv

\newcommand{\upperRomannumeral}[1]{\uppercase\expandafter{\romannumeral#1}}

 \usepackage{graphicx}
\usepackage[pagewise]{lineno}
\usepackage{stackengine}
\usepackage{stmaryrd}
\usepackage[normalem]{ulem} \usepackage{titlesec}
\usepackage{natbib}
	\let\cite=\citep
	\usepackage{fancyhdr}
	\usepackage[hidelinks,colorlinks,citecolor=blue,pagebackref]{hyperref}
	\usepackage[capitalise,noabbrev]{cleveref}
\usepackage{rotating}

\usepackage{tocloft}
\advance\cftsecnumwidth 0.5em\relax
\advance\cftsubsecindent 0.5em\relax
\advance\cftsubsecnumwidth 0.5em\relax

\usepackage{xr}
\newcommand{\constranab}[2]{\cite[{#1\ref{0-#2}}]{joshi-anabelomorphy}}

\newcommand{\constrone}[2]{\cite[{#1\ref{I-#2}}]{joshi-teich}}
\newcommand{\constrtwo}[2]{\cite[{#1\ref{II-#2}}]{joshi-teich-estimates}}
\newcommand{\constrtwoh}[2]{\cite[{#1\ref{II5-#2}}]{joshi-teich-def}}

\newcommand{\constrfour}[2]{\cite[{#1\ref{IV-#2}}]{joshi-teich-abc}}

\newtoggle{draft}
\togglefalse{draft}

\theoremstyle{theorem}
\newtheorem{theorem}[equation]{Theorem}      \newtheorem{theoremdef}[equation]{Theorem-Definition}
\newtheorem{lemma}[equation]{Lemma}          \newtheorem{corollary}[equation]{Corollary}  \newtheorem{proposition}[equation]{Proposition}
\newtheorem{prop*}[equation]{Proposition}

\newtheorem{conj}[equation]{Conjecture}

\theoremstyle{remark}

\theoremsymbol{\ensuremath{\bullet}}

\theoremstyle{definition}
\newtheorem{remark}[equation]{Remark}
\theoremsymbol{\ensuremath{\bullet}}
\newtheorem*{rmk*}{Remark}
\theoremsymbol{\ensuremath{\bullet}}
\newtheorem{defn}[equation]{Definition}
\theoremsymbol{\ensuremath{\bullet}}
\newtheorem*{defns*}{Definition}
\theoremsymbol{\ensuremath{\bullet}}

\newcommand{\bthmdef}{\begin{theoremdef}}
\newcommand{\ethmdef}{\end{theoremdef}}

\let\oldproofname=\proofname
\renewcommand{\proofname}{{\bfseries\rmfamily\textup{\oldproofname}}}

\titleformat{\subsection}[runin]{\normalfont\bfseries}{\ssep\thesubsection}{.5em}{}[{\ \ }]
\titlespacing{\subsection}{0pt}{1.5ex plus .1ex minus .2ex}{0pt}
\titleformat{\subsubsection}[runin]{\normalfont\bfseries}{\ssep\thesubsubsection}{.5em}{}[{\ \ }]
\titlespacing{\subsubsection}{0pt}{1.5ex plus .1ex minus .2ex}{0pt}

\newcommand{\nwss}{\numberwithin{equation}{subsection}}
\newcommand{\nwsss}{\numberwithin{equation}{subsubsection}}

\newcommand{\subpara}{\subsection{}\nwss}
\newcommand{\subparat}[1]{\subsection{#1}\nwss}
\newcommand{\subsubpara}{\subsubsection{}\nwsss}
\newcommand{\subsubparat}[1]{\subsubsection{#1}\nwsss}

\newcommand{\ssep}{\S\,}

\crefname{section}{§}{§§}
\crefname{subsection}{§}{§§}
\crefname{subsubsection}{§}{§§}

\usepackage{microtype}
\usepackage[margin=15pt,font=small,labelfont={bf,sf}]{caption}

\let\into=\hookrightarrow
\let\isom=\simeq

\let\tensor=\otimes
\newcommand{\A}{\mathscr{A}}
\newcommand{\abs}[1]{\left\vert#1\right\vert}

\newcommand{\bF}{{\bar{F}}}
\newcommand{\bQ}{{\bar{\Q}}}

\newcommand{\C}{{\mathbb C}}

\newcommand{\Ext}{{\rm Ext}\,}
\newcommand{\F}{{\mathbb F}}
\newcommand{\frob}{\mathfrak{Frob}}
\newcommand{\gal}{{\rm Gal}}
\newcommand{\Gm}{\mathbb{G}_m}

\newcommand{\N}{\mathscr{N}}

\newcommand{\Q}{{\mathbb Q}}
\newcommand{\R}{{\mathbb R}}

\newcommand{\spec}{{\rm Spec}}
\newcommand{\Spec}{{\rm Spec}}

\newcommand{\Z}{{\mathbb Z}}

\renewcommand{\O}{{\mathscr O}}
\renewcommand{\P}{{\mathbb P}}
\renewcommand{\wp}{{\mathfrak p}}

\newcommand{\Aut}{{\rm Aut}}
\newcommand{\fm}{{\mathfrak{M}}}

\newcommand{\invlim}{\varprojlim}
\let\fm=\fa
\newcommand{\mapright}[1]{{\xymatrix{{}\ar[r]^{#1}&{}}}}

\renewcommand{\bpro}{\begin{proposition}}
	\renewcommand{\epro}{\end{proposition}}
\renewcommand{\bcon}{\begin{conj}}
	\renewcommand{\econ}{\end{conj}}

\newcommand{\ebh}{\widehat{\bE}}
\newcommand{\lbh}{\widehat{\bL}}

\renewcommand{\fm}{\mathfrak{m}}

\newtcolorbox[auto counter,crefname={rosetta Stone Fragment}{rosetta Stone Fragements}]{rosetta}[2][]{colback=white,coltitle=black,colframe=white!25!brown,fonttitle=\bfseries,
	title=Rosetta Stone Fragment  \thetcbcounter:  for \joshiros\ and \iutthr  #2,#1}

\newtcolorbox[auto counter,crefname={boxedcontent}{boxedcontent}]{boxedcontent}[2][]{colback=white,coltitle=black,colframe=white!25!brown,fonttitle=\bfseries,
	title= #2,#1}

\title{Construction of Arithmetic Teichmuller Spaces III: \\ A `Rosetta Stone' \\
	and a proof of Mochizuki's Corollary 3.12\\ 
	\textcolor{red}{{\large\bf Preliminary version for comments}}	
}
\author{Kirti Joshi}

\begin{document}
	\maketitle
\epigraphwidth0.65\textwidth
\epigraph{ It is not enough to be in the right place at the right time. 
You should also have an open mind at the right time.}{Paul Erdos \cite{spechter-erdos-book}}

\newcommand{\iut}{\cite{mochizuki-iut1, mochizuki-iut2, mochizuki-iut3,mochizuki-iut4}}
\newcommand{\iutabc}{\cite{mochizuki-iut1, mochizuki-iut2, mochizuki-iut3}}
\newcommand{\topics}{\cite{mochizuki-topics1,mochizuki-topics2,mochizuki-topics3}}
\newcommand{\joshiros}{\cite{joshi-teich,joshi-untilts,joshi-teich-def}}

\begin{abstract}
	This is a continuation of my work on Arithmetic Teichmuller Spaces. This paper establishes a number of important results  including (1) a proof Mochizuki's Corollary 3.12  (2)  a `Rosetta Stone' for a parallel reading \iut\ and \joshiros, and (3) a proof that Mochizuki's gluing of Hodge-Theaters, Frobenioids along prime-strips \`a la \iut\ is naturally provided by  the existence of  Arithmetic Teichmuller Spaces.
\end{abstract}

\lhead{}

\iftoggle{draft}{\pagewiselinenumbers}{\relax}
\newcommand{\act}{\curvearrowright}
\newcommand{\lmp}{{\Pi\act\Ot}}
\newcommand{\lmpi}{{\lmp}_{\int}}
\newcommand{\lmpf}{\lmp_F}
\newcommand{\Om}{\O^{\times\mu}}
\newcommand{\Omf}{\O^{\times\mu}_{\bF}}
\renewcommand{\N}{\mathbb{N}}
\newcommand{\yoga}{Yoga}
\newcommand{\gl}[1]{{\rm GL}(#1)}
\newcommand{\bK}{\overline{K}}
\newcommand{\reptrip}{\rho:G_K\to\gl V}
\newcommand{\reptripp}[1]{\rho\circ\alpha:G_{\ifstrempty{#1}{K}{{#1}}}\to\gl V}

\newcommand{\benumlab}{\begin{enumerate}[label={{\bf(\arabic{*})}}]}
\newcommand{\benumlabresume}{\begin{enumerate}[resume,label={{\bf(\arabic{*})}}]}
\newcommand{\benumlabstart}[1]{\begin{enumerate}[start={#1},label={{\bf(\arabic{*})}}]}

\newcommand{\ord}{\mathop{\rm ord}\nolimits}	
\newcommand{\kcs}{K^\circledast}
\newcommand{\lcs}{L^\circledast}
\renewcommand{\A}{\mathbb{A}}
\newcommand{\bfq}{\bar{\mathbb{F}}_q}
\newcommand{\tripod}{\P^1-\{0,1728,\infty\}}

\newcommand{\vseq}[2]{{#1}_1,\ldots,{#1}_{#2}}
\newcommand{\anab}[4]{\left({#1},\{#3 \}\right)\anabelmap\left({#2},\{#4 \}\right)}

\newcommand{\gln}{{\rm GL}_n}
\newcommand{\glo}[1]{{\rm GL}_1(#1)}
\newcommand{\glt}[1]{{\rm GL_2}(#1)}

\newcommand{\linv}{\mathfrak{L}}
\newcommand{\bedef}{\begin{defn}}
\newcommand{\eedef}{\end{defn}}
\newcommand{\bdefn}{\begin{defn}}
\newcommand{\edefn}{\end{defn}}
\renewcommand{\act}[1][]{\overset{#1}{\curvearrowright}}
\newcommand{\bfx}{\overline{F(X)}}
\newcommand{\anabelmap}{\leftrightsquigarrow}
\newcommand{\ban}[1][G]{\mathscr{B}({#1})}
\newcommand{\pit}{\Pi^{temp}}
 
 \newcommand{\bL}{\overline{L}}
 \newcommand{\bkm}{\bK_M}
 \newcommand{\vbk}{v_{\bK}}
 \newcommand{\vbkm}{v_{\bkm}}
\newcommand{\ocs}{\O^\circledast}
\newcommand{\ot}{\O^\triangleright}
\newcommand{\ocsk}{\ocs_K}
\newcommand{\otk}{\ot_K}
\newcommand{\ok}{\O_K}
\newcommand{\oko}{\O_K^1}
\newcommand{\oks}{\ok^*}
\newcommand{\Qpb}{\overline{\Q}_p}
\newcommand{\Qpbh}{\widehat{\overline{\Q}}_p}
\newcommand{\tr}{\triangleright}
\newcommand{\ocpt}{\O_{\C_p}^\tr}
\newcommand{\ocpf}{\O_{\C_p}^\flat}
\newcommand{\sG}{\mathscr{G}}
\newcommand{\sX}{\mathscr{X}}
\newcommand{\sxfe}{\sX_{F,E}}
\newcommand{\sxfep}{\sX_{F,E'}}
\newcommand{\sxcpte}{\sX_{\cpt,E}}
\newcommand{\sxcptep}{\sX_{\cpt,E'}}
\newcommand{\loglt}{\log_{\sG}}
\newcommand{\fc}{\mathfrak{t}}
\newcommand{\ku}{K_u}
\newcommand{\kup}{\ku'}
\newcommand{\kt}{\tilde{K}}
\newcommand{\sGpf}{\sG(\O_K)^{pf}}
\newcommand{\hgm}{\widehat{\mathbb{G}}_m}
\newcommand{\bE}{\overline{E}}
\newcommand{\sY}{\mathscr{Y}}
\newcommand{\syfe}{\mathscr{Y}_{F,E}}
\newcommand{\syfep}{\mathscr{Y}_{F,E'}}
\newcommand{\syfqp}[1]{\mathscr{Y}_{\cptl{#1},\Q_p}}
\newcommand{\syfqpe}[1]{\mathscr{Y}_{{#1},E}}
\newcommand{\syfqpep}[1]{\mathscr{Y}_{{#1},E'}}
\newcommand{\fJ}{\mathfrak{J}}
\newcommand{\fM}{\mathfrak{M}}
\newcommand{\locvar}{local arithmetic-geometric anabelian variation of fundamental group of $X/E$ at $\wp$}
\newcommand{\fjxep}{\fJ(X,E,\wp)}
\newcommand{\fjxe}{\fJ(X,E)}
\newcommand{\fpc}[1]{\widehat{{\overline{\F_p(({#1}))}}}}
\newcommand{\cpt}{\C_p^\flat}
\newcommand{\cpti}{\C_{p_i}^\flat}
\newcommand{\cptl}[1]{\C_{p,{#1}}^\flat}
\newcommand{\fja}[1]{\fJ^{\rm arith}({#1})}
\newcommand{\ainfe}{A_{\inf,E}(\O_F)}
\renewcommand{\ainfe}{W_{\O_E}(\O_F)}
\newcommand{\gmh}{\widehat{\mathbb{G}}_m}
\newcommand{\sE}{\mathscr{E}}
\newcommand{\bpi}{B^{\varphi=\pi}}
\newcommand{\bpip}{B^{\varphi=p}}
\newcommand{\onto}{\twoheadrightarrow}

\newcommand{\cpmax}{{\C_p^{\rm max}}}
\newcommand{\xan}{X^{an}}
\newcommand{\yan}{Y^{an}}
\newcommand{\bPi}{\overline{\Pi}}
\newcommand{\bPit}{\bPi^{\rm{\scriptscriptstyle temp}}}
\newcommand{\Pit}{\Pi^{\rm{\scriptscriptstyle temp}}}
\renewcommand{\pit}[1]{\Pi^{\scriptscriptstyle temp}_{#1}}
\newcommand{\pitk}[2]{\Pi^{\scriptscriptstyle temp}_{#1;#2}}
\newcommand{\pio}[1]{\pi_1({#1})}
\newcommand{\fTeich}{\widetilde{\fJ(X/L)}}
\newcommand{\vphi}{\varphi}
\newcommand{\sgt}{\widetilde{\sG}}
\newcommand{\sxqp}{\mathscr{X}_{\cpt,\Q_p}}
\newcommand{\syQp}{\mathscr{Y}_{\cpt,\Q_p}}

\newpage
{\hypersetup{linkcolor=Mahogany}\sffamily
\tableofcontents
}

\newcommand{\mywork}[1]{\textcolor{teal}{#1}}

\togglefalse{draft}
\newcommand{\FF}{\cite{fargues-fontaine}}
\iftoggle{draft}{\pagewiselinenumbers}{\relax}

\newcommand{\attportion}{Sections~\cref{se:number-field-case}, \cref{se:construct-att}, \cref{se:relation-to-iut}, \cref{se:self-similarity} and \cref{se:applications-elliptic}}

\newcommand{\Pib}{\overline{\Pi}}
\newcommand{\four}{Sections~\cref{se:grothendieck-conj}, \cref{se:untilts-of-Pi}, and \cref{se:riemann-surfaces}}

\newcommand{\omu}{\O_{\bQ_p}^{\mu}}
\newcommand{\lmod}{L_{\rm mod}}
\newcommand{\ttheta}{\widetilde{\Theta}}
\newcommand{\tthetaj}[1]{\ttheta_{\scriptscriptstyle{Joshi},#1}}
\newcommand{\tthetam}[1]{\ttheta_{\scriptscriptstyle{Mochizuki},#1}}
\newcommand{\moccor}{\cite[Corollary 3.12]{mochizuki-iut3}}
\newcommand{\thetam}{{\ttheta}_{\scriptscriptstyle{Mochizuki}}}
\newcommand{\thetamp}{{\ttheta}_{\scriptscriptstyle{Mochizuki},p}}
\newcommand{\thetaj}{{{\ttheta}_{\scriptscriptstyle{Joshi}}}}
\newcommand{\thetajp}{{{\ttheta}_{\scriptscriptstyle{Joshi},p}}}
\newcommand{\thetajpq}{{{\ttheta}_{\scriptscriptstyle{Joshi},p}^?}}
\newcommand{\sM}{\mathscr{M}}
\newcommand{\pib}{\overline{\Pi}}
\newcommand{\bN}{\mathbb{N}}
\newcommand{\sD}{\mathscr{D}}
\newcommand{\sF}{\mathscr{F}}
\newcommand{\sL}{\mathscr{L}}

\newcommand{\bdrp}{B_{dR}^+}
\newcommand{\bdr}{B_{dR}}

\newcommand{\iutthr}{\cite{mochizuki-iut1,mochizuki-iut2, mochizuki-iut3}}

\newcommand{\thetajpp}{{\widehat{\Theta}}_{\scriptscriptstyle{Joshi},p}}
\newcommand{\thetaja}{{\widehat{\Theta}}_{Joshi}}
\newcommand{\thetajpph}{{\widehat{\widehat{\Theta}}}_{\scriptscriptstyle{Joshi},p}}
\newcommand{\thetajppa}{{\widehat{\widehat{\Theta}}}_{\scriptscriptstyle{Joshi}}}

\newcommand{\ells}{{\ell^*}}

\definecolor{darkmidnightblue}{rgb}{0.0, 0.2, 0.4}
\definecolor{carmine}{rgb}{0.59, 0.0, 0.09}
\newtcbox{\mybox}[1][red]{on line,
	arc=0pt,outer arc=0pt,colback=white,colframe=darkmidnightblue,
	boxsep=0pt,left=0pt,right=0pt,top=2pt,bottom=2pt,
	boxrule=0pt,leftrule=1pt, rightrule=1pt,bottomrule=1pt,toprule=1pt}
\newcommand{\tmb}[1]{\mybox{#1}}
\newcommand{\present}{the present series of papers (\cite{joshi-teich,joshi-untilts,joshi-teich-estimates,joshi-teich-def} and this paper)}

\newcommand{\bcris}{B_{cris}}

\newcommand{\sP}{\mathscr{P}}
\newcommand{\spc}{\sP^{Teich}}
\newcommand{\pitop}[1]{\pi_1^{top}(#1)}
\newcommand{\sppi}{\sP(\Pi\hookleftarrow\pib)}
\newcommand{\sppim}{\sP_{\scriptscriptstyle{Mochizuki}(\Pi\hookleftarrow\pib)}}
\newcommand{\sppij}{\sP_{\scriptscriptstyle{Joshi}}(\Pi\hookleftarrow\pib)}

\newcommand{\hol}[3]{\mathfrak{hol}_{#1}(#3)_{#2}}
\newcommand{\holt}[2]{\mathfrak{hol}(#1)_{#2}}
\newcommand{\xdia}{X^\lozenge}
\newcommand{\holdia}{{\mathfrak{Hol}(X/E)}}
\newcommand{\spd}{{\rm Spd}}
\newcommand{\perf}{{\rm Perf}}
\newcommand{\perfld}{{\rm PerfFld}}
\newcommand{\spa}{{\rm Spa}}
\newcommand{\bnddsub}{{\tiny\circ}}

\newcommand{\xad}{X^{ad}}

\newcommand{\sppimold}{\sP_{\scriptscriptstyle{Mochizuki}(\Pi\hookleftarrow\pib)}}
\renewcommand{\sppim}{\sP_{\scriptscriptstyle{Mochizuki}'(\Pi\hookleftarrow\pib)}}

\newcommand{\qprojqp}[1]{\mathcal{QProj}_{#1}}
\newcommand{\projqp}[1]{\mathcal{Proj}_{#1}}

\newcommand{\chxfe}{{\rm Ch^1(\sxfe)}}
\newcommand{\divxfe}{{\rm Div}(\sxfe)}

\newcommand{\tsigma}{\tilde{\sigma}}

\newcommand{\ismg}{{\mathbf{Ism}(G)}}

\newcommand{\okbt}{\O_{\bE;K}^\triangleright}

\newcommand{\Ob}{\overline{\O}} 

\newcommand{\divp}{Div_{+}}

\newcommand{\divt}{\widetilde{{\rm Div_{+}}}}
\newcommand{\pdivt}{\widetilde{{\rm PDiv}}}

\newcommand{\fjxlp}{\fJ(X,L_\wp)}

\newcommand{\bbvl}{\mathbb{V}_L}
\newcommand{\bbvlp}{\mathbb{V}_{L'}}
\newcommand{\bbvlmod}{\mathbb{V}_{\lmod}}
\newcommand{\bbvlmodgood}{\mathbb{V}_{\lmod}^{good}}
\newcommand{\bbvlgood}{\mathbb{V}_{L}^{good}}
\newcommand{\bbvlpgood}{\mathbb{V}_{L'}^{good}}
\newcommand{\bbvlmodoss}{\mathbb{V}_{\lmod}^{odd,ss}}
\newcommand{\bbvloss}{\mathbb{V}_{L}^{odd,ss}}
\newcommand{\bbvlposs}{\mathbb{V}_{L'}^{odd,ss}}
\newcommand{\ubblv}{\underline{\mathbb{V}}}
\newcommand{\ubblvgood}{\underline{\mathbb{V}}^{good}}
\newcommand{\ubblvoss}{\underline{\mathbb{V}}^{odd,ss}}
\newcommand{\ubblvossp}{\underline{\mathbb{V}}^{odd,ss}_p}

\newcommand{\ul}[1]{\underline{#1}}

\newcommand{\breveB}{{\breve{B}}}
\newcommand{\breveBten}{{\breve{B}}^\otimes}
\newcommand{\breveBb}{{\breve{\mathbf{B}}}}
\newcommand{\breveBbten}{{\breve{\mathbf{B}}^\otimes}}

\newcommand{\tSigma}{\widetilde{\Sigma}}

\newcommand{\bigB}{\bigotimes_{E'\anabelmap E} \left(\oplus_{j=1}^\ells B_{E'} \right)}

\newcommand{\thetajh}{{\widehat{\Theta}}_{\scriptscriptstyle{Joshi}}}

\newcommand{\sI}{\mathscr{I}}
\newcommand{\utheta}{U_\Theta}

\newcommand{\logsh}[1]{\mathscr{I}({#1})}
\newcommand{\logshq}[1]{\mathscr{I}^{\Q}({#1})}
\newcommand{\logshqp}[1]{\mathscr{I}^{\Q_p}({#1})}
\newcommand{\logshqpv}[1]{\mathscr{I}^{\Q_{p_v}}({#1})}
\newcommand{\logshe}{\logsh{E}}
\newcommand{\logsheq}{\logshq{E}}
\newcommand{\logsheplv}{\logsh{L_v}}
\newcommand{\logshepqlv}{\logshq{L_v}}

\newcommand{\bS}{\mathbb{S}}

\newcommand{\sIp}{\sI_p}
\newcommand{\sIpq}{{\sI}_p^{\Q_p}}
\newcommand{\sIm}{\mathscrbf{I}_{\scriptscriptstyle{Mochizuki}}}
\newcommand{\sImj}{\mathscrbf{I}_{\scriptscriptstyle{Joshi}}}
\newcommand{\sImq}{\mathscrbf{I}_{\scriptscriptstyle{Mochizuki}}^{\Q}}
\newcommand{\sImjq}{\mathscrbf{I}_{\scriptscriptstyle{Joshi}}^{\Q}}
\newcommand{\tsIm}{\widetilde{\mathscrbf{I}}_{\scriptscriptstyle{Mochizuki}}}
\newcommand{\tsImj}{\widetilde{\mathscrbf{I}}_{\scriptscriptstyle{Joshi}}}
\newcommand{\tsImq}{\widetilde{\mathscrbf{I}}_{\scriptscriptstyle{Mochizuki}}^{\Q}}
\newcommand{\tsImjq}{\widetilde{\mathscrbf{I}}_{\scriptscriptstyle{Joshi}}^{\Q}}
\newcommand{\bsIm}{\boldsymbol{\mathscrbf{I}}_{\scriptscriptstyle{Mochizuki}}}
\newcommand{\bsImj}{\boldsymbol{\mathscrbf{I}}_{\scriptscriptstyle{Joshi}}}
\newcommand{\bsImq}{\boldsymbol{\mathscrbf{I}}_{\scriptscriptstyle{Mochizuki}}^{\Q}}
\newcommand{\bsImjq}{\boldsymbol{\mathscrbf{I}}_{\scriptscriptstyle{Joshi}}^{\Q}}
\newcommand{\Vol}{{\rm Vol}}
\newcommand{\mulog}[2]{{\logVol}_{#1}({#2})}

\newcommand{\EB}{{}^EB} 
\newcommand{\flog}{\mathfrak{log}}

\newcommand{\arith}[1]{\mathfrak{arith}(#1)}
\newcommand{\arithl}{\arith{L}}
\newcommand{\adel}[1]{\mathfrak{adel}(#1)}
\newcommand{\adell}{\adel{L}}
\newcommand{\by}{{\bf y}}
\newcommand{\uby}[1]{{\underline{#1}}}
\newcommand{\ubyz}{{\underline{\by}}_0}
\newcommand{\bz}{\mathbf{z}}

\newcommand{\syflwwp}{\sY_{\C_{p_w}^\flat,L_w'}}
\newcommand{\syflwp}{\sY_{\cpt,L_w'}}
\newcommand{\sxflwo}{\sX_{\cpt,L_w'}}
\newcommand{\syflv}{\sY_{\lvbht,L_v}}
\newcommand{\sxflv}{\sX_{\lvbht,L_v}}
\newcommand{\bsY}{\mathscrbf{Y}}
\newcommand{\bsX}{\mathscrbf{X}}
\newcommand{\yadl}{\bsY_L}
\newcommand{\yadlmod}{\bsY_{\lmod}}
\newcommand{\yadq}{\bsY_{\Q}}
\newcommand{\xadl}{\bsX_L}
\newcommand{\xadq}{\bsX_{\Q}}
\newcommand{\yadlpoint}{\{(L_\wp\into K_{\wp}, K_\wp\isom \cpt)\}_{\wp\in\vlnon}}
\newcommand{\xadlpoint}{\{(L_\wp\into K_{\wp}, K_\wp\isom \cpt)\}_{\wp\in\vlnon}}

\newcommand{\yadlp}{\bsY_{L'}}
\newcommand{\yadlmodp}{\bsY_{\lmod'}}
\newcommand{\xadlp}{\bsX_{L'}}

\newcommand{\Blp}{\mathbb{B}_{L'}}
\newcommand{\Btlp}{\widetilde{\mathbb{B}}_{L'}}
\newcommand{\Blpp}{\mathbb{B}^+_{L'}}
\newcommand{\Btlpp}{\widetilde{\mathbb{B}}^+_{L'}}

\newcommand{\bsI}{\mathscrbf{I}}
\newcommand{\thetajti}{{\thetaj}^{\tilde{\mathscrbf{I}}}}
\newcommand{\thetaji}{{\thetaj}^{\mathscrbf{I}}}
\newcommand{\thetamti}{{\thetam}^{\tilde{\mathscrbf{I}}}}
\newcommand{\thetami}{{\thetam}^{\mathscrbf{I}}}

\newcommand{\fragments}{\ref{tab:rosetta-stone1}, \ref{tab:rosetta-stone2}, \ref{tab:rosetta-stone3}, \ref{tab:rosetta-stone5}, \ref{tab:rosetta-stone4}}
\newcommand{\rosettastone}{Rosetta Stone Fragments~\fragments}
\setstretch{1.1} 

\newcommand{\extciteI}[1]{\cite[{I-\ref{#1}}]{joshi-teich}}
\newcommand{\bvarphi}{\boldsymbol{\varphi}}

\addtocontents{toc}{\protect\setcounter{tocdepth}{3}}
\newpage
\section{Introduction}\label{se:intro}
\subparat{The goals of this paper}\label{intro:goals} \textcolor{cyan}{[For a quick outline of results see \cref{intro:outline}.]} This paper is a continuation of \cite{joshi-teich-def} which lays the foundations of global (i.e. over a number field) theory of Arithmetic Teichmuller Spaces. Familiarity with the results of \cite{joshi-teich,joshi-untilts,joshi-teich-estimates,joshi-teich-def} is essential in reading this paper.   The present paper has two goals:
\benumlab
\item To provide a `Rosetta Stone'  for a parallel reading of \iutabc\ and \present. This is established in \cref{se:intro-rosetta-stone}  with \rosettastone\ and (also  see \cref{se:frobenioids}, \cref{se:perfect-frob-perfectoid}). 
\item The second goal of this paper is to establish the assertion known as (Shinichi) Mochizuki's Corollary~3.12. This is introduced in \cref{intro:moccor1} while the (global) geometric case is discussed in \cref{intro:geom-moccor}.  [Details appear in (\cref{se:construction-of-thetaj-and-thetam}, \cref{se:fundamental-estimate-joshi}, \cref{se:mochizuki-construction-thetam}) and the geometric case (\cref{se:appendix-geom-case}).]
\eenum
The present paper deals exclusively with the (global) arithmetic case of the said corollary.  This requires  formulating an Arithmetic Teichmuller Theory for Number Fields which is detailed in \cite{joshi-teich-def} (Mochizuki asserts a less precise version of such a theory).  

In \cite{mochizuki-iut4}, Mochizuki claims that the arithmetic case, \moccor, leads to the proof of several important Diophantine conjectures, including the $abc$-Conjecture. This is taken up in \cite{joshi-teich-abc}. 

A detailed introduction to the present paper follows.

\subparat{Classical Teichmuller Theory as a model for the Arithmetic Case}\label{intro:classical-teich} 
A familiarity with classical Teichmuller Theory is highly recommended for  understanding the arithmetic case of \moccor\ established in \present. 

From an algebro-geometric perspective, classical Teichmuller Theory may be understood as providing a natural variation of the periods of Riemann surfaces which complements and refines the variation of periods constructed by Riemann (for explicit examples see \cref{se:appendix-geom-case}). 

Likewise, the theory of \present\ and \iut\ (via \rosettastone), both  arise from the existence of natural variations of $p$-adic periods, at all primes $p$ simultaneously, of  a fixed variety over a number field. For an explicit example of such a variation see \constrtwoh{\ssep}{se:appendix2}. 

So the perspective of \present\ (and \iut) is to create, in the arithmetic case, a theory which parallels the classical Teichmuller Theory. [Mochizuki's methods do not extend to higher genus curves or to higher dimensions. Specifically, the discussion in \iut,  is restricted to elliptic curves definable over a fixed number field. So my approach has an added advantage.]

\subparat{Mochizuki's Key Principle of Inter-Universality: The role of arbitrary geometric base-points}\label{intro:geom-base} 
According to Mochizuki,  his \textit{Key Principle of Inter-Universality} \cite[\ssep I3, Page 25]{mochizuki-iut1}, lies at the foundation of the theory which he has called \textit{Inter-Universal Teichmuller Theory}. In brief, Mochizuki's Key Principal of Inter-Universality requires one to work with arbitrary geometric base-points for tempered fundamental groups.

To put it succinctly: \textit{Mochizuki's Key Principle of Inter-Universality i.e. working with arbitrary geometric base-points is an anabelian proxy for deformations (of arithmetic).} My work  provides the most natural way of quantifying this statement (this is discussed below in \cref{intro:deform1}, \cref{intro:deform2}). 

As noted in \cite{joshi-untilts, joshi-teich-quest}, formulating the theory using Mochizuki's Key Principle of Inter-Universality i.e. using arbitrary geometric base-points necessarily implies that (arbitrary) algebraically closed perfectoid fields enter \iut\ (and  one cannot limit one-self to $\C_p$ (for all primes $p$)). This necessarily forces one to take the path I have taken in \present.  \textcolor{cyan}{[A detailed discussion of Mochizuki's principle (and how it appears in my theory) is given in \cref{ss:geom-base-points}.]}

It should be noted that, no discussion of the centrality of the role of arbitrary geometric base-points in \iut\ (required by Mochizuki's key principle) appears in any of \cite{fesenko-iut}, \cite{fucheng}, \cite{scholze-stix},  \cite{yamashita}, \cite{dupuy2020statement,dupuy2020probabilistic}, \cite{scholze-review}, \cite{saidi-review-iut}.

\subparat{Number Fields are Deformable}\label{intro:deform1} Mochizuki's presentation of his ideas in \iut\ is quite labyrinthine.
Hence, unfortunately, even the sharpest of readers can loose  sight of the profoundly original idea within \iut. So let me explicate it for the readers: for the past hundred and fifty years, since the works of Ernst Kummer and Richard Dedekind, number theorists have viewed a fixed number field as a static background in which the drama of Diophantine Geometry unfolds. Mochizuki asserts that this static background itself is a dynamic, pliable, deformable entity i.e. the arithmetic of a fixed number field itself is bendable. [Mochizuki offers no convincing proof of this claim in \iut. In \cite{joshi-teich-def}, I have established this claim quite robustly and independently by demonstrating  the existence of topological deformations of the  arithmetic of a fixed number field (and hence my use of the term bendable is quite appropriate). So at this juncture, Mochizuki's assertion itself is true.] 

Moreover, Mochizuki posits that this bendability of a number field may be harnessed to establish one of the most sought-after, fundamental conjectures in Diophantine Geometry. That is why I believe Mochizuki's idea is  pioneering and truly remarkable.

\subparat{Quantification of Deformations of Number Fields}\label{intro:deform2} While Mochizuki has asserted that his approach may be understood as proceeding by dismantling, distorting and reassembling of the arithmetic of a fixed number field \cite[Page 12]{mochizuki-iut1},  nevertheless, a precise mathematical quantification of what exactly this means  has remained elusive in \iut\ and its numerous expositions. That is why a clear understanding of its claims has not emerged in the decade following the initial public release of \iut\ and skepticism regarding \iut\ has prevailed and solidified with the appearance of \cite{scholze-stix} (my discussion of this is in \cite{joshi-teich-quest,joshi-teich-summary-comments}). 

In \cite{joshi-teich-def}, I have established the most precise version of Mochizuki's claim. Notably, I demonstrate that primes of a fixed number field provide, in a precise sense, multiplicative coordinates for the number field and that  the claimed distortions arise from general changes of these prime coordinates in a manner reminiscent of the Virasoro Uniformization  in Conformal Field Theory \cite{beilinson88}; and by   \constrtwoh{Theorem }{th:consequence-inequivalent-arithmeticoids} these changes also distort all local analytic geometries at all primes. 

Importantly, my work makes it clear that non-trivial distortions of the fixed number field (the theory literally permits topologically bending the rules of arithmetic in the precise sense of \cite{joshi-teich-def}) and   $p$-adic metrics (at all primes),  do arise and are  mediated by these general changes of prime coordinates (readers familiar with the role of Virasoro Uniformization in Conformal Field Theory will recognize why my usage of the term `bendability of the number field is quite apt). [A less precise version of this appears as Mochizuki's Indeterminacy Ind2 in \iut. Unfortunately a precise demonstration of the existence and quantification of the properties of Ind2 does not seem possible without invoking my approach. Without a precise quantification, it is not possible to say, for example, when two arithmetic holomorphic structures differ by an Indeterminacy of type Ind2. I also demonstrate precisely how Mochizuki's indeterminacy Ind3 arises. [My discussion, with proofs, of how Mochizuki's Three indeterminacies arise is in \cref{ss:Mochizuki-indeterminacies}.] These are precisely the sort of assertions which are needed for a rigorous and satisfactory formulation (and proof) of \moccor.]

\subparat{The role of the Absolute Grothendieck Conjecture}\label{intro:abs-groth} Let me say a few words which will clarify the role of the absolute Grothendieck conjecture in my theory (and Mochizuki's theory). Understanding the consequences of the existence of realms where one has validity and invalidity of this conjecture turns out to be central in arriving at a clear understanding of my theory and Mochizuki's theory.  However, discussion of the failure  of this conjecture (i.e. the existence of realms of its invalidity) does not occur in \iut\ at all as it stays focused on validity of the conjecture. 

First let me discuss the realm where one has the validity of this conjecture. Mochizuki in \cite[Corollary 2.3]{mochizuki07-cuspidalizations-proper-hyperbolic} established that the absolute Grothendieck Conjecture is true for incomplete hyperbolic curves over $p$-adic fields arising from hyperbolic curves definable over number fields. As noted in \cite{joshi-untilts} this validity assertion is equivalent to the assertion that the arithmetic Teichmuller space of such a curve is connected and the field of moduli of this curve is its field of definition. [This is analogous to the theorem that a Shimura variety is defined over a suitable reflex field which is the minimal field of definition of its points definable over a number field.].

Now let me discuss two realms where the absolute Grothendieck conjecture fails to hold and why both the instances of its failure are central to the two theories. It is well-known   that the absolute Grothendieck conjecture fails for hyperbolic Riemann surfaces and this failure arises from the existence of Teichmuller and moduli spaces of Riemann surfaces. Pursuing this idea in the $p$-adic context, in \cite{joshi-teich,joshi-teich-summary-comments}, I established that the absolute Grothendieck conjecture is false in the category of Berkovich spaces over algebraically closed, complete (rank one) valued fields containing $\Q_p$ and that this signals the existence of an Arithmetic Teichmuller Theory (this is one of the principal results of \cite{joshi-teich,joshi-untilts}). 

The second realm where the absolute Grothendieck conjecture fails was observed in \cite{joshi-gconj} where I established that the absolute Grothendieck conjecture fails  for a complete Fargues-Fontaine curve $\sxcpte$ for suitable $p$-adic fields (and all primes $p$). By the main theorem of \cite{joshi-teich}, one may view $\sxcpte$ as an approximate  parameter space of arithmetic Holomorphic structures. Since the idea of \cite{mochizuki-iut3,mochizuki-iut4} is to average over arithmetic holomorphic structures, considering all Fargues-Fontaine curves $\sxcptep$ with \'etale fundamental groups isomorphic to $\pi_1^{et}(\sxcpte)\isom G_E$ (for all primes $p$) ensures that this averaging of \cite{mochizuki-iut3,mochizuki-iut4} takes place over the largest possible set of arithmetic holomorphic structures permitted by the theory! In other words, Mochizuki's Indeterminacy Ind1 (at all primes) arises from the failure of the absolute Grothendieck conjecture for Fargues-Fontaine curves \cite{joshi-gconj} (at all primes). This is how Mochizuki's Indeterminacy Ind1 emerges from the point of view of \present. [A detailed discussion of Mochizuki's Indeterminacies is in \cref{se:intro-rosetta-stone}.]

However, working with $\sxcpte$ has one major drawback: residue fields of the closed classical points of $\sxcpte$ do not come equipped with natural valuations. This is why one must work with incomplete Fargues-Fontaine curves $\sY_{\cpt,E}$ whose closed classical points provide residue fields equipped with valuations \cite{fargues-fontaine}. This important point, together with the central observation of \cite[Appendix]{joshi-teich-estimates} that the Frobenius of $\sY_{\cpt,E}$ has all the properties Mochizuki asserts for his $\flog$-link, establishes unequivocally that the definition of arithmetic holomorphic structures made in \cite{joshi-untilts} is the most natural one even for \iut. 

The next central point in \cite{joshi-teich-estimates} is the constructions of certain correspondences on $\sY_{\cpt,E}$ which have all the properties of Mochizuki's $\Theta_{gau}$-Link and notably I demonstrate precisely how Mochizuki's asserted valuation scaling property (of a $\Theta_{gau}$-Link) and especially how a certain scaling factor $j^2$ (for $j=1,2,\ldots,\frac{\ell-1}{2}$) appears non-trivially (the global arithmetic consequence of this factor is established in \Cref{cor:val-norm-theta}--this assertion is  \cite[Theorem A(ii)]{mochizuki-iut1}). The  non-triviality of this factor is a key sticky point in \cite{mochizuki-iut3}  and especially in proof of the main theorem of \cite{mochizuki-iut4} because of which its existence is denied in \cite[Section 2.2]{scholze-stix}. My proof of the non-triviality of this factor is \Cref{th:adelic-theta-link} (the local case is detailed in \cite{joshi-teich-estimates}) and its global arithmetic consequences are given by \Cref{cor:val-norm-theta}. 

\subparat{Global aspects: The central role of global arithmetic of number fields}\label{intro:global-aspects} 
Let me say that the theory of the present series of papers is a global arithmetic theory and I want to point out how these global aspects play a central role in this paper. There are several important new ingredients needed for the global version of the discussion of the previous paragraph.  One key ingredient is the theory of (global)  Arithmetic Teichmuller Spaces of Number Fields described in \cite{joshi-teich-def} which provides the non-trivial deformations of any fixed number field; several other new ingredients which are required is detailed here. Let me briefly describe the global aspects. 

One key feature of the global theory of \cite{joshi-teich-def} is the fact that each avatar $\arithl_\by$ of a fixed number field $L$  defines a hyperplane (see \constrtwoh{\ssep}{ss:period-map-prod-form}) given by the validity of the product formula for the number field $$H_\by\subset \bigoplus_{v\in\mathbb{V}_L} \R.$$ 

While in \cite{mochizuki-iut3}, Mochizuki has indicated the key role of the product formula through the appearance of global realified Frobenioids in his theory, a key innovation, \constrtwoh{Theorem }{th:hyperplane} and \constrtwoh{Remark }{rem:hyperplane}, of \cite{joshi-teich-def} is that the product formula for the normalized arithmeticoid $\arithl_\by^{nor}$ defines a new arithmetic period mapping  
\begin{equation}\begin{tikzcd}
\yadl\arrow[r, "\by\mapsto H_\by"] & \P(\oplus_{v\in\mathbb{V}_L}\R)
\end{tikzcd}
\end{equation}
On one hand, the valuations on untilts parameterized by a Fargues-Fontaine curve at some arbitrary prime (e.g. $\syQp$) cannot be all simultaneously normalized and on the other hand, local changes of valuations are globally constrained by the validity of the product formula. These two observations imply that the above arithmetic period mapping is a highly non-trivial global arithmetic period mapping in which an arithmetic avatar $\arithl_\by$ provides a hyperplane $H_\by$ (as above) which moves under symmetries (such as the actions of the absolute Galois group, global Frobenius and $L^*$) established in \constrtwoh{Theorem }{th:galois-action-on-adelic-ff}. [Let me remark that \cite[Subsection 2.2]{scholze-stix} argue incorrectly that this hyperplane $ H_\by$ does not move--this incorrect argument is an important point in their (incorrect) conclusions regarding \iut. For a detailed analysis of \cite{scholze-stix}, \cite{scholze-review} see \cite{joshi-report}.] 

This period mapping combined with Mochizuki's $\Theta_{gau}$-links (encapsulated in the natural constructions of \cref{ss:Mochizuki-theta-gau-link} and \constrtwo{\ssep}{se:ansatz}) which scale valuations (via  \Cref{th:hyperplane-theta-lnk} and \Cref{re:hyperplane-theta-lnk}) at primes of semi-stable reduction of an elliptic curve, also force renormalizations of valuations at primes of good or non-semistable reduction so that the product formula continues to hold (\Cref{cor:val-norm-theta}).  That these phenomenon surrounding these global product formula hyperplanes do occur is asserted without adequate proof in \cite{mochizuki-iut3} (\Cref{re:hyperplane-theta-lnk} provides a detailed discussion of this point). However, the \textit{global-product-formula-as-a-period-mapping} optik of \cite{joshi-teich-def} is my own innovation beyond \iutthr.

Another important point is that  the aforementioned (global) action of $L^*\act\yadl$ (given by \constrtwoh{Theorem }{th:galois-action-on-adelic-ff}), which plays an important role in the aforementioned period mapping is also incorporated in the construction of Mochizuki's Adelic Ansatz given in \cref{se:adelic-ansatz}.  Hence through this, the action $L^*\act\yadl$  also intervenes in the construction of the theta-values locii constructed here through the stability properties of Mochizuki's Adelic Ansatz (\Cref{def:mochizuki-adelic-ansatz}) with respect to this global action (\Cref{th:adelic-theta-link}).  The effect of this action in the Diophantine contexts of heights is given in \constrtwoh{\ssep}{se:heights}. To put it succinctly, \textit{global properties of the number field $L'$ are woven into the construction of Mochizuki's Global Ansatz $\tSigma_{L'}$ (through \Cref{th:adelic-theta-link}) and hence into the construction of all the theta-values sets constructed in this paper.}

Let me also remark that the existence of Initial Theta Data (see \cref{se:theta-data}) which is vital to the constructions of the theta-values locii in \cref{se:construction-of-thetaj-and-thetam} and hence to \moccor\ is a global arithmetic fact (which cannot be proved locally). While this existence is detailed in \cite{joshi-teich-abc} (Mochizuki's proof of this existence is in \cite{mochizuki-iut4}), it is central to the present paper. My approach to constructing (global) theta-values locii is qualitatively the same as that adopted in \cite[Theorem 3.11]{mochizuki-iut3}.

\textcolor{red}{At any rate, these are some of the many ways in which the theory of this paper and \cite{joshi-teich-def} is a global theory.} 

\subparat{Mochizuki's Corollary 3.12}\label{intro:moccor1} Now to the main topic of this paper, namely \moccor\ and its proof (\cref{se:construction-of-thetaj-and-thetam}, \cref{se:fundamental-estimate-joshi}, \cref{se:mochizuki-construction-thetam}). Since the claim of \iut\ is so remarkable, its ideas rather original, a complete proof of \moccor, is surely needed for anyone wanting to be persuaded of the validity of the central assertion of \cite{mochizuki-iut4}. Unfortunately, the approach to this corollary in \iut\ is so convoluted that many readers have arrived at a limited understanding of what the corollary claims and what methods are used  for arriving at its conclusion. Mochizuki's \cite{mochizuki-gaussian} is readable but offers  rudimentary and unworkable analogies while \cite{mochizuki-essential-logic} deflects the main mathematical existential question, raised by \cite{scholze-stix} and myself (as noted earlier, in \cite{joshi-untilts-2020,joshi-teich,joshi-untilts}, I have successfully addressed the question of the existence of arithmetic holomorphic structures head-on by constructing them explicitly). 

\subparat{The geometric case of Mochizuki's Corollary 3.12}\label{intro:geom-moccor}
As is well-known, Szpiro formulated  the Arithmetic Szpiro Conjecture based on his discovery (and proof) of the Geometric Szpiro Conjecture (\cite{szpiro81}, \cite{szpiro1991}). Since then several proofs of the Geometric Szpiro Conjecture have appeared. The proofs  \cite{bogomolov00} and \cite{zhang01} are Teichmuller theoretic in their approach. In \cref{se:appendix-geom-case}, I provide a discussion, from the point of view of this paper and of \cite{mochizuki-iut3}, of the proofs geometric Szpiro inequality due to \cite{bogomolov00}, \cite{zhang01}, exhibiting, in these proofs the existence of such avatars of a Riemann surface and related symmetries--especially the role of a global Frobenius morphism. Notably, in \cref{th:geometric-case-of-moccor}, I  prove the geometric case  of \cite[Corollary 3.12]{mochizuki-iut3}  in the context of Geometric Szpiro Inequality as established in \cite{bogomolov00} and \cite{zhang01}. Readers are strongly advised to read formulation  of the geometric case of \moccor\ given in \cref{se:appendix-geom-case} before proceeding with the arithmetic case which occupies the bulk of this paper.

\subparat{Many versions of Mochizuki's Corollary 3.12}\label{intro:moccor2} My important discovery is that, in fact, there are many versions of \moccor\ (Mochizuki considers a specific one). The formulation of any version of \moccor\ requires the construction of the respective theta-values sets--the main ones constructed here  in \cref{se:construction-of-thetaj-and-thetam}, \cref{se:mochizuki-construction-thetam} are $\thetaj^{\Blp}$, $\thetami$, $\thetaji$ (and their varianFts). 

It should be emphasized that the central difficulty of establishing \moccor\ is the  construction of any of these theta-values sets with a certain list of properties. On the other hand,   \cite{scholze-stix} have doubted that such a set could possibly exist and the conclusions of \cite{scholze-stix} are centrally dependent on their non-existence assertion.  \textit{However}, due to \present\ existence of such sets is no longer in doubt. 

Each of these sets  is the collation of theta-values arising from a certain  set, dubbed  \textit{Mochizuki's (global) Ansatz},  consisting of certain tuples of distinct arithmetic holomorphic structures  (each tuple  is a $\Theta_{gau}$-Link in the sense of \iut). The decorations ${\Blp}, \bsI$ refer to the codomain used for obtaining theta-values. To be sure, the construction of each of these sets requires the existence of all of the following

\benumlab
\item the choice of a global arithmetic datum  \textit{Initial Theta-Data} (known by this name in \iutabc).
\item The existence of a plurality of
\begin{enumerate}
\item  arithmetic holomorphic structures (at each prime), 
\item $\flog$-Links, 
\item $\Theta$-Links, 
\item global deformations of a number field,
\end{enumerate}
\item and the establishment of the symmetry properties of these data.
\item a common codomain for collating theta-values provided by {\bf(1)-(3)}.
\end{enumerate} 
[From my perspective, the existence of these data is the construction of Arithmetic Teichmuller Spaces and a precise quantification of their properties documented in \present.] 

Importantly, \cite[Theorem 3.11]{mochizuki-iut3} asserts the construction of $\thetami$--subject to the existence of {\bf(1)-(4)} above. The  set $\thetami$ is referred to therein as the \textit{multi-radial representation of theta-values} i.e. in the parlance of \cite[Theorem 3.11]{mochizuki-iut3}, each of the sets $\thetaj^{\Blp},\thetaji,\thetami$ are  multi-radial representations of theta-values (with respect to their respective theta-value codomains). Here, $\bsI$ is Mochizuki's tensor packet codomain for theta-values. 
[Beware that the superscript $\bsI$ in $\thetami,\thetaji$ refers to $\bsIm$ and $\bsImj$ respectively and these two codomains are slightly different versions of each other, my version being simpler of the two. In \cite[\ssep 3]{mochizuki-iut3}, $\bsIm$ is called the ``tensor product log-shell.''] There are variants $\thetamti,\thetajti$ in which  $\bsI$ is replaced by a product codomain which is simpler and adequate for obtaining lower bound given by \moccor\ ($\Blp$ in $\thetaj^{\Blp}$ refers to a product codomain and it also comes with its tensor product version). An important point is that $\Blp$ is a ring  (while $\bsI$ is not) and so in this ring one can multiply theta-values arising from different arithmetic holomorphic structures. This allows me to demonstrate that  Mochizuki's tensor product structure (detailed in \cite[\ssep 3]{mochizuki-iut3}) arises essentially from the product structure in this ring. However, the tensor packet codomain $\bsI$ is better suited for tighter upper bounds on the theta-values sets and  is indispensable for \cite{joshi-teich-abc} and \cite{mochizuki-iut4}.  As I demonstrate here, the choice of any of these theta-values sets leads to a formulation of \moccor\ specific to this choice and hence one obtains several distinct avatars  of  \moccor\ (including Mochizuki's version). In each case, the proof of the fundamental estimate of \moccor\ proceeds along the lines of the proto-type of \moccor\ which I have established in \cite{joshi-teich-estimates}. In \cite[Appendix I]{joshi-teich-def}, I establish the geometric case of this corollary.  At this point, using \rosettastone\ (and the above data constructed in \present) one may also obtain a proof of \moccor\ along the lines of \cite{mochizuki-iut3}. Notably, I demonstrate that Mochizuki's gluing of Hodge-Theaters and Frobenioids along prime-strips also arises quite naturally from my point of view.

\subparat{Why this work is needed: The incompleteness of \iutthr}\label{intro:incomplete} 
The proof of \moccor, is incomplete  as the proof rests on  establishing the existence of a plurality of arithmetic holomorphic structures (and their symmetry properties). Here, I use the term plurality in the sense of being logically distinct. Thus for example, classical Teichmuller Theory asserts the existence of a plurality of complex structures on a fixed topological surface. One needs the existence of arithmetic holomorphic structures at all primes of a number field, and deformations of a number field itself. These central points are claimed, but not established in \iutthr. This was first pointed out in \cite{scholze-stix}, and later, more robustly in my works. On the other hand, \cite{scholze-stix} (and \cite{scholze-review}) ignore Mochiuzki's Key Principle of Inter-Universality (see \cref{intro:geom-base}) and also arrived at the incorrect conclusion that distinct arithmetic holomorphic structures cannot  exist at all. That assertion of \cite{scholze-stix,scholze-review} has been  disproved in \cite{joshi-untilts-2020,joshi-teich,joshi-untilts}, by constructing large families of distinct arithmetic holomorphic structures and establishing the existence of deformations of a fixed number field in \cite{joshi-teich-def} (the existence of such deformations is also asserted by Mochizuki). [For a detailed analysis of \cite{scholze-stix}, \cite{scholze-review} see \cite{joshi-report}.] My approach to arithmetic holomorphic structures, detailed in \cite{joshi-teich,joshi-untilts,joshi-teich-def}, puts the theory on par with the classical notion of complex holomorphic structures and provides arithmetic holomorphic structures in Mochizuki's group theoretic sense. Notably,  Mochizuki's group theoretic arithmetic holomorphic structures may be distinguished by actual Berkovich analytic (i.e. holomorphic) structures and one obtains Mochizuki's group theoretic arithmetic holomorphic structures by applying the (tempered) fundamental group functor to the relevant Berkovich analytic spaces. My theory may also be applied to a number field and this provides a deformation theory of number fields (\cite{joshi-teich-def}), whose existence is also asserted in \iut.

\subparat{On \cite{mochizuki-essential-logic}}\label{intro:essen-logic} Mochizuki  in his rejoinder \cite{mochizuki-essential-logic} to \cite{scholze-stix}  asserted that their argument is incorrect. However, the arguments presented in \cite{mochizuki-essential-logic} do not convincingly establish the said existence but rather mandates the  existence of distinct arithmetic holomorphic structures. [Read in the classical context of Riemann surfaces, this mandate argument is tantamount to  mandating that distinct complex structures exist on a topological surface.]

My point is that this mandate argument of \cite{mochizuki-essential-logic}, is 
\benumlab
\item wholly inadequate because  \iut\ also requires a way of distinguishing arithmetic holomorphic structures and precisely proving how these transform under symmetries of the theory (this is a key point needed in the construction of the set which underlies \moccor); 
\item also  completely unnecessary, since distinct arithmetic holomorphic structures \textit{do} exist for canonical reasons \cite{joshi-teich,joshi-untilts}. 
\eenum 
Notably, \cite{mochizuki-essential-logic} suggests the existence of distinct arithmetic holomorphic structures but it is not clearly established in \iutthr, nor can it be quantitatively established using the framework of \iutthr\ or \topics.

\subparat{Other articles on \iut}\label{intro:other} Apart from the present series of papers, several other articles on the subject of \iut\ have appeared:

\cite{fesenko-iut}, \cite{fucheng}, \cite{scholze-stix}, \cite{yamashita},  \cite{dupuy2020statement,dupuy2020probabilistic}, \cite{scholze-review}, \cite{saidi-review-iut}.  

Let me say a few words about these articles: \cite{fesenko-iut}, \cite{fucheng}, and \cite{yamashita} are surveys and they essentially repeat \iut\ verbatim on all its key points. In \cite{dupuy2020statement,dupuy2020probabilistic}, they do not state the inequality \moccor\ as a theorem, but work relative to Mochizuki's assertion. The reviews \cite{scholze-stix}, \cite{scholze-review} have been discussed by me in \cite{joshi-untilts}, \cite{joshi-teich-quest} and in the present paper. 

It should be noted that except for a passing reference to Mochizuki's philosophy of such structures in \cite{yamashita}, the term `arithmetic holomorphic structure' (coined and used by Mochizuki himself) does not appear in  any of  \cite{fesenko-iut},  \cite{fucheng}, \cite{scholze-stix}, \cite{yamashita}, \cite{dupuy2020statement,dupuy2020probabilistic}, \cite{scholze-review}, \cite{saidi-review-iut}. 

But importantly, as Mochizuki asserts in  \cite[Page 24]{mochizuki-iut1}, arithmetic holomorphic structures are central to \iut\ and  \present\ is the only one which treats arithmetic holomorphic structures in detail and establishes their existence in the context of \iut.

\subparat{Two points of departure from \iutthr}\label{intro:differences} Apart from the construction of many theta-value sets, there are two important points of departure between my treatment of \moccor\ and that of \cite{mochizuki-iut3}. 

First point arises from my observation that Mochizuki's rational log-shell can be naturally identified with a suitable Bloch-Kato subspace, of the form $H^1_e(-,-)$, which consists of ``Frobenius trivial, crystalline'' Galois cohomology classes. [The reason for the quotes is a bit technical and as follows: such classes correspond to working with $B_e=\bcris^{\varphi=1}$ rather than with $B$, but beware that since construction of $\bcris$ requires fixing an algebraic closure of $\Q_p$, a choice of $\bQ_p$ is tacit in writing this equality, while in the present theory $\bQ_p$ moves with an arithmetic holomorphic structure (the construction of $B$ does not require a choice of $\bQ_p$). Hence this equality is an isomorphism in the context of \present. Because of this, the  relevant  Frobenius action on the Galois cohomology is by a non-identity automorphism (in general).] The identification of this cohomology with Mochizuki's rational log-shell is via the inverse of the relevant Bloch-Kato exponential i.e. with a relevant Bloch-Kato logarithm. Mochizuki's integral log-shell is the image of the relevant integral subgroup of this rational cohomology. This is the  genesis of Mochizuki's integral log-shell and one now has a natural reason why it is relevant in the first place.   [To understand why ``Frobenius triviality'' is relevant in the context of the arithmetic Szpiro Conjecture, see \cite[Appendix I]{joshi-teich-def} where this property appears in  the proof of the Geometric Szpiro Conjecture of \cite{bogomolov00,zhang01}.] 

However, Theta-values do not live in this subspace. Hence,  I construct ``Frobenius trivial, crystalline''  galois cohomology classes arising from a Tate elliptic curve (for any choice of a Tate parameter).  Mochizuki's approach to circumventing the fact that the theta-values do not live in the relevant subspace  is to use a multiplicative action  \cite[Figure 3.1]{mochizuki-iut3} of theta values  on the relevant local components of $\bsIm$ for defining $\thetami$. Using the inverse of the relevant Bloch-Kato exponential, one recovers Mochizuki's multiplicative approach from mine. 

The second important point of departure is that I work with products of norms (i.e. absolute values) in $\Q_p$-Banach spaces (in the case of $\thetami,\thetaji$, these spaces are finite dimensional over $\Q_p$) while Mochizuki works with volumes in tensor product of finite dimensional $\Q_p$-vector spaces. In the $p$-adic context, translation invariance of $p$-adic Haar measure (on a $p$-adic field $E$) means volumes are related to absolute values (i.e. $\Vol_E({\alpha+\beta\O_E})=\abs{\beta}_E$). This observation can be used to calculate (suitably weighted) volumes of tensor product regions of the type $\beta_1\O_{E_1}\tensor_{\Z_p}\beta_2\O_{E_2}\tensor\cdots \tensor_{\Z_p}\beta_n\O_{E_n}$, which are contained in the set $\thetami$ considered in \moccor,  by means of the cross-norm property of tensor product norm,  available in the context of the present paper, asserts that $$\abs{v\tensor w}_{V\tensor W}=\abs{v}_V\cdot\abs{w}_W$$
and allowing one to relate volumes  to absolute values (at least for the tensor product regions of the above type). This allows me to translate my inequalities (established using products of absolute values) to Mochizuki's tensor product volume inequality asserted in \moccor. On the other hand, the point of \cite{mochizuki-iut4} is that the volume of the region occurring in \moccor\ may be bounded from above by more direct means (complementing the lower bound given by \moccor). Let me also remark, at this juncture, that  a proof of \moccor\ may also be obtained by combining \cite{mochizuki-iut3} and \present, and \rosettastone. However, that is not a path I have chosen to follow here.

\subparat{An outline of this paper}\label{intro:outline} Here is a brief outline of this paper. A parallel reading of \present\ is strongly recommended. In \cref{se:holomorphoids}, I setup the basic frame work of holomorphoids (expanding on \cite{joshi-teich-def}) and demonstrate the existence of distinct scheme theories (this is asserted in \iutthr) and specialize quickly to holomorphoids of elliptic curves over number fields. In \cref{se:theta-data},  I introduce the global data which is called an Initial Theta Data in \cite{mochizuki-iut1}. \cref{se:adelic-ansatz} is dedicated to the construction of certain correspondences which are global in nature (local constructions are in \cite{joshi-teich-estimates}) and referred to in \iutthr\ as $\Theta_{gau}$-Links. In \cref{se:local-preliminaries}, I recall some rings of $p$-adic Hodge Theory. \cref{se:intro-rosetta-stone} provides a `Rosetta Stone' (\rosettastone) for translating between \iut\ and \present. In \cref{se:construction-of-thetaj-and-thetam}, \cref{se:mochizuki-construction-thetam}, I provide a detailed  construction of the many flavors $\ttheta_{\scriptscriptstyle{Joshi}}^{\Blp}$, $\thetami$ and $\thetaji$ of the fundamental theta-values set required for \moccor. In \cref{se:fundamental-estimate-joshi}, the fundamental estimate i.e. \moccor\ is established for my version $\ttheta_{\scriptscriptstyle{Joshi}}^{\Blp}$ of this set and in  \cref{se:mochizuki-construction-thetam}, the fundamental estimate is established for Mochizuki's version $\thetami$  and also $\thetaji$ (this is \Cref{cor:cor312}). In \cref{se:frobenioids}, I establish the claims of \iutthr\ regarding Frobenioids and Hodge-Theaters especially the existence of many distinct Hodge-Theaters of each type considered in \iutthr. In \cref{se:perfect-frob-perfectoid}, I provide a comparison of perfect Frobenioids and perfectoids. In \cref{se:appendix-geom-case}, I provide a formulation (with proof)  of the geometric case of \moccor\ in \Cref{th:geometric-case-of-moccor}.

\subparat{Acknowledgments}\label{intro:ack}
I thank Taylor Dupuy for reminding me to clarify the relationship between holomorphic hulls in the sense of \cite{mochizuki-iut3} and convex closures which I use here; Mochizuki has also reminded me of the same point in his  \href{https://www.kurims.kyoto-u.ac.jp/~motizuki/Report\%20on\%20a\%20certain\%20series\%20of\%20preprints\%20(2024-03).pdf}{recent comments}. This is addressed in \Cref{pr:convexes-and-hulls}. I have also addressed   another of Mochizuki's concern about clarifying the global aspects of my theory (see \cref{intro:global-aspects} for details). I thank Peter Scholze for some conversations surrounding this paper.

\addtocontents{toc}{\protect\setcounter{tocdepth}{3}}
\section{Holomorphoids of a variety over a number field}\label{se:holomorphoids}
\nwss

\newcommand{\V}{\mathbb{V}}
\newcommand{\vq}{\V_{\Q}}
\newcommand{\vl}{\V_L}
\newcommand{\vnl}{\V_L^{non}}
\newcommand{\vlnon}{\vnl}
\newcommand{\vlarch}{\mathbb{V}_L^{arc}}
\newcommand{\vlp}{\V_{L'}}
\newcommand{\vlpp}{\V_{L',p}}
\newcommand{\vnlp}{\V_{L'}^{non}}
\newcommand{\vlnonp}{\vnlp}
\newcommand{\vlarchp}{\mathbb{V}_L^{'arc}}
\newcommand{\fA}{\mathfrak{A}}
\newcommand{\lvbh}{\widehat{\overline{L}}_v}
\newcommand{\lvbht}{{\widehat{\overline{L}}_v}^\flat}

In this paper, I will freely use the results of \cite{joshi-teich,joshi-teich-estimates,joshi-untilts,joshi-teich-def} and familiarity with the ideas of \cite{joshi-formal-groups,joshi-anabelomorphy} will also be useful. Familiarity with classical Teichmuller Theory (for example the basic texts \cite{imayoshi-book} or \cite{nag-book} or  \cite{lehto-book}) is absolutely essential in understanding the above papers--especially an understanding of moduli of curves will not be adequate. The terms \textit{anabelomorphism}, \textit{amphoric} and \textit{unamphoric}, which were introduced in \constranab{\ssep}{ss:anab-defs} and, will be used here throughout. 

\subparat{Holomorphoids} The global arithmetic input for the present theory is developed in \cite{joshi-teich-def} and that is a starting point for the present paper. \textcolor{red}{However, I will not recall the definition of Arithmetic Teichmuller Spaces given in \cite{joshi-untilts,joshi-teich-def,joshi-teich} but work with points of such spaces and these (points) will be precisely defined below and called (global/local) holomorphoids.}

Let $L$ be a number field with no real embeddings (this is a simplifying assumption and not a technical restriction). I will use the theory of arithmeticoids of a number field developed in \cite{joshi-teich-def}.  For the definitions, notations and the basic properties of  arithmeticoids, the reader is referred to \cite{joshi-teich-def}. Especially, parallel reading of \cite{joshi-teich-def} is highly recommended. The adelic Fargues-Fontaine curves $\yadl$ were introduced in  \constrtwoh{\ssep}{se:adelic-ff-curves} and these curves and their properties will also be used here throughout.

I will also assume throughout this paper that all arithmeticoids $\arithl_\by$ are normalized i.e. $\arithl_\by=\arithl_\by^{nor}$  \constrtwoh{\ssep}{ss:normalization} i.e. the valuations on each component of $\by$ are normalized so that the  Artin-Whaples product formula (\cite{lang-diophantine})  holds for each $x\in L^*$:
\be\label{eq:prod-formula} 
\prod_{v\in\V_L}\abs{x}_{K_v}=1.
\ee

An important point demonstrated in \constrtwoh{Theorem }{th:hyperplane} is that $\yadl$ is equipped with a certain global arithmetic period mapping which is given by the validity of the product formula  \eqref{eq:prod-formula} in each arithmeticoid $\arithl_\by^{nor}$. As pointed out in \constrtwoh{Remark }{rem:hyperplane}, the global arithmetic period mapping is non-constant. The validity of the product formula in a normalized arithmeticoid and the global arithmetic period mapping are central in understanding the global arithmetic aspects of this theory. The relationship between the main construction of this paper and the global arithmetic period mapping of  \constrtwoh{Theorem }{th:hyperplane} is established here in \Cref{th:hyperplane-theta-lnk}. 

The following definition is motivated by \constrone{Definition }{def:arith-hol-strs} (same as) \constrone{Definition }{def:arith-hol-space-local}, \cite{joshi-teich-def}, and  holomorphoid defined in it is a point of the global Arithmetic Teichmuller Space \cite{joshi-teich,joshi-untilts,joshi-teich-def} and each local holomorphoid defined below  is a point of the local Arithmetic Teichmuller Space of \constrone{Definition }{def:arith-hol-space-local} (see \Cref{le:point-lemma}). 
\begin{defn}\label{def:holomorphoids}\
A \textit{holomorphoid of a scheme over a number field}, denoted $(X,L\into \arith{L}_\by)$ (or simply $(X,\arith{L}_\by)$), consists of
an arithmeticoid $\arithl_\by$ of a number field $L$ (with $\by\in\yadl$)  (i.e. one given is a choice of an arithmeticoid $\by\in\yadl$), and a scheme $X/L$ over $L$ (viewed as a scheme $X/\arithl_\by$ i.e. as a scheme over $L$ referencing the arithmetic of $L$ given by $\arithl_\by$ (see \constrtwoh{Theorem }{th:existence-inequivalent-arithmeticoids}). Write $\by=(y_v)_{v\in\vl}\in\yadl$, 
then for each $v\in \vl$, one is also given a morphism of analytic spaces 
\be \label{eq:geom-base-point-holomorphoid}
\sM(K_{y_v})\to \xan_{L_v}.
\ee
Morphisms of holomorphoids may be defined in the obvious way (\cite{joshi-teich-def}). Owing to \eqref{eq:geom-base-point-holomorphoid}, one may think of a holomorphoid as a pointed object and I will use the term \textit{pointed holomorphoid} when one wants to emphasize this aspect of a holomorphoid.
One may also consider holomorphoids $(X,L\into \arith{L}_\by)$ with $\by\in \yadl^{max}$ (resp. $\by\in \yadl^{\R}$)  which are constructed in \constrtwoh{Definition }{def:adelic-ff-curves-maximally-complete} (resp. \constrtwoh{Definition }{def:adelic-ff-curves-realified}).
\end{defn}

\brem\label{re:why-ymax-etc}\ 
\benumlab
\item As will be seen in \Cref{le:point-lemma}, each local  holomorphoid given by a holomorphoid $(X,\arithl_\by)$ gives a point (i.e. an object) of the local Arithmetic Teichmuller Space $\fJ(X,L_v)$ constructed in \constrone{Definition }{def:arith-hol-space-local}, \constrone{Definition }{def:arith-hol-space-local}. But the category $\fJ(X,L_v)$ considered in \cite{joshi-teich} is bigger than the category of holomorphoids considered here. 
\item Similarly, each (global) holomorphoid $(X,\arithl_\by)$ gives a point (i.e. an object) of the adelic Teichmuller space $\widetilde{\fJ(X,L)}$ of \cite{joshi-teich,joshi-untilts}, but again the category considered there is bigger than what is being considered here. 
\item Each $\by$ is a point of the space $\yadl$ constructed in \constrtwoh{\ssep}{se:adelic-ff-curves} which refers to a specific avatar or deformation of the arithmetic of $L$.
\item The usefulness of working with holomorphoids arising from $\yadl^{max}$ and $\yadl^\R$ is given by \constrtwoh{Proposition }{pro:working-ymax}, \constrtwoh{Proposition }{pro:working-yR} and the theory of maximally complete fields \cite{kaplansky42}, \cite{poonen93}.
\eenum
\erem

\begin{defn}\label{def:geom-base-point}\
In the notation and terminology of \Cref{def:holomorphoids}, let $(X,L\into \arith{L}_\by)$, with $\by=(y_v)_{v\in\vl}\in\yadl$, be a holomorphoid of $X/L$.  Then
\benumlab
\item  for each $v\in\vl$,  the datum of the untilt $y_v=(L_v\into K_{y_v},K_{y_v}^\flat\isom \lvbht)$ provides the triple $$\left(\xan_{L_{v}},\xan_{K_{y_v}},\xan_{K_{y_v}}\to\xan_{L_{v}}\right),$$ 
called the \textit{holomorphic structure at $v$ of the holomorphoid $(X,L\into \arith{L}_\by)$} or the \textit{local holomorphic structure at $v$ of the holomorphoid  $(X,L\into \arith{L}_\by)$}, and which
	consists of the Berkovich analytic spaces $\xan_{K_{y_v}}$, $\xan_{L_v}$ and the base change morphism of analytic spaces $\xan_{K_{y_v}}\to\xan_{L_v}$. The construction of holomorphic structure
	$$(X,L\into \arith{L}_\by) \mapsto \left(\xan_{L_{v}},\xan_{K_{y_v}},\xan_{K_{y_v}}\to\xan_{L_{v}}\right)$$
	is clearly functorial in the holomorphoid $(X,L\into \arith{L}_\by)$.
	\item Let $(X,L\into \arith{L}_\by)$ be a holomorphoid of $X/L$. Then the collection of morphisms of analytic spaces  $$\left\{ \sM(K_{y_v})\to \xan_{L_v} : v\in\vl \right\}$$  will be called the \textit{geometric base-point of the holomorphoid $(X,L\into \arith{L}_\by)$} and for each  $v\in\vl$ the morphism $$\sM(K_{y_v})\to \xan_{L_v}$$ will be called the \textit{geometric base-point of $(X,L\into \arith{L}_\by)$ at $v$} or the
	 \textit{$K_{y_v}$-geometric base-point of $\xan_{L_v}$}. The construction 
	$$(X,L\into \arith{L}_\by) \mapsto \left(\xan_{K_{y_v}}\to\xan_{L_{v}}\right)$$ is again functorial in $(X,L\into \arith{L}_\by)$.
	\item Let $(X,L\into \arith{L}_\by)$ be a holomorphoid. Then the construction
	$$(X,L\into \arith{L}_\by) \mapsto \by\in \yadl$$
	is also clearly functorial in $(X,L\into \arith{L}_\by)$.
	\item The following notational conventions will be in force throughout. I will  often write $\xan_{K_v}$ instead of $\xan_{K_{y_v}}$ and often write $$X/\arith{L}_\by$$ in place of $(X,L\into \arith{L}_\by)$. 
	\item Sometimes, I will use the more suggestive notation $$\holt{X/L}{\by}$$ to denote the holomorphoid of $X/L$ given by the arithmeticoid $\arith{L}_\by$. If $\by=(y_v)_{v\in\vl}$, then I will use the suggestive notation $$\hol{}{y_{v}}{X/L}$$ to denote the corresponding local holomorphoid at $v$. \eenum
\end{defn}

\blem\label{le:point-lemma} 
Let $(X,L\into \arith{L})$ be a holomorphoid of $X/L$. Let $v\in\vl$. Then
the local holomorphoid $\hol{}{y_{v}}{X/L}$ at $v$ is a point (i.e. an object) of the local Arithmetic Teichmuller Space $\fJ(X/L_v)$ defined in \constrone{Definition }{def:arith-hol-space-local} and hence each holomorphoid $(X,L\into \arith{L})$ provides a unique object of the adelic Arithmetic Teichmuller Space $\fJ(X/L)=\prod_{v\in\vl} \fJ(X/L_v)$.
\elem
\bp 
This is immediate from \constrone{Definition }{def:arith-hol-space-local} and \cite{joshi-teich}. The explicit correspondence is given (using \Cref{def:holomorphoids}) by the rule
$$(X,\arithl_\by)\mapsto \left(X/L_v,y_v=(L_v\into K_{y_v},K_{y_v}^\flat\isom \lvbht), \sM(K_{y_v})\to \xan_{L_v}\right)_{v\in\vl}\in \fJ(X/L)$$
with $\by=(y_v)_{v\in\vl}$.
\ep

\brem For connecting with \iut, let me remark that the passage to tempered fundamental groups is facilitated by  geometric base-points at $v\in\vl$ provided by any holomorphoid $\holt{X/L}{\by}$ of $X/L$  via \Cref{def:geom-base-point}{\bf(2)}. One can also compute the geometric tempered fundamental group by picking a morphism $\sM(K_{y_v})\to\xan_{L_v}$ which factors as $\sM(K_{y_v})\to\xan_{K_{y_v}}\to\xan_{L_v}$.
\erem

\brem 
Let me point out the analogy with classical Teichmuller Theory. Let $\Sigma$ be a compact Riemann surface of genus $g\geq 1$. Let $T_{\Sigma}$ be its Teichmuller space. Let $\abs{\Sigma}$ be the topological space underlying $\Sigma$.  Then a holomorphoid of $\abs{\Sigma}$ is a point  $\Sigma'\in T_{\Sigma}$ of the Teichmuller space of $\Sigma$. Notably, a holomorphoid of $\abs{\Sigma}$ comes equipped with its own complex analytic space $\Sigma'$ with $\abs{\Sigma'}=\abs{\Sigma}$ i.e. a complex (holomorphic) structure on $\abs{\Sigma}$. As is seen from \Cref{le:point-lemma}, in this analogy, each  (global) holomorphoid (resp. local holomorphoid) defined in \Cref{def:holomorphoids} is a point of the global Arithmetic Teichmuller Space (resp. local Arithmetic Teichmuller Space) of \cite{joshi-teich,joshi-untilts}.
\erem

\brem\label{rem:diamonds-elaborated} \textcolor{red}{This is an elaboration of \constrone{Remark }{re:diamonds} and  will be useful to readers familiar with the Theory of Diamonds \cite{scholze-diamonds}}. Let $X/L$ be as in \Cref{def:holomorphoids}.   Let $\by\in\yadl$. Let $v\in\vlnon$ and let $p_v$ be the residue characteristic of $L_v$.  Let $X^\Diamond_{L_v}$ be the diamond associated by \cite[Def. 15.5]{scholze-diamonds} with the analytic adic space $X_{L_v}^{ad}$ associated to $X/L_v$ and let $\abs{X^\Diamond_{L_v}}$ denote the underlying topological space of $X^\Diamond_{L_v}$. Then by \constrone{Remark }{re:diamonds} the set of equivalence classes of local holomorphoids of $X/L_v$   i.e. the set
$$\left\{ \hol{}{y_v}{X/L} \right\}/\sim$$
(essentially this set is the local Arithmetic Teichmuller Space of \constrone{Definition }{def:arith-hol-space-local}) is formally similar and closely related with the set $X^\Diamond_{L_v}(\C_{p_v}^\flat)$ \cite[Def. 15.5]{scholze-diamonds}. The projection $\left\{\hol{}{y_v}{X/L} \right\}\to y_v$ corresponds to the projection 
$$\abs{X^\Diamond_{L_v}}\to \abs{{\rm Spd}(L_v)}.$$
The adelic Arithmetic Teichmuller Space of \cite{joshi-teich,joshi-untilts,joshi-teich-def} deals with the following sort of object:
$$\prod_{v\in\vl} \abs{X^\Diamond_{L_v}} \to \prod_{v\in\vl} \abs{{\rm Spd}(L_v)},$$
which is roughly corresponds to the projection $$\left\{\holt{X/L}{\by} :\by\in\yadl \right\}\to \yadl$$ given by $\holt{X/L}{\by}\mapsto\by\in\yadl$. At the moment there is no good global geometric theory of objects of the form $$\prod_{v\in\vl} X^\Diamond_{L_v}.$$ In the present series of papers one considers this type of object subject to the global symmetries considered in \cite{joshi-teich-def} acting on the base.  This complicated object is the type of dynamical arena of Arithmetic Teichmuller Spaces considered in the present series of papers.
\erem

\subparat{Category of holomorphoids} If $\by,\by'\in \yadl$ and if $(X,L\into \arith{L}_{\by})$ and $(X',L\into \arith{L}_{\by'})$ are holomophoids of schemes $X,X'$ over $L$, then one can talk about morphisms between holomorphoids $(X',L\into \arith{L}_{\by})\to (X,L\into \arith{L}_{\by'})$ in the obvious way. So one has a natural category of holomorphoids of schemes over  arithmeticoids of $L$. This category comes equipped with the projection $(X,L\into \arith{L}_{\by})\mapsto\by$.

\subparat{Holomorphoids and the existence of distinct scheme theories}\label{ss:scheme-theories} [\textit{The purpose of this section is to establish Mochizuki's assertion \cite[\ssep I3]{mochizuki-iut1} about existence of distinct scheme theories.}]  Fix an arithmeticoid $\arithl_\by$. Note that, given a holomorphoid $(X,\arithl_\by)$, one can consider the pair $$(X,\arithl_\by)^{sch}$$ obtained by forgetting the geometric base-point datum (given by \Cref{def:geom-base-point}) of the holomorphoid $(X,\arithl_\by)$. The meaning of this notation $(X,\arithl_\by)^{sch}$ is as follows:   $(X,\arithl_\by)^{sch}$ is the scheme $X$ viewed as a scheme over $L$ with the arithmetic of $L$ specified by $\arithl_\by$ (in the sense of \cite{joshi-teich-def}).  I will refer to $(X,\arithl_\by)^{sch}$ as a scheme over the arithmeticoid $\arithl_\by$. 

\textit{One way to think about $(X,\arithl_\by)^{sch}$ is that it is a scheme over $L$ equipped with its private version of the arithmetic of $L$ (and $\bL,\{ L_v,\bL_v\}_{v\in\vl}$) using which many familiar quantities  and structures  associated to $X$ (such as Galois groups, fundamental groups) may be computed. (The holomorphoid $(X,L\into\arithl_\by)$, which gives to $(X,\arithl_\by)^{sch}$,  also comes equipped with the same  private copy of arithmetic.} 

For a fixed $\arithl_\by$, the category, denoted $\mathcal{Sch}/\arith{L}_{\by}$,  with objects $(X,\arithl_\by)^{sch}$ (and natural morphisms) is equivalent to the category of schemes over $L$ (by  forgetting $\arithl_\by$).  I will refer to $\mathcal{Sch}/\arith{L}_{\by}$ as category of schemes over the arithmeticoid $\arith{L}_\by$.

Note that in general, for $\by\neq\by'$ in $\yadl$,    there may be no canonical or natural equivalence between the categories $\mathcal{Sch}/\arith{L}_{\by}$ and $\mathcal{Sch}/\arith{L}_{\by'}$ which respects the structures presented by the arithmeticoids $\arith{L}_{\by}$ and $\arith{L}_{\by'}$.   
More precisely:
 
By \constrtwoh{Theorem }{th:existence-inequivalent-arithmeticoids},  two arithmeticoids $\arith{L}_{\by}$ and $\arith{L}_{\by'}$ may not even be (topologically) equivalent,  one may think of $\left\{\mathcal{Sch}/\arith{L}_{\by}\right\}_{\by\in\yadl}$ as providing distinct or  (topologically) inequivalent models of the category $\mathcal{Sch}/L$ of schemes   over $L$. [This type of assertion occurs without proof in \cite[\ssep I2]{mochizuki-iut1}. Notably, \iutthr\ frequently refers to distinct miniature models of scheme theory, but omits a precise quantification of these notions.] If one forgets the information of the relevant arithmeticoids, then these categories are, of course, all equivalent.  

Let me summarize this discussion in the following way:
\bthm\label{th:exist-cats-of-schemes-over-L} 
Each member of $$\left\{\mathcal{Sch}/\arith{L}_{\by}\right\}_{\by\in\yadl}$$
is a category equivalent, via the forgetful functor, to the category of $\mathcal{Sch}/L$. But given two distinct arithmeticoids $\arith{L}_{\by_1}$ and $\arith{L}_{\by_2}$, there may be no equivalence between $\mathcal{Sch}/\arith{L}_{\by_1}$ and $\mathcal{Sch}/\arith{L}_{\by_2}$ which  is compatible with the two arithmeticoids. 
\ethm

\brem\label{re:cat-sch-holmorphoid}
In the context of \Cref{th:exist-cats-of-schemes-over-L}, one would in fact like to consider a category whose class of objects is
$$\boldsymbol{\mathcal{Sch}}_{L}=\left\{(X,\arithl_\by)^{sch}: (X,\arithl_\by)^{sch}\in \mathcal{Sch}/\arith{L}_{\by}\text{ and } {\by\in\yadl}\right\}$$ 
with morphisms defined in the obvious way and view this as a fibered category 
\be\label{eq:proj-to-arithmeticoid} \begin{tikzcd}
\boldsymbol{\mathcal{Sch}}_L\ar[rr, "{(X,\arithl_\by)^{sch}\mapsto\by}"] &&\yadl
\end{tikzcd}.
\ee 
Thus the  fiber over $\by$ is $\mathcal{Sch}/\arith{L}_{\by}$. Symmetries  acting on $\yadl$, detailed in  \cite{joshi-teich-def},  move the fibers around.  If one wants to work with schemes (as \iut\ does) then one must work with the fibered category $\boldsymbol{\mathcal{Sch}}_L\to\yadl$ together with symmetries of the base category $\yadl$ acting on objects of $\boldsymbol{\mathcal{Sch}}_L$. Presently, I do not know of any general formulation of the theory of fibered categories equipped with symmetry actions given by their respective base categories.
\erem

\brem 
Let me remark that for $\dim(X)=1$ and the genus of $X$ is $g=1$ (this is the case of interest in \iut), then in this case, the notion of holomorphoids is a more precise version of Mochizuki's (decorated) Hodge-Theaters. More precisely, given a  holomorphoid $\holt{X/L}{\by}$, one may construct from it a  decorated Hodge-Theater (for each decoration for Hodge-Theaters used in \iut) (for more on Hodge-Theaters see \cref{se:frobenioids}). In the context of Mochizuki's Hodge-Theaters (with a fixed decoration), the role of an arithmeticoid $\arithl_\by$ is played by a global Frobenioid of $L$ and the role of the functor considered in \cref{eq:proj-to-arithmeticoid} is given by the association $$\left\{ H : H \text{ a fixed decoration Hodge Theater}\right\} \to \text{the global Frobenioid of } L \text{ given by }H.$$ For details see \cref{ss:hodge-theater}.
\erem

\subparat{Holomorphoids of elliptic curves}\label{ss:elliptic-curve-assumptions} From now on, I will assume that $C$ is an elliptic curve over an arithmeticoid $\arith{L}$ of $L$ and write $X=C-\{O\}$ where $O\in C(L)$ is the origin of the group law. Observe that $X/\arith{L}$ is an hyperbolic curve of topological type $(g,n)=(1,1)$. I will refer to $C/\arith{L}$ (resp. $X/\arith{L}$) a holomorphoid of the elliptic curve $C/L$ (resp. a holomorphoid of the canonically punctured elliptic curve $X/L$).

To apply the results of my theory (of \cite{joshi-teich,joshi-untilts,joshi-teich-def}) to the context of \iut, the following will be useful:
\bthm\label{th:holomorphoids-elliptic-case} 
Let   $\holt{C/L}{\by}$ and $\holt{X/L}{\by}$ be (pointed) holomorphoids of $C/L$ and $X/L$ respectively.  Then one has (using the geometric base-points of these holomorphoids to compute absolute Galois groups and various fundamental groups): 
\benumlab
\item $X/L$ is a hyperbolic curve of strict Belyi Type.
\item One has for each $v\in\vl$ 
\benum
\item a surjection of topological groups
$$\pit{X/L_v;K_v} \to G_{L_v,K_v}.$$
\item Let $\bL_v\subset K_v$ be the algebraic closure of $L_v\subset K_v$. Then one also has the monoid  of non-zero elements of $\bL_v$ equipped with natural action of Galois:
$$\pit{X/L_v,K_v}\to  G_{L_v,K_v}\act \bL_v^*,$$ and,
\item the monoid of units of $\O_{\bL_v}^*$ of $\O_{\bL_v}$ equipped with its natural action of Galois:  
$$\pit{X/L_v,K_v}\to  G_{L_v,K_v}\act \O_{\bL_v}^*,$$ 
and,
\item the monoid  $\O_{\bL_v}^\triangleright=\O_{\bL_v}-\{0\}$ of non-zero elements of $\O_{\bL_v}$ equipped with natural action of Galois:
$$\pit{X/L_v,K_v}\to  G_{L_v,K_v}\act \O_{\bL_v}^\triangleright.$$
\eenum
\item If $X/\arith{L}_{\by_1}$ and $X/\arith{L}_{\by_2}$ are two distinct holomorphoids of $X/L$ then the objects provided by {\bf(1)--(2)} corresponding to $X/\arith{L}_{\by_1}$ and $X/\arith{L}_{\by_2}$  are all abstractly isomorphic but are equipped with distinct arithmetic labels and their arithmetic and geometric properties are quantifiably distinct from each other.
\eenum
\ethm
\bp 
This is an elementary consequence of the results of \cite{joshi-teich,joshi-untilts,joshi-teich-def}.
\ep

\brem Following remarks will be useful:
\benumlab
\item  As is established in \cref{se:intro-rosetta-stone}, \cref{se:frobenioids}, all other objects of \iutthr\ such as Frobenioids, Hodge-Theaters (of all types), prime-strips etc. can be obtained from every holomorphoid $X/\arith{L}_\by$ considered here.
\item
In \cite{mochizuki-essential-logic}, Mochizuki has argued that the logic of his theory (\iutthr) requires treating two isomorphs as being distinct i.e. the logic of the theory is to work with Isomorph A {\bf and} Isomorph B as opposed to Isomorph A {\bf or} Isomorph B. The assertion \Cref{th:holomorphoids-elliptic-case}{\bf(3)} (proved in the present paper) provides a clear mathematical reason for doing so, namely, these objects (arising from $X/\arith{L}_{\by_1}$ and $X/\arith{L}_{\by_2}$) while being (abstractly) isomorphic, they can be treated as distinct because they do indeed have distinct arithmetic and geometric properties. 
\item To put it clearly, in my approach, one works with $X/\arith{L}_{\by_1}$ {\bf and} $X/\arith{L}_{\by_2}$ for natural reasons (they are quantifiably distinct, but isomorphic, objects in the theory). No realignment of logic is needed to validate my theory. Moreover, my work provides distinct objects even for \iutthr. Readers should recognize that my work  provides the most natural and mathematically precise formulation of the theory created in \iutthr\ and is fully compatible with the mathematical portions of \cite{mochizuki-essential-logic}.
\eenum
\erem

\newcommand{\uBB}{\u{\mathbb{B}}}
\newcommand{\uBBbreve}{\u{\breve{\mathbb{B}}}}

\section{Fixing the  Initial Theta Data \`a la Mochizuki}\label{se:theta-data}
Let me remark that the main theorems of this paper will require the existence of global arithmetic data called Initial Theta Data (as in the title of this section). The existence of the Initial Theta Data is established in \constrfour{Theorem }{th:existence} (and in \cite{mochizuki-iut4}), so all  the principal results of this paper will require the existence of the Initial Theta Data  whose definition and properties are detailed in this section. \textit{Importantly, the existence of the Initial Theta Data (\constrfour{Theorem }{th:existence}) is a global assertion i.e. requiring one to work over a number field and it cannot be established by local means.} 

From now on, the following assumptions will be in force.  These assumptions essentially corresponds to fixing the \textit{Initial Theta Data} in \cite{mochizuki-iut1}:
\subparat{Basic Hypotheses}\label{ss:theta-data-fixing}
\benumlab
\item Let $L/\Q$ be a number field such that $L$ has no real embeddings.  
\item Let $X/L$ be an hyperbolic curve of Strict Belyi Type  obtained from a geometrically connected, elliptic curve $C/L$ by setting $X=(C-\{O\})/L$, where $O\in C(L)$ is the origin of the group law of $C$. 
\item Fix an algebraic closure $\bL$ of $L$ and write $G_L=\gal(\bL/L)$.  
\item Let $\lmod$ be the field of moduli of $X/L$. 
\eenum
If one needs to emphasize the field one is working with, then I will write $X_{L}, X_{\lmod}$, $X_{L'}$ etc. Frequently one will need to work with all the three simultaneously.

\subparat{Notational conventions} I will also fix some notation regarding valuations of $L$ and $\lmod$:
\benumlabresume
\item Let $\bbvl$ (resp. $\bbvlmod$) be the sets of equivalence classes of non-archimedean valuations of $L$ (resp. $\lmod$). 
\item Let $\bbvlmodoss$ be the set of primes of $\lmod$ with odd residue characteristics at which $X$ has bad, semi-stable reduction i.e. split multiplicative reduction.
\item Let $\bbvlmodgood=\bbvlmod-\bbvlmodoss$.
\item One can also define similar subsets of valuations of $L$. 
\eenum

The following lemma is immediate from the definitions above:
\blem\label{le:reduction-good-bad}
With the above assumptions on $X$ and $L$, and notational conventions one has
\benumlab
\item  $w\in\bbvlmodgood$ if and only if $w|2$ or $X$ has good reduction at $w$ or  has bad additive  reduction at $w$. 
\item For any $w|2$ lying in $\bbvlmodgood$, $X$ can have good or bad (additive or multiplicative) reduction at $w$.
\eenum
\elem

\subparat{Initial Theta-data}\label{ss:theta-data-fixing2} Now to the assumptions of \ssep\ref{ss:theta-data-fixing}, let me add the following. Assume in addition to the assumptions of \ref{ss:theta-data-fixing} that
\benumlabstart{9}
\item $L/\lmod$ is a Galois extension,
\item The field $L$ contains all points of order $2\cdot 3=6$ of $C/L$ i.e. $L(C[6](\bL))=L$.
\item Let $\ell\geq 5$ be a prime number, write $\ells=\frac{\ell-1}{2}$, and assume that  $\ell$ does not divide the degree  $[L:\lmod]$ of $L/\lmod$.
\item Suppose that the image of the mod $\ell$ representation $$\rho_{C/L;\ell}:G_L\to \Aut_{\Z/\ell}(C[\Z/\ell])=GL_2(\Z/\ell)$$  contains $SL_2(\Z/\ell)$.  
\item Let $L'$ be the fixed field of the kernel of the homomorphism $\rho_{C/L;\ell}:G_L \to GL_2(\Z/\ell)$.
\item Fix a bijection $\bbvlmod\isom \ubblv\subset \bbvlp$.
\eenum
\brem Let $v|p$ be any prime of $L$ lying over some prime number $p$ and let $L_v$ be the completion of $L$ at $v$.  By \constrone{Theorem }{th:field-of-def},  for every prime $p$ and for every prime $v\in \bbvl$ with  $v|p$, every local holomorphoid of $X/L$ at $v$ has the field of definition $\lmod$ (notably, the field of definition is independent of $p$ and $v|p$). [Similarly for any $w|p$ in $L'$ one has the category $\fJ(X/L_w')_{\cpt}$ of \cite{joshi-teich,joshi-untilts} whose every object has $\lmod$ as its field of definition.]
\erem

\subparat{More notation}\label{ss:val-sets-defn}
\begin{align*}
\ubblv^{non}& =\ubblv\cap\bbvlp^{non}, \text{ and},\\
\ubblv^{arc}&=\ubblv\cap\bbvlp^{arc};\\
\ubblvgood&=\ubblv\cap \bbvlpgood, \text{ and}, \\ 
\ubblvoss&=\ubblv\cap \bbvlposs.
\end{align*}
Similarly (for $L$)
\begin{align*}
\ubblv^{non}_L& =\ubblv\cap\bbvl^{non}, \text{ and},\\
\ubblv^{arc}_L&=\ubblv\cap\bbvl^{arc};\\
\ubblvgood_L&=\ubblv\cap \bbvlgood, \text{ and}, \\ 
\ubblvoss_L&=\ubblv\cap \bbvloss.
\end{align*}
For a rational prime $p$, let 
\begin{align*}
\ubblv_p&=\{ w\in \ubblv: w|p\},\\
\ubblvgood_p&=\{ w\in \ubblvgood: w|p\},\\
\ubblvossp&=\{ w\in \ubblvoss: w|p\}.
\end{align*}

\newcommand{\inithtdata}{\cref{ss:theta-data-fixing}, \cref{ss:theta-data-fixing2}}
\newcommand{\assumptions}{\cref{ss:elliptic-curve-assumptions}, \inithtdata}

\subsubparat{Existence of local torsion points of order $2\ell$} An important consequence of \cref{ss:theta-data-fixing} and \cref{ss:theta-data-fixing2} is the following. Let $w\in \ubblvoss$ be a prime of $L'$ lying over a prime $v\in\bbvloss$ of $L$. Then $C/L'_{w}$  has a non-trivial point of exact order $2\cdot\ell$ defined over $L'_{w}$ i.e. $$C[2\cdot \ell](L'_{w})\neq0.$$
In particular, if $w|v$ for a prime $v|p$ of $L$, then  the Tate parameter of $C/L_v$ has an $2\ell^{th}$-root in $L'_{w}$. \textcolor{red}{The construction of the theta-values locii in \cref{se:construction-of-thetaj-and-thetam}, \cref{se:mochizuki-construction-thetam} will require the existence of a point of exact order $2\ell$.}

\subsubparat{Elementary bound}
For use later on, when one wants to vary $C/L$, let me record the following elementary consequence of the assumptions \assumptions.
\blem\label{le:bounded-degree-theta-data} 
Let $L$ be a number field, $\bL$ be a fixed algebraic closure of $L$. For any $0<A\in\R$, let $$L^{\leq A}\subset \bL$$ be the smallest subfield containing all finite extensions $L'/L$ of degrees $$[L':L]\leq A.$$ Let $\ell$ be a fixed odd prime number. Suppose that $C/L$ is an elliptic curve  equipped with data \cref{ss:theta-data-fixing}. Then as $C/L$ varies over $\sM_{1,1}(L)$ the data \cref{ss:theta-data-fixing} for $C/L$ provide a family of finite extensions $L'/L$ (as above) (i.e.  $L'$ is the fixed field of $\ker(\rho_{C/L;\ell})$)  which are contained in $L^{\leq A}$ with $A=|GL_2(\Z/\ell)|\ll \ell^4$.
\elem
\bp 
Note that $L^{\leq A}/L$ is an extension of infinite degree. The proof is clear. The degree of the fixed field of $\ker(\rho_{C/L;\ell})$ is bounded by the order of the Galois group of $L'/L$ and by construction of $L'/L$, one has  $\gal(L'/L)\subset GL_2(\Z/\ell)$ and as $\ell$ is fixed, the latter group has bounded order $A=|GL_2(\Z/\ell)|\ll \ell^4$.
\ep

\subparat{Existence of Initial Theta Data}
Let me say this clearly: \textcolor{red}{the existence of Initial Theta Data \constrfour{Theorem }{th:existence} is a global assertion i.e. requiring one to work over a number field and with compactly bounded subsets of algebraic points and this existence cannot be proved by $p$-adic local means (even for all primes). This existence theorem is necessary for the central constructions and results of the present paper (and \cite{mochizuki-iut3}) and hence also for \cite{joshi-teich-abc} (and \cite{mochizuki-iut4}).}
The existence of Initial Theta Data  satisfying \assumptions, \inithtdata\ and some additional properties is a global arithmetic assertion established in \constrfour{Theorem }{th:existence} (resp. in \cite[Corollary 2.2]{mochizuki-iut4}). For ease of comparison between my papers and Mochizuki's, I have preserved Mochizuki's content wise order of appearance of results: this means that results which appear in \cite{mochizuki-iut3} appear in this paper and results of \cite{mochizuki-iut4} appear in \cite{joshi-teich-abc}. This does not necessarily reflect the logical dependence of results in \cite{mochizuki-iut3} on the existence of Initial Theta Data which is proved in \cite[Corollary 2.2]{mochizuki-iut4} (note that the key Diophantine estimate \cite[Theorem 1.10]{mochizuki-iut4} also depends, among other critical assumptions, on the existence of Initial Theta Data). In what follows, I will assume the existence of Initial Theta Data satisfying \assumptions, \inithtdata\ for the rest of this paper. 

\section{Construction of Mochizuki's Adelic Ansatz and Adelic Theta-links}\label{se:adelic-ansatz}
In this section I want to provide the construction of the adelic version of certain correspondences (introduced in \cite{joshi-teich-estimates,joshi-teich-def}), which in the terminology of \iutthr, are called $\Theta_{gau}$-links. The construction provided here is the adelic version of the fundamental constructions of \cite{joshi-teich-estimates} where I gave the construction of these objects for one fixed prime of a number field.

\subparat{The variant $\yadl'$ of $\yadl$} Let $L$ be a number field, which will be assumed to have no real embeddings. Let $C/L$ and $X/L$ be chosen as above and fix an  Initial Theta Data \assumptions. Let $L'\supset L$ be the number field provided by the  Initial Theta Data. Since $L$ has no real embeddings, neither does $L'$. Let $\yadlp$ be the adelic curve defined in \constrtwoh{\ssep}{se:adelic-ff-curves} using the field $L'$. In fact it is more convenient to work with the variant $\yadlp'$, defined in \constrtwoh{Remark }{rem:yadl-variant}, of $\yadl$ (the latter is defined in \constrtwoh{Definition }{def:adelic-ff-curves}). 

This means that $\yadlp'$ is the topological space
\be \yadlp'=\prod_{p} \prod_{w |p} \abs{\sY_{\cpt,L_w'}},\ee
the product over $p$ runs over all rational primes including $p=\infty$. This is possible because, for every rational prime $p$, one has an inclusion of valued fields $\widehat{\bar{L}}_w' \subset \C_p$ for every $w|p$.

I will define the adelic version of Mochizuki's Ansatz (see \cite{joshi-teich-estimates} for the local version). Mochizuki's Adelic Ansatz $\tSigma_{L'}$ is a subset 
\be\label{eq:adelic-ansatz}
\tSigma_{L'}\subset \left(\yadlp'\right)^\ells=\overbrace{\yadlp'\times\yadlp'\times\cdots\times\yadlp'}^{\ells\text{ factors}}.
\ee constructed using the  Initial Theta Data \inithtdata.

\subparat{Mochizuki's Adelic Ansatz using $\yadlp'$}\label{ss:mochizuki-adelic-ansatz}
\newcommand{\sGal}[1]{\mathbf{G}_{#1}}
For $w\in\vlp$ and $y_w\in\syflwp$, let $K_{y_w}$ be the residue field of $y_w$ in $\syflwp$ and let $p_w$ be the residue characteristic of $K_{y_w}$ (this means for $w\in\vlnonp$,  $K_{y_w}^\flat$ has characteristic $p_w>0$, and for $w\in\vlarchp$, $p_w=\infty$).

Recall from \constrtwoh{Theorem }{th:galois-action-on-adelic-ff} that $\yadlp$ is equipped with the following three actions: $\sGal{L'}\act\yadlp$, ${L'}^*\act\yadlp$, and $\bvarphi\act\yadlp$. This provides the following actions on $\yadlp^\ells$.
\bdefn\label{def:actions-on-tuples}Let $L'\supset L$ be number fields. Then one has the following actions on $\yadlp^\ells$.
\benumlab
\item The diagonal action of the group $\sGal{L'}\act\yadlp^\ells$: given by
$$(g=(g_w)_{w\in\vlp},\bz=(\by_1,\by_2,\ldots,\by_\ells))\in\sGal{L'}\times\yadlp'^\ells\mapsto g\cdot\bz=(g\cdot\by_1,g\cdot\by_2,\ldots,g\cdot\by_\ells)).$$ 
\item The diagonal action $\bvarphi^\Z\act\yadlp^\ells$ of the group $\bvarphi^\Z$ generated by the global Frobenius given by
$$\bvarphi(\bz)=(\bvarphi(\by_1),\bvarphi(\by_2),\ldots,\bvarphi(\by_\ells)).$$
\item The action ${L'}^*\act\yadlp^\ells$ of ${L'}^*$ and hence of ${L'}^*\supset L^*$ given by
$$(x,\bz=(\by_1,\by_2,\ldots,\by_\ells))\mapsto x\cdot\bz=(x\cdot\by_1,x^{2^2}\cdot \by_2,\ldots,x^{\ells^2}\cdot\by_\ells)$$ in which the action of ${L'}^*$  on the $j^{th}$ factor of $\yadlp^{\ells}$ is by $x^{j^2}$ acting on this factor.
\eenum
\edefn

Let $\bz=(\by_1,\by_2,\ldots,\by_\ells)\in(\yadlp')^\ells$ be an arbitrary point $(\yadlp')^\ells$. Then $\bz$ can be thought of as $$\bz=(\bz_w)_{w\in\vlp},$$
i.e. each $\bz$ may be thought of as a collection of $\ells$-tuples $\bz_w\in\syflwwp^\ells$ for each $w\in\vlp$.

The set $\tSigma_{L'}\subset  (\yadlp')^\ells$ is given as follows. 
\bdefn\label{def:mochizuki-adelic-ansatz} 
Mochizuki's Adelic Ansatz $\tSigma_{L'}$ is the subset of $(\yadlp)^\ells$ defined as follows: $\bz\in\tSigma_{L'}$ if and only if the following condition holds true: 
\be\label{eq:adelic-ansatz2}\bz_w=\begin{cases} (y_w,y_w,y_w,\ldots,y_w)\in \syflwp^\ells & \text{ if } w \in \vlp-\ubblvoss \text{ and} \\
						(y_{w,1},y_{w,2},\ldots,y_{w,\ells}) \in \tSigma_{\cpt,L_w'} & \text{ if } w\in \ubblvoss.
\end{cases}
\ee
where, for $w\in\bbvlposs$, $$\tSigma_{\C_{p_w}^\flat,L_w'}\subset \syflwwp^\ells$$ is a variant of the local version of Mochizuki's Ansatz (constructed in \cite{joshi-teich-estimates}) and the construction of $\tSigma_{\C_{p_w}^\flat,L_w'}$ is detailed in  the next paragraph \cref{ss:ansatz-local-def}.
\edefn
\subsubparat{Mochizuki's Adelic Ansatz: local aspects \`a la \cite{joshi-teich-estimates}}\label{ss:ansatz-local-def} Fix a prime $w\in\bbvlposs$ write $p=p_w$ for its residue characteristic and write  $E_w'=L'_w$ (for notational continuity with \cite{joshi-teich-estimates}). I will construct \textit{Mochizuki's Ansatz} \be\tSigma_{\cpt,L_w'}\subset \syflwp^\ells.\ee The definition of this set given below is a variant of the definition made in \cite{joshi-teich-estimates} which uses $L_{w}'$ in a crucial way.   

Consider the continuous mapping (with finite fibers) (by \cite[Proposition 2.3.20]{fargues-fontaine}) 
\be\label{eq:y-proj} \abs{\sY_{\cpt,E_w'}}\to \abs{\sY_{\cpt,\Q_p}}\ee 
\be\sY_{\cpt,E_w'}\ni(E_w'\into K,K^\flat\isom\cpt) \mapsto (K,K^\flat\isom \cpt)\in\syQp\ee
given by forgetting  $E_w'$ in the isometric embeddings $\Q_p\subset E_w'\into K$. 

Let 
\be\tilde{\Sigma}_{\cpt,E_w'}\subset \sY_{\cpt,E_w'}^\ells\ee be the inverse image, under the above morphism \eqref{eq:y-proj}, of Mochizuki's primitive ansatz \be \tilde{\Sigma}_{\cpt}\subset \syQp^\ells\ee which is defined in \cite[\ssep 6]{joshi-teich-estimates}. This completes the construction of Mochizuki's Adelic Ansatz $\tSigma_{L'}$ \eqref{eq:adelic-ansatz2}.

\bdefn\label{dfn:primitive-ansatz}
I will call the set $\tilde{\Sigma}_{\cpt,E_w'}=\tilde{\Sigma}_{\C_{p_w},E_w'}$ \textit{Mochizuki's local Ansatz in $\sY_{\C_{p_w},E_w'}^\ells=\sY_{\C_{p_w}^\flat,L_w'}^\ells$}. 
\edefn

\subsubparat{Fundamental Properties of Mochizuki's Adelic Ansatz}\label{sss:fun-prop-adelic-ansatz} This paragraph details the fundamental properties of $\tSigma_{L'}$ which will be used throughout this paper. 
\bthm\label{th:adelic-theta-link} 
Let $L$ be a number field, let $C/L$ and $X/L$ be as above (\cref{ss:elliptic-curve-assumptions}) and an  Initial Theta Data \inithtdata\ be fixed. Then the set $\tSigma_{L'}$ constructed above \eqref{eq:adelic-ansatz2} has the following properties. Let $\sGal{L'}$ be the group defined in \constrtwoh{\ssep}{ss:GalL-defined} (for $L'$). Let $\tSigma_{L'}$ be the set constructed above. 
\benumlab
\item  The set $\tSigma_{L'}$ is stable under the natural action of $\sGal{L'}$ on $\left(\yadlp'\right)^\ells$  which is given using  the natural action of $\sGal{L'}$  on $\yadlp'$ (constructed in \constrtwoh{Theorem }{th:galois-action-on-adelic-ff}).
\item The set $\tSigma_{L'}$ is stable under the natural action of the (global) Frobenius $\boldsymbol{\varphi}$ on $\left(\yadlp'\right)^\ells$ given using the natural action of the (global) Frobenius $\boldsymbol{\varphi}$ on $\yadlp'$. 
\item Hence $\tSigma_{L'}$ is also stable under the natural action of $(L')^*$ on $\yadlp$ constructed in \constrtwoh{Theorem }{th:galois-action-on-adelic-ff}.
\item Let $(\by'_1,\by_2',\ldots,\by_\ells')\in \tSigma_{L'}$, write, for $j=1,2,\ldots,\ells$, $\by_j'=(y'_{j,w})_{w\in\vlp}$.  For each $w\in\ubblvoss$. Let $$(y'_{w,1},y'_{w,2},\ldots,y'_{w,\ells})\in \sY_{\cpt,L_w'}^\ells$$ be the $w$-component of $(\by_1',\by_2',\ldots,\by_\ells')$. Let $K_{w,1}',K_{w,2}',\ldots,K_{w,\ells}'$ be the residue fields of $y'_{w,1},y'_{w,2},\ldots,y'_{w,\ells}$  in $\sY_{\cpt,L_w'}$ respectively. Then one has the following relationship between their valuations
\be\label{eq:valuation-scaling-relationship} 
\abs{-}_{K_{y_{w,j}'}}=\abs{-}_{K_{y_{w,1}'}}^{j^2} \text{ for }j=1,2,\ldots,\ells.
\ee
\item Moreover for each $w\in\vlp-\ubblvoss$ one has the following relationship between their valuations:
\be 
\abs{-}_{K_{y_{w,j}'}}=\abs{-}_{K_{y_{w,1}'}} \text{ for }j=1,2,\ldots,\ells.
\ee
\eenum
\ethm

\bp 
The first assertion is a consequence of \constrtwo{Proposition }{pr:prim-ansatz-galois-stable} and \constrtwo{Theorem }{pr:lift-vals}. The second assertion {\bf(2)} is immediate from \constrtwo{Proposition }{pr:frob-stable-cannonical} and \constrtwoh{Theorem }{th:galois-action-on-adelic-ff}.  

Let me give a proof of {\bf(3)}. This can be proved along the lines of \constrtwo{Theorem }{pr:lift-vals}, but let me deduce the present case as a consequence of \constrtwo{Theorem }{pr:lift-vals}. To see this, let $w|p$. Then one has a  quasi-finite morphism  \eqref{eq:y-proj} $\sY_{\cpt,L_w'}\to\syQp$. Let $y_{p,j}\in\syQp$ be the image of $y_{w,j}'$ under this morphism and let  $(y_{p,1},y_{p,2},\ldots,y_{p,\ells})\in\tSigma_{\cpt}\subset \syQp^\ells$ be the corresponding point of $\syQp^\ells$. 

In \constrtwo{Theorem }{pr:lift-vals}, I have shown that \eqref{eq:valuation-scaling-relationship} holds for the residue fields $K_{y_{p,1}}, K_{y_{p,2}}, \ldots,K_{y_{p,\ells}}$ of $y_{p,1},y_{p,2},\ldots,y_{p,\ells}$. By  \cite[Proposition 2.3.20(3)]{fargues-fontaine} one has an isomorphism of local rings of $y_{p,j}$ and $y'_{w,j}$ for each $j$. Hence \eqref{eq:valuation-scaling-relationship} also holds for the residue fields $K_{y_{w,1}'}, K_{y_{w,2}'}, \ldots,K_{y_{w,\ells'}}$ of $y'_{w,1},y'_{w,2},\ldots,y'_{w,\ells}$. 

The last assertion is immediate from the definition of $\tSigma_{L'}$.

This completes the proof.
\ep

The scaling relationship established in \Cref{th:adelic-theta-link} (at the primes of semi-stable reduction) has fundamental consequences for global arithmetic. Mochizuki asserts the content of the next corollary as \cite[Theorem A(ii)]{mochizuki-iut1}:

\bcor\label{cor:val-norm-theta} Suppose that the notations and assumptions of \Cref{th:adelic-theta-link} remain in force. Then  the  non-trivial scaling relationship established by \Cref{th:adelic-theta-link}{\bf(3)} at primes $w\in\ubblvoss$, forces that suitable (re)normalization of valuations must be introduced at primes $w\in\vlp-\ubblvoss$ so that the product formula to continue to hold for each normalized arithmeticoid $\arithl_{\by'_j}^{nor}$ given by $(\by'_1,\by_2',\ldots,\by_\ells')$.
\ecor
\bp 
This is clear from \Cref{th:adelic-theta-link} and the global constraint placed by  the requirement that the product formula holds.
\ep

\subsubparat{Mochizuki's Ansatz is the set of $\Theta_{gau}$-Links in \iut}\label{sss:ansatz-theta-link}
\bdefn\label{def:theta-link}
Let $(\by'_1,\by_2',\ldots,\by_\ells')\in \tSigma_{L'}$. Then the assignment $$\by_1'\mapsto (\by'_1,\by_2',\ldots,\by_\ells')\in \tSigma_{L'}$$ will be called a \textit{$\Theta_{gau}$-Link}. Sometimes I will write this as 
$$\tSigma_{L'}\ni(\by'_1,\by_2',\ldots,\by_\ells')\mapsto\by_{\ells}'\in\yadlp'.$$ Thus $\tSigma_{L'}$ is the set of $\Theta_{gau}$-links (more precisely $\Theta_{gau}$-Links in $\yadlp'$ or $\yadlp'$-$\Theta_{gau}$-Links). One may also consider $\Theta_{gau}$-Links using variants $\yadlp^{'max}$ and $\yadlp^{'\R}$ of $\yadlp'$ considered in \cite{joshi-teich-def}.
\edefn
\brem\  
\benumlab 
\item Thus a $\Theta_{gau}$-Link is simply a point of Mochizuki's Adelic Ansatz $\tSigma_{L'}$ (\cref{def:mochizuki-adelic-ansatz}).
\item Mochizuki's version of $\Theta_{gau}$-Link is introduced in \cite{mochizuki-iut3} and Mochizuki asserts in \cite[Proposition 4.1(iv), Page 344]{mochizuki-iut3} that it has the valuation scaling property established here as \Cref{th:adelic-theta-link}{\bf(3)}. This property is central to \moccor\ and \cite{mochizuki-iut4}. [Mochizuki's $\Theta_{gau}$-Link was first established rigorously in \cite{joshi-teich-estimates} and  the global case is \Cref{th:adelic-theta-link}{\bf(3)}.]
\item By \Cref{pr:theta-div-correspondence}, one may think of Mochizuki's Adelic Ansatz as a divisorial correspondence  on the adelic Fargues-Fontaine curve $\yadlp'$ (see \constrtwoh{\ssep}{ss:correspondences}). 
\eenum
\erem

The following is useful in the context of \iut:
\bcor 
Mochizuki's Adelic Ansatz $\tSigma_{L'}$ is a metrisable space with the topology and the metric induced from $\yadlp^\ells$. 
\ecor
\bp 
This is immediate from \constrtwoh{Theorem }{th:distance-bet-arith} as $\yadlp$ is a metrisable topological space.
\ep

\brem 
In particular, for a choice of a metric on $\yadlp$,  it makes sense to talk about the distance between two distinct $\Theta_{gau}$-Links.
\erem

\subparat{Mochizuki's Ansatz arises from a correspondence on $\yadlp'$}
The following proposition is based on the notion of correspondences on $\yadl$ considered in  \constrtwoh{\ssep}{ss:correspondences}. The assertion is the following:
\bpro\label{pr:theta-div-correspondence} 
Each point $$\bz\in\tSigma_{L'}=(\by'_1,\by_2',\ldots,\by_\ells')\in \tSigma_{L'}$$ defines a divisorial correspondence 
$$\by_1'\mapsto \sum_{j=1}^\ells(\by_j')=(\by'_1)+(\by_2')+\cdots+(\by_\ells')$$
on $\yadlp$  of degree $\ells$
\epro
\bp 
This is clear from the construction of $\tSigma_{L'}$ given above and \constrtwoh{\ssep}{ss:correspondences}.
\ep

\subparat{The choice of a standard point of $\yadl'$}\label{ss:stand-point} I will choose a point $\by_0\in\yadlp$ (resp. $\by_0'\in \yadlp'$), which will be referred to as the \textit{standard point}, and whose construction is detailed, for the field $L$, in \constrtwoh{\ssep}{pa:stand-arithmeticoid}  but here I work with $\yadlp$ and $\yadlp'$ respectively. To understand what this choice entails let me set out some notation. 

Let $w\in\vlp$ be a prime of $L'$ with $w|p$. Let $L_w'$ be the completion of $L'$ at $w$. The choice of $\by_0=(y_w)_{\in\vlp}$ (resp. $\by_0'=(y_w')_{\in\vlp}\in \yadlp'$), is such that each $y_w\in\sY_{{\lbh}^{'\flat}_w,L_w'}$ (resp. $y_w'\in\sY_{\C_{p}^\flat,L_w'})$ is  a closed classical point in the fiber of the canonical quotient morphism $$\sY_{{\lbh}^{'\flat}_w,L_w'}\to\sX_{{\lbh}^{'\flat}_w,L_w'} (\textit{ resp. } \sY_{\cpt,L_w'}\to\sX_{{\cpt},L_w'})$$ such that the image of this chosen point under the composite morphism $$\sY_{{\lbh}^{'\flat}_w,L_w'}\to\sX_{{\lbh}^{'\flat}_w,L_w'}\to\sxqp (\textit{ resp. } \sY_{\cpt,L_w'}\to\sX_{{\cpt},L_w'}).$$ is the canonical point of $\sX_{\lbh^{'\flat}_w,\Q_p} (resp. \sxqp$). (this makes sense for all $w\in\vlp$ by \cite{joshi-teich-def}).  

The  residue field of the local ring of $y_w$ (resp. $y_w'$) in $\sY_{{\lbh}^{'\flat}_w,L_w'}$ (resp. in $\sY_{\cpt,L_w'}$) is an  algebraically closed, perfectoid field  which can be naturally identified with $\lbh_w$ (resp. naturally identified with $\C_p$) and  equipped with a natural action of $G_{L_w}$. The above discussion is summarized in following:

\bpro 
Let $\by_0'=(y_w')_{\in\vlp}$ be a standard point of $\yadlp'$. Then for each $w\in\vlnonp$, the image of $y_w'$ under the canonical quotient morphism $\yadlp'\to \xadlp'$ is the canonical point of $\xadlp'$.
\epro

\brem Let me remark that there are countably infinite possible choices for the standard point $\by_0$ (resp. $\by_0'$). This is because the set of valuations $\vl$ of any number field $L$ is countable and for each $v\in \vl$, the set of possible choices of closed classical points $y_v$ (chosen as above) is countable hence the set of $\by_0=(y_v)_{v\in\vl}$ as above is countable. However it is important to recognize that  the choice of any of the $y_v$'s, and hence of $\by_0$, is not unique. Any one of the countably many choices for $\by_0$ are equally valid choices (for doing arithmetic and geometry). But, as has been established in \cite{joshi-teich,joshi-teich-estimates,joshi-teich-def,joshi-untilts}, each choice  sees $p$-adic arithmetic (and geometry) slightly differently relative to other possible choices and so one can average over such choices.
\erem

\subparat{The choice of a standard point of Mochizuki's Ansatz}\label{ss:stand-theta-link} Now I will choose a standard point of $\bz_\Theta\in \tSigma_{L'}$ which one may think of as a standard $\Theta$-Link (for the purposes of \iutthr). Let $\bz_\Theta=(\by_1,\by_2,\ldots,\by_{\ells})\in \tSigma_{L'}$ be chosen such that $\by_{\ells}=\by_0'\in\yadlp'$  is the point standard point of $\yadlp'$ chosen in \cref{ss:stand-point}.

\newcommand{\sV}{\mathcal{V}}
\newcommand{\adeloid}{\mathcal{A\!d\!e\!l\!o\!i\!d}}
\newcommand{\ideloid}{\mathcal{I\!d\!e\!l\!o\!i\!d}}
\newcommand{\sH}{\mathscr{H}}
\subsection{Relationship with the global arithmetic period mapping of $\yadl$}\label{ss:theta-hyperplanes}\nwss
I want to discuss the relationship between properties of Mochizuki's Ansatz constructed above, especially \Cref{th:adelic-theta-link}, \Cref{cor:val-norm-theta}, and the global arithmetic period mapping of $\yadlp$ constructed in \constrtwoh{\ssep}{ss:period-map-prod-form}. This relationship  (\Cref{th:adelic-theta-link} and \Cref{th:hyperplane-theta-lnk}) is an important, but undemonstrated point in \cite{mochizuki-iut3}, and  the absence of proof in \iutthr, led to the denial of its existence in \cite[Section 2.2]{scholze-stix}.
\bthm\label{th:hyperplane-theta-lnk} 
Let $L$ be a number field, fixed as above.  Let $\bz=(\by_1,\by_2,\ldots,\by_\ells)\in\tSigma_{L'}$. Let $(\arithl_{\by_1},\ldots,\arithl_{\by_\ells})$ be the tuple of arithmeticoids provided by $\bz$. Consider the $\R$-vector space considered in  \constrtwoh{\ssep}{ss:period-map-prod-form} $$\sV_{L,j}=\bigoplus_{v\in\vl} \R_{v,j},$$
where each $\R_{v,j}=\R$ is given its usual absolute value.
Then 
\benumlab
\item  for $j=1,2,\ldots,\ells$, each the ideloid $\ideloid_{\arithl_{\by_j}}$ provides a homomorphism of groups
$$\log_{\arithl_{\by_j}}: \ideloid_{\arithl_{\by_j}}  \to \sV_{L,j}$$ given by
$$ \ideloid_{\arith{L}_{\by_j}} \ni (x_v)_{v\in \vl} \mapsto  (\log\abs{x_v}_{K_{y_v}})_{v\in\V_L}.$$
\item Each normalized arithmeticoid $\arithl_\by$ provides a degree homomorphism
$$ \ideloid_{\arith{L}_{\by_j}} \to \sV_{L,j} \to \R$$
given by 
$$ (x_v)_{v\in \vl}\mapsto (\log\abs{x_v}_{K_{y_v}})_{v\in\V_L}\mapsto  \sum_{v\in\vl} \log(\abs{x_v}_{K_{y_v}})=\deg_{\by_{\by_j}}((x_v)_{v\in\vl}).$$
\item As a consequence, each normalized $\arith{L}_{\by}^{nor}$ provides a hyperplane $H_{\by_j}\subset \sV_{L,j}$ given by the logarithm of product formula equation \eqref{eq:prod-formula} for  $L\into\arithl_{\by_j}^{nor}$:
$$H_{\by_j}:\deg_{\by_j}((x)_{v\in\V_L})=\sum_{v\in\V_L}\log\abs{x}_{K_{y_{j,v}}}=0\qquad (\text{for all }x\in L^*).$$
\item The scaling of valuations by the factor $j^2$ at primes of semi-stable reduction which occurs via \Cref{th:adelic-theta-link}, together with the requirement that the global product formula holds for $\arithl_{\by_j}^{nor}$, means that the local valuations (given by an arithmeticoid $\arithl_{\by_j}$) at primes of non-semistable reduction must also be suitably (re)normalized so that the product formula holds (this is given to us by \Cref{cor:val-norm-theta}).
\item As the product formula \eqref{eq:prod-formula} holds for each normalized arithmeticoid $\arithl_{\by_j}^{nor}$, the hyperplane $H_{\by_j}$ is stable under multiplication by elements of $L^*$.  
\item Let $\P(\sV_L)$ be the projective space of hyperplanes in $\sV_{L}$.  Then one has a natural and non-constant function defined on Mochizuki's Adelic Ansatz (\ssep\ref{ss:mochizuki-adelic-ansatz}):
$$\yadlp^\ells\supseteq \tSigma_{L'}\to \prod_{j=1}^\ells\P(\sV_{L,j})$$
given by $\tSigma_{L'}\ni\bz=(\by_1,\ldots,\by_\ells)\mapsto (H_{\by_1},\ldots, H_{\by_\ells})$. 
\item Let $\P(\tSigma_{L'})\subset \prod_{j=1}^\ells\P(\sV_{L,j})$ be the image of $\tSigma_{L'}$ under this function. Then $\P(\tSigma_{L'})$ is equipped with a natural action of all the groups and symmetries given by \Cref{th:ansatz-arithmeticoids} which act up on $\tSigma_{L'}$.
\eenum
\ethm
\bp 
This is clear from  \constrtwoh{Theorem }{th:hyperplane},  and the properties of $\tSigma_{L'}$ established so far and detailed in \Cref{th:adelic-theta-link}, \Cref{cor:val-norm-theta}.
\ep

\brem\label{re:hyperplane-theta-lnk}\newcommand{\sFrob}{\mathcal{F\!r\!o\!b}}\ 
\benumlab
\item I expect that the above period mapping is continuous for a suitable topology on the codomain, but I hope to take  this up in a separate paper.
\item Mochizuki's discussion of the role of product formulas is in \cite[Remark 3.9.6]{mochizuki-iut3}, \cite[Remark 5.2.1]{mochizuki-iut1}, \cite[Page 71]{mochizuki-gaussian}. 
\item It is worth noting that in the context of \cite[Corollary 3.12]{mochizuki-iut3} Mochizuki considers (without a transparent proof)  a version of this construction by means of the theory of realified Frobenioid $\sFrob(L)^\R$ of $L$ (see my discussion of this Frobenioid is in \ssep\ref{se:frobenioids}). Specifically Mochizuki works with the tuple of realified Frobenioids $(\sFrob(L)_1^\R,\ldots, \sFrob(L)_\ells^\R)$ which provides, by the same degree mapping, the hyperplanes considered here.
\item The existence of the phenomenon and properties established by \Cref{th:hyperplane-theta-lnk} is vital to establishing  \cite[Corollary 3.12]{mochizuki-iut3}. But this type of proposition is claimed without proof in \cite{mochizuki-iut3}, \cite{fucheng}, \cite{yamashita}. Notably, because of the scaling property \Cref{th:adelic-theta-link},  the hyperplane $H_{\by_{j}}$ comes equipped with a scaling factor of $j^2$ (for appropriate coordinate $v\in\V_L$ depending on the semi-stable elliptic curve as in \cite{joshi-teich-rosetta}). The lack of proof of this property in \cite{mochizuki-iut3} leads to the assertions of impossibility of this phenomenon in \cite[Section 2.2, Page 10]{scholze-stix}.
\eenum
\erem

\section{Some local preliminaries}\label{se:local-preliminaries}
I will recall a few local preliminaries and constructions which are not covered in \cite{joshi-teich-estimates} but which will be used in the remainder of this paper.

Let $p$ be a prime. Let $\bQ_p\supset \Q_p$ be an algebraic closure of $\Q_p$. Let $\C_p$ be the completion of $\bQ_p$, and equip $\C_p$ with its canonical action of the absolute Galois group $G_{\Q_p}$ of $\Q_p$ (computed using $\bQ_p$). The following result will be used throughout the rest of the paper.

\subparat{Basic rings of $p$-adic Hodge Theory}

For any choice of an algebraic closure $\bQ_p$ of $\Q_p$ one has the construction of $p$-adic period rings of \cite{fontaine94a}. Let $\bdrp$ be the ring of $p$-adic periods constructed in \cite{fontaine94a} for the (Galois group, field) pair $(G_{\Q_p},\Q_p)$ and let $\bdr$ be its quotient field. By the results of several mathematicians, notably, Fontaine, Messing, Faltings and Tsuji, $\bdr$ is the ring of de Rham periods of varieties finite type, and  defined over finite extensions of $\Q_p$ contained in $\bQ_p$. The construction of $\bdr$ is functorial in the sense of \cite[1.5.6]{fontaine94a}.

\subparat{The ring $B_p=B_{\cpt,\Q_p}$} Let $B_p=B_{\cpt,\Q_p}$  be the ring constructed in \cite[Chapitre 2]{fargues-fontaine}; if $p$ is clear from the context then, I will write $B=B_p$. 
By \cite[Th\'eor\`eme 6.5.2]{fargues-fontaine} and \cite[10.1.1]{fargues-fontaine}), there exists an element $t\in B$ which generates the maximal ideal of the canonical point of $\sxqp$. From the point of view of modern $p$-adic Hodge Theory (as detailed in \cite{fargues-fontaine}), the ring $B[\frac{1}{t}]$ is also quite natural.  Let $B^+_p=B_{\cpt,\Q_p}^+$ be the ring constructed in \cite[1.10]{fargues-fontaine}; if $p$ is clear from the context then, I will write $B^+=B_p^+$.

\subsubparat{The ring $B_p^+$} Let me also remark that for the ring $B_p=B=B_{\cpt,\Q_p}$ one has $B\supset B^+\supset W(\O_{\cpt})[1/p]$ and $W(\O_{\cpt})[1/p]\supset W({\overline{\F}_p})[1/p]=\Q_p^{unram}$ (this is the maximal unramified extension of $\Q_p$) and so $B_p$ and $B_p^+$ is naturally a module for $\Q_p^{unram}$. Hence $B$ and $B^+$ are  $E_0$-algebras for any finite unramified extension $E_0/\Q_p$. This remark will be used in what follows.

\subsubparat{The rings $\bdrp,\bdr$}  By \cite[Th\'eor\`eme 6.5.2]{fargues-fontaine}, the ring $\bdrp$ is the completion of the local ring of $\sxqp$ at the canonical point of $\sxqp$ and hence the ring $\bdrp$ is a discrete valuation ring. Moreover the maximal ideal of this local ring  is generated by the element $t$  and $\bdr=\bdrp[\frac{1}{t}]$. 

\brem On the other hand  $B$ has uncountably many maximal ideals corresponding to closed classical points of $\syQp$ and my construction of $\thetaj$ and $\thetam$ require using the fact that $B$ has infinitely many maximal ideals. \erem

\subsubparat{The ring $B_E$} Let $E$ be a $p$-adic field and let $E_0\subset E$ be its maximal unramified subfield.  By \cite[Chap 1, 1.6.2 and Chap. 10]{fargues-fontaine}, the ring $B_E$ associated to the datum $(F=\cpt,E)$ is the ring $B_{E}=B_{\cpt,E}$.

\blem\label{le:ring-BE} Let $E$ be a $p$-adic field and let $E_0\subset E$ be the maximal unramified subfield of $E$. Then one has a natural isomorphism
\be 
B_E\isom B_{\cpt,\Q_p}\tensor_{E_0} E.
\ee
\elem
\bp 
This is \cite[Proposition 1.6.9]{fargues-fontaine}.
\ep
\subsubparat{The rings $B,\bdr$} The ring $B=B_{\cpt,\Q_p}$ is related to $\bdr,\bdrp$ as follows.
 For a $p$-adic field $E$, let $E_0\subset E$ denote the maximal unramified subfield of $E$. Then by \cite[Proposition 10.2.7]{fargues-fontaine} one has an injection (for each finite extension $E/\Q_p$ contained in $\bQ_p$).
\be\label{eq:B-in-BdR} 
B_E=B\tensor_{E_0} E\into \bdr.
\ee

In particular
\be\label{eq:B-in-BdR2} 
B=B_p\into \bdr.
\ee

\newcommand{\tB}{\tilde{B}}
\subsubparat{The ring $B\tensor_{\Q_{p}}E$}\label{pa:ring-BE-var} It will be useful to consider a ring which is bigger than the ring $B_E$ (\Cref{le:ring-BE}) which is better suited for calculations. Consider \be\tB=B\tensor_{\Q_p} E\ee then by usual transitivity properties of tensor products with respect to the inclusions $\Q_p\subset E_0\subset E$ and so one can write this ring as
\be 
B\tensor_{\Q_p} E=(B\tensor_{\Q_p}E_0)\tensor_{E_0}E
\ee
and since $B$ is also an $E_0$ algebra, 
\be\label{eq:tilde-B-prop}
B\tensor_{\Q_{p}}E_0=\overbrace{B\oplus \cdots \oplus B}^{[E_0:\Q_p] factors}.
\ee
Hence 
\be 
\tB_E=\overbrace{B_E\oplus \cdots \oplus B_E}^{[E_0:\Q_p] factors}.
\ee
At any rate, $\tB_E=B\tensor_{\Q_p}E$ is naturally a finite and free $B=B_{\cpt,\Q_p}$-module through this identification. 

One also has the natural diagonal embedding 
\be B_E\into \tB_E=B\tensor_{\Q_p} E= B_E^{\oplus [E_0:\Q_p]}.
\ee

\subsubparat{The ring $B_e$}\label{ss:Be} Let me also introduce the ring \cite[Th\'eor\`eme 6.5.2]{fargues-fontaine} 
\be B_e=B^+[1/t]^{\vphi=1}=B[1/t]^{\vphi=1}=\bcris^{\vphi=1}.
\ee
By \cite[Th\'eor\`eme 6.5.2]{fargues-fontaine} $B_e$ can be identified with the affine coordinate ring of the curve $\sxqp-\{x_{can}\}$ i.e. $\sxqp-\{x_{can}\}=\spec(B_e)$.
One also has, by \cite[Theorem 6.5.2]{fargues-fontaine},  the inclusions 
\be \xymatrix{
B^+\ar@{^{(}->}[r] & \bcris^+\ar@{^{(}->}[r] &  \bcris^+[1/t]=\bcris\ar@{^{(}->}[r] &\bdr.}\ee
These inclusions will allow one to view the theta-value locus as living in each of these period rings. Direct construction of the theta-value locus in $\bcris$ or $\bdr$ does not seem to be possible.

\brem  By \cite{fontaine94a,fontaine94b}, the rings $\bdr\supset \bcris$ depend on the choice of a fixed algebraic closure $\Q_p$. On the other hand, the rings $B,B_e,B_E$ of \cite{fargues-fontaine} do not depend on such a choice.
\erem

\subparat{Enlarging subsets of $B_{E}^\ells$ using $\Aut(G_{E})$-action, Frobenius action and convex closures}\label{ss:enlargements-by-anabelomorphy} Let $E$ be a $p$-adic field.
Given a subset $S$ of the ring $B_{E}$, one can enlarge it in many natural ways. Let me review some of these enlargements which will be applied to obtain the final versions of the theta-values locii in various related rings.

Suppose $S\subset B_{E}$ (the considerations of this paragraph will be applied with $E=E_w'$ for each relevant $w\in\ubblvossp$). I will enlarge this set further as follows: 
\subsubparat{$\Aut(G_E)$-enlargements}
Let $S\subset B_{E}$ or $B_{E}^\ells$ , consider the natural action of $\Aut(G_{E})$ on the data $G_{E}\act B_{E}$ by pre-composing the action of $G_{E}$ on $B_{E}$ with $\sigma\in\Aut(G_{E})$. More precisely, if $\sigma\in\Aut(G_{E})$ then consider the pre-composed or twisted action 
\be 
G_{E}\mapright{\sigma} G_{E}\act B_{E}
\ee
and let $S^{\sigma}$ be the image of $S$ under this new action of $G_{E}$ on $B_{E}$. This consideration can be applied verbatim to the action $G_{E}$-action on $B_{E}^\ells$ and subsets in this context. So if $S\subset B_{E}^\ells$ is any subset  then one simply  takes the union 
\be 
\bigcup_{\sigma\in \Aut(G_{E})} S^\sigma\subset B_{E}^\ells.
\ee
This provides a \textit{Galois  enlargement} of any subset  $S\subset B_{E}^\ells$. 

\brem\ 
\benumlab
\item By \cite{joshi-gconj}, the Fargues-Fontaine curve $\sxfe$ is not uniquely identified by its \'etale fundamental group--this group is isomorphic to the absolute Galois group $G_E$ (this is the failure of the absolute Grothendieck conjecture for Fargues-Fontaine curves which is established in \cite{joshi-gconj}).
\item In terminology of \iut, taking unions over the images under automorphisms of $\Aut(G_{E})$ is Mochizuki's way of applying the indeterminacy Ind1 in \cite[Theorem 3.11]{mochizuki-iut3} (my discussion of this is in \cref{ss:Mochizuki-indeterminacies}).
\item In this sense Mochizuki's indeterminacy Ind1 arises from the failure of the absolute Grothendieck conjecture for Fargues-Fontaine curves established in \cite{joshi-gconj}.
\eenum
\erem
\subsubparat{Frobenius stablization enlargement}
Another possible enlargement which can be applied is to replace $S\subset B_E$ by \be S\mapsto \bigcup_{n\in \Z}\vphi^n(S),\ee where $\vphi$ is the Frobenius morphism of $B_{E}$. This obviously applies to subsets of $B_E^\ells$ and variants. This is the \textit{Frobenius enlargement} of $S$.
\subsubparat{Convex closure enlargement} Finally one can replace such a Galois and Frobenius   enlarged set $S$ by the  closure of its convex hull   $\overline{S}$  in $B_{E}$ with respect to the Frechet topology i.e. one can replace $S$ by the smallest, closed, convex subset containing $S$. This set is closed and convex (see \cite[Chapter I]{schneider-book} or \cite{schikhof-book1} for basic properties of convex sets in $p$-adic Banach spaces).

\brem 
My version of \moccor\ and Mochizuki's version are both about lower bounds and it is important to recognize that enlargements (i.e. adding elements to the set $S$) do not affect lower bounds on the supremum of norms of elements of $S$. On the other hand, on a suitably enlarged set, an upper bound may be relatively easier to obtain. This is the rationale for establishing \moccor.
\erem

\newcommand{\thetajq}{\thetaj^?}

\section{Construction of $\thetaj^?$}\label{se:construction-of-thetaj-and-thetam}
Let $L$, $C/L$ and $X/L$ be as assumed in   \cref{ss:elliptic-curve-assumptions} and \inithtdata.  In this section, subject to these assumptions, I provide the construction of the theta-values locii $\thetaj^?$ (for suitable choice of decorations $?$) and Mochizuki's set $\thetami$ (and its variant $\thetaji$).  The sets $\thetaji$ and $\thetami$ are needed for the statements of the fundamental estimate i.e the statements of \moccor\ in the respective theories. The construction of $\thetami$ is the central assertion of \cite[Theorem 3.11]{mochizuki-iut3}. 

One of the central difficulties which many readers have faced in ascertaining the veracity of \iutthr\ has been the lack of a quantitative demonstration of the existence of distinct arithmetic holomorphic structures required for the construction of the set $\thetami$ which plays a fundamental role in \moccor. 

This difficulty has been overcome through the present series of papers: my demonstration of the existence of arithmetic holomorphic structures is in \cite{joshi-teich,joshi-untilts,joshi-teich-def} and my construction of $\thetami$ is in \cref{se:mochizuki-construction-thetam}. This section provides a construction of a simpler subset $\thetaji$ which is adequate for establishing \moccor.
 
From my point of view the existence of these locii  is a consequence of the existence of arithmetic Teichmuller spaces and its various properties established in \cite{joshi-teich,joshi-untilts,joshi-teich-def,joshi-teich-estimates} (and the present paper). Mochizuki's approach to the construction of $\thetami$ (modulo the existence of arithmetic holomorphic structures and their properties) is detailed in \cite[Theorem 3.11]{mochizuki-iut3} and its proof.

The core idea of the construction of these sets is the same as that of \iutabc--namely one collates, in a suitably chosen ambient set, some arithmetic information (Tate parameters or more generally theta-values) arising from many different holomorphoids of $X/L$. Depending on the choice of ambient set used to collate the theta values, one obtains the version of the locus of that flavor. My  construction of Mochizuki's set $\thetami$ is detailed in \cref{se:mochizuki-construction-thetam}. It should be remarked that the verb \textit{to collate} is not used in the context of \cite[Theorem 3.11, Corollary 3.12]{mochizuki-iut3}, but in my opinion it clarifies immensely what \cite[Theorem 3.11]{mochizuki-iut3} is all about. In fact, this verb is never used in \iut, \cite{mochizuki-gaussian} and \cite{mochizuki-essential-logic} at all.

The approach adopted here is based on \cite{joshi-teich-estimates} but, like \iutthr, uses the (global)  Initial Theta Data \inithtdata\  in a crucial way and is a more elaborate one than the local prototype of \cite{joshi-teich-estimates} which provides a natural demonstration of the  structures claimed in \cite{mochizuki-iut3}. A reading of that paper (\cite{joshi-teich-estimates}) is highly recommended along with this one.

\subparat{Descent of holomorphoid data to $\lmod$}
The following proposition is useful in understanding how one may descend holomorphoid data from $L'$ to $\lmod$ using \inithtdata. 
\bthm\label{th:lmod-arith-descent}
Let $C,X,L,L'$ be as in \assumptions. Let $\holt{X/L'}{\by}$ be a holomorphoid of $X/L'$ and let  $\arith{\by}$ be the arithmeticoid of $L'$ given by this holomorphoid.  Then 
\benumlab
\item the bijection (see \inithtdata) $$\V_{\lmod}\isom \underline{\V}\subset \V_{L'}$$ together with the arithmeticoid $\by$ provides us with an arithmeticoid $\arith{\lmod}_{\underline{\by}}$  with ${\underline{\by}}\in\yadlmod$ naturally given using $\by$ and the above bijection.
\item  In particular, every holomorphoid $\holt{X/L'}{\by}$ of $X/L'$ provides us an arithmeticoid $$\arith{\lmod}_{\ul{\by}}$$ of $\lmod$.
\eenum
\ethm
\bp To prove all the assertions, it is enough to provide a natural construction of $\underline{\by}\in\yadlmod$ from $\by\in\yadlp$ given by the data of the holomorphoid $\holt{X/L}{\by}$ and the bijection $\V_{\lmod}\isom \underline{\V}\subset \V_{L'}$ provided by \assumptions. By definition (\constrtwoh{\ssep}{se:adelic-ff-curves}) one has that $$\yadlp=\prod_{w\in\vlp} \abs{\syflwp}.$$ Now $\underline{\by}\in\yadlmod$ is the point obtained by projection of $\by=(y'_w)_{w\in\vlp}$ on to the coordinates corresponding to the ``index'' set identified by $\underline{\V}\isom\V_{\lmod}$ using the natural (continuous) mapping \cite[Proposition 2.3.20]{fargues-fontaine} $$f_{w|v}:\abs{\sY_{\cpt,L_w'}}\to \abs{\sY_{\cpt,{\lmod}_{,v}}} \qquad (\text{ with } w|v|p)$$ 
i.e. define $$\ul{\by}=(f_{w|v}(y_w'))_{w\in\ul{\V}}\in\yadlmod'.$$
Then $\by\mapsto \uby{\by}\in\yadlmod$. This completes the construction of $\ul{\by}$.
\ep

\subsubpara\label{ss:underline-notation-convention} From now on, and without further mention, $\arith{\lmod}_{\ul{\by}}$ will mean the arithmeticoid of $\lmod$ constructed by \Cref{th:lmod-arith-descent} from the arithmeticoid $\arith{L'}_{\by}$ and in general underlined version of objects will refer to objects arising from arithmeticoid $\arith{\lmod}_{\ul{\by}}$ constructed using \Cref{th:lmod-arith-descent}.

\subparat{The definition of $\thetajq$ and its variants}
The construction of the set $\thetajq$ is  adelic and it will be constructed here by constructing subsets $\thetajpq$ for every prime $p$ and then taking product over all rational primes $p$. The construction of the adelic set $\thetajq$ (and of each $\thetajpq$) requires fixing an  Initial Theta Data \inithtdata. This is also required  in \iutthr\ for the construction of $\thetami$. 

\subsubparat{Local components of $\thetajp$}\label{ss:definition-of-thetaj-using-ansatz}.
The construction of $\thetajp$ is based on the local construction of \cite{joshi-teich-estimates}. The important advantage of the method of \cite{joshi-teich-estimates}  is that for any prime $w|p$,  it provides  theta-values locii $$\thetaj^{B_{?}}$$ using a wide class of rings $B_{?}$ used in $p$-adic Hodge Theory, for example 
$${B_{?}}=B, B[1/t], B_I\ (I\subset [0,1]\subset \R), \bcris\text{ and } \bdr,$$ (see \Cref{th:local-theta-val-locus}).  At the same time, by choosing a suitable Galois cohomology group (for each $w|p$ and for each rational prime $p$) to serve as a receiver of theta-values, viewed as cohomology classes via \cite[Proposition 1.4]{mochizuki-theta}, one  establishes the construction of $\thetam$ along the lines of \iutthr\ using my canonical approach to arithmetic holomorphic structures. 

On the other hand working with these rings of $p$-adic Hodge Theory provides yet another distinct advantage: one can multiply (and add) together theta-values arising from distinct arithmetic holomorphic structures for distinct primes $w$ lying over $p$. This bundling of theta-values arising from distinct arithmetic holomorphic structures and for distinct primes $w|p$ (for each fixed rational prime $p$) also requires fixing an Initial theta-data \inithtdata. In \iutthr, Mochizuki works with a Galois cohomology group (for each prime $w|p$ and each prime $p$) and so there is no natural ring structure to multiply theta-values at all. To circumvent this problem, Mochizuki forces multiplication of theta-values by means of his theory of tensor packets \cite[Section 3]{mochizuki-iut3} (to achieve this bundling of theta-values).

In particular, my approach to the construction of $\thetajq$ also demonstrates that  Mochizuki's framework of tensor-packets, arises far more naturally from the theory of Arithmetic Teichmuller spaces as documented in \present.

\subparat{Mochizuki's Adelic Ansatz and arithmeticoids} Let $L$, $C/L$, and $X/L$ be as in \cref{ss:elliptic-curve-assumptions}. Fix an Initial theta-data \inithtdata. Let $$\bz=(\by_1,\ldots,\by_\ells)\in\tSigma_{L'}$$ be a point of Mochizuki's Adelic Ansatz \cref{se:adelic-ansatz} i.e. $\bz$ is a $\Theta_{gau}$-link (\cref{def:theta-link}). Then by \cite{joshi-teich-def} $\bz$ provides an $\ells$-tuple of arithmeticoids of $$(\arith{L'}_{\by_1},\arith{L'}_{\by_2}, \ldots, \arith{L'}_{\by_\ells})$$ and hence one obtains an $\ells$-tuple of holomorphoids of $$(X/\arith{L'}_{\by_1},X/\arith{L'}_{\by_2}, \ldots, X/\arith{L'}_{\by_\ells}).$$ of $X/L$. The following is a fundamental consequence of \cite{joshi-teich-def}:
\bthm\label{th:ansatz-arithmeticoids} For any $\bz=(\by_1,\ldots,\by_\ells)\in\tSigma_{L'}$, 
the holomorphoids  $$X/\arith{L'}_{\by_1},X/\arith{L'}_{\by_2}, \ldots, X/\arith{L'}_{\by_\ells}$$ are all distinct and need not even provide isomorphic pairs of analytic spaces at any prime $w\in\vlnonp$.
\ethm
\bp 
This is immediate from \constrtwoh{Theorem }{th:consequence-inequivalent-arithmeticoids}.
\ep

\subparat{The collation of theta-values} Let $\bz\in(\by_1,\ldots,\by_\ells)\in\tSigma_{L'}$ (in Mochizuki's parlance, each $\bz$ is a $\Theta$-Link). For   the given $X,L,L'$ (as in \inithtdata), \Cref{th:ansatz-arithmeticoids} provides, for each $\bz\in\tSigma_{L'}$ and for each $w\in\ubblvoss$, the tuple of Tate parameters 
\be\label{eq:tuple-tate-params} \left(q_{w,\arith{\by_1}},q_{w,\arith{\by_2}},\ldots,q_{w,\arith{\by_\ells}}\right).
\ee

The basic idea behind the construction of $\thetaj^?$ and $\thetam^?$ is that one wants to collate, for the given $X,L,L'$ (as in \inithtdata), and for each $\bz\in\tSigma_{L'}$ and for each $w\in\ubblvoss$, the tuple of Tate parameters \eqref{eq:tuple-tate-params}
 in a suitable container set (here $q_{w,\arith{\by_j}}$ is the Tate parameter of $X/L_w'$ computed in the arithmeticoid $\arith{\by_j}$). 

The main problem here is that, even for any single $w|p$, the individual elements $q_{w,\arith{\by_1}}$, $q_{w,\arith{\by_2}}$, $\ldots$ , $q_{w,\arith{\by_\ells}}$ of the tuple given in eq.~\eqref{eq:tuple-tate-params}  do not live in a natural common set where these elements may be compared with each other.   \textit{This is where the choice of a common set $?$ (in the superscripts of $\thetam^?$, $\thetajq$) comes into picture.} In \cite{joshi-teich,joshi-teich-estimates}, I recognized that one may choose fundamental ring of $p$-adic Hodge theory as a common location for comparing theta-values and there I introduced therein the theory of lifting of theta-values (this is recalled here briefly in \cref{ss:main-def-theta-values-locus-joshi}).

Mochizuki's choice of his common set  in \cite{mochizuki-iut3} is based on his theory of Galois Theoretic Theta-evaluation \cite[Proposition 1.4]{mochizuki-theta} (applied to the arithmeticoid $\arith{\by_j}$) and then observing that the relevant Galois cohomology group provided by any arithmeticoid  is amphoric (\constrtwoh{Proposition }{pr:arith-isom-gal-cohom}). The method of collation of cohomology classes (in the cohomology of a fixed arithmeticoid) is established here in \cref{ss:collation-of-classes}. Hence the methods developed in \present\ is better adapted for providing a uniform way at arriving at both the constructions.

\subsubpara Let $\bz\in(\by_1,\ldots,\by_\ells)\in\tSigma_{L'}$. Fix a prime number $p$. By \Cref{th:adelic-theta-link}, for each $w\in\ubblvossp$ (this means $w|p$), let $(y_{w,1},\ldots,y_{w,\ells})$ be the $w$-component of $\bz$. Let $E_w'=L_{w}'$ be the completion of $L'$ at $w$. I will use $E_w'$ instead of $L_w'$ for notational continuity with \cite{joshi-teich-estimates} and \cite{fargues-fontaine}. 
I will now apply the construction of the theta-value locus detailed in \cite[\ssep 7]{joshi-teich-estimates} for each $w\in\ubblvossp$. 

The construction given in \cite{joshi-teich-estimates} works with the ring $B=B_p$. Here I work with the ring $B_{E'_w}$. So some changes to the constructions of \cite{joshi-teich-estimates} are needed. Let me detail them here. For simplicity of notation (as $w$ is fixed for this discussion) write 
$$(y_1',y_2',\ldots,y_{\ells}')=(y_{w,1},\ldots,y_{w,\ells})$$
i.e. $w$ is suppressed from the subscript until such time as one needs it and again for simplicity of notation write $E=E_w'=L_w'$.

Following \Cref{th:ansatz-arithmeticoids} (or \cite{joshi-teich-estimates}), one obtains from an $\ells$-tuple of local holomorphoids 
\be 
(\hol{}{y_{1}'}{X/L_w'}, \hol{}{y_{2}'}{X/L_w'}, \ldots, \hol{}{y_{\ells}'}{X/L_w'})
\ee
of holomorphoids of $X/E_w'$ lying over  $(y'_1,y'_2,\ldots,y'_\ells)\in \tSigma_{\cpt}^{E_w'}$ as our $\Theta_{gau}$-Link.

In the next section, I will define theta-values as in \cite[\ssep 7]{joshi-teich-estimates} but consider lifts to  $B_{E_w'}$ as follows.

\subsubparat{Lifting theta-values}\label{ss:main-def-theta-values-locus-joshi} Let me briefly recall results established in \constrtwo{\ssep}{se:construction-ttheta} which deals with lifting of theta-values to the ring $B$ and other rings of this type. To avoid notational clutter, let $E$ be a $p$-adic field (in applications as above it $E=E_w'$ as above). Here one wants to work with the ring $B_E$ but the procedure is the same. Fix a Lubin-Tate formal group $\sG/\O_{E}$. Let $\pi\in\O_E$ be a uniformizer of $\O_E$.

Let $y'\in\sY_{\cpt,E}$ be a closed classical point and let $\hol{X/E}{y'}{\cpt}=(X/E,(E\into K,K^\flat\isom \cpt), *_K:\sM(K)\to \xan_E)$ be a holomorphoid of $X/E$. Let $\xi=\xi_1$  be one of the theta values (see \cite[\ssep 6]{joshi-teich-estimates}) in $K$. 

Then by \constrtwo{Theorem }{thm:theta-pilot-object-appears} one has a Teichmuller lift $[z]\in W_{\O_E}(\O_{\cpt})$ with the property that for any $\rho\in (0,1]\subset \R$, one has $$\abs{[z]}_\rho=\abs{z}_{\cpt}=\abs{\xi}_K.$$

Since $W_{\O_{E}}(\O_{\cpt})\subset B_E$, $[z]$ provides a lift of $\xi$ to $B_E$. Such a lift is not uniquely defined and so one considers  special lift of $\xi$ to $B$ and one can add to $[z]$ any element of the Tate-module $T(\sG(\O_K))\subset B_E^{\vphi=\pi}$ (\cite[Chapitre 4]{fargues-fontaine}).

The Tate module $T(\sG(\O_K))\subset B^{\vphi=\pi}$ is a free $\O_E$-submodule of rank one which depends on $y'$ \cite[Chapitre 4]{fargues-fontaine}. Let $t_{K_{y'}}\in B_E$ be its generator. Then the lifts of $\xi$ to $B_E$ of  interest to us are the lifts of the form
\be\label{eq:simple-lift} \Xi=[z]+\O_E\cdot t_{K_{y'}}.
\ee

A lift \eqref{eq:simple-lift} $\Xi$ originating from the holomorphoid $\hol{X/E}{y'}{\cpt}$  will be labeled simply by $$\Xi^z_{K_{y'}}$$ instead of $$\Xi_{\hol{}{y'}{X/L_w'}}$$ and this will often be contracted to $\Xi_{K_{y'}}$.  

\brem   Note that if $[z],[z']$ are two lifts of $\xi\in K_{y'}$  to $B_E$ then $\abs{z}_{\cpt}=\abs{z'}_{\cpt}$ by \cite[Proposition 2.2.16]{fargues-fontaine}. Hence the absolute values of $z,z'\in\cpt$ for any pair of Teichmuller lifts $[z],[z']$ of $\xi\in K$ are independent of the choice of the Teichmuller lifts. 
\erem
 
\newcommand{\bxi}{\boldsymbol{\Xi}}
\subsubparat{Admissible lifts of theta-values} 
With this preparation I am ready to define the set of admissible lifts of theta-values, or more simply an admissible theta-values-lift, to $B_E^\ells$.

Let $(y_1',y_2',\ldots,y'_\ells)\in\tSigma_{\cpt}^E\subset \sY_{\cpt,E}^\ells$ be a tuple of points in Mochizuki's Ansatz $\tSigma_{\cpt}^E$ for $B_E$. Let $$i_1:E\into K_{y_1'}, i_2:E\into K_{y_2'}, \ldots, i_\ells:E\into K_{y_\ells'}$$
be isometric embeddings of $E$ in the residue fields of the point $(y_1',y_2',\ldots,y'_\ells)$.

\bdefn\label{def:additive-theta-pilot-object}
An \textit{admissible lift of theta-values at $w\in\ubblvossp$} or an \textit{admissible theta-values-lift at $w\in\ubblvossp$} is the tuple (in $B_E^\ells$):
\be\label{eq:additive-theta-pilot-object} \bxi_{\lambda,\bz,w}^{z_1,\ldots,z_\ells}=(\Xi^{z_1}_{\lambda,K_{y_1'}}:=[z_1]+i_1(\lambda)t_{K_{y_1'}},\ldots, \Xi^{z_\ells}_{\lambda,K_{y_\ells'}}:=[z_\ells]+i_\ells(\lambda)t_{K_{y_\ells'}})\in B_E^\ells,\ee
where  $$\lambda \in \O_E;$$
and $$[z_1],\ldots,[z_\ells]$$ are Teichmuller lifts to $B_E$ of the theta values $$\xi_{1;K_{y_1'}}, \ldots, \xi_{1;K_{y_1'}} $$
  from the residue fields $$K_{y_1'},\ldots,K_{y_\ells'}$$
   of $(y_1',y_2',\ldots,y'_\ells)$ to the ring $B_E$ under the natural homomorphism
   $$B_E\to B_E/\fm_{y'_j}=K_{y_j'} \text{ for } j=1,\ldots,\ells$$ 
   and where $t_{y_j'}\in \fm_{y'_j}$. 
   The method of lifts are described earlier and also detailed in \constrtwo{\ssep}{se:lifting-values-to-B}.
\edefn 

Often I will suppress the $z_1,\ldots,z_\ells$ from the notation $\bxi_{\lambda,\bz,w}^{z_1,\ldots,z_\ells}$ and write 
$$\bxi_{\lambda,\bz,w}$$
for the right hand side of \eqref{eq:additive-theta-pilot-object}, and often I will also suppress $\lambda$ and write 
$$\bxi_{\bz,w}$$
for the right hand side of \eqref{eq:additive-theta-pilot-object}. In particular, I will write $\bxi_{\lambda,\bz,w}^{z_1,\ldots,z_\ells}$ \eqref{eq:additive-theta-pilot-object} simply as 
$$\bxi=\bxi_{\bz,w}=(\Xi_{K_{y_1'}},\ldots, \Xi_{K_{y_\ells'}})\in B_E^\ells.$$

\brem 
Let me remark that each $\bxi_{\lambda,\bz,w}^{z_1,\ldots,z_\ells}$ is the analog of Mochizuki's $\Theta$-Pilot Object in \cite[Definition 3.8(i)]{mochizuki-iut3}, \cite[Page 72]{mochizuki-gaussian}. 
\erem

\subpara\label{ss:pilot-obj-components-in-gen} If $w\in\vlp-\ubblvoss$, then one defines
\be 
\bxi_{\bz,w}=([1],\ldots,[1]) \in B_{E}^\ells.
\ee 
Here $[1]\in W(\O_{\cpt})\subset B_E$ is the Teichmuller element.
The following lemma is trivial:
\blem\label{le:triv-lem}
 For each $\rho\in(0,1)\subset\R$ one has $$\abs{[1]}_{B_E,\rho}=1.$$
\elem
\bp
Clear from  \cite[Chapitre 1, D\'efinition 1.4.1]{fargues-fontaine}.
\ep

\subparat{Admissible lifts and Tate parameters} The connection between admissible theta-values lifts and Tate parameters is given by the following:
\bpro\label{pr:admissible-lifts-tate-params} 
Let $X,L,L'$ be as in \cref{ss:elliptic-curve-assumptions}, \inithtdata. Let $\bz\in\tSigma_{L'}$, and let $w\in\ubblvossp$ for some prime $p$. Let $\bxi_{\bz,w}=(\Xi_{\lambda,K_{y_1'}},\ldots, \Xi_{\lambda,K_{y_\ells'}})$ be an admissible theta-values-lift at $w$. For $j=1,\ldots,\ells$ let $$\psi_j:B_E\to K_{y_{j}}$$ be the natural quotient homomorphism arising by taking quotient modulo the (principal) maximal ideal corresponding to the closed classical point $y_j'\in\sY_{\cpt,E_w'}$. Then for $j=1,\ldots,\ells$ one has 
$$\psi_j(\Xi_{\lambda,K_{y_1'}})=\psi_j([z_{w,j}])=\xi_{1,K_{y_j'}}=q_{X/E_w',y_j'}\in K_{y_j'} ,$$
where is $q_{X/E_w',y_j'}$ the Tate parameter of $X/E_w'$ in the arithmetic holomorphic structure provided by $y_j'$.
\epro
\bp 
Let $\fm_j\subset B_E$ be the maximal ideal corresponding to the point $y_j'$. Then by \cite[Chap. IV and Chap. VI]{fargues-fontaine} $\fm_j=(t_{K_{y_j'}})$ is a principal ideal and this generator maps to zero under $B_E\to B_E/\fm_j=K_{y_j'}$. Hence 
$$\psi_j(\Xi_{\lambda,K_{y_1'}})=\psi_j([z_{w,j}]).$$
The assertion now follows as $[z_{w,j}]$ is a Teichmuller lift of $\xi_{1,K_{y_j'}}=q_{X/E_w',y_j'}$.
\ep

\subsubparat{The theta-values locus in $\thetaj^{B_E}$}\label{pa:basic-thetaj-def} Now I am ready to define the basic theta-values locus in $B_E^\ells$. 
The theta-values set  constructed  below (as well as in \iut) is  the union of all the admissible theta-values-lifts. 

\bdefn\label{de:basic-theta-val-locus} Let $X,L,L'$ be as in \cref{ss:elliptic-curve-assumptions}, \inithtdata. For $w\in\vlnonp$, write $E=E_w'=L_w'$ for notational simplicity. Then the \textit{theta-values locus at (a fixed) $w\in\vlnon$}  is the subset $\thetaj^{B_E}\subset B_E^\ells$ is defined as follows: 
\benumlab
\item if  $w\in\ubblvossp$, then $\thetaj^{B_E}\subset B_E^\ells$ consists of all admissible theta-values lifts  $\bxi_{\lambda,\bz,w}^{z_1,\ldots,z_\ells}$ at a fixed $w\in\ubblvossp$ \eqref{eq:additive-theta-pilot-object} where $\bz$ runs over $\tSigma_{L'}$. Explicitly this is the set:
\begin{multline}\label{eq:basic-thetaj-def} \thetaj^{B_E}=\left\{ \bxi_{\lambda,\bz,w}^{z_1,\ldots,z_\ells}= ([z_1]+i_1(\lambda)\cdot t_{K_{y_1'}},\ldots, [z_\ells]+i_\ells(\lambda)\cdot t_{K_{y_\ells'}})\in B_E^\ells: \right.\\
\qquad\qquad \text{ for all } [z_1], \ldots,[z_\ells] \text{ lifting } \xi_{1,K_{y_1'}},\ldots,\xi_{1,K_{y_\ells'}} \text{ to } B_E, \\
\text{ for all }\lambda\in\O_E, \text{ and, } \\
\left. \qquad\qquad \qquad  \qquad\qquad \qquad \qquad\qquad \text{ for all } \bz \in\tSigma_{L'} \ (\text{here } \bz_w=(y_1',y_2',\ldots,y'_\ells) \in\tSigma_{\cpt}^{E})   \right\}.\hfill
\end{multline}
\item If $w\in\vlp-\ubblvoss$, then by the definition of $\bxi_{\bz,w}$ for $w\in\vlp-\ubblvoss$ one has
\be 
\thetaj^{B_E}=\left\{\bxi_{\bz,w}: \bz\in\tSigma_{L'}\right\}=\left\{([1],\ldots,[1]) \right\}\in B_{E}^\ells.
\ee
\eenum 
In particular, $\thetaj^{B_E}\subset B_E^\ells$ is defined for all $w\in\vlp$.
\edefn

\subparat{Enlargements of subsets of $\thetaj^{B_E}$} Now one applies the enlargement procedures \cref{ss:enlargements-by-anabelomorphy} to the set $\thetaj^{B_E}$ defined in \Cref{de:basic-theta-val-locus} to obtain an enlargement of $\thetaj^{B_E}$. \textit{This  enlarged set will again be denoted by $\thetaj^{B_E}$.} The enlarged set $\thetaj^{B_E}$ will be  called the \emph{basic theta-values locus for the ring $B_E$}. Let me record the construction of this set in the following theorem:

\bthm\label{th:enlarged-theta-value-locus-local} Let $C/L,X/L$ and $L'$ be as in \inithtdata. Then for each $w\in\ubblvossp$ and with $E=E_w'=L_w'$ one has a basic theta-values locus $\thetaj^{B_{E_w'}}=\thetaj^{B_{L_w'}}$ is 
\benumlab
\item Galois stable for $G_{E_w'}$,
\item $\vphi$-stable for where $\vphi$ is the Frobenius of $B_{E_w'}$,
\item and also $\Aut(G_{E_w'})$-stable in the sense of earlier definitions.
\eenum
\ethm 
\bp 
The proof is clear by construction and the definitions. 
\ep

\subsubpara{Variants of $\thetaj^{B_E}$} Let me sketch the construction of some useful variants of $\thetaj^{B_E}$. Let $I\subset [0,1]\subset \R$ be an interval and let $B_I$ be the ring defined in \cite[D\'efinition 1.6.2]{fargues-fontaine}. The ring $B$ is obtained as the ring $B_I$ for $I=(0,1)$ i.e. $B=B_{(0,1)}$ by \cite[D\'efinition 1.6.2]{fargues-fontaine}. By imitating the above procedures, one can also construct  theta-values locii $\thetaj^{B_I}$ with values in the rings $B_I$ for each interval $I\subset [0,1]$.  

Notably, one also has a naturally defined set $$\thetaj^{B}$$ (\constrtwo{\ssep}{se:construction-ttheta}) with values in $B$ instead of $B_E$.

Now let $E_0\subset E$ be the maximal absolutely unramified subfield of $E$. Then one can also work with $B_I\tensor_{E_0}E$ instead of $B_I$ and construct the set $$\thetaj^{B_I\tensor_{E_0}E}.$$  

Let $t\in B$ be the element constructed in \cite[10.1.1]{fargues-fontaine}.
Then one can also  construct theta-values locii $$\thetaj^{B[1/t]}$$ and $$\thetaj^{B_E[1/t]}$$  with values in  $(B[1/t])^\ells$ and  $(B_E[1/t])^\ells$ respectively.  

For future reference, let me also discuss the variant of the ring $B_E$ discussed in \cref{pa:ring-BE-var}. Let $E_0\subset E$ be the maximal, absolutely unramified subfield of $E$. Let $\tB_E=B\tensor_{\Q_p} E$ be the ring defined in \cref{pa:ring-BE-var}. By \eqref{eq:tilde-B-prop} one has an embedding
$$\tB_E = \overbrace{B_E\oplus\cdots \oplus B_E}^{[E_0:\Q_p] copies} = B_E^{[E_0:\Q_p]}.$$
Hence one may also construct theta-value locus
$$\thetaj^{\tB_E}\subset \tB_E^\ells.$$
One can also consider $\tB_E[1/t]$ and construct the theta-values locus
$$\thetaj^{\tB_E[1/t]}\subset \left(\tB_E[1/t]\right)^\ells.$$

Let me record these constructions in the following theorem:
\bthm\label{th:local-theta-val-locus}
Let $L, C/L,X/L$ and $L'$ be as in \cref{ss:elliptic-curve-assumptions}, \inithtdata. Let $w\in\ubblvoss$ and write $E=E_w'=L_w'$. Let $B_?$ be any one of the following rings (\cite{fargues-fontaine}) $$\left\{ B,B[1/t],B_E,\tB_E,B_E[1/t],\tB_E[1/t] \right\}\cup\left\{B_I,B_I\tensor_{E_0} E: I\subset [0,1]\subset \R \right\}.$$ Then one can construct the theta-values-locus $$\thetaj^{B_?} \subset B_{?}^\ells$$
satisfying properties similar to those established above.
\ethm

\brem 
Note that the rings in the second set have more complicated behavior with respect to Frobenius but a formulation of the properties in this case is not difficult.
\erem

\subparat{Embedding $\thetaj^{B_E}$ in $\bdr$} The rings $B_?$ are related to other rings of classical $p$-adic Hodge Theory (\cite{fontaine94a}) in various ways. Using such relationships one can embed the theta-values locus $\thetaj^{B_?}$ in the conventional rings of classical $p$-adic Hodge Theory. Notably one ha the following:

\bthm\label{th:local-theta-val-in-bdr} 
Let $L, C/L,X/L$ and $L'$ be as in \cref{ss:elliptic-curve-assumptions}, \inithtdata. Let $w\in\ubblvoss$ and write $E=E_w'=L_w'$. One has an embedding
\be \thetaj^{B_E}\subset B_E^\ells\into \bdr^\ells.\ee
\ethm
\bp 
This is immediate from the fact \eqref{eq:B-in-BdR} that $B_E\into \bdr$.
\ep

\bdefn\label{eq:def-theta-bdr}
Let $L, C/L,X/L$ and $L'$ and $E$ be as in the previous theorem.
Write $$\thetaj^{\bdr}$$ for the image of $\thetaj^{B_E}$ in $\bdr^\ells$.
\edefn

\subparat{Non-degeneracy of $\thetaj^{B_E}$} The following lemma is useful in understanding $\thetaj^{B_E}$:
\blem 
Let $\fm_1,\fm_2,\ldots,\fm_{\ells}$ be closed maximal ideals of $B_E$. Let $$\tau:B_E^\ells\to \prod_{j=1}^\ells B_E/\fm_j$$
be the  homomorphism given using the natural surjection $\tau_j: B_E\to B_E/\fm_j$ for each $j=1,\ldots,\ells$.
Then the image $$\tau(\thetaj^{B_E})\neq(0,0,\ldots,0).$$ 
\elem
\bp 
It suffices to note that $\thetaj^{B_E}$ contains elements of the form
$$\bxi_{0,\bz,w}^{z_1,\ldots,z_\ells}=([z_1],[z_2],\ldots,[z_\ells])$$
and by \cite[Lemma 2.2.13]{fargues-fontaine}, it follows that the image of any Teichmuller lift $[z]\in W(\O_{\cpt})\tensor_{\Z_p}\O_E\subset B_E$ is non-zero in $B_E/\fm_j$ for each $j=1,\ldots,\ells$, hence $$\tau([z_1],\ldots,[z_\ells])\neq(0,0,\ldots,0)\in \prod_{j=1}^\ells B_E/\fm_j.$$
This completes the proof.
\ep

\subsubparat{The adelic ring $\Blp$ of \cite{joshi-teich-def}}
Let $\Blp$ (resp. $\Blpp$) be the adelic topological ring constructed in \constrtwoh{Definition }{def:adelic-B-ring}  considered here for the field $L'\supset L$.

Explicitly $\Blp$ is the topological ring
\be
\Blp = \prod_{w\in \vlp} B_{L'_w}.
\ee

Let me also define a variant, denoted $\Btlp$, of the topological ring $\Blp$ using the $\Q_p$-Fr\'echet algebra $\tB_{L_w'}=B\tensor_{\Q_p}L_w'$ considered in \cref{pa:ring-BE-var}. Let $\Btlp$ be the topological ring
\be  
\Btlp =   \prod_{w\in \vlp} \left(\tB_{L'_w}\right)= \prod_{w\in \vlp} \left(B\tensor_{\Q_p}{L'_w}\right).
\ee

\subparat{The definition of the adelic theta-values locus $\thetaj^{\Blp}$}
Let $L,C/L,X/L$ and $L'$ be as in \cref{ss:elliptic-curve-assumptions}, \inithtdata.
Let $\Blp$ (resp. $\Blpp$) be the adelic topological ring constructed in \constrtwoh{Definition }{def:adelic-B-ring} considered here for the field $L'\supset L$.

Explicitly $\Blp$ is the topological ring
\be
\Blp = \prod_{w\in \vlp} B_{L'_w}.
\ee

\bdefn\label{ss:adelic-theta-val-locus-in-Blp} The \textit{adelic theta-values locus  with values in $\Blp$, denoted by $\thetaj^{\Blp}$,} is the set of all   $$\bxi_{\bz}=(\bxi_{\bz,w})_{w\in\vlp} \in \Blp^\ells,$$ 
for each $\bz\in\tSigma_{L'}$. Hence 
\be 
\thetaj^{\Blp}=\left\{ (\bxi_{\bz,w})_{w\in\vlp} \in \Blp^\ells : \bz\in \tSigma_{L'}  \right\} \subset \Blp^\ells.
\ee
\edefn 

By \cref{th:local-theta-val-locus} and \eqref{eq:tilde-B-prop} one can also embed $\thetaj^{\Blp}\subset \Btlp^\ells$ and I will write
\be 
\thetaj^{\Btlp} 
\ee
for the image of $\thetaj^{\Blp}$ in  $\Btlp^\ells$.

\subsubparat{Elements $\bxi_\bz$ as Mochizuki's $\Theta$-pilot objects}\label{rem:theta-pilot-object-defined}
Let me remark that each admissible theta-values-lift $\bxi_\bz$  is the additive analog of  Mochizuki's \textit{$\Theta$-pilot object} (see \cite[Definition 3.8]{mochizuki-iut3}).  Especially, Mochizuki's theta-values set $\thetam$, according to \cite[Corollary 3.12]{mochizuki-iut3}, is  the ``holomorphic hull'' of the collection of all possible $\Theta$-pilot objects- for the definition ``holomorphic hulls'' see \cite[Remark 3.9.5(ii)]{mochizuki-iut3}.  

In \cref{se:mochizuki-construction-thetam}, I will demonstrate that each $\bxi_\bz\in\thetaj^{\Blp}$ has a Galois cohomological manifestation and this cohomological object will be a $\Theta$-Pilot in the sense of \cite{mochizuki-iut3}.

\subsubparat{Descending to $\lmod$}
Assume $L,C/L,X/L,\ell,L'$ are as in \cref{ss:elliptic-curve-assumptions}, \inithtdata. The curve $X/L$ is in fact definable over $\lmod$. Hence it will be convenient to provide a theta-values locus which does in fact live over $\lmod$. This is accomplished by means of the bijection $\V_{\lmod}\isom\ubblv\subset \vlp$. By taking the image of the theta-values locii under the projection constructed in the proof of \Cref{th:lmod-arith-descent} and the trace morphism $\Blp\to \mathbb{B}_{\lmod}$ given by taking, for each unique prime   $w\in\ubblv$ of $L'$ lying over $v\in\V_{\lmod}$ given by the bijection $\ubblv\isom \V_{\lmod}$, the trace morphism $B_{L_w'}\to B_{L_{mod,v}}$,  one obtains the image of theta-values locii in $\mathbb{B}_{\lmod}^\ells$ which I will refer to as the \textit{descent of the theta-value locii to $\lmod$}.

\section{Proof of Mochizuki's Corollary~3.12 for $\thetaj^{\Blp}$}\label{se:fundamental-estimate-joshi}
\subparat{} Now the stage is set for proving my version of the fundamental estimate known as of Mochizuki's Corollary 3.12. A prototype version of which was established in \constrtwo{Theorem }{th:main}. The version proved here works in the number field setting and exhibits all features of \moccor\ except for the fact that I work with the ring $\Blp$ defined above. Mochizuki's version works with Galois cohomology groups and will be established in \ssep\ref{se:mochizuki-construction-thetam}.

Throughout this section, the assumptions and notational conventions of \cref{ss:elliptic-curve-assumptions} and \inithtdata\ will be in force.

\subpara The following definition will be used throughout. 
\bdefn\label{def:size-of-pilot-objects}
Let $\rho\in (0,1]\subset \R$. Let $\bz\in\tSigma_{L'}$, let $$\bxi=(\bxi_{\bz,w})_{w\in\vlp}\in\thetaj^{\Blp}.$$ 
For each $w\in\vlp$, $\bxi_{\bz,w}\in B_{L_w'}^\ells$ is an $\ells$-tuple
$$\bxi_{\bz,w}=(\Xi_1,\Xi_2,\ldots,\Xi_\ells)\in B_{L_w'}^\ells.$$  Define the \textit{$(B_{L'_w},\rho)$-size of $\bxi_{\bz,w}$} to be the product
\be\abs{\bxi_{\bz,w}}_{B_{L'_w},\rho}=\abs{\Xi_1}_{B_{L'_w,\rho}}\cdot \abs{\Xi_2}_{B_{L'_w,\rho}}\cdots\cdot \abs{\Xi_\ells}_{B_{L'_w,\rho}}\in \R.\ee
be the product of absolute values with respect to the norm  $\abs{-}_{B_{L'_w},\rho}$ on $B_{L'_w}$ given by \cite[D\'efinition 1.10.2]{fargues-fontaine}.
Define the \textit{$(\Blp,\rho)$-size of $\bxi_{\bz}$} to be the product
\be \abs{\bxi_{\bz}}_{\Blp,\rho}=\prod_{w\in\vlp}\abs{\bxi_{\bz,w}}_{B_{L'_w},\rho}.
\ee
\edefn

\blem
For all $\rho\in(0,1]\subset \R$ and for all $\bz\in\tSigma_{L'}$,  
$$\abs{\bxi_{\bz}}_{\Blp,\rho}<\infty.$$
\elem
\bp 
By the definition of $\bxi_{\bz}$, and by \Cref{le:triv-lem}, for all $w\in\vlp-\ubblvoss$, $\bxi_{\bz,w}$ has absolute value equal to  one. Hence this product has a finite number of terms $\neq 1$ and so is finite.
\ep

\brem  
It will be often convenient to group the product definition $\abs{\bxi_{\bz}}_{\Blp,\rho}$ as follows 
\be
\abs{\bxi_{\bz}}_{\Blp,\rho}=\prod_{p\in\vq}\prod_{w\in\ubblvossp}\abs{\bxi_{\bz,w}}_{B_{L'_w},\rho}.
\ee
\erem

\bdefn 
Let $L,C/L,X/L$ and $L'$ be as in \cref{ss:elliptic-curve-assumptions}, \inithtdata.
For any $\rho\in(0,1]\subset \R$ let
\be 
\abs{\thetaj^{\Blp}}_{\Blp,\rho}=\sup\left\{\abs{\bxi}_{\Blp,\rho}: \bxi\in \thetaj^{\Blp}\right\}.
\ee
and let
\be 
\abs{\thetaj^{\Blp}}_{\Blp}=\sup\left\{\abs{\bxi}_{\Blp,\rho}: \bxi\in \thetaj^{\Blp}\text{ and } \rho\in(0,1]\subset \R \right\}.
\ee
\edefn

\brem 
This was considered in \cite{joshi-teich-estimates}  for a fixed $p$-adic field $E$ (and hence a fixed rational prime $p$).  The above definition obviously deals with the global case.
\erem

\subparat{The fundamental theta-values estimate} The following fundamental estimate is the global version of  \constrtwo{Theorem }{th:main} and is my version of \moccor:
\bthm\label{th:moccor} 
Let $C/L$ be an elliptic curve and suppose $\ell$ is an odd prime such that \cref{ss:elliptic-curve-assumptions}, \inithtdata\  hold for $C,L,\ell$. Let $L'/L$ be the finite extension of $L$ defined in \inithtdata. For any $w\in\ubblvoss$, let $p_w$ be the rational prime lying below $w$. Choose a normalization of the $p_w$-adic valuation on, $\C_{p_w}$, a completed algebraic closure of $\Q_{p_w}$ to be such that for any uniformizer $\pi_w\in \O_{L_w'}$ one has $$\abs{\pi_w}_{\C_{p_w}}=p_w^{-1}.$$  Let $\thetaj^{\Blp}$ be the adelic theta-values locus  for the ring $\Blp$ (defined in \cref{ss:adelic-theta-val-locus-in-Blp}). 
Then 
\[
\abs{\thetaj^{\Blp}}_{\Blp}\geq  \prod_{w\in\ubblvoss} \abs{q_{w}^{1/2\ell}}^\ells_{\C_{p_w}}
\]
\ethm
\bp  The set of such rational primes $p$ is finite (by \inithtdata, $\ubblvossp\neq\emptyset$ if and only if  $p$ is an odd rational prime such that for some $w|p$ the curve $X$ has split multiplicative reduction), hence the product in the assertion is finite.

Recall from \cref{ss:stand-point} and \cref{ss:stand-theta-link} the construction of the standard point $\bz_{\Theta}\in\tSigma_{L'}$. I will write $\bz=\bz_{\Theta}$ for brevity of notation. By construction, the set $\thetaj^{\Blp}$ contains elements of the  $\bxi_{\bz_\Theta}=\bxi_{\bz}$. More precisely,  $\thetaj^{\Blp}$ contains elements of the from
$$(\bxi^{\alpha_{w,1}, \alpha_{w,2}, \ldots,\alpha_{w,\ells}}_{\lambda,\bz,w})_{w\in\vlp}.$$
Again for brevity of notation write 
$$\bbalpha_w=(\alpha_{w,1}, \alpha_{w,2}, \ldots,\alpha_{w,\ells})$$
and 
$$\bbalpha=(\alpha_{w,1}, \alpha_{w,2}, \ldots,\alpha_{w,\ells})_{w\in\vlp}\in \prod_{w\in\vlp} \left(W(\O_{\C_{p_w}^\flat})\tensor \O_{L_w'}\right)\subset \Blp.$$
Then $\thetaj^{\Blp}$ contains elements of the form 
$$\bxi_{\lambda,\bz}^\bbalpha$$ and especially it contains elements of the form
$$\bxi_{0,\bz}^\bbalpha.$$
By the definition of $\abs{\thetaj^{\Blp}}$ one obtains
\be 
\abs{\thetaj^{\Blp}} \geq \abs{\bxi_{0,\bz}^\bbalpha}_{\Blp}.
\ee
Hence it is sufficient to find a lower bound for $\abs{\bxi_{0,\bz}^\bbalpha}_{\Blp}$ for some choice of $\bbalpha$.

As $\bxi_{0,\bz}^\bbalpha=(\bxi_{0,\bz,w}^{\bbalpha_w})_{w\in\vlp}$ and for $w\not\in\ubblvossp$ and for any $\rho\in[0,1]\subset \R$ one has 
$$\abs{\bxi_{0,\bz,w}^{\bbalpha_w}}_{B_{L_w'},\rho}=1.$$
So it will be sufficient to consider $\abs{\bxi_{0,\bz,w}^{\bbalpha_w}}_{B_{L'w}}$ for $w\in\ubblvoss$.

Let us write $$\bz=\bz_{\Theta}=(\bz_1,\bz_2,\ldots,\bz_\ells)=(y_{w,1}',y_{w,2}',\ldots,y_{w,\ells}')_{w\in\vlp}\in \left(\yadlp'\right)^\ells.$$ Then
by definition, for $w\in\ubblvoss$, $\bbalpha_w$ is an $\ells$-tuple of Teichmuller lifts of the tuple $(q_{w,j})_{j=1,\ldots,\ells}$ consisting of the Tate parameters of the $\ells$-tuple of holomorphoids $(C/L_{w}',y'_{w,j}\in \sY_{\C_{p}^\flat,L_w'})_{j=1,\ldots,\ells}$ of $C/L'$, with the Tate parameters  computed in the $\ells$-tuple of residue fields $(K_{y_{w,j}})_{j=1,\ldots,\ells}$ of $(y'_{w,j})_{j=1,\ldots,\ells}$.

Moreover $\bxi_{0,\bz,w}^\bbalpha$ runs through the local (i.e. $w$) component  $\thetaj^{B_{L_w'}}$ of the adelic theta-values locus $\thetaj^{\Blp}$ with values in $\Blp$. Hence one can use the method of proof of \constrtwo{Theorem }{th:main} for each $w\in\ubblvoss$ to establish that, for $\ell\gg 0$, one has
$$\abs{\bxi_{0,\bz,w}^{\bbalpha_w}}_{B_{L_w'},\rho}>\abs{q^{1/2\ell}_w}_{\C_{p_w}}^\ells.$$

Hence by passage to the product one obtains
$$\abs{\bxi_{0,\bz}^{\bbalpha}}_{\Blp,\rho}=\prod_{w\in \vlp} \abs{\bxi_{0,\bz,w}^{\bbalpha_w}}_{B_{L_w'},\rho} \geq \prod_{w\in\ubblvoss} \left(\abs{q^{1/2\ell}_w}_{\C_{p_w}}^\ells\right),$$
and so by definition $$\abs{\thetaj^{\Blp}}_{\Blp}\geq \prod_{w\in\ubblvoss} \left(\abs{q^{1/2\ell}_w}_{\C_{p_w}}^\ells\right).$$
This proves the assertion.
\ep

\subparat{Some corollaries of the fundamental estimates}\label{pa:precompact} Let me first make the following definition.
A subset $\Phi=\prod_{w\in\vlp}\Phi_w \subset \Blp^\ells$  is said to be \textit{a special precompact subset} if 
\benumlab
\item $\Phi$ contains the lifts of theta values arising from special point $\bz_{\Theta}$ constructed in \Cref{th:moccor}, i.e. $\bxi_{\bz_\Theta}\in\Phi$, and
\item for every $w\in\ubblvoss$, the factor $\Phi_w\subset B_{E'_w}^\ells$ is precompact i.e. the  closure of $\Phi_w$  in $B_{E'_w}^\ells$ is compact (for the Fr\'echet-topology) and for all other primes $w\in\vlp$, $\Phi_w=\{1\}$.
\eenum
\bcor With the hypothesis of the \Cref{th:moccor}.  Suppose $\Phi\subset \thetaj^{\Blp}$ is a special precompact subset of   $\thetaj^{\Blp}$ in the sense of \cref{pa:precompact},  then one has 
\[ 
\infty >\abs{\Phi}_{\Blp} >\prod_{w\in\ubblvossp\neq \emptyset} \abs{q_{w}^{1/2\ell}}^\ells_{\C_{p_w}}.
\]
\ecor
\bp
Let $\Phi\subset\thetaj^{\Blp}$ be a special precompact subset. As any finite extension of $\Q_p$ is locally compact, special precompactness of $\Phi$ is equivalent to $\Phi_w$ being a bounded subset of $B_{E'_w}^\ells$ for all $w\in\vlp$ (for the Fr\'echet topology) and as the set of $\ubblvoss$ is finite, so $\sup_{0<\rho<1}\abs{\Phi}_{\Blp,\rho}$ is evidently  bounded and hence the second assertion is established.
\ep

\bcor 
In the notation of \Cref{th:moccor}, and with $\rho=1$, one has
\[
1\geq \abs{\thetaj^{\Blp}}_{\Blp,1}\geq  \prod_{w\in\ubblvoss} \abs{q_{w}^{1/2\ell}}_{\C_{p_w}}.
\]
\ecor
\bp 
The upper bound follows from the fact that for each $w\in\ubblvoss$, and $\rho=1$, one can apply \constrtwo{Proposition }{pr:trivial-upper-bound-on-ttheta}. The lower bound is \Cref{th:moccor}.
\ep

\subparat{Prime Bundling rings $\breveB_{L',p}$ and $\breveBb_{L',p}$}\label{ss:bundling-rings}
One important aspect of working with theta-value locii with values in a ring is that one can use the multiplicative structure of the ring to multiply theta-values arising from distinct arithmeticoids. For example in $\Blp$ it makes perfect sense to talk about 
$$\bxi_{\bz}\cdot \bxi_{\bz'} \in \Blp,$$
as well as 
$$\bxi_{\bz}+ \bxi_{\bz'} \in \Blp.$$
Similarly one can also multiply  the local components at two distinct primes $w_1,w_2\in\ubblvossp$ lying over a fixed rational prime:
$$(\bxi_{\bz,w_1}\tensor 1)\cdot(1\tensor\bxi_{\bz,w_2})\in B_{L_{w_1}'}\tensor_{\Q_p}B_{L_{w_2}'},$$ 
as well as
$$\bxi_{\bz,w_1}\tensor 1+1\tensor\bxi_{\bz,w_2}\in B_{L_{w_1}'}\tensor_{\Q_p}B_{L_{w_2}'}.$$
This observation allows me to demonstrate the tensor packet structure of \cite[Section 3]{mochizuki-iut3} which plays in the main theorem of \cite{mochizuki-iut4}. I call this \textit{bundling of theta-values}.  The purpose of this bundling is to arrive at a minimal theta-values locus (this becomes important in providing tight bounds) and this together with the fact that $\underline{\V}\isom \V_{\lmod}$, one can expect to descend these data back to $\lmod$.

In this section I will describe this bundling structure explicitly. One needs to work with several Fargues-Fontaine rings $B_E$ (for varying $E$) simultaneously and leads to some auxiliary rings which I call \textit{bundling rings}. The definition of these bundling rings is the main content of this section.

The notations and assumptions of \cref{ss:elliptic-curve-assumptions}, \inithtdata\ will strictly be in force. Following \iutthr, I will work with all the primes $w|p$ in $\ubblv$ simultaneously.

\subsubparat{More bundling rings} For simplicity, I will write 
\begin{align}
\vlpp&=\{w\in\vlp: w|p \},\\
\ubblv_p&=\left\{w\in\ubblv: w|p \right\}, \text{ and}\\
\ubblvossp&=\left\{w\in\ubblvoss: w|p \right\}.
\end{align}
 The local components, $\left\{\thetaj^{B_{E'_w}}:w\in\ubblvossp\right\}$, of the theta-torsion values locus $\thetaj^{\Blp}$,  can be viewed as subsets of the following direct product of rings:
\be 
\left(\bigoplus_{w\in\ubblvossp} B_{E_w'}\right)^\ells=\prod_{j=1}^\ells \left(\bigoplus_{w\in\ubblvossp}B_{E_w'}\right)
\ee 
where for notational continuity I will write $E_w'=L'_w$. 

More generally, let
\be 
\breveB_{L',p}=\left(\bigoplus_{w\in\vlpp}B_{E_w'}\right).
\ee

The right hand side obviously depends on $L'$ and hence the notation includes $L'$ dependence. If $L'$ is not in doubt, I will simply write $\breveB_p=\breveB_{L',p}$.

Then one can write the ring $\Blp$ as 
\be\Blp=\prod_{p\in\vq} \breveB_{L',p}.\ee

\brem For $w|p$, one can think of elements of  $\breveB_p^\ells$ as $\ells$-tuples of elements of $\breveB_p$ indexed by $j=1,\ldots,\ells$. Each element of $\breveB_p$ is an element of the $B_p$-algebra in the parenthesis and so one can think of elements of $\breveB_p$ as having coordinates indexed by  the set
\be 
\left\{E_{w}': w\in\vlpp \right\}.
\ee 
\erem

\blem 
For each $p$, $\breveB_{L',p}$ is a finite module over $B_p$.
\elem
\bp 
For $w|p$,  one has $B_{E_w'}=B_p\tensor_{E_{w,0}'}E_w'$ where $E_{w,0}'\subset E_w'$ is the maximal absolutely unramified subfield of $E_w'$.  So $\breveB_{L',p}$ is a $B_p$-algebra through this identification and for each $w|p$ this $B_{E_w'}$ is a finite module over $B_p$ and so the assertion follows.
\ep

\subsubparat{The ring $\breveBb_p$} To understand my next definition, let me provide some motivation. Suppose for some prime $p$, each $w\in\vlpp$ the absolute residue degree, $f(w|p)$, of $E_w'/\Q_p$ is one i.e. $f(w|p)=1$ for all $w\in\vlpp$. Then the maximal unramified subfield $E_{w,0}'\subset E_w'$ satisfies
$$E_{w,0}'=\Q_p \text{ for all } w\in\vlpp.$$
Then we may write 
$$\breveB_{L',p}=B\tensor_{\Q_p}\left(\bigoplus_{w\in \vlpp}E_w'\right).$$

In general this condition on the absolute residue degree does not hold. This presents some difficulty in understanding $\breveB_{L',p}$ transparently and motivates my next definition.

For simplicity let me write:
\be\label{eq:def-breveBB} 
\breveBb_p=\bigoplus_{w\in\ubblvossp} \left(B\tensor_{\Q_p}{E_w'}\right),
\ee
Then using the fact that $B_E\into B\tensor_{\Q_p}E$ one has
\be\label{eq:def-breveB-extended} 
\breveB_p^\ells=\left(\bigoplus_{w\in\ubblvossp} B_{E_w'}\right)^\ells\into \breveBb_p^\ells
\ee
The theta-values locus obviously embeds in $\breveBb_p^\ells$ by means of the diagonal embedding. 

\blem\label{le:breveB-extended-simplified}
One has the simplification (this is  the reason for working with $\breveBb_p$ as opposed to $\breveB_p$):
$$
\breveBb_p=B\tensor_{\Q_p}\left( \bigoplus_{w\in\ubblvossp} {E_w'}\right)
$$
\elem
\bp
The proof is immediate from the definition \eqref{eq:def-breveBB}.
\ep

\subsubparat{Bundling theta-values: the locus $\thetajp^{\breveB_p}$} Now let 
\be\thetajp^{\breveB_p}=\bigoplus_{w\in \ubblvossp} \thetajp^{E_w'}\into \breveB_p=\bigoplus_{w\in\ubblvossp} B_{E_w'}\ee
Thus $\thetajp^{\breveB_p}$ assembles together  the collection of all possible lifts (as constrained in \constrtwo{\ssep}{se:construction-ttheta}) viewed as lifts in $B_{E_{w}'} \subset{\breveB_p}$  for each $w\in\ubblvossp$.

\subsubparat{theta-values locii in variants of the $B$-rings}  For each prime $p$, one can construct bundling rings using other rings considered here:
$$B_p, B_p[1/t], B_{e}, \breveB_p,  \breveB_p[1/t], \breveBb_p, \breveBb_p[1/t], \bdr,$$
 and hence bundled sets such as $\thetajp^{\breveB_p}$ can be constructed using  these rings.  
I will indicate the choice of the ring for the construction by using the notation for the ring in the superscript, for example,  as $$\thetajp^{B_p}, \thetajp^{\breveB_p} \ etc.$$ 

Thus one has constructed the theta-values locii
\be 
\thetajp^{B},\thetajp^{B[1/t]},\thetajp^{\breveB}, \thetajp^{\breveB[1/t]},\thetajp^{\breveBb},\thetajp^{\breveBb[1/t]}, \thetajp^{B_e},\cdots
\ee

So my point is that the locii of the sort constructed by Mochizuki can be constructed as soon as one has the Arithmetic Teichmuller spaces of \cite{joshi-teich}, \cite{joshi-teich-estimates} and the present paper.
At any rate one has  proved the following theorem:
\bthm\label{th:theta-values-loci-in-B-rings-and-adelic-version} 
Fix an elliptic curve $C/L$ as in \cref{ss:theta-data-fixing}. 
Then using the choice of the theta function as in \cite[Proposition 1.6]{mochizuki-theta} or \constrtwo{\ssep}{se:theta-values}, \constrtwo{\ssep}{se:ansatz} at each $w|p$ and a choice of a ring, say $B_?$, from the set $$\left\{B=B_p, B[1/t],  \breveB,  \breveB[1/t], \breveBb, \breveBb[1/t], B_e, \bdr \right\}\bigcup \left\{B_I: I\subset [0,1]\subset \R{\vphantom{\breveBb[1/t]}}\right\}$$ and the corresponding bundling ring, one has a naturally defined theta-values locus viewed in bundling rings (stated here for $B_E$):
$$\thetajp^\breveB\subset \breveB^\ells_p.$$
By passage to the product over all primes $p$:
$$\thetaj^{\Blp}=\prod_p\thetajp\subset \prod_{p} \breveB_p^\ells=\Blp^\ells.$$
 Using the diagonal inclusion $\breveB_p\into \breveBb_p^\ells$ considered above, one can consider the theta-value locus $$\thetaj^{\breveBb}\subset \prod_p \breveBb_p^\ells.$$
\ethm

\subparat{Norms on tensor product of Banach spaces}
Familiarity with the archimedean theory of tensor products of Banach spaces is useful but not essential--for the archimedean theory see \cite{kaniuth-banach-algebras-book} or \cite{diestel-grothendieck-tensor-products-book}. I will use the construction of tensor product of Banach spaces over $\Q_p$ for which \cite[Chap IV, \ssep 17]{schneider-book} is a convenient reference. To apply these results (of \cite{schneider-book}) one needs to remember that the Banach spaces $B_E$ over $E$ and as the $p$-adic field $E$ is discretely valued, $E$ is a spherically complete field. 

\subsubparat{Basic examples} Let $\{V_\alpha\}$ be a finite multi-set of (1) finite dimensional $\Q_p$ vector spaces or (2) $\Q_p$-Banach spaces of infinite dimensions. Two collections of such will be important to us. The first  is the collection of $p$-adic fields which occur above
\be\{V_\alpha \}=\{E_w':w\in\ubblvoss, w|p\}\ee
viewed as a collection of finite dimensional $\Q_p$ vector spaces which are complete with respect to their usual $p$-adic norms $\abs{-}_{E_w'}$ and hence are Banach spaces with these norms. Second is the  collection  
\be\{V_\alpha \}=\{B_{E_w'}:w\in\ubblvoss, w|p\}\ee
of infinite dimensional Fr\'echet algebras over $\Q_p$ equipped with Fr\'echet norms $\abs{-}_\rho$ for every $\rho\in[0,1]\subset \R$ constructed in \cite{fargues-fontaine} and hence can be treated as  Banach spaces with respect to any of the Fr\'echet norms.

\subsubparat{Tensor products and the cross-norm property}\label{ss:cross-norm-prop} Let $E$ be a $p$-adic field for rational prime $p$  or let  $E=\R$ or $\C$ i.e. the archimedean case is allowed in this discussion.

In general, there are many distinct tensor products of a pair of Banach spaces over $E$. Let $(V,\abs{-}_V),(W,\abs{-}_W)$ be Banach-spaces over $E$. Suppose one has a tensor product $(V\tensor_E W,\abs{-}_{V\tensor_W})$ of $(V,\abs{-}_V)$ and $(W,\abs{-}_W)$.  Then the following property, if true, of this tensor product will be called the \textit{cross-norm property} (of this tensor product):
\be\label{eq:cross-norm-prop0} 
\abs{v\tensor w}_{V\tensor_{E} W}=\abs{v}_V\cdot \abs{w}_W\qquad \forall v\in V,w\in W.
\ee

In the archimedean case, by \cite{grothendieck-memoir} or \cite{diestel-grothendieck-tensor-products-book}, any tensor product of Banach spaces (over $E=\R,\C$) sandwiched between the injective tensor product and the projective tensor product has the cross-norm property \eqref{eq:cross-norm-prop0}. More precisely, in the archimedean case, by  \cite{grothendieck-memoir} or \cite{diestel-grothendieck-tensor-products-book}, the projective tensor product is the greatest cross-norm, while the injective tensor product is the smallest cross-norm. 

Now assume $E$ is non-archimedean. I will not recall here the definition of the projective tensor product $\abs{-}_{V\tensor_{E} W}$ given in \cite[Chapter 2, 2.1.7]{bosch-non-arch-analysis} or \cite[Chapter 4, \ssep 17]{schneider-book}. 

Then one has cross-norm property \eqref{eq:cross-norm-prop0} for the projective tensor Banach spaces over the field $E$ (\cite[Proposition 17.4]{schneider-book}). Note that since the $p$-adic field $E$ is complete and discretely valued field, so it is  spherically complete \cite[\ssep 1, Lemma 1.6]{schneider-book}. 

As a consequence of spherical completeness of $E$, one has by \cite[Theorem 1]{vanderput1967} that the injective and projective tensor products coincide  and so one may speak of the  unique tensor product. [One of the difficulties of working with fields which are not spherically complete is that a Banach space may have no non-zero continuous linear forms on it; on the other hand, from the point of view of $p$-adic integration (and measures), the existence of linear forms is desirable in any case.] Moreover under spherical completeness of $E$, one has additionally the nice properties \cite[Proposition 17.8, 17.9]{schneider-book} of the tensor product. This remark applies to  all  Banach spaces over $p$-adic fields considered here and the following summarizes this discussion:
\bthm\label{th:cross-norm-prop} 
Let $E$ be a $p$-adic field or let $E=\R,\C$.  
\benumlab
\item If $E=\R,\C$, then the projective and the injective tensor products have the cross-norm property \eqref{eq:cross-norm-prop0} (and so does any norm sandwiched between these two norms) (\cite{grothendieck-memoir} or \cite[Theorem 1.1.3 ]{diestel-grothendieck-tensor-products-book}). 
\item If $E$ is a $p$-adic field, then projective tensor product and the injective tensor products coincide (and one can speak of the tensor product of Banach spaces over a $p$-adic field $E$) (\cite[Theorem 1]{vanderput1967}) and,
\item  the cross-norm property
 \eqref{eq:cross-norm-prop0} holds for  any pair of Banach spaces over any $p$-adic field $E$ . 
 \item Especially, for $E=\Q_p$, one has
$$
\abs{v\tensor w}_{V\tensor_{\Q_p} W}=\abs{v}_V\cdot \abs{w}_W.
$$
\item In particular, for $p$-adic fields $E_1,E_2\supseteq E$, the tensor product $B_{E_1}{\tensor}_{E} B_{E_2}$ has the  cross-norm property for any choice of norms on each of its factors.
\eenum
\ethm

\subsubparat{Passage to completed tensor products} While passage to completed tensor product is unnecessary for what I will do here, it is worth remarking that by \cite[Chapter 4, Proposition 17.6]{schneider-book}, if $V,W$ are Fr\'echet spaces over $\Q_p$ then completed tensor product $V\widehat{\tensor}_{\Q_p} W$ is a Fr\'echet space. In particular, for $p$-adic fields $E_1,E_2$, the tensor product $B_{E_1}\widehat{\tensor}_{\Q_p} B_{E_2}$ is a Fr\'echet space.

\subsubparat{Tensor products again} The discussion of this paragraph provides an analog of the above equality for the example of Banach algebras  \[\{B_{E}:E\text{ a p-adic  field}\}\]
equipped with their multiplicative norms i.e. valuations $$\{ \abs{-}_{B_E;\rho}\} \text{ for } \rho\in [0,1]\subset \R$$
constructed in \cite[Chapitre 2]{fargues-fontaine}. The following lemma is a standard consequence of valuation theory and is easily established:
\blem\label{le:tensor-product-lemma-1} 
Let $E_1,E_2$ be two $p$-adic fields, let $\rho_1,\rho_2\in [0,1]\subset \R$. Write $B_1=B_{E_1},B_2=B_{E_2}$ and equip them with the multiplicative norms $\abs{-}_{B_i;\rho_i}$ for $i=1,2$. Then the algebra $B_1\tensor_{\Q_p} B_2$ is equipped with a valuation i.e. a multiplicative norm $\abs{-}_{B_1\tensor B_2;\rho_1,\rho_2}$ such that for every $b_1\in B_1$ and for every $b_2\in B_2$ one has
\begin{multline*}
\abs{b_1\tensor b_2}_{B_1\tensor B_2;\rho_1,\rho_2}=\abs{(b_1\tensor1)\cdot(1\tensor b_2)}_{B_1\tensor B_2;\rho_1,\rho_2}=\\ \qquad\qquad =\abs{(b_1\tensor1)}_{B_1\tensor B_2;\rho_1,\rho_2}\cdot \abs{(1\tensor b_2)}_{B_1\tensor B_2;\rho_1,\rho_2}=\abs{b_1}_{B_1;\rho_1}\abs{b_2}_{B_2;\rho_2}.
\end{multline*}
especially
\[ 
\abs{b_1\tensor b_2}_{B_1\tensor B_2;\rho_1,\rho_2}=\abs{b_1}_{B_1;\rho_1}\abs{b_2}_{B_2;\rho_2}.
\]
\elem
I will refer to $\abs{-}_{B_1\tensor B_2;\rho_1,\rho_2}$ as the \textit{tensor product norm} of $(B_1,\abs{-}_{B_1;\rho_1})$ and $(B_2,\abs{-}_{B_2;\rho_2})$.

\subsubparat{Another variant} Here is a Banach theoretic version of the above lemma which I could not find in existing literature and of which I will provide a proof.
\blem\label{le:tensor-product-lemma-2} 
Let $E_1,E_2$ be two $p$-adic fields, let $\rho_1,\rho_2\in [0,1]\subset \R$. Write $B_1=B_{E_1},B_2=B_{E_2}$ and equip them with the multiplicative norms $\abs{-}_{B_i;\rho_i}$ for $i=1,2$. Then the algebra $B_1\tensor_{\Q_p} B_2$ is equipped with  a multiplicative norm $\abs{-}_{B_1\tensor B_2;\rho_1,\rho_2}$ such that for every $b_1\in B_1$ and for every $b_2\in B_2$ one has
\begin{multline*}
\abs{b_1\tensor b_2}_{B_1\tensor B_2;\rho_1,\rho_2}=\abs{(b_1\tensor1)\cdot(1\tensor b_2)}_{B_1\tensor B_2;\rho_1,\rho_2}=\\ \qquad\qquad =\abs{(b_1\tensor1)}_{B_1\tensor B_2;\rho_1,\rho_2}\cdot \abs{(1\tensor b_2)}_{B_1\tensor B_2;\rho_1,\rho_2}=\abs{b_1}_{B_1;\rho_1}\abs{b_2}_{B_2;\rho_2}.
\end{multline*}
especially
\[ 
\abs{b_1\tensor b_2}_{B_1\tensor B_2;\rho_1,\rho_2}=\abs{b_1}_{B_1;\rho_1}\abs{b_2}_{B_2;\rho_2}.
\]
\elem
\bp 
By \cite{lang-algebra}, there is a natural way to define an algebra structure on the tensor product $B_1\tensor_{\Q_p}B_2$ compatible with the algebra structures of each of the factors, namely
$$(b_1\tensor b_2)\cdot (c_1\tensor c_2)=(b_1c_1)\tensor (b_2c_2).$$
Then $$\abs{(b_1\tensor b_2)\cdot (c_1\tensor c_2)}_{B_1\tensor_{\Q_p}B_2}=\abs{(b_1c_1)\tensor (b_2c_2)}_{B_1\tensor_{\Q_p}B_2}$$ and by \cite[Chap 4, Proposition 17.4]{schneider-book}
$$\abs{(b_1c_1)\tensor (b_2c_2)}=\abs{b_1c_1}_{B_1}\abs{b_2c_2}_{B_2}$$
and as the norms on $B_1,B_2$ are multiplicative, one gets 
$$\abs{(b_1c_1)\tensor (b_2c_2)}=\abs{b_1c_1}_{B_1}\abs{b_2c_2}_{B_2}=\abs{b_1}\abs{c_1}\abs{b_2}\abs{c_2}=\abs{b_1\tensor c_1}\abs{b_2\tensor c_2}$$
and therefore
$$\abs{(b_1c_1)\tensor (b_2c_2)}=\abs{b_1\tensor c_1}\abs{b_2\tensor c_2},$$
and the assertion of the lemma is immediate from this formula by considering the special case
$(b_1\tensor 1)(1\tensor b_2)=b_1\tensor b_2$.  This proves the assertion.
\ep

\subsubparat{The lower bound lemma} The proof of the   following lemma is now self-evident:
\blem\label{le:lower-bound-lemma} 
Suppose one is given a finite collection of $\Q_p$-Banach spaces $\{V_\alpha \}_{\alpha}$ and suppose that one is also given an element $m_\alpha\in V_{\alpha}$ for each $V_{\alpha}$. Then for the projective tensor product norm on $\bigotimes_\alpha V_\alpha$, one has
\be 
\abs{\bigotimes_{\alpha} m_\alpha}_{\otimes_\alpha V_\alpha} = \prod_{\alpha}\abs{m_\alpha}_{V_\alpha}.
\ee
In particular, if $S\subset \bigotimes_{\alpha} V_\alpha$ is a subset containing some pure tensor $\bigotimes_\alpha m_\alpha$ then one has
\be  
\sup\left\{ \abs{s}_{\otimes_\alpha V_\alpha}: s\in S \right\} \geq \prod_{\alpha}\abs{m_\alpha}_{V_\alpha}.
\ee
\elem

\subsubparat{Adelic considerations} The above consideration will now be adelized as follows. Suppose for every prime $p$ one has a finite collection of $\Q_p$-Banach spaces $\{V_{\alpha,p}\}$ (the indexing set for the $\alpha$ may depend on $p$ but I am suppressing the notation for the indexing set) and suppose $S_p\subseteq \bigotimes V_{\alpha,p}$ is a subset. Then consider the product
\[
S=\prod_p S_p\subset \prod_p\bigotimes_{\alpha} V_{\alpha,p}.
\]
Define 
\[
\abs{S_p}=\sup\{ \abs{s_p}_{\tensor_\alpha V_{\alpha,p}}: {s_p\in S_p}\},
\]
and define (with the convention that $\infty\cdot \infty=\infty, 0\cdot\infty=\infty$)
\[
\abs{S}=\prod_p\abs{S_p}\in \R\cup\{\infty\}.
\]

\subsubparat{Basic examples} The principal cases of interest here are \[ 
S_p=\thetajp\] and \[S_p=\thetamp\] for all rational primes $p$ lying below the primes in $\ubblvoss$ and for the remaining primes $S_p=\{1\}$ or $S_p$ is the unit ball in the $p$-component. In this case, the sets $S_p$,  satisfy the additional condition that
$\abs{S_p}=1$   for all but finitely many primes $p$  and for the remaining finitely many primes, one expects $\abs{S_p}<\infty$. So $\abs{S}<\infty$ in the principal case of interest.

\subsubparat{Sets of adelic and bounded adelic types}  The above example  leads to the following definition. 
\benumlab
\item I will say that $S$ as above is a \textit{set of adelic type} if
$\abs{S_p}=1$  for all but finitely many primes $p$, and \\
\item I will say that $S$ is a \textit{set of bounded adelic type} if $S$ is a set of adelic type and $\abs{S_p}<\infty$ for every prime $p$. 
\eenum

The following lemma is a consequence of the  construction of the theta-values sets:
\blem 
The set $\thetaj^{\Blp}$ and its  variants for different choices of the $B$-rings are of adelic type.
\elem

\subparat{The fundamental estimate for $\thetaj^{\breveBten}$ and $\thetaj^{\breveBbten}$}
\subsubsection{The tensor product rings $\breveBten$ and $\breveBbten$}\nwsss Let me now provide a proof of the fundamental estimate using the above detailed theory of tensor products. 
To keep the discussion simple, I will work with tensor product of $\Q_p$-F\'rechet  vectors spaces and, for any rational prime $p$ such that $\ubblvossp\neq\emptyset$, introduce  the rings \be\breveBten_{L',p}=\bigotimes_{w\in\ubblvossp\neq\emptyset}B_{E_w'},\ee
and the ring
\be\breveBbten_{L',p}=\bigotimes_{w\in\ubblvossp\neq\emptyset}(B_p\tensor_{\Q_{p}}{E_w'}).\ee
\textcolor{red}{Strictly speaking, I should work  with completed tensor product $\hat{\otimes}$, but I will avoid doing this for simplicity.}  

\subsubsection{The theta-values locii $\thetaj^{\breveBten}$ and $\thetaj^{\breveBbten}$}\nwsss For each $p$ such that $\ubblvossp\neq\emptyset$, one has the homomorphism (given by the construction of tensor products):
$$\prod_{w\in\ubblvossp} B_{E_w'}\to \breveBten_{L',p}$$
which maps an element $(x_w)_{w\in\ubblvossp}\mapsto \bigotimes_{w\in\ubblvossp} x_w$
and 
$$\prod_{w\in\ubblvossp} (B_p\tensor_{\Q_p}{E_w'})\to \breveBbten_{L',p}$$
which maps an element $(x_w)_{w\in\ubblvossp}\mapsto \bigotimes_{w\in\ubblvossp} x_w$.
This can be applied to lifts of theta-values in $\prod_{w\in\ubblvossp} B_{E_w'}$ and $\prod_{w\in\ubblvossp} (B_p\tensor_{\Q_p}{E_w'})$ respectively. In particular one may take the image of the theta-values locii $\prod_{w\in\ubblvossp}\thetaj^{B_{E_w'}}$ in $\left(\breveBten_{L',p}\right)^\ells$ and similarly in $\left(\breveBbten_{L',p}\right)^\ells$. This gives us theta-values locii 
\be \thetaj^{\breveBten}\subset \prod_{p,\ubblvossp\neq\emptyset}\left(\breveBten_{L',p}\right)^\ells\ee and 
\be \thetaj^{\breveBbten}\subset \prod_{p,\ubblvossp\neq\emptyset}\left(\breveBbten_{L',p}\right)^\ells \ee respectively. Note that both the products are over a finite set of rational primes $p$. One may also extend the definition of these rings primes lying over all rational primes. 

\subsubsection{The fundamental estimate}\nwsss
\bthm\label{th:moccor-tensor-product1} 
Let $C/L$ be an elliptic curve and suppose $\ell$ is an odd prime such that \cref{ss:elliptic-curve-assumptions}, \inithtdata\  hold for $C,L,\ell$. Let $L'/L$ be the finite extension of $L$ defined in \inithtdata. For any $w\in\ubblvoss$, let $p_w$ be the rational prime lying below $w$. Choose a normalization of the $p_w$-adic valuation on, $\C_{p_w}$, a completed algebraic closure of $\Q_{p_w}$ to be such that for any uniformizer $\pi_w\in \O_{L_w'}$ one has $$\abs{\pi_w}_{\C_{p_w}}=p_w^{-1}.$$  Let $\thetaj^{\Blp}$ be the adelic theta-values locus  for the ring $\Blp$ (defined in \cref{ss:adelic-theta-val-locus-in-Blp}). 
Then 
\[
\abs{\thetaj^{\breveBten}}_{\breveBten}\geq  \prod_{w\in\ubblvoss} \abs{q_{w}^{1/2\ell}}^\ells_{\C_{p_w}},
\]
and a similar assertion for $\abs{\thetaj^{\breveBbten}}_{\breveBbten}$.
\ethm

\bp
By \constrtwo{Theorem }{th:main}, for each $\ubblvossp\neq\emptyset$, the set $\thetaj^{\breveB}{}_{,p}$ contains a pure tensor $\tensor_{w|p}Z_w$ such that $$\abs{Z_w}_{B_{E'_w},\rho}\geq \abs{q_{w}^{1/2\ell}}_{\C_{p_w}},$$ with $Z_w$ a theta-Pilot object for $X/E_w'$ as defined earlier  and by \constrtwo{Theorem }{th:main} and its proof.

Now the properties of tensor products $\tensor_{w|p}B_{E'_w}$ established in \Cref{le:tensor-product-lemma-1}, \Cref{le:tensor-product-lemma-2}  show that
\be\abs{\tensor_{w |p}Z_w}_{B_{E'_w},\rho}=\prod_{w |p}\abs{Z_w}_{B_{E'_w},\rho}.\ee

By \constrtwo{Theorem }{th:main}, for each $0<\rho<1$ and each $p$ for which $\ubblvossp\neq\emptyset$ and each $w\in\ubblvossp$  one has the lower bound 
\be 
\abs{Z_w}_{B_{E'_w},\rho}\geq \abs{q_{w}^{1/2\ell}}_{\C_{p_w}}.
\ee

Thus  by taking products over all $p$ for which $\ubblvossp\neq\emptyset$  and the extending the product (trivially) to all $p$ and all primes $w|p$ in $L'$  one gets:
\be\prod_{p}\abs{\tensor_{w |p}Z_w}_{B_{E'_w},\rho}=\prod_{p}\prod_{w\in\ubblvossp\neq\emptyset}\abs{Z_w}_{B_{E'_w},\rho}\geq \prod_{p}\prod_{w\in\ubblvossp\neq\emptyset}\abs{q_{w}^{1/2\ell}}_{\C_{p_w}}.\ee

Now the claim of Theorem~\ref{th:moccor} follows from \Cref{le:lower-bound-lemma}. This proves the assertion.
\ep

\subparat{The theta-values locii in  $\breveBb_p\tensor_{B_p}\breveBb_p\tensor_{B_p} \cdots \tensor_{B_p} \breveBb_p$ and in $\breveBb_p$}
One may further exploit multiplicative properties one has available in the situation at hand but which is not utilized in the construction of $\thetaj^{\breveBten}\subset \left(\breveBbten\right)^\ells$ given above. This construction does not explicitly demonstrate the tensor packet structure which is asserted in \cite[Section 3]{mochizuki-iut3}. I will now demonstrate this as follows.

\subsubparat{Multiplication of lifts of theta-values} 
Since $\breveBb_p$ is a ring, one can also multiply elements in it. More precisely given a tuple of elements  $(\breve{b}_1,\ldots,\breve{b}_\ells)\in \breveBb_p^\ells$ one can also multiply them in $\breveBb_p$ i.e. consider the mapping
$$\breveBb_p^\ells \to \breveBb_p,$$
given by
$$(\breve{b}_1,\ldots,\breve{b}_\ells)\mapsto \breve{b}_1\cdot \breve{b}_2\cdots \breve{b}_\ells.$$
This is multi-linear in each coordinate. This multi-linear mapping factors though the tensor product. In this case it is convenient to treat  $\breveBb_p$ as a finite $B_p$-algebra and consider the natural factorization 
\be 
\begin{tikzcd}
	\breveBb_p^\ells \ar[r]\ar[dr]& \breveBb_p\tensor_{B_p}\breveBb_p\tensor_{B_p} \cdots \tensor_{B_p} \breveBb_p\ar[d]\\
	&   \breveBb_p.
\end{tikzcd}
\ee
where the first mapping is $(b_1,\ldots, b_\ells)\mapsto b_1\tensor b_2\tensor \cdots\tensor b_\ells$ and the second is the natural mapping $b_1\tensor b_2\tensor \cdots\tensor b_\ells\mapsto b_1\cdot b_2\cdot \cdots\cdot b_\ells$.

The advantage of these considerations becomes clear from the following elementary decomposition lemma:

\blem\label{le:direct-sum-decomp2} One has the decomposition of the $B_p$-algebra: 
$$ \breveBb_p\tensor_{B_p}\breveBb_p\tensor_{B_p} \cdots \tensor_{B_p} \breveBb_p=B_p\tensor_{\Q_p} \left(\bigoplus_\alpha E_\alpha''\right)
$$
where the sum on the right runs over a finite multi-set of finite extensions $E_\alpha''$ of $\Q_p$ contained in $\bQ_p$. 
\elem
\bp 
This is immediate from the definition of $\breveBb_p$ and properties of tensor products. Indeed, since 
\be
\breveBb_p=B_{\cpt,\Q_p}\tensor_{\Q_p} \left(\bigoplus_{w|p,w\in \bbvlmod}E_w\right)
\ee 
and as the tensor product $\breveBb_p\tensor_{B_p}\breveBb_p$ reduces to computing tensor products of a finite collection of $p$-adic fields,
and so the assertion of the lemma follows.
\ep

\brem
In \cite[Section 3]{mochizuki-iut3}, especially  \cite[Remark 3.1.1, Remark 3.9.1]{mochizuki-iut3},  and \cite{mochizuki-iut4}, Mochizuki forces the tensor product structure and arrives essentially at the above lemma.
Especially note that appearance of finite number of $p$-adic fields $E_\alpha''$ in \Cref{le:direct-sum-decomp2} and (hence also in \Cref{th:theta-value-locus-tensor-product} below) should be compared to \cite[\ssep 3, Remark 3.1.1]{mochizuki-iut3} and \cite{mochizuki-iut4}.
\erem

\subsubsection{} Now one may define theta-locii with the above structures as follows. Let 
\be
\thetajpp \subset \breveBb_p\tensor_{B_p}\breveBb_p\tensor_{B_p} \cdots \tensor_{B_p} \breveBb_p\isom  B_p\tensor_{\Q_p} \left(\bigoplus_\alpha E_\alpha''\right)
\ee
be the image of 
\be \thetajp\subset \breveBb_p^\ells \to \breveBb_p\tensor_{B_p}\breveBb_p\tensor_{B_p} \cdots \tensor_{B_p} \breveBb_p
\ee under the above mapping and let
\be
\thetajpph \subset \breveBb_p
\ee
be the image of 
\be \thetajp\subset \breveBb_p^\ells \to  \breveBb_p
\ee under the multiplication mapping.
In particular, by \Cref{le:direct-sum-decomp2} one has
\be 
\thetajpp\subset B_p\tensor_{\Q_p} \left(\bigoplus_\alpha E_\alpha''\right),
\ee
and $\thetajpph$ is the image of $\thetajpp$ under the vertical arrow in the above diagram.

\subsubpara This leads to the following adelic theta-values locus which should be compared with the appearance of the decomposition in \cite{mochizuki-iut3}: except for $B_p$ factor on the right, the locus in \moccor\ is of this type:

\bthm\label{th:theta-value-locus-tensor-product} The adelic theta-values locus for the ring $\breveBb_p$ (for each prime $p$), is denoted by $\thetajh^{\breveBb}$ and is the naturally defined locus constructed above: 
$$
{\thetajh}^\breveBb=\prod_p \thetajpp^{\breveBb} \subset \prod_p \left(\breveBb_p\tensor_{B_p}\breveBb_p\tensor_{B_p} \cdots \tensor_{B_p} \breveBb_p\right) = \prod_p \left(B_p\tensor_{\Q_p}\left(\bigoplus_\alpha E_{p,\alpha}''\right)\right),
$$
where the finite collections of fields $E_{p,\alpha}''$ appearing on the extreme right is the collection of fields appearing in \Cref{le:direct-sum-decomp2}. Similarly, one has an adelic theta-values locus $\thetajppa$ using the ring $\breveBb_p$ (defined above):
$$ 
{\thetajppa}^{\breveBb}\subset \prod_p \breveBb_p.
$$
One also has lower bounds on these locii given by \Cref{th:moccor-tensor-product1}.
\ethm

\brem 
In the context of \moccor\ and \constrtwo{Theorem }{th:main} one wants to work with a collection of norms given using $B$ on $\breveBb$ and one wants to evaluate $$\sup\left\{\abs{z}_{\breveBb}:z\in{\thetajh}^\breveBb  \right\}$$ and also $$\sup\left\{\abs{z}_{\breveBb}:z\in{\thetajppa}^\breveBb  \right\},$$
since ${\thetajppa}^\breveBb$ is the image ${\thetajh}^\breveBb$ so one has
$$\sup\left\{\abs{z}_{\breveBb}:z\in{\thetajh}^\breveBb  \right\}\geq \sup\left\{\abs{z}_{\breveBb}:z\in{\thetajppa}^\breveBb  \right\}.$$
Notably, since $$\abs{v\tensor w}_{V\tensor W}=\abs{v}_V\cdot\abs{w}_W$$
an element of $z\in{\thetajh}^\breveBb$ also provides a lower bound on the supremum for ${\thetajppa}^\breveBb$.

On the other hand, the above inequality suggests that the passage to the tensor product version ${\thetajppa}^\breveBb$ should be expected to provide tighter upper bounds! 
\erem

\subparat{theta-values locii and compactly bounded subsets of $\sM_{1,1}(\bQ)$} This paragraph is not needed in the main proofs. Let me remark that in the above constructions, and in \Cref{th:theta-value-locus-tensor-product}, I have suppressed from the notation ${\thetajh}^\breveBb$ the dependence of this set on $X/L$ (and $\ell$). I will ignore $\ell$ for the moment and write ${\thetajh}^\breveBb(X/L)$ to indicate dependence on $X/L$. Let $\sM_{1,1}(L)$ be the $L$-valued points of moduli stack $\sM_{1,1}$ of $1$-pointed elliptic curves. 

Then it is tempting to expect that the following makes sense: let
\be 
\thetajh^{\breveBb}(\sM_{1,1}(L))=\bigcup_{X/L\in\sM_{1,1}(L)} {\thetajh}^{\breveBb}(X/L).
\ee
and expect that 
\be\label{eq:naive-theta-bound}
\abs{\thetajh^{\breveBb}(\sM_{1,1}(L))}_{\breveBb}<\infty.
\ee 

Especially if this is true  then any upper bound for it can serve as a uniform upper bound for ${\thetajh}^{\breveBb}(X/L)$ for every $X/L\in\sM_{1,1}(L)$.

\textcolor{red}{However eq.~\eqref{eq:naive-theta-bound} is too optimistic.} Fixing initial theta data as in \cref{ss:theta-data-fixing} means one is considering elliptic curves over $L$ which acquire semi-stable reduction over extensions of $L'/L$, \textit{however} the field $L'$ will vary as the  curve varies in moduli in general (the  degree $L'/L$ remains bounded by \Cref{le:bounded-degree-theta-data}). 

On the other hand such bounds may exist for arbitrary choice of compactly bounded subsets of $\sM_{1,1}(\bQ)$ (see \cite{mochizuki-general-pos} for definition of compactly bounded subsets). Let $\bQ^{\leq d}$ be the union of all finite extensions of $\Q$ (all contained in $\bQ$) of degrees bounded by the $d$ (specified by for example by \Cref{le:bounded-degree-theta-data}).

For $\mathscr{K}\subset \sM_{1,1}(\bQ)$ a compactly bounded subset supported at a finite set of rational primes $\Sigma$ one can expect to define 
$$\thetajh^{\breveBb}(\mathscr{K}\cap \sM_{1,1}(\bQ^{\leq d}))$$ and one can expect to prove that 
\be\label{eq:naive-theta-bound2} \abs{\thetajh^{\breveBb}(\mathscr{K}\cap\sM_{1,1}(\bQ^{\leq d}))}_{\breveBb}<\infty.\ee
Very roughly this is  the  idea behind the computations of \cite[Corollary 2.2]{mochizuki-iut4}

\section{A `Rosetta Stone' for  \joshiros\ and \iutthr}\label{se:intro-rosetta-stone}
Let me now provide a fundamental `Rosetta Stone' (\rosettastone) for the readers, so that they are able to translate objects and concepts between  the theories of \joshiros\ and \iutthr. 

Once a reader has assimilated and understood the contents of the \rosettastone, hopefully, the readers will be persuaded that there exists corresponding objects in \joshiros\ which allows one to establish the proofs of the relevant claims of \iutthr\ with complete clarity.

\subparat{How to use the Rosetta Stone for reading \iutthr} The \rosettastone\ must be read along with my discussion of other objects of \iutthr\ such as Frobenioids and Hodge Theaters provided here in  \cref{se:frobenioids}. My work together with the (\rosettastone) provides an adequate toolkit for the reader to be able to navigate \iutthr. [The existence of distinct such objects is asserted without proof in \iutthr.]

\subsubparat{Notation} Notations and assumptions of \assumptions\ will strictly be in force for full compatibility with \iutthr. This means $L$ is a number field with no real embeddings, $C/L$ is an elliptic curve and $X=C-\{O\}/L$ is the hyperbolic curve (of strict Belyi Type by Mochizuki's definition) of topological type $(1,1)$. Notably, the isomorphism class $[X/L]$ of $X/L$ defines a moduli point in the moduli stack $\sM_{1,1}$ and the field of definition of this moduli point is $\lmod$ by \assumptions. 

This data, by the results of \present, provides the category  $\fJ(X/L)$ of the adelic arithmetic Teichmuller space associated to $X/L$ \cite{joshi-teich,joshi-untilts}. Each object of this category, by means of \constrtwoh{Proposition }{pr:arith-datum-frob}, provides an arithmeticoid $\arith{L}$ of $L$. As is demonstrated in \cite{joshi-teich-def}, $\arith{L}$ is a deformation of the arithmetic of number field $L$. Working with holomorphoid $(X,\arithl)$ (\Cref{def:holomorphoids}) allows us to view $X$ viewed as a scheme for the version of arithmetic of $L$ provided by the particular arithmeticoid $\arithl$. As a scheme $X/\arith{L}$ is, of course isomorphic to $X/L$.  For the properties of arithmeticoids see \constrtwoh{\ssep}{se:arithmeticoid-adeloid-frobenioid}. To be sure, a holomorphoid provides a global arithmetic structure and by \Cref{def:geom-base-point}, it also provides local arithmetic holomorphic structures at all primes $v\in\vl$. For the properties of arithmetic holomorphic structures and arithmetic teichmuller spaces see \cite{joshi-teich,joshi-teich-def, joshi-untilts}.

\subparat{\Cref{tab:rosetta-stone1}: Arithmetic Holomorphic Structures}
An idea central to both \present\ and especially \iutthr, is the notion of arithmetic holomorphic structures. 
\subsubparat{Arithmetic Holomorphic structures in \cite{joshi-teich,joshi-untilts,joshi-teich-def}}\label{ss:ahs-1}
I will take the definition of arithmetic holomorphic structures as given in \constrone{Definition }{def:arith-hol-strs}, \cite[Definition 4.1]{joshi-untilts}  (which sets out the local theory at each prime of $L$) and demonstrates that this provides all corresponding  local data of \iutthr. [Existence of such structures is not demonstrated in \iutthr.] The global theory of arithmetic holomorphic structures is detailed in \cite{joshi-untilts,joshi-teich-def}.

\subsubparat{The role of geometric base-points in \iut}\label{ss:geom-base-points}\nwsss
Mochizuki's \iut, is founded upon \textit{Mochizuki's Key Principle of Inter-Universality},  which is stated in \cite[\ssep I3, Page 25]{mochizuki-iut1}. According to this principle, Mochizuki's \iut\ requires one  to work with arbitrary geometric base-points. Furthermore, as Mochizuki reminds us \cite[Page 26]{mochizuki-iut1} (and in many other places in \iut\ and \cite{mochizuki-gaussian}), the domains and codomains of all the key operations of his theory, namely $\flog$-links and $\Theta_{gau}$-links, refer to distinct geometric base-points.  

The reader may ask: why is it important to work with arbitrary geometric base-points? \textit{Working with arbitrary geometric base-points is the anabelian proxy for deformations of algebraic varieties.} In \iut\ the idea is to use arbitrary geometric base-points (in the theory of tempered fundamental groups) as giving deformations of arithmetic. However, the exact mechanism by which one obtains such deformations was not established in \iut\ but has been precisely quantified in \cite{joshi-untilts, joshi-teich,joshi-teich-def}. Let me note that distinct geometric base-points also play an important role in the Chabauty-Kim Theory of \cite{kim2005} (i.e. usage of arbitrary geometric base-points is important in other Diophantine contexts).

My definition of arithmetic holomorphic structures \constrone{Definition }{def:arith-hol-strs} (and also see \Cref{def:geom-base-point}) includes geometric base-points and hence my approach is fully compatible with the approach taken in \iut. 

The following table compares the notion of geometric base-points in the theory of \'etale and tempered fundamental groups.
\begin{boxedcontent}[label=tab:geom-base-point,grow to right by=1cm,grow to left by=1.5cm]{Comparison of geometric base-points for \'etale and tempered fundamental groups}
\centering{\small
		\begin{tabular}{|p{1in}|p{2.7in}|p{2.7in}|}\hline 
			&	\'etale fund. group	& tempered fund. group  \\ 
			\hline 
			Input datum & $$(X/\Q_p, *_K:\Spec(K)\to X)$$ where $X/\Q_p$ is a scheme, $K\supset \Q_p$ is an alg. closed field, and  $*_K:\Spec(K)\to X$ is a $K$-geom. base point of $X$	 & $$(\xan_{\Q_p}, *_K:\sM(K)\to \xan_{\Q_p})$$ where $\xan_{\Q_p}$ is an analytic space over $\Q_p$, $K\supset \Q_p$ an alg. closed, perfectoid  field, and $*_K:\sM(K)\to \xan_{\Q_p}$ is a  $K$-geom. base point of $\xan_{\Q_p}$. \\ 
			\hline  
			& isom. class of alg. closed fields $K\supset \Q_p$ is determined by  the cardinality of $K$ &  isom. class of alg. closed, perfectoid fields $K\supset \Q_p$ is not determined by  cardinality of $K$.\\
			\hline
		\end{tabular} 
	}
\end{boxedcontent}
\clearpage

\brem\ 
\benumlab
\item While Mochizuki recognized the centrality of this principle for \iut, Mochizuki's approach to his theory is group theoretic and stays focused on  utilizing the group theory surrounding the tempered fundamental group as the key tool in building his theory, however,  his approach ignores the rich arithmetic, algebra and geometry which arbitrary geometric base-points for tempered fundamental groups provide (as can be seen in \present). This is an important difference in our respective viewpoints.
\item Let me remark that the anabelian reconstruction techniques of \topics\ are of considerable importance in \iut. However, \cite[Theorem 1.9]{mochizuki-topics3} (which is the main theorem of that paper) does not determine the geometric base-point datum while, in accordance with Mochizuki's Key Principle of Inter-Universality, this (undetermined) datum is of central importance in \iut.
\item Let me also say that, as far as I am aware, though Mochizuki's statement of this principle has appeared in all versions of \cite{mochizuki-iut1} (for example it is included in the hard copy of RIMS-Preprint from 2012),  neither this principle nor  the centrality of the role of arbitrary geometric base-points (in \iut) finds any mention in   \cite{fesenko-iut},  \cite{fucheng}, \cite{scholze-stix}, \cite{yamashita}, \cite{dupuy2020statement,dupuy2020probabilistic}, \cite{scholze-review}, \cite{saidi-review-iut}. 
\eenum
\erem

\subsubparat{Arithmetic Holomorphic Structures in \iutthr}
One notable problem with the claims of \iutthr\ is that the existence of distinct arithmetic holomorphic structures is not established in it and my reading of \cite{mochizuki-essential-logic} is that this issue cannot be resolved by methods of \iutthr. My central observation is that arithmetic holomorphic structures provided here in \cref{ss:ahs-1} also provide arithmetic holomorphic structures in \iutthr. \textit{My comparison of arithmetic holomorphic structures in the two theories is given in \Cref{tab:rosetta-stone1}.}

\subsubparat{Fundamental Dichotomy of an Arithmetic Holomorphic Structure}\label{sss:fund-dichotomy} Arithmetic holomorphic structures in both the theories (\iutthr\ and \joshiros)  consists of two parts: 
\benumlab
\item The first portion is tuned into arithmetic and topology of the field, and geometry of the variety while
\item  the second portion tracks independent changes of valuations and hence requires recording the value group.
\eenum
I call this  \textit{the fundamental dichotomy of an Arithmetic Holomorphic Structure}. Importantly, arithmetic holomorphic structures in both the theories exhibit such a dichotomy. The non-triviality of the theories rests on exhibiting the existence of plurality of such data. This is a  central problem for \iutthr\ which is not resolved by the artificial ``resetting of logic'' recommended in \cite{mochizuki-essential-logic}; on the other hand the theory of  \present\ provides a natural plurality of such data from the very beginning and these data can be applied to exhibit the corresponding data in \iutthr.

\subsubpara{} Let $X/\arith{L}_\by$ be a holomorphoid of $X/L$. In \cite{joshi-teich,joshi-untilts} the fundamental dichotomy arises from the data of the arithmeticoid $\arith{L}_\by$. Let $$\by=(K_v,K_v^\flat\isom \C^\flat_{p_v})_{v\in\vlp}\in\yadlp'$$ provided by $\arith{L}_\by$. [Recall that by \constrtwoh{Proposition }{pr:arith-datum-frob},  each arithmetic holomorphic structure in the sense of \cite{joshi-untilts} provides, the datum of the point $\by$.] 

The portion of the datum  given by the perfectoid field $K_v$ (with an embedding of $L_v$) provides the geometric data of the analytic spaces $(X/L_v,X/K_v, \xan_{K_v}\to \xan_{L_v})$ and these provide  preferred isomorphs of the tempered fundamental group of $X/L_v$, the geometric tempered fundamental group of $X/K_v$ and the absolute galois group of $L_v$ computed for the algebraic closure of $L_v$ in $K_v$ (see \constrone{Theorem }{thm:arith-hol-strs}).

The second portion is the data of the isometry $K_v^\flat\isom \C^\flat_{p_v}$ which records valuation specific information of $L_v\subset \bL_v\subset K_v$  in one fixed value group, namely $\abs{\C^\flat_{p_v}}$, even as the untilt  $(K_v,K_v^\flat\isom \C^\flat_{p_v})_{v\in\vlp}$ moves (as $\by$ moves).

\subparat{Mochizuki's Dichotomy: \'Etale-pictures and Frobenius-pictures via \cite{joshi-teich}}\label{ss:mochizuki-dichotomy}
Mochizuki has a similar dichotomy for arithmetic holomorphic structures. The two portions of Mochizuki's dichotomy for arithmetic holomorphic structures  are termed in \iutthr\ as
\benumlab
\item an ``\'etale-picture'' (see \cite[\ssep 3, Corollary 3.9]{mochizuki-iut1}) and, 
\item a ``Frobenius-pictures'' (see \cite[\ssep 3, Corollary 3.8]{mochizuki-iut1}).
\eenum

The relationship between Dichotomy of Arithmetic Holomorphic Structures {sss:fund-dichotomy}  and  Mochizuki's dichotomy of \'etale picture (or \'etale like data) and Frobenius-picture or (Frobenius-like data) is described for each prime $v\in\vl$ separately. Hence fix a $v\in\vl$ and let $\hol{X/L_v}{y_v}{F_v}$ be an arithmetic holomorphic structure on $X/L_v$. Let  $y_v=(K_v,i_{K_v}:K_v^\flat\isom F_v)$ be point of $\sY_{F_v,L_v}$  provided by $\hol{X/L_v}{y_v}{F_v}$. Then I claim that  Mochizuki's Dichotomy of \'etale-picture and Frobenius-pictures arises from the Dichotomy of Arith Holomorphic Structure given by the datum of the local holomorphoid $\hol{X/L_v}{y}{F_v}$ as by \cref{pr:mochizuki-etale-picture} and \cref{pr:mochizuki-frobenius-picture} and is summarized in \Cref{tab:rosetta-stone1}:

\subsubparat{Mochizuki's \'Etale Picture or \'etale-like data}\label{ss:etale-like-data}\nwsss
Mochizuki's \'etale-picture consists of the group theoretic data provided by the  proposition below:
\bpro\label{pr:mochizuki-etale-picture}
	Each holomorphoid   $\holt{X/L}{\by}$ of $X/L$ provides,  for each $v\in\vl$,  Mochizuki's \'etale-like picture, labeled by the geometric data of the (local) holomorphoid of $\hol{X/L_v}{y_v}{F_v}$ from which it arises and  consists of the following data:
	\benumlab
	\item the algebraically closed perfectoid field $K_v$ which contains an isometrically embedded $\Q_{p_v}$ and hence $K_v$ provides an algebraic closure of any finite extension $L_v$ of $\Q_p$. Thus $K_v$ provides a preferred algebraic closure $L_v\subset \bL_{v;K_v}\subset K_v$ of $L_v$ and hence provides a preferred isomorph $G_{L_v;K_v}=\gal(\bL_v/L_v)\isom G_{L_v}$ of absolute Galois group of $L_v$;
	\item  the Galois category  $\mathrm{Temp}_{X/L_v;K_v}$  (as defined in \cite[Chap. III, 2.1.4]{andre-book})  of tempered coverings of $\xan_{L_v}$ which are defined over finite extensions of $L_v$ contained in $K_v$;
	\item and a preferred isomorph $\pit{X/L_v;K_v}$ of the tempered fundamental group of $X/L_v$ (see \constrone{Theorem }{thm:arith-hol-strs}).
	\item Importantly, in my theory, this data appears with the distinguishable label given by the holomorphoid  $\hol{X/L_v}{y_v}{F_v}$ from which this ``\'etale-picture'' arises.
	\eenum
\epro
\bp 
The proof is clear.
\ep

\brem 
There are other \'etale like data in \iutthr, these may be also obtained in a similar manner as above.
\erem

\subsubparat{Frobenius-picture or Frobenius-like data}\label{ss:frobenius-like-data}
\nwsss
Now I claim that $\hol{X/L_v}{y_v}{F_v}$ also provides Mochizuki's Frobenius-picture or Frobenius-like data. 

\bpro\label{pr:mochizuki-frobenius-picture}
	Let $\holt{X/L}{\by}$ be a holomorphoid of $X/L$ and let $\hol{X/L_v}{y_v}{F_v}$ be the local holomorphoid of $X/L_v$ (for any $v\in\vl$) given by $\holt{X/L}{\by}$. Then
	\benumlab
	\item the isometry $i_{K_v}:K_v^\flat\isom F_v$ identifies the value group of $K_v$ with the value group of $F_v$ \cite{scholze12-perfectoid-ihes};
	\item in particular for all anabelomorphic holomorphoids of $X/L_v$ in the category $\fJ(X/L_v)_{F_v}$, one has a fixed value group--namely the value group of $F_v$, to compare valuations.
	\item let $\bL_v\subset K_v$ be the algebraic closure of $L_v$ in $K_v$ and $\O_{\bL_v}^\triangleright=\O_{\bL_v}^\triangleright-\{0\}\subset \O_{K_v}^\triangleright-\{0\}$ be the multiplicative monoid of the ring of integers $\O_{\bL_v}\subset \bL_v\subset K_v$. Then $\O_{\bL_v}^\triangleright$ is a Frobenius-like (monoidal) data of  \iutthr,\cite{mochizuki-gaussian}.
	\eenum
\epro
\bp 
This is clear.
\ep
The value group information provided by this proposition is  Mochizuki's Frobenius picture \cite{mochizuki-iut1}. 

\brem Let me remark that the definition of $K^\flat$ \cite{scholze12-perfectoid-ihes} shows how Frobenius i.e. $p^{th}$-powers play a central role in the construction of $K^\flat$. So the term Frobenius-like or Frobenius-picture is quite reasonable in this context.  One important point to be noted is that even if $K$ is fixed, there are many inequivalent untilts, $(K,K^\flat \isom F)$, of $F$ and hence there are many inequivalent Frobenius-like pictures. There are other types of Frobenius-like data in \iut. These are also provided by the datum $\hol{X/E}{y}{F}$. 
\erem

\subsubparat{Existence of \'Etale and Frobenius Pictures} Since the geometric data of holomorphoids of $X/L$ provide  \'etale-pictures and Frobenius-pictures of \iutthr\ one has now established the following:
\bthm
	Each holomorphoid $\holt{X/L}{\by}$ of $X/L$, provides, for each $v\in\vl$, an arithmetic holomorphic structure in the sense of \cite{joshi-untilts} and also  provides Mochizuki's Dichotomy of an arithmetic holomorphic structure: namely an ``\'Etale picture'' and a ``Frobenius picture'' in the sense of \iutthr, which is indexed by the holomorphoid $\holt{X/L}{\by}$ which it arises from. In particular, there exist many inequivalent \'etale and Frobenius pictures in \iutthr.
\ethm

\brem I will not use the terminology of \'etale and Frobenius pictures, but the above results make it clear that my approach does indeed have all the properties claimed by Mochizuki in \iutthr\ where he has asserted the existence of such inequivalent data without proofs.
\erem

\subsubparat{$L^*$-action on \'Etale and Frobenius pictures}
There is natural action of $L^*$ on the set of pairs consisting of (\'etale,Frobenius)-pictures and the justification of the existence of this action is inadequate (in my opinion) in  \iutthr, though Mochizuki asserts the principal the property of this action \cite[Remark 1.2.3(ii)]{mochizuki-iut3}, namely that if $\alpha\in L^*$ acts trivially, then $\alpha$ is a root of unity. 

In \constrtwoh{Theorem }{th:galois-action-on-adelic-ff}, I demonstrate the existence of the natural multiplicative action of $L^*\act \yadl$. From the above description of (\'etale,Frobenius)-pictures, one sees at once, that $L^*\act \yadl$ provides the required action on the sets of such pairs (arising, as above, from arithmeticoids in $\yadl$). The fundamental property of this action alluded to in the above paragraph is proved in \constrtwoh{Theorem }{th:fundamental-property-frobenius}.

\subsubparat{\Cref{tab:rosetta-stone1}}
Before exhibiting \Cref{tab:rosetta-stone1}, let me recall some of Mochizuki's notation \cite[Figure I1.2, Page 9]{mochizuki-iut1} (Mochizuki's base number field is $F$ while mine is $L$): 

\begin{itemize}
	\item $v$ is a prime of a fixed number $L$,
	\item $L_v$ is its completion at $v$, 
	\item $\O_v$ is its ring of integers, 
	\item $\bL_v$ is an algebraic closure of $L_v$ and 
	\item $\bL_v^\times=\bL_v^*=\bL_v-\{0\}$;
\item $\Ob_v\subset \bL_v$ is the integral closure of $\O_v$ in $\bL_v$;
	\item  $\Ob_v^\times\subset \bL_v^\times$ is group of units;
\item $\Pi_v$ is the tempered fundamental group of a fixed hyperbolic curve $X/L_v$. 
\end{itemize}
The notation for the last column of \Cref{tab:rosetta-stone1} is as declared in the present paper.

\clearpage
\newpage
\numberwithin{table}{subsubsection}
\begin{rosetta}[label=tab:rosetta-stone1,grow to right by=1cm,grow to left by=2cm]{}
\centering
	\begin{tabular}{|p{1in}|p{2in}|p{3.5in}|}\hline 
	Theory/Object	& \iutthr & \joshiros \\ 
	\hline 
	\'etale picture	&group-like objects $$\Pi_v, G_v, \overline{\Pi}_v\into \Pi_v\twoheadrightarrow G_v$$ and their ``abstract isomorphs''& each primary geometric  object $$\displaystyle{(\xan_{L_v}, (K_v\supset E,K_v^\flat\isom F_v), *_K:\sM(K_v)\to \xan_{L_v})}$$   provides all the  group-like objects $$\pit{X/L_v;K_v}, G_{L_v;K_v}, \pit{X/K_v}\into \pit{X/L_v}\onto G_{L_v;K_v},$$ corresponding group-like objects are all isomorphic. The datum $$(*_K:\sM(K_v)\to \xan_{L_v})$$ is the datum of a geometric base-point.\\ 
	\hline 
	Frobenius picture	& valuation theoretic objects:
	$$G_v\act\bL_{v}^* $$
	$$G_v\act\Ob_v^*,$$
	 $$G_v\act\Ob_v^\triangleright$$ and abstract isomorphs of all of these  & primary valuation theoretic object is the given isomorphism $$K_v^\flat\isom F_v^\flat,$$ one also has Mochizuki-style Frobenius/valuation theoretic objects:
	 $$G_{L_v;K_v}\act \bL_{v;K_v}^* \subset {K_v}^*$$
	 $$G_{L_v;K_v}\act \O_{\bL_{v;K_v}}^* \subset \O_{K_v}^*$$
	  $$G_{L_v;K_v}\act \O_{\bL_{v;K_v}}^\triangleright \subset \O_{K_v}^\triangleright$$  given by the geometric data of holomorphoids. \\ 
	\hline 
\end{tabular} 
\tcblower 
	Dichotomy of Arithmetic Holomorphic structures in \joshiros\ and \iutthr.
\end{rosetta}
\clearpage
\newpage

\subsubparat{\Cref{tab:rosetta-stone2}}\numberwithin{table}{subsubsection}
The table given below (\Cref{tab:rosetta-stone2}) will be useful in understanding the role of holomorphoids in \iutthr.
\begin{rosetta}[label=tab:rosetta-stone2,grow to right by=1cm,grow to left by=0.75cm]{}
	\begin{tabular}{|p{2in}|p{2in}|p{2in}|}
		\hline 
		\iutthr	& \joshiros & notation in \joshiros\\ 
		\hline 
		mono-analytization at $v$	& working with a specific (local) holomorphoid of $X/L_v$ \cref{se:holomorphoids} & $\holt{X/L_v}{y_v}$ \\ 
		\hline 
		coric objects, structures at $v$	& working with $\varphi_v$-invariant objects, structures &  \\ 
		\hline 
		mono-analytic prime-strip $\sF_{\scriptscriptstyle{Mochizuki},v}^?$  (at $v$) & the prime-strip $\sF_{\scriptscriptstyle{Joshi},v}^?$  given by a specific (local) holomorphoid of $X/L_v$ \cref{ss:prime-strips}  & $\sF_{\scriptscriptstyle{Joshi}}^?(\holt{X/L_v}{y_v})$  \\ 
		\hline 
		mono-analytic log-shell at $v$	& the log-shell given by a specific (local) holomorphoid of $X/L_v$ & \\ 
		\hline
		global coricity & corresponds to working with $\yadlp/\bvarphi^\Z=\xadlp$ or its variants such as $\yadlp/L^{'*}$ & \\
		\hline 
	\end{tabular}
\tcblower
Mono-analyticizations, coric structures of \iutthr\ via holomorphoids of \joshiros.
\end{rosetta}
\subparat{Mochizuki's Diagram of coverings}\nwss
Let $\arith{L}_{\by}$, with $\by=(y_v)_{v\in\vl}$, be an arithmeticoid of $L$ and suppose one has holomorphoids of $C/L,X/L,L'$ be as given by \cref{ss:elliptic-curve-assumptions}, \inithtdata\ arising from $\arith{L}_{\by}$.  Then one has a certain diagram of coverings constructed using $\arith{L}_{\by}$ (which will not be recalled here). The constructions of these coverings is detailed in \cite{mochizuki-theta} and also in \cite[\ssep1 and \ssep2]{mochizuki-iut1}. 
\blem 
Let $\holt{X/L}{\by}$, with $\by=(y_v)_{v\in\vl}$, be a holomorphoid of $X/L$ and suppose that \cref{ss:elliptic-curve-assumptions}, \inithtdata\ hold. Then 
\benumlab
\item one has a construction of various finite \'etale and tempered coverings of $X/L_v$ (for each $v\in\vlnon$) detailed in \cite[\ssep 1, Page 42]{mochizuki-iut1}, \cite[Definition 3.1(d), Page 68]{mochizuki-iut1} and \cite[Definition 3.1(e), Page 69]{mochizuki-iut1} and moreover, 
\item the isomorphism class of the diagram of topological groups, obtained by applying the \'etale (resp. tempered fundamental group) functor to this diagram of coverings obtained using $\holt{X/L}{\by}$, is independent of the choice of the arithmeticoid $\arith{L}_{\by}$ provided by the holomorphoid $\holt{X/L}{\by}$.
\eenum
\elem
\bp 
The construction of coverings being referred to here and detailed in \cite[\ssep 1, Page 42]{mochizuki-iut1}, \cite[Definition 3.1(d), Page 68]{mochizuki-iut1} and \cite[Definition 3.1(e), Page 69]{mochizuki-iut1} is quite geometric and so there are no essential difficulties in establishing these coverings in the presence of the arithmeticoid $\arith{L}_{\by}$. As is established in Mochizuki's paper, the isomorphism class of the diagram of groups is determined from the \'etale fundamental group (resp. the tempered fundamental group) of $X/L$ (resp. $X/L_v$) and its stack quotients by $\{\pm1\}$.
\ep
\subparat{\Cref{tab:rosetta-stone3}: Prime-strips}\label{ss:rosetta-stone}
\nwss
In  \Cref{tab:rosetta-stone3}, I will consider holomorphoids  $(X/E,X/K)\in\fjxe$ with the restriction $K^\flat=\cpt$.  This restriction corresponds to working with the category $\fjxe_{\cpt}$ of  \constrone{\ssep}{se:construct-att} and places cardinality restriction on $K$ (because $K$ and $K^\flat$ have the same cardinality). Cardinality restrictions are missing in \iutthr.  [By \cite{scholze12-perfectoid-ihes}, fixing the cardinality of the tilt $K^\flat$ to be that of $\cpt$ implies the perfectoid field  $K$ provided by any untilt of $\cpt$ has the cardinality of $\cpt$. The latter has the cardinality of $\C_p$ which has the cardinality of the continuum.]
This is cardinality restriction poses no problems in the context of \iutthr\ (the theory of \joshiros\ is far more general).

In \Cref{tab:rosetta-stone3}, second column adheres to Mochizuki's notation \cite[Figure I1.2, Page 9]{mochizuki-iut1} (his base field is $F$ while mine is $L$):

\begin{itemize}
	\item $v$ is a prime of a fixed number $L$,
	\item $L_v$ is its completion at $v$, 
	\item $\O_v$ is its ring of integers, 
	\item $\bL_v$ is an algebraic closure of $L_v$ and 
	\item $\bL_v^\times=\bL_v^*=\bL_v-\{0\}$;
	\item $\mu(\bL_v)\subset \bL^\times_v$ is the group of roots of unity in $\bL_v$.
	\item $\Ob_v\subset \bL_v$ is the integral closure of $\O_v$ in $\bL_v$;
	\item  $\Ob_v^\times\subset \bL_v^\times$ is group of units;
	\item  $\Ob_v^{\times\mu}=\Ob_v^{\times}/\mu(\bL)$;
	\item $\Pi_v$ is the tempered fundamental group of a fixed hyperbolic curve $X/L_v$. 
\end{itemize}
The notation for the last column of the table  $L_v, \bL_v$ is clear and rest of the notation is as declared in the present paper.

\newcommand{\cptv}{\C_{p_v}^\flat}

\newpage
\begin{minipage}{\textwidth}
\begin{rosetta}[label=tab:rosetta-stone3,grow to right by=1cm,grow to left by=2cm]{}\renewcommand{\arraystretch}{1.5}	
	\begin{tabular}{|c|c|c|}
		\hline 
		\rule[-1ex]{0pt}{2.5ex} Mochizuki's Notation & Object in \cite{mochizuki-iut3} & Object in \cite{joshi-teich} \\ 
		\hline 
		\rule[-1ex]{0pt}{2.5ex} $F$ & $F$ base field & $L$ base field \\ 
		\hline 
		\rule[-1ex]{0pt}{2.5ex} $\bF$ & $\bF$ alg. clos. & $\bL$  {\tiny algebraic closure  provided by an arithmeticoid} \\ 
		\hline 	
		\rule[-1ex]{0pt}{2.5ex} $F_v$ & $F_v$ completion at v & $L_v$ completion at v \\ 
		\hline 
		\rule[-1ex]{0pt}{2.5ex} $\bF_v$ & $\bF_v$ alg. clos. & $\bL_v\subset K_v$ {\tiny algebraic closure in the $v$-component $K_v$ provided by an arithmeticoid} \\ 
		\hline 
		\rule[-1ex]{0pt}{2.5ex} $\sD$ & $\Pi_v$ & $\pit{(X/L_v,X/K_v)} (\isom \pit{(X/L_v,X/\C_{p_v})})$ \\ 
		\hline 
		\rule[-1ex]{0pt}{2.5ex} $\sD^\vdash$ & $G_v$ & $G_{L_v;K_v} (\isom G_{L_v;K_v}) \qquad (\bL_v\subset K_v)$ \\ 
		\hline 
		\rule[-1ex]{0pt}{2.5ex} $\sF$ & $\Pi_v\act \Ob_{v}$ & $\pit{X/L_{v};K_v}\act \O_{\bL_{v}}\qquad (\bL_v\subset K_v)$  \\ 
		\hline 
		\rule[-1ex]{0pt}{2.5ex} $\sF^\vdash$ & $G_v\act\Ob_v^\times\times q_v^\bN$ & $G_{L_v;K_v}\act\O_{\bL_v}^\times\times q_{(X/L_v,X/K_v)}^\bN \qquad (\bL_v\subset K_v)$\\ 
		\hline 
		\rule[-1ex]{0pt}{2.5ex} $\sF^{\vdash\times}$  & $G_v\act\Ob_v^{\times}$ & $G_{L_v;K_v}\act\Ob_v^{\times} \qquad (\bL_v\subset K_v)$ \\ 
		\hline 
		\rule[-1ex]{0pt}{2.5ex} $\sF^{\vdash\times\mu}$ & $G_v\act\Ob_v^{\times\mu}$ & $G_{L_v;K_v}\act \O_{\bL_v}^{\times\mu} \qquad (\bL_v\subset K_v)$ \\ 
		\hline 
		\rule[-1ex]{0pt}{2.5ex} $\sF^{\vdash\blacktriangleright\times\mu}$ & $G_v\act\Ob_v^{\times\mu}\times q^\bN$  & $G_{L_v}\act \O_{\bL_v}^{\times\mu} \times q_{(X/L_v,X/K_v)}^\bN$ \\ 
		\hline 
		\rule[-1ex]{0pt}{2.5ex} $\sF^{\vdash\blacktriangleright}$  & $G_v\act q^\bN$ & $G_{L_v;K_v}\act q_{(X/L_v,X/K_v)}^\bN \qquad (\bL_v\subset K_v)$ \\ 
		\hline 
		\rule[-1ex]{0pt}{2.5ex} $\sF^{\vdash\perp}$  & $G_v\act \Ob_v^{\mu_{2\ell}}\times q^\bN$ & $G_{L_v;K_v}\act \O_{\bL_v}^{\mu_{2\ell}}\times q_{(X/L_v,X/K_v)}^\bN \qquad (\bL_v\subset K_v)$ \\ 
		\hline 
\hline 
		\multicolumn{3}{|c|}{Holomorphoid Notation For Prime-strips}\\
		\hline
		\rule[-1ex]{0pt}{2.5ex} $\sF^?$ &  $\sF_{\scriptscriptstyle{Mochizuki}}^{?}$ & $\sF_{\scriptscriptstyle{Joshi}}^{?}(\hol{X/L_v}{y_v}{\cptv})$ \\ 
		\hline

	\end{tabular}
\end{rosetta}		
\end{minipage}
\vfill\eject

\subparat{Primes-strips \`a la \iutthr\ from \joshiros}\label{ss:prime-strips} A prime-strip (\cite[\ssep I1, Page 9]{mochizuki-iut1}) in \iutthr\ consists of certain monoidal data equipped with the action of the local Galois group. While the notion of prime-strips  of \cite[\ssep I1, Page 9]{mochizuki-iut1} is not used in the theory of \joshiros, prime-strips are certainly available in the theory of \joshiros\ and there is a correspondence between prime-strips of \iutthr\  and those of \joshiros. This correspondence between prime-strips of the two theories may be transcribed using \Cref{tab:rosetta-stone3} (reading it from left to right) from the corresponding object in the rightmost column (which provides objects in the theory of \present). For example:

A prime-strip $\sF^{?}$ arising from a datum of the  point $(X/L_v,X/K_v)$ of the Arithmetic Teichmuller Space $\fJ(X,L_v)_{\cpt}$   would be notated as
\be\sF^{?}_{\scriptscriptstyle{Joshi}}(X/L_v,X/K_v)\ee to indicate the datum of the point of the Arithmetic Teichmuller Space giving rise to this particular prime-strip. So 
where Mochizuki \cite[\ssep I1, Page 9]{mochizuki-iut1} writes 
\be\label{eq:mochizuki-prime-strip}
\sF^{\vdash\times\mu}=G_v\act\O_{\bL_v}^{\times\mu}
\ee
I will write 
\be\sF^{\vdash\times\mu}_{\scriptscriptstyle{Joshi}}(X/L_v,X/K_v)=G_{L_v}\act \O_{\bL_v}^{\times\mu},\qquad L_v\subset \bL_v \subset K_v\ee
where $G_{L_v}$ is the absolute Galois group of $\bL_v/L_v$ and where $\bL_v\subset K_v$ is the algebraic closure of $\bL_v$ in $K_v$.

Of course,  if one is dealing with more than one holomorphoid simultaneously, then one can add specific notation and write 
$$\sF_{\scriptscriptstyle{Joshi}}^{?}(\hol{X/L_v}{y_v}{F_v})$$
instead.

Similarly 
\be\sF^{\vdash\blacktriangleright\times\mu}_{\scriptscriptstyle{Mochizuki}} = G_v\act\O_{\bL_v}^{\times\mu}\times q^\bN_{X/L_v,X/K_v},\qquad \bL_v\subset \bL_v \subset K_v\ee
of \iutthr, will be
\be\label{eq:joshi-prime-strip}\sF^{\vdash\blacktriangleright\times\mu}_{\scriptscriptstyle{Joshi}}(X/L_v,X/K_v)=G_{L_v}\act \O_{\bL_v}^{\times\mu} \times q_{(X/L_v,X/K_v)}^\bN,
\ee
where $G_{L_v}$ is as before, $\bL_v$ is as before and where $q_{(Y/E',Y/K)}$ is the Tate quasi-period of the projective Tate elliptic curve $C/L_v$ underlying $X/L_v$ (with the period being computed in the analytic space $C^{an}_{K_v}$. 

\subsubparat{Realified Prime-strips a la Mochizuki from \joshiros}\nwsss
Let $v\in\vl$ and let $\holt{X/L_v}{y_v}=(X/L_v,X/K_{y_v})$ be a holomorphoid of $X/L_v$. 
Recall that by \cite{fargues-fontaine}, the data of the point $y_v\in\sY_{\cpt,L_v}$ provides a valuation $\abs{-}_{K_{y_v}}$ on the residue field $K_{y_v}$ and hence a valuation $\abs{-}_{\bL_v}$ on the algebraic closure $\bL_v\subset K_{y_v}$ of $L_v$ contained in $K_{y_v}$. 

For  datum of this holomorphoid, and for each prime-strip $$\sF_{\scriptscriptstyle{Joshi}}^?(\holt{X/L_v}{y_v})$$ given in \Cref{tab:rosetta-stone3} arising from this datum, one also has  a \textit{realified prime-strip} \be\label{eq:realified-prime-strip}\sF_{\scriptscriptstyle{Joshi}}^{|?}(\holt{X/L_v}{y_v})\ee in the sense of \iutthr. My notation for realified prime-strips follows Mochizuki's convention  \cite[Table I1.2 on page 9]{mochizuki-iut1} which puts an extra vertical bar $|?$ at the start of the decoration $?$ to denote the corresponding realified version of that prime-strip. 

Thus for example, the realified prime-strip corresponding to the prime-strip $$\sF^{\vdash\blacktriangleright\times\mu}_{\scriptscriptstyle{Joshi}}(X/L_v,X/K_{y_v})=G_{L_v}\act \O_{\bL_v}^{\times\mu} \times q_{(X/L_v,X/K_{y_v})}^\bN$$ will be denoted by
$$\sF^{\Vdash\blacktriangleright\times\mu}_{\scriptscriptstyle{Joshi}}(X/L_v,X/K_{y_v})$$  and consists of the data 
$$(\sF^{\vdash\blacktriangleright\times\mu}_{\scriptscriptstyle{Joshi}}(X/L_v,X/K_{y_v}), \abs{-}_{K_{y_v}}:K_{y_v}\to \R)$$ consisting of the prime-strip $\sF^{\Vdash\blacktriangleright\times\mu}_{\scriptscriptstyle{Joshi}}(X/L_v,X/K_{y_v})$ plus the valuation $$\abs{-}_{K_{y_v}}:K_{y_v}\to \R $$ given by \cite{fargues-fontaine} which allows us to compute absolute values of elements of $\bL_v\subset K_{y_v}$ and in particular calculate the absolute value of the Tate parameter $$q_{(X/L_v,X/K_{y_v})}\in K_{y_v}.$$ 

\brem Mochizuki's method of keeping track of valuations is through the notion of realified Frobenioids. By my discussion of Frobenioids in \cref{ss:frobenioid-non-archimedean}, one can use the given valuation $(K_{y_v},\abs{-}_{K_{y_v}}\to \R)$ to obtain the relevant realified frobenioid.
\erem

These observations are summarized in the following proposition.
\bpro\label{pr:prime-strips-constructed} 
Each holomorphoid $\holt{X/L}{\by}$ with $\by=\{y_v\}_{v\in\vl}$ provides (for each $v\in\vl$) an isomorph, $$\sF_{\scriptscriptstyle{Joshi}}^?(X/L_v)_{y_v},$$ of Mochizuki's prime-strip $$\sF_{Mochizuki,L_v}^?$$ (as in \Cref{tab:rosetta-stone3}) and also an isomorph $$\sF_{\scriptscriptstyle{Joshi}}^{|?}(X/L_v)_{y_v}$$ of the corresponding realified prime-strip $$\sF_{Mochizuki,L_v}^{|?}$$ (as in \Cref{tab:rosetta-stone3}).
In particular, there are many geometrically and arithmetically distinguishable prime-strips in \iutthr.
\epro

\brem In particular, \cref{pr:prime-strips-constructed} provides clearly distinguishable geometric labels for the prime-strips. This means that there exist many distinct (but abstractly isomorphic) prime-strips in \iut-- this is needed in the proofs of \cite[Theorem 3.11]{mochizuki-iut3} and \moccor.
\erem
\subparat{The fundamental Prime-Strips of \cite{joshi-teich,joshi-untilts}}\nwss
In \present, I do not use the notion of prime-strips but it is certainly possible to define  prime-strips as indicated in \Cref{tab:rosetta-stone3}. But \Cref{tab:rosetta-stone5} shows that there is however a more natural and even canonical object which serves this  role and is given by 
\be\label{eq:canonical-prime-strip}
G_{L_v}\act \hgm(\O_{\cpt}) \isom G_{L_v}\act B^{\vphi=p}.
\ee
(the isomorphism is given by \cite[Proposition 4.4.6]{fargues-fontaine}) which I shall call the \textit{canonical prime-strip at $p$} provided by the theory of  \present.
By \Cref{tab:rosetta-stone5}, there are natural, and equally important, variants  of \eqref{eq:canonical-prime-strip}:
\begin{align}
G_{L_v}\act \abs{\syQp},\\
\shortintertext{and}
G_{L_v}\act \abs{\sxqp}.
\end{align}
\brem
The fundamental prime-strip at $p$ \eqref{eq:canonical-prime-strip} provides a number of advantages over Mochizuki's definition of prime-strips \eqref{eq:mochizuki-prime-strip}. Notably, I demonstrate (in \Cref{th:gluing-prime-strips}) that the set of prime-strips in \iutthr\ is equipped with a natural action of $G_{L_v}$ which corresponds to the galois action on the set of arithmetic holomorphic structures via its action on $\sY_{\cpt,L_v}$ i.e. the natural operations of the theory permute the prime strips of \iutthr. This is an  important point, even for the proofs of \cite{mochizuki-iut3}, however,  it is not established in \iutthr. 
\erem

\subsubparat{\Cref{tab:rosetta-stone5}}
The next \Cref{tab:rosetta-stone5} provides a detailed comparison of prime-strips in the two theories.

\begin{minipage}{\textwidth}
\begin{rosetta}[label=tab:rosetta-stone5,grow to right by=1cm,grow to left by=2cm]{}
	\renewcommand{\arraystretch}{3}	
	\begin{center}
	\begin{tabular}{|c|c|}
		\hline 
		\multicolumn{2}{|c|}{Fundamental Prime strips}\\
		\hline 
		\multicolumn{2}{|c|}
		{$	
		\begin{rcases*}
			\left\{G_v\act \O_{L_v}^*:v\in\vlnon \right\}  &  \\ 
\left\{G_v\act \O_{L_v}^*/\mu(\bL_v):v\in\vlnon \right\}  &
		\end{rcases*}\leftrightsquigarrow\left\{G_v\act\hgm(\O_{\cpt}) :v\in\vlnon \right\}
		$}
	\\
	\hline 
	\multicolumn{2}{|c|}{Mochizuki style versions of the Fundamental strips}\\
	\hline 
$	
		\left\{G_v\act \O_{L_v}^*:v\in\vlnon \right\} $ &  $\leftrightsquigarrow\left\{G_v\act\hgm(\O_{\cpt})/\Z_p^* :v\in\vlnon \right\}$ \\ 
		\hline 
		$\left\{G_v\act \O_{L_v}^*/\mu(\bL_v):v\in\vlnon \right\}  $  &
		$\leftrightsquigarrow\left\{G_v\act\hgm(\O_{\cpt})/\Q_p^* :v\in\vlnon \right\}$\\
	\hline 
	\multicolumn{2}{|c|}{Mochizuki style mono-analytic versions of Fundamental strips}\\
	\hline 
$	
	\left\{G_v\act \O_{L_v}^*:v\in\vlnon \right\} $ &  $\left\{G_{v;K_{y_v}}\act \O_{L_v;K_{y_v}}^*:v\in\vlnon \right\}_{\by\in\yadl}$ \\ 
	\hline 
	$\left\{G_v\act \O_{L_v}^*/\mu(\bL_v):v\in\vlnon \right\}  $  &
	$\left\{G_{v;K_{y_v}}\act \O_{L_v;K_{y_v}}^*/\mu(\bL_{v;K_{y_v}}):v\in\vlnon \right\}_{\by\in\yadl}$
	\\		
	\hline 
	\multicolumn{2}{|c|}{Geometric versions of the Fundamental strips}\\
	\hline
	\textcolor{red}{no geometric versions available in IUTT} & \textcolor{blue}{geometric versions are listed below} \\
	\hline
	   $\scriptstyle{\textcolor{red}{serves\ as\ a\ proxy:\ }}$ \ 
	   $\left\{G_v\act \O_{L_v}^*:v\in\vlnon \right\}$ & $\left\{G_v\act \syflv:v\in\vlnon \right\}$ \\ 
		\hline 
		$\scriptstyle{\textcolor{red}{serves\ as\ a\ proxy:\ }}$ \ 	 $\left\{G_v\act \O_{L_v}^*/\mu(\bL_v):v\in\vlnon \right\}$ & $\left\{G_v\act \sxflv:v\in\vlnon \right\}$ \\ 
		\hline 
		\end{tabular} 
	\end{center}
	\tcblower
	Fundamental Prime-strips in the two theories
\end{rosetta}
\end{minipage}
\clearpage
\newpage

\subparat{Gluing  along prime-strips \`a la Mochizuki via \joshiros}\label{ss:gluing-prime-strips} 
Now let me prove that the analog of Mochizuki's gluing of prime-strips \cite[Page 23]{mochizuki-iut1} exists in the theory of present series of papers because of \Cref{th:gluing-prime-strips}. 

Let me describe Mochizuki's gluing of two  prime-strips of the type:
\be
\sF_{\scriptscriptstyle{Mochizuki}}^{\vdash\blacktriangleright\times\mu}=G_v\act\Ob_{v}^{\times\mu}\times q^\bN
\ee
the core idea of \iutthr, especially see \cite[Page 23]{mochizuki-iut1}, is that group theoretic data of \iutthr\ i.e. two isomorphs $(\Pi_1,\pib_1)$ and $(\Pi_2,\pib_2)$ of $(\Pi,\pib)$ are glued along an isomorphism of the copies of $$\sF_{\scriptscriptstyle{Mochizuki}}^{\vdash\blacktriangleright\times\mu}=G_v\act\Ob_{v}^{\times\mu}\times q^\bN.$$  
provided by each of these pairs of fundamental groups.
\brem
	\iutthr\ attempts to accomplish this gluing \cite[Page 23]{mochizuki-iut1} without explicitly exhibiting the existence of distinct primes-strips, i.e. prime-strips in \iutthr\ do not have any arithmetically (or any natural or Teichmuller Theoretic) distinguishable labels.	Mochizuki's gluing claim  \cite[Page 23]{mochizuki-iut1}  is established here in \Cref{th:gluing-prime-strips} by means of the theory developed in  \present\ and applies to context of \iutthr. 
\erem

Prime-strips are local objects so in the next theorem I work with a fixed prime $v\in\vl$.

\bthm\label{th:gluing-prime-strips} 
Let the notation be as in \cref{ss:elliptic-curve-assumptions}. Let $\by=(y_v)_{v\in\vl}\in\yadl$ and let $$\sigma=(\sigma_v)_{v\in\vl}:\yadl\to\yadl$$ be an isomorphism of $\yadl$ given by any of the actions on $\yadl$ established in \constrtwoh{Theorem }{th:galois-action-on-adelic-ff} and \constrtwoh{Corollary }{cor:global-frobenius}. Let $v\in\vl$. Let $\sigma(\by)=(\sigma_v(y_v))_{v\in\vl}\in\yadl$ be the image of $\by$ under $\sigma$. Let $\sigma_v(y_v)=(L_v^\sigma\into K_{\sigma_v(y_v)}, K_{\sigma_v(y_v)}^\flat\isom \cptv)\in \sY_{\cpt,L_v}$ be the untilt of $\cptv$ provided by $\sigma_v(y_v)$.  Let $(X/L_v,X/K_{y_v})$ (resp. $(X/L_{v}^\sigma,X/K_{\sigma(y_v)})$) be the pair of analytic spaces provided by the local holomorphoid $\hol{X/L_v}{y_v}{\C_{p_v}^\flat}$ (resp. $\hol{X/L^\sigma_v}{\sigma_v(y_v)}{\C_{p_v}^\flat}$)--here $L_v^\sigma$ is isomorphic to $L_v$, but the embeddings $L_v\into K_v$ and $L_v^\sigma\into K_{\sigma_v(y_v)}$ may not be topologically comparable, and hence the distinct notation). Let $\bL_v\subset K_v$ (resp. ${\bL}^\sigma_v\subset K_{\sigma_v(y_v)}$ be the algebraic closure of $L_v$ (resp. $L_v^\sigma$) in the respective (perfectoid) overfields.
Then one has an isomorphism (non-unique)
$$\pit{(X/L_v,X/K_{y_v})}\isom \pit{(X/L_v^\sigma,X/K_{\sigma(y_v)})},$$
and any such isomorphism provides an isomorphism
$$G_{L_v}\isom G_{L_v^\sigma},$$
and an isomorphism of each of the  prime-strips $\sF_{\scriptscriptstyle{Joshi}}^{?}$ in \Cref{tab:rosetta-stone3}. For example:
\be\xymatrix{
	\sF^{\vdash\times\mu}_{\scriptscriptstyle{Joshi}}(X/L_v,X/K_{y_v})\ar@2{-}[d]& &\sF^{\vdash\times\mu}_{\scriptscriptstyle{Joshi}}(X/L^\sigma_v,X/K_{\sigma_v(y_v)})\ar@2{-}[d]\\
	G_{L_v}\act \O_{\bL_v}^{\times\mu} \ar[rr]^\isom & & G_{L_v^\sigma}\act \O_{{\bL}^\sigma_v}^{\times\mu}}\ee and an isomorphism
\be\xymatrix{\sF^{\vdash\blacktriangleright\times\mu}_{\scriptscriptstyle{Joshi}}(X/L_v,X/K_{y_v})\ar@2{-}[d]& & \sF^{\vdash\blacktriangleright\times\mu}_{\scriptscriptstyle{Joshi}}(X/L^\sigma_v,X/K_{\sigma_v(y_v)})\ar@2{-}[d]\\
	G_{L_v}\act  \O_{\bL_v}^{\times\mu} \times q_{(X/L_v,X/K_v)}^\bN \ar[rr]^\isom &   & G_{L_v^\sigma}\act \O_{{\bL}^\sigma_v}^{\times\mu} \times q_{(X^\sigma/L_v^\sigma,X/K_{\sigma_v(y_v)})}^\bN.}\ee
\ethm

\brem
	In	 the parlance of \iutthr, \cite{mochizuki-gaussian}, the isomorphisms of prime-strips given by \Cref{th:gluing-prime-strips} can be stated as follows:  one has two distinct copies  of Mochizuki's prime-strips ${}^1\sF_{\scriptscriptstyle{Mochizuki}}^{?}$ and ${}^2\sF_{\scriptscriptstyle{Mochizuki}}^{?}$ which are identified via an abstract isomorphism	
	\be
	\xymatrix{ 
		{}^1\sF_{\scriptscriptstyle{Mochizuki}}^{?} \ar@2{-}[d] \ar[r]^\isom &  {}^2\sF_{\scriptscriptstyle{Mochizuki}}^{?} \ar@2{-}[d]  \\
		\sF^{\vdash\blacktriangleright\times\mu}_{\scriptscriptstyle{Joshi}}(X/L_v,X/K_{y_v}) &   \sF^{\vdash\blacktriangleright\times\mu}_{\scriptscriptstyle{Joshi}}(X/L^\sigma_v,X/K_{\sigma_v(y_v)}).
	}
	\ee
	However, \textcolor{red}{let me caution the reader that the two sides of this diagram make separate contributions to the set whose construction is asserted in \moccor. So this gluing does not really represent the creation of a new geometric object by gluing (as asserted in \cite[Example 2.4.7]{mochizuki-essential-logic}). That is why I have disagreed with claims of \cite[\ssep 3.4]{mochizuki-essential-logic} and \cite[Example 2.4.7]{mochizuki-essential-logic}. More precisely, for a workable Teichmuller Theory of the sort one wants to establish here, one must demonstrate that both sides  of the above diagram exist for intrinsic reasons and hence the analogy with projective line  \cite[Example 2.4.7]{mochizuki-essential-logic} is highly misleading.}
\erem

\subparat{Mochizuki's log-links}\label{ss:log-link-indentified}\nwss
The local theory of Mochizuki's  $\flog$-links  is detailed in \cite[Theorem 10.15.1]{joshi-teich-estimates} \constrtwo{Theorem }{th:flog-links}. Notably, it identifies the $\flog$-link at a fixed prime $v\in\vl$ explicitly. That assertion provides the  following global definition of Mochizuki's $\flog$ from the point of view of \present:
\bdefn
	Let $X/L$ be as defined in \cref{ss:elliptic-curve-assumptions} (or more generally, let $X$ be a geometrically connected, quasi-projective variety over a number field $L$ with no real embeddings). Let $\arith{L}_\by$ be an arithmeticoid of $L$. Let $\boldsymbol{\varphi}:\yadl\to \yadl$ be the global Frobenius morphism. Then the assignment 
	$$X/\arith{L}_\by\mapsto X/\arith{L}_{\boldsymbol{\varphi}(\by)}$$
	will be called \textit{Mochizuki's $\flog$-link} for $X/L$. 
\edefn

\brem\ 
	\benumlab	
	\item As established in \constrtwo{\ssep}{se:log-links}, Mochizuki's diagram of iterated $\flog$-links \cite[Remark 1.1.1]{mochizuki-iut3} corresponds to the diagram
	$$\begin{tikzcd}
	\cdots & X/\arith{L}_\by \arrow[d,rightsquigarrow] & X/\arith{L}_{\bvarphi(\by)} \arrow[d,rightsquigarrow] & X/\arith{L}_{\bvarphi^2(\by)} \arrow[d,rightsquigarrow] & \cdots \\
	\cdots & \by \arrow{r}{\bvarphi}[swap]{\flog} & \bvarphi(\by) \arrow{r}{\bvarphi}[swap]{\flog} & \bvarphi^2(\by)  & \cdots
	\end{tikzcd}$$
	\item I have demonstrated in \constrtwoh{\ssep }{se:heights} that heights are dependent on the choice of arithmeticoids used for computing them. Let $x\in X(\bL)$ be a closed point. As demonstrated in \constrtwoh{\ssep }{se:heights}, the (logarithmic) height of $x$:  $$h(x)_{\sL,\by},$$ computed using some ample line  bundle $\sL$ on $X$, and  with respect to the arithmeticoid $\arith{L}_{\by}$ can be compared with the height $$h(x)_{\sL,\boldsymbol{\varphi}(\by)}$$
	computed using $\arith{L}_{\boldsymbol{\varphi}(\by)}$. And one can even compare $h(x)_{\sL,\by}$ with $$h(x)_{\sL,\boldsymbol{\varphi}(\by)}, h(x)_{\sL,\boldsymbol{\varphi^2}(\by)}, h(x)_{\sL,\boldsymbol{\varphi^3}(\by)},\cdots.$$
	\eenum
\erem

\subparat{Mochizuki's $\Theta_{gau}$-Links}\label{ss:Mochizuki-theta-gau-link}\nwsss
My approach to Mochizuki's theta-links (specifically the $\Theta_{gau}$-links) is detailed in  \cref{def:mochizuki-adelic-ansatz}, \cref{ss:ansatz-local-def}, \Cref{th:adelic-theta-link} and \cite{joshi-teich-estimates}  and hence will not be repeated here. In \cite[Appendix]{joshi-teich-estimates}, I establish the existence of theta-links in the context of proofs of the Geometric Szpiro Inequality \cite{bogomolov00}, \cite{zhang01}.

\subparat{Mochizuki's ``Three Indeterminacies'' \cite[Theorem 3.11]{mochizuki-iut3}}\label{ss:Mochizuki-indeterminacies}
As I have pointed out in \cite{joshi-teich-summary-comments}, from the point of view of  classical Teichmuller spaces and Arithmetic Teichmuller spaces, if one holds the (tempered) fundamental group fixed, then the holomorphic  structure (i.e. the holomorphoid) giving rise to the said group, is necessarily indeterminate i.e the mapping 
$$\hol{}{y_v}{X/L_v}\mapsto \pit{X/L_v}$$
does not determine the holomorphoid i.e. the arithmetic holomorphic structure uniquely.

More precisely, if the topological isomorphism class of the fundamental group $\pit{X/L_v}$  is held fixed, then a holomorphoid  $\hol{}{y_v}{X/L_v}$ giving rise to this group is necessarily an indeterminate quantity as there are many distinct holomorphoids $\hol{}{y_v}{X/L_v}$ (in the classical case, there are many distinct complex structures) which provide the said group. 

\subsubpara{} So in my view, Mochizuki's Three Indeterminacies (see \cite[Theorem 3.11]{mochizuki-iut3}) which underlie his theory, are necessarily consequences of the existence of Arithmetic Teichmuller Spaces. Let me illustrate how this remark allows one to arrive at Mochizuki's three indeterminacies. I am listing them as they appear in the theory in the order they get applied in \cite[Theorem 3.11 and Corollary 3.12]{mochizuki-iut3} (the perspective is my own, Mochizuki's perspective is discussed in \cite[Theorem 3.11]{mochizuki-iut3}).

\begin{description}
	\item[Ind3] In \cite[Theorem A (aka Theorem 3.11), Page 425]{mochizuki-iut3}, this indeterminacy arises from the upper semi-compatibility of $\flog$-Links. From my point of view, this indeterminacy arises by \constrtwo{Theorem }{th:flog-kummer-correspondence} (reader may recall that in \constrtwo{\ssep}{se:log-links}, I have established that working with the Frobenius orbit of a point on $\syflv$ corresponds to working with a column of $\flog$-links). To understand this,	let $\by_0$ be the standard point \cref{ss:stand-point}. This choice is not canonical and one can also use $\bvarphi(\by_0)$ as a standard point. Write $\by_0=(y_v)_{v\in\vlp}$. One wants to understand how the choice of $\by_0$ over $\bvarphi(\by_0)$ affects arithmetic holomorphic structures and related Galois cohomology calculations. The  local aspect of this, at each $v$, is explicated in  \constrtwo{Theorem }{th:flog-kummer-correspondence}. To recall what \constrtwo{Theorem }{th:flog-kummer-correspondence} provides,  let $p_v$ be the residue characteristic of $v$. Then at each $v$, the choice of $y_v$ over $\varphi_v(y_v)$ means that the homomorphism
	$$\Z_{p_v}(1)\to B_{\cptv,\Q_{p_v}}$$
	which send a generator to the element $t\in B_{\cptv,\Q_{p_v}}$
	is uniquely determined only up to  a $p$-multiple. 
	\item[Ind2] This indeterminacy arises as a consequence of \constrone{Theorem }{th:main3},  \constrtwoh{Theorem }{th:galois-action-on-adelic-ff} and is the indeterminacy of the point $\by\in\yadlp$ which provides the arithmetic holomorphic structure at each $v\in\vlp$. One can write this Mochizuki style as indeterminacy with respect to $\Q_p$-linear isomorphisms $\sigma:\hgm(\O_{\C_{p_v}^\flat})\to \hgm(\O_{\C_{p_v}^\flat})$ (\constrone{Theorem }{th:main3} and \constrtwoh{Theorem }{th:galois-action-on-adelic-ff} which demonstrate how such an isomorphism changes arithmetic holomorphic structures). 
	
	In other words, Mochizuki's Ind2 is the indeterminacy  with respect to the choice of a point of the space of arithmetic holomorphic structures provided by the relevant Fargues-Fontaine curve $\sxqp$. The precise demonstration of the existence of Ind2 is an \textit{important contribution of \cite{joshi-teich}}. As remarked in \constrone{\ssep}{pa:virasoro-unif-remark} and \constrtwoh{Remark }{rm:global-frob-remark}, this indeterminacy is the arithmetic analog, at all primes, of the Virasoro action in the Geometric Langlands Program \cite{beilinson00b}.
	
	\item[Ind1] \textcolor{red}{This indeterminacy is what I have called anabelomorphy in \cite{joshi-anabelomorphy}. Mochizuki's contributed \constranab{\ssep}{ss:anabelomorphy-and-iut} to  clarify the role of anabelomorphy in \iut\  and \cite{joshi-anabelomorphy} provides other arithmetic consequences of anabelomorphy.} The  precise origins of Ind1 are quite arithmetic-geometric and not group theoretic as one might be led to believe because of the group theoretic methods of \iutthr. This indeterminacy Ind1 arises as follows. By \joshiros, one may view the Fargues-Fontaine curve $\sX_{\cptv,L_v}$ (and especially $\sY_{\cpt,L_v}$) as a sort of moduli for arithmetic holomorphic structures. By \cite[Theorem 2.1]{joshi-gconj}, this moduli is not uniquely determined by the absolute Galois group $$G_{L_v}=\pi_1^{et}(\sX_{\cptv,L_v}).$$ So one must consider all Fargues-Fontaine curves $\sX$ with $\pi_1^{et}(\sX)=G_{L_v}$ as providing arithmetic holomorphic structures which must be included to ensure the largest collection of arithmetic holomorphic structures over which one calculates averages in \moccor. \textit{One way to understand this is that in \joshiros, one may think of $\sX_{\cptv,L_v}$ as a parameter space for arithmetic holomorphic structures up to an action of Frobenius, then  $\sX_{\cptv,L_v}$ is itself not uniquely identified by its \'etale fundamental group!} [Also see \constrtwoh{Remark }{rm:global-frob-remark}.] The following proposition is useful:
\end{description}

\bpro\label{le:finite-anab-classes} 
	Let $E$ be a $p$-adic field contained in an algebraic closure $\bQ_p$ of $\Q_p$. Then there are only finitely many $p$-adic fields $E'\subset\bQ_p$ such that $E'$ is anabelomorphic to $E$. In particular, for a fixed geometric point $\bQ_p$ (fixed for computing all the relevant fundamental groups), there are only finitely many isomorphism classes of Fargues-Fontaine curves $\sX_{\cpt,E'}$ which are anabelomorphic to $\sX_{\cpt,E}$ i.e  such that 
	$$\pi_1^{et}(\sX_{\cpt,E'}) = G_E.$$
\epro
\bp 
By \cite{joshi-anabelomorphy}, if $E$ is a $p$-adic field, then its absolute degree $[E:\Q_p]$ is an amphoric quantity i.e. for any $p$-adic field $E'$ which is anabelomorphic to $E$, one has $[E:\Q_p]=[E':\Q_p]$. By standard results--for example Krasner's Lemma, there are only a finitely many subfields of $\bQ_p$ of a bounded degree over $\Q_p$. So the assertion follows from the main theorem of \cite{joshi-gconj}.
\ep

\subsubparat{Canonical Interpretation of Mochizuki's Indeterminacies} One may describe these indeterminacies in the following canonical way:
\begin{description}
	\item[Crystalline Indeterminacy or Mochizuki's Ind3] The homomorphism $\Z_p(1)\to B\to \bcris$ is determined uniquely only up to a $p$-multiple. From my point of view, this indeterminacy corresponds to working with Arithmetic Holomorphic Structures arising from the Frobenius orbit of a point in the Fargues-Fontaine curve $\syflv$ at all primes $v$.
	\item[Inderminacy of Arith. Hol. Strs. or Mochizuki's Ind2] For each $v\in\vl$, the local holomorphoid of $X/L_v$ is not uniquely determined by $\pit{X/L_v}$.
	\item[Anabelomorphy at all primes or Mochizuki's Ind1] This arises from the fact  that the isomorphism class of $G_{L_v}$ (for each prime $v$) does not uniquely determine $L_v$ in its isomorphism class.  This was already clarified by Mochizuki who contributed \cite[Paragraph 1.6]{joshi-anabelomorphy}. My work provides a more canonical view of this: this indeterminacy arises from the fact \cite{joshi-gconj} that the parameter space $\sxfe$ for arithmetic holomorphic structures is not uniquely determined by it \'etale fundamental group.
\end{description}
\subsubparat{\Cref{tab:rosetta-stone4}}
Finally one has  \Cref{tab:rosetta-stone4} which provides a translation between several objects central to \moccor.

\vskip0.5in

\begin{minipage}{\textwidth}
	\begin{rosetta}[label=tab:rosetta-stone4]{}
\renewcommand{\arraystretch}{1.5}
		\begin{center}	
		\begin{tabular}{|p{0.4\textwidth}|p{0.4\textwidth}|}
			\hline 
			\hfill	\iutthr \hfill & \hfill\joshiros\hfill \\ 
			\hline
			$\flog$-Link & $X/\arith{L}_{\by}\to X/\arith{L}_{\bvarphi(\by)}$ \\ 
			\hline
			$\Theta_{gau}$-Link & $\left\{\bz_w\in\tSigma_{\C_{p_w}^\flat,L_w'}\right\}_{w\in\ubblvossp}$ \\ 
			\hline
			distinct global and local Kummer Theories & global and local Kummer Theory provided by each global  holomorphoid $X/\arith{L}_\by$\\ 
			\hline
			local Isom-Groups $\left\{G_v\act \bL_v^{*\mu} \right\}_{v\in\vl}$ & $\left\{G_{L_v}\act B_{p_v}^{\varphi=p_v} \right\}_{v\in\vl}$ \\ 
			\hline
			distinct Hodge Theaters $\mathcal{HT}^?$ of each type $?$  & Hodge Theater $\mathcal{HT}^?_{\by}$ of each type $?$ provided by each arithmeticoid $\arith{L}_{\by}$\\
			\hline
			A frobenioid of the number field $L$ & A frobenioid given by each arithmeticoid $\arith{L}_{\by}$\\
			\hline
		\end{tabular}
	\end{center}
	\end{rosetta}
\end{minipage}
\clearpage

\section{Construction of $\thetam^\bsI$ and $\thetaj^\bsI$ using Mochizuki's tensor packets of $\log$-shells}\label{se:mochizuki-construction-thetam}
Mochizuki's construction of $\thetam^\bsI$ requires constructing his codomain $\bsIm\subset \bsImq$) for receiving theta-values.  This codomain $\bsIm$ (resp. $\bsImq$) is constructed in \cite[Section 3]{mochizuki-iut3} using Galois cohomology (Mochizuki's notation is some what more complicated than mine). It is also possible to work with somewhat simpler codomains $\bsImj\subset \bsImjq$ constructed in Mochizuki's style. The locus theta-values locus $\thetaji$  is somewhat simpler than $\thetami$ and is adequate for establishing \moccor. The construction of $\bsIm$ (resp. $\bsImq, \bsImj, \bsImjq$) is an adelic construction and carried out separately at each prime and requires some additional preparation.  Once one has constructed the codomain $\bsIm$ (resp. $\bsImq,\bsImj, \bsImjq$), one may consider theta-values locii in $\bsIm$ (resp. $\bsImq,\bsImj, \bsImjq$). 

\subparat{Bloch-Kato subspaces in Galois cohomology}\label{ss:log-shells} \textcolor{red}{Readers are advised to read \constrtwoh{\ssep}{se:arith-and-galois-cohomology} before reading this section.} As established  there, Galois cohomology of a number field depends on the choice of an arithmeticoid. This point will be used throughout.

Mochizuki's construction of the codomain $$\bsIm$$ for constructing $\thetami$ is based on his notion of log-shells \cite[Definition 5.4(iii)]{mochizuki-topics3}. Here I will provide a more intrinsic description of log-shells  than is done in \cite{joshi-anabelomorphy} and is based on Bloch-Kato subspaces in Galois cohomology. To do this I recall the definitions of Bloch-Kato subspaces in Galois cohomology.

\subsubparat{The non-archimedean case}
Let $E\supset \Q_p$ be a $p$-adic field. Fix an algebraic closure $\bE$ of $E$, let $G_E=\gal(\bE/E)$ and let $\ebh$ be its completion equipped with its canonical action of $G_E$. By  well-known classical results  the prime $p$, and the topological groups $\O_E^*\subset E^*$ are amphoric (see \cite[Theorem 3.3]{joshi-anabelomorphy} or \cite{jarden79}). 

My next point is important. By \cite[Proposition 1.2.1]{mochizuki04} or \cite[Proposition 4.2]{hoshi-topics} one sees that the triple $G_E\act \mu(\bE)$ consisting of  the absolute Galois group $G_E$, the subgroup $\mu(\bE)\subset \bE^*$  consisting of the roots of unity in $\bE$ equipped with the natural action of $G_E$, is amphoric and any isomorphism of topological groups $$G_{E_1}\isom G_{E_2}$$ of local fields $E_1,E_2$ (with chosen algebraic closures $\bE_1,\bE_2$) provides an isomorphism of the corresponding triples $$\left(G_{E_1}\act \mu(\bE_1)\right) \isom \left(G_{E_2}\act \mu(\bE_2)\right)$$
and hence an isomorphism
$$\left(G_{E_1}\act \hat{\Z}(1)_{E_1}\right) \isom \left(G_{E_2}\act \hat{\Z}(1)_{E_2}\right).$$ This should help understanding how one approaches Galois cohomology from the point of Arithmetic Teichmuller Spaces (for more details see \cite{joshi-teich-def}).

Now to my third point. One has an isomorphism of $\Q_p$-vector space (\cite{perrin-riou1994}): 
\be\Ext^1_{G_E}(\Q_p(0),\Q_p(1)) = H^1(G_E,\Q_p(1)).\ee

\brem By \cite[Proposition 1.4]{mochizuki-theta}, Mochizuki's Galois theoretic theta classes live in $\Ext^1_{G_E}(\Q_p(0),\Q_p(1))=H^1(G_{E},\Q_p(1))$. The archimedean case is discussed below. 
\erem 

Recall further that one has the decomposition of $\Q_p$-vector spaces (see \cite{perrin-riou1994}):
\be
H^1(G_E,\Q_p(1))  \isom H^1_f(G_E,\Q_p(1))\oplus \Q_p \isom E\oplus \Q_p.
\ee
This isomorphism is explicitly given as $q\mapsto (\log_E(q),v_E(q))$, so this isomorphism tracks the valuation of $q$ separately.

Following \cite{bloch90,nekovar93,perrin-riou1994}, write
\begin{align*}
H_f^1(G_E,\Q_p(1)) &=\ker(H^1(G_E,\Q_p(1))\to H^1(G_E,\bcris)),\\
H_{st}^1(G_E,\Q_p(1)) &=\ker(H^1(G_E,\Q_p(1))\to H^1(G_E,B_{st})),\\
H_g^1(G_E,\Q_p(1)) &=\ker(H^1(G_E,\Q_p(1))\to H^1(G_E,\bdr)),\\
H_e^1(G_E,\Q_p(1)) &=\ker(H^1(G_E,\Q_p(1))\to H^1(G_E,\bcris^{\varphi=1})),
\end{align*}
and  for $*\in\{e,f,st,g \}$ define the integral versions as
$$H_*^1(G_E,\Z_p(1)) =H^1(G_E,\Z_p(1))\cap H^1_*(G_E,\Q_p(1)).$$

The first three subspaces are the subspaces corresponding respectively to subspaces consisting of crytalline, semi-stable and de Rham extensions of $\Q_p(1)$ by $\Q_p(0)$  and $H_e^1(G_E,\Q_p(1))$ 
consists of crystalline, Frobenius trivial extensions of $\Q_p(1)$ by $\Q_p(0)$ contained  in $H^1_f(G_E,\Q_p(1))$.
From \cite{bloch90,nekovar93,perrin-riou1994} one further has:
 \be\label{eq:fontaine-subspaces-identified} H_e^1(G_E,\Q_p(1))=H_f^1(G_E,\Q_p(1))\subsetneqq H^1_{st}(G_E,\Q_p(1))= H^1_g(G_E,\Q_p(1))= H^1(G_E,\Q_p(1)).\ee
The first equality will play an important role in what follows.

\subsubparat{The archimedean case} [\textcolor{red}{It will be useful to read my discussion of Schottky uniformization in \ssep\ref{se:appendix-geom-case}.} 
Mochizuki does not use the Schottky uniformization approach for his archimedean calculations, but my view is that working with Schottky uniformization at the archimedean primes sets the archimedean and the non-archimedean primes of split-multiplicative reduction on par because both have the same type of uniformization i.e. the Schottky uniformization at archimedean primes corresponds to Tate uniformization at non-archimedean primes of split multiplicative reduction. 

Following \ssep\ref{se:appendix-geom-case},  for any archimedean prime $v\in\vlarch$ of a number field $L$ (assuming as in \inithtdata, that $L$ has  no real embeddings), one has $L_v\isom \C$ and I will write:
\be H^1(L_v,\Z(1)):=Ext^1_{\Z-MHS}(\Z(0),\Z(1))\isom L_v^*\isom \C^*.\ee

As explained in \constrtwoh{\ssep}{ss:hodge-str-schottky} this is the archimedean analog of Galois cohomology (for $v\in\vlnon$) $H^1(G_{L_v},\Q_p(1))$ considered above. 

To keep the theory at archimedean primes further aligned with the theory at non-archimedean primes, I define the \textit{archimedean Frobenius-invariant Bloch-Kato  subspace}  to be the multiplicative submonoid of $\C^*$ given by
\be 
H^1_{e}(\C,\Z(1))=\left\{ q\in Ext^1_{\Z-MHS}(\Z(0),\Z(1)): 0<\abs{q}_{\C}\leq 1\right\} \subset H^1(\C,\Z(1)).
\ee

Notably, if $q$ is a Schottky parameter of an elliptic curve over $\C$, then $$q\in H^1_{e}(\C,\Z(1))\subset H^1(\C,\Z(1))=Ext^1_{\Z-MHS}(\Z(0),\Z(1)).$$

\subparat{Mochizuki's Log-shells from the point of view of \cite{joshi-teich}}\label{ss:log-shells-joshi}
\subsubpara 
\textit{Mochizuki's log-shell} used in \cite{mochizuki-iut3} should be viewed as a rule on the arithmetic Teichmuller space $\fjxe$:
\bdefn\label{def:mono-analytic-logshell1}
Suppose $\holt{X/L}{\by}$ is a holomorphoid of $X/L$ and $v\in\vlnon$ is a prime of $L$ with $v|p$ and $\holt{X/L_v}{y_v}$ is the local holomorphoid of $X/L$ at $v$ given by $\holt{X/L}{\by}$. Then
the \textit{Mochizuki's log-shell given by $\holt{X/L}{\by}$ at $v$}   is the $\Z_p$-module given as follows
\be\label{eq:mono-analytic-log-shell} 
\logsh{(X/L_v,\xan/K_v)}=\begin{cases}\frac{1}{p}\log(\O_{L_v}^*)  \subset  \frac{1}{p}\log(\O_{K_v}^*) & \text{ if } p\geq 3,\\
						 \frac{1}{4}\log(\O_{L_v}^*) \subset \frac{1}{4}\log(\O_{K_v}^*) & \text{ if } p=2.
						\end{cases}
\ee
The \textit{Mochizuki's $\Q_p$-log-shell given by $\holt{X/L}{\by}$ at $v$} is the $\Q_p$-vector space 
 \be\logshqp{(X/L_v,\yan/K_v)}=\logsh{(X/L_v,\yan/K_v)}\tensor_{\Z_p}\Q_p.\ee
\edefn

To connect with \cite{mochizuki-iut3}, one uses the following:
\blem\label{le:mono-analytic-logsh} 
The log-shell $\logsh{(X/L_v,\xan/K_v)}$ (\Cref{def:mono-analytic-logshell1}) is the \textit{mono-analytic log-shell} \cite[Proposition 1.2(vi)]{mochizuki-iut3} given by the arithmetic holomorphic structure $\holt{X/L}{\by}$.
\elem
\bp 
This is immediate from Mochizuki's definition of mono-analytic log-shells (for example \cite[Proposition 1.2(vi)]{mochizuki-iut3}) and the definition of arithmetic holomorphic structures in \cite{joshi-teich,joshi-untilts}.
\ep

An important property of the $\Q_p$-logshells is following identification with the Bloch-Kato subspaces:
\bpro\label{pr:fontaine-subspaces-identified} For any $v\in\vlnon$ one has
\benumlab 
\item $\logsh{(X/L_v,\xan/K_v)}$  is a $\pit{X/L_v;K_v}$-amphoric  (and hence also $G_{L_v;K_v}$-amphoric) $\Z_p$-module.
\item $\logshqp{(X/L_v,\xan/K_v)}$ and $H^1_e(G_{L_v;K_v},\Q_p(1))$ are $\pit{X/L_v;K_v}$-amphoric  (and hence also $G_{L_v;K_v}$-amphoric) $\Q_p$ vector spaces.
\item One has  $$\logsh{(X/L_v,\xan/K_v)}=H^1_f(G_{L_v;K_v},\Z_p(1)) = H^1_e(G_{L_v;K_v},\Z_p(1)),$$ and
$$\logshqp{(X/L_v,\xan/K_v)}=H^1_f(G_{L_v;K_v},\Q_p(1)) = H^1_e(G_{L_v;K_v},\Q_p(1)).$$
\item In particular, $\logsh{(X/L_v,\xan/K_v)}$ (resp.  $\logshqp{(X/L_v,\xan/K_v)}$) consists of integral (resp. rational) Frobenius-invariant Galois cohomology classes.
\eenum
\epro 
\bp 
The first and the second assertions are immediate from the amphoricity of $\O_{L_v}^*$ and of the prime number $p$ (see \constranab{Theorem }{th:third-fun-anab}, \constranab{Theorem }{th:anab-fontaine-sub} and \constranab{Remark }{re:log-shell-and-fontaine-subspace}). The third is immediate from \eqref{eq:fontaine-subspaces-identified} which also yields the amphoricity of $H^1_e(G_{L_v;K_v},\Q_p(1))$. The fourth assertion is a consequence of the third.
\ep

\brem These assertion are important in the context of \iutthr: Thanks to the explicit identification of the Frobenius morphism $\bvarphi:\yadlp'\to\yadlp'$ of $\yadlp'$ and $\flog$-Links in \cref{ss:log-link-indentified} one sees that  $H^1_e(G_{L_v;K_v},\Q_p(1))$ which serves as a container for such cohomology classes of interest is stable under iterates of $\flog$-Links.  This stability property \eqref{pr:fontaine-subspaces-identified}  underlies Mochizuki's  choice of $\logshqp{(X/L_v,\xan/K_v)}=H^1_e(G_{L_v;K_v},\Q_p(1))$ as the foundation for constructing the codomain for defining the theta-values locus in \cite[Section 3]{mochizuki-iut3}. [This proposition is not established in \iutthr.] This codomain is constructed in the next few paragraphs.
\erem

\newcommand{\bsi}{\boldsymbol{\mathscrbf{I}}}
\newcommand{\bsiq}{\boldsymbol{\mathscrbf{I}^{\Q}}}
\newcommand{\blogsh}[1]{\bsi(#1)}
\newcommand{\blogshq}[1]{\bsi^{\Q}(#1)}

\subsubparat{Adelic log-shells} It is convenient to work adelically at this point. 
\bdefn 
Let $\holt{X/L}{\by}$ be a holomorphoid of $X/L$. Then \textit{Mochizuki's adelic log-shell} is the topological group (with product topology) given by: 
\be \blogsh{\holt{X/L}{\by}} =\prod_{v\in\vl} \logsh{\holt{X/L_v}{y_v}},\ee
and \textit{Mochizuki's rational adelic log-shell} is given by
\be \blogshq{\holt{X/L}{\by}} =\prod_{v\in\vl} \logshqpv{\holt{X/L_v}{y_v}}.\ee
\edefn
\bpro\
One has the following identification of Mochizuki's adelic log-shell and Mochizuki's rational adelic log-shell:
\be
\begin{tikzcd}
\logsh{\holt{X/L}{\by}} \ar[r,equal]\ar[d,hook] & H^1(\arith{L}_{\by},\Z(1))\ar[d,hook] \\
\logshq{\holt{X/L}{\by}} \ar[r,equal] & H^1_e(\arith{L}_{\by},\Q(1)).	
\end{tikzcd}
\ee
\epro
\bp 
This is immediate from the definition of the adelic cohomology of an arithmeticoid given in \constrtwoh{\ssep}{se:arith-and-galois-cohomology} and \Cref{pr:fontaine-subspaces-identified}.
\ep

\subparat{The Banach space $\frac{1}{t}B^{\vphi=p}$ as the log-shell for \cite{joshi-teich}}
The theory  developed in \present\ one has a more natural definition of a log-shell. This paragraph provides the detail of this assertion.
\subsubparat{The subspace $B^{\vphi=p}$} Let me begin by recalling from \cite{fargues-fontaine}:
\bthm\label{th:B-is-gen-in-deg-one}\,
\benumlab
\item One has an isomorphism of $\Q_p$-Banach spaces 
$$\begin{tikzcd}
\hgm(\O_{\cpt}) = 1+\fm_{{\cpt}} \ar [r, "\isom"] & B^{\vphi=p}  
\end{tikzcd}$$
given by $1+x\mapsto \log([1+x])$.
\item The ring $B$ is generated as a $\Q_p$-algebra by $B^{\vphi=p}$.
\eenum
\ethm
\bp 
These are established in \cite[Proposition 4.4.6]{fargues-fontaine} and \cite[Th\'eor\`eme 6.2.1]{fargues-fontaine} respectively.
\ep

\subsubpara{} It is precisely because of the \Cref{th:B-is-gen-in-deg-one} and Mochizuki's definition of the log-shell that I have asserted in \cite{joshi-teich-summary-comments} that $B^{\vphi=p}$ is the analog of Mochizuki's log-shell in \cite{joshi-teich,joshi-teich-estimates,joshi-teich-summary-comments}. Here is a more precise construction.
\bdefn
For any prime $p<\infty$, let 
\be 
\sI_{joshi,p}=\frac{1}{t}B^{\vphi=p}.
\ee
Then $\sI_{joshi,p}$ will be referred to as  the log-shell for $B[1/t]$. 
\edefn 

From  \cref{ss:Be} recall the ring $B_e=B[{1}/{t}]^{\vphi=1}$,
that is, $B_e$ is the $\Q_p$-algebra of Frobenius-invariant elements of $B[{1}/{t}]$ (generated by $\frac{1}{t}B^{\vphi=p}$).

The geometric meaning of $\sI_{joshi,p}$  is explicated by the following.
\blem\label{le:frob-tab-log-shell}\  
\benumlab
\item $$\sI_{joshi,p}=\frac{\log(1+\fm_{{\cpt}})}{t}.$$
\item For all $n\in \Z$ one has
$$\vphi^n(\sI_{joshi,p})\subset \sI_{joshi,p}.$$
\item One has $$
\sI_{joshi,p}=\frac{1}{t}B^{\vphi=p}=B_e^{deg\leq 1}.$$
\item 
Notably, the log-shell $\sI_{joshi,p}$ consists of functions on $\sxqp$ of degree $\leq 1$ which are inariant under the Frobenius morphism of $B[1/t]$.
\eenum
\elem
\bp 
The first assertion is immediate from \Cref{th:B-is-gen-in-deg-one} and the definition of $\sI_{joshi,p}$.
It is enough to prove the second assertion for $n=1$. First note that if $\frac{b}{t}\in \frac{1}{t}B^{\vphi=p}$ then \be\vphi(\frac{b}{t})=\frac{\vphi(b)}{\vphi(t)}=\frac{p\cdot b}{p\cdot t}=\frac{b}{t}\ee and so $\frac{b}{t}$ is a $\vphi$ invariant or Frobenius fixed element of $B$  (also see \cite[Th\'eor\`eme 6.5.2]{fargues-fontaine}). Thirdly, by  \cite{fargues-fontaine}, $\sxqp$ can be obtained as $$\sxqp={\rm Proj}(B_e)$$
and one has the equality
\be 
\sI_{joshi,p}=\frac{1}{t}B^{\vphi=p}=B_e^{deg\leq 1}.
\ee
 Hence the assertion is clear.
\ep

\brem Stability of $\sI_{joshi,p}$ under iterates of the Frobenius morphism provided by \Cref{le:frob-tab-log-shell} corresponds to the stability of Mochizuki's log-shells under iterates of Mochizuki's $\flog$-links (see \cite[Remark 1.1.1, Remark 1.2.2(iii)]{mochizuki-iut3} and \cite[Appendix]{joshi-teich-estimates}). 
\erem

\subparat{Construction of Mochizuki's Tensor Packet codomain $\bsIm$ (resp. $\bsImq$) for $\ttheta_{\scriptscriptstyle{Mochizuki}}^\bsI$}\label{ss:tensor-packet-codomain}
Let me briefly describe Mochizuki's choice of codomain or container for theta-values. Mochizuki's formalism of tensor packets is discussed in \cite[\ssep 3]{mochizuki-iut3}.  \textcolor{red}{Assumptions \cref{ss:elliptic-curve-assumptions}, \cref{ss:theta-data-fixing} will now be strictly in force.}

In \iutthr, Mochizuki bundles theta-values arising primes $w|p$ by means of tensor products of log-shells. This is  the analog in \iutthr\ of the formalism of bundling rings  considered earlier in \ssep\ref{ss:bundling-rings}. The formalism of bundling rings allows us to use the ring structure to multiply the theta-values at primes $w|p$ in a common ring. The formalism of tensor products and tensor packets \cite[Section 3]{mochizuki-iut3} provide this functionality in \cite{mochizuki-iut3}.

\subsubpara{} Let $\by_0=(y_w)_{w\in\vlp}$ be the chosen standard point (as chosen in \cref{ss:stand-point}). By \Cref{th:lmod-arith} one obtains a standard point $\ubyz\in\yadlmod'$ and hence an arithmeticoid $\arith{\lmod}_{\ubyz}$ of $\lmod$. This will be used in what follows.

Let $\ubyz=(z_w')_{w\in\ubblv}$ (recall that the bijection $\ubblv\isom \V_{\lmod}$ is fixed by \inithtdata). Mochizuki constructs his codomain for theta-values using tensor product, for all for $w|p$, of the log-shells  
\be \logsh{\holt{X/L}{z_w'}}=H^1_e(G_{L_w';K_{z_w'}},\Q_p(1))=H^1_e(G_{L_w';K_{z_w'}},\Q_p(1))
\ee  (see \cref{ss:log-shells}). 
With the choice of $\by_0$ and the induced choice $\ubyz$ given by the identification $\ubblv\isom\V_{\lmod}$ given by \inithtdata, I will drop the reference to the arithmeticoid $\arith{\lmod}_{\ubyz}$ in writing these groups and simply write
$$\sI(L_w')=\logsh{\holt{X/L}{z_w'}}.$$
The following auxiliary notation will also be used for brevity.
Let $\ubblv_p$ denote the set of primes $w\in L'$  such that $w\in\ubblv$ and $w|p$. Write
\begin{align}
	\sIp(L') & =\bigoplus_{w\in\ubblv_p}\logsh{L_w'}\\
	\sIpq(L') & =\bigoplus_{w\in\ubblv_p}L_w'.
\end{align}
where $\sIpq(L')$ is the $\Q_p$-vector space generated by each direct summand of $\sIp(L')$. Note that the product is over primes of $\ubblv$ and not over all primes of $\vlp$. If $L'$ is clear from the context, I will simply write $\sIp$, $\sIpq$ etc.

\subsubpara{} Of primary importance is the case of rational primes $p$ for which one has $\ubblvossp\neq\emptyset$. Since finite direct products and finite direct sums coincide for abelian groups and vector spaces, it will be convenient to write both of these multiplicatively instead of Mochizuki's method of writing additively. Direct products  are better for adelic situations in any case:
\be\sIp=\prod_{w\in\ubblv_p}\logsh{L_w'}\ee
and write 
\be\sIpq=\prod_{w\in\ubblv_p}\logshqp{L_w'}.\ee

The following definition will be useful in understanding $\prod_p \sIp$ in terms of the cohomology of the standard arithmeticoid $\arith{L'}_{\by_0}$ chosen in \cref{ss:tensor-packet-codomain}:

\bdefn 
Let $\V_{L'}\supset \ubblv\isom \V_{\lmod}$ be the bijection given by \inithtdata. Let 
$$\ul H^1_e(\arith{L_{\by_0}'},\Z(1))=\prod_{w\in \ubblv} H_e^1(\arith{L'_{\by_{0}}},\Z(1)), $$
and similarly define $\ul H^1_e(\arith{L_{\by_0}'},\Q(1))$.
\edefn

With this notation, one has the following identification:

\bpro\label{pr:galois-cohomology-section} 
Let $\arith{L'}_{\by_0}$ be  the standard arithmeticoid. Then for each choice of a bijection $\ubblv\isom \V_{\lmod}$ in \inithtdata,  one has 
\benumlab
\item the equality:
$$\ul H^1_e(\arith{L'}_{\by_0},\Q(1))=\prod_{p}\sIpq=\prod_p\prod_{w\in\ubblv_p} \sI^{\Q_p}(L_w').$$
\item One also has the restriction homomorphism
$$H^1_e(\arith{\lmod}_{\ubyz},\Q(1))\mapright{res}  \ul H^1_e(\arith{L'}_{\by_0},\Q(1)), $$
which on each local factor is the restriction homomorphism on Galois cohomology.
\item One also has the corestriction homomorphism
$$\ul H^1_e(\arith{L'}_{\by_0},\Q(1))\mapright{cores} H^1_e(\arith{\lmod}_{\ubyz},\Q(1))$$
which on each local factor is the corestriction homomorphism on Galois cohomology.
\item The composite of restriction and corestriction is multiplication by the local degree $[L_w':{\lmod}_{,v}]$ for each $\ubblv\ni w|v\in\V_{\lmod}$ i.e. $(q_v)\mapsto cores(res(q_v))$ is the mapping $(q_v)\mapsto (q_v^{[L_w':{\lmod}_{,v}]})$.
\eenum
\epro
\bp 
The proof is clear from the constructions, the definitions of these objects and standard properties of Galois cohomology.
\ep

\subsubpara{} 
Mochizuki's method in \iutthr\ relies on anabelian reconstruction theory and the proof requires recovering a lot of information about $X/L$ purely from the tempered fundamental groups of $X/L_v$ for $v\in\vl$. This necessitates tracking a lot more information than is needed for my methods above. This is not required in my approach. Instead of appealing to Mochizuki's Anabelian Reconstruction \topics, one can use the Rosetta-Stone \cref{ss:rosetta-stone} to translate the required objects from my theory and Mochizuki's Reconstruction Theory is not needed. Nevertheless  to be fully compatible with \iutthr, I will work with Mochizuki's setup for defining the codomain for his theta-values locus.

\subsubpara{}\label{ss:extended-theta-link} I will follow \cite[Section 3]{mochizuki-iut3} closely for the most part. It will be convenient to work with the graph of the correspondence which gives rise to Mochizuki's Adelic Ansatz \cref{se:adelic-ansatz} instead of working with Mochizuki's Adelic Ansatz. This means that instead of working with $$\bz=(\bz_1,\ldots,\bz_\ells)\in\tSigma_{L'}$$ one works with the $\ells+1$-tuple $$(\bz_0,\bz_1,\ldots,\bz_\ells)$$ where $$\bz_0\mapsto(\bz_1,\ldots,\bz_\ells)$$ is the correspondence on $\yadl$ giving rise to $\tSigma_{L'}$. Notably, if $\bz_{\Theta}$ is the standard point \cref{ss:stand-theta-link} of $\tSigma_{L'}$, then one has the $\ells+1$-tuple $$\bz_{\Theta}=(\by_0',\by_1',\ldots,\by_\ells')$$ associated with this point.

\subsubparat{Mochizuki's Procession of prime-strips} Following \cite[Section 3]{mochizuki-iut3}, for $0\leq j \leq \ells$, let $\bS_{j+1}=\{0,1,\ldots, j\}$  so that each set $\bS_{j+1}$ has exactly $j+1$ elements. Similarly define  $\bS_{j+1}^\divideontimes=\{1,\ldots, j\}$  so that each set $\bS_{j+1}^\divideontimes=\{1,\ldots, j\}$ has exactly $j$ elements. 

I will think of these sets as giving functions (with values in finite sets) on the set of $\ells+1$ tuples rising from $\tSigma_{L'}$. This is done  as follows. For each $\ells+1$ tuple $\bz=(\bz_0,\bz_1,\ldots,\bz_\ells)$ of arithmeticoids as above, one has, for each $j=0,1,\ldots,\ells$, the sets
 $$\bS_{j}(\bz)=\left\{ \bz_0, \bz_1, \ldots, \bz_{j} \right\}$$
 and 
 $$\bS_{j}^\divideontimes(\bz)=\left\{ \bz_1, \ldots, \bz_{j} \right\}.$$
 Notably 
 $$\bS_{\ells}(\bz)=\left\{ \bz_0, \bz_1, \ldots, \bz_{\ells} \right\}$$
 and
  $$\bS_{\ells}^\divideontimes(\bz)=\left\{ \bz_1, \ldots, \bz_{\ells} \right\}.$$ The following is an elementary consequence of this observation

\blem 
For each $\ells+1$-tuple  $\bz=(\bz_0,\bz_1,\ldots,\bz_\ells)$ (as constructed above) and arising from $\bz=(\bz_1,\ldots,\bz_\ells) \in \tSigma_{L'}$. Then one has the obvious chain of inclusions of sets of arithmeticoids
$$\bS_0(\bz)=\left\{\bz_0\right\} \subset \bS_1(\bz)=\left\{\bz_0,\bz_1\right\}\subset \cdots \subset \bS_{\ells}(\bz)=\left\{\bz_0,\bz_1,\ldots,\bz_\ells\right\},$$
and
$$\bS_1^\divideontimes(\bz)=\left\{\bz_1\right\}\subset \bS_2^\divideontimes(\bz)=\left\{\bz_1, \bz_2\right\}\subset \cdots \subset \bS_{\ells}^\divideontimes(\bz)=\left\{\bz_1,\ldots,\bz_\ells\right\}.$$
Applying Mochizuki's ``prime-strip'' constructions $\sF_{\scriptscriptstyle{Mochizuki}}^?$ of \cref{ss:prime-strips}, \cref{pr:prime-strips-constructed} and Rosetta Stone Fragment~\ref{tab:rosetta-stone3}, to these chains of subsets,  one obtains a procession of prime-strips in the sense of \cite[Theorem 3.11(i)]{mochizuki-iut3}.
\elem 

\brem Mochizuki's choice of the sets $A=\bS_{j+1}$ for $j=1,\ldots,\ells$ is dictated by his reconstruction theoretic approach in which Mochizuki reconstructs a specific instance of $X/L_w'$ from anabelian data. Notably in \cite{mochizuki-iut1}, $\bS_{\ells}$ is in bijection with the set of cusps of a certain tempered covering constructed from $X/L_v$ at a prime $v\in\vlnonp$. From the point of view of \present, anabelian reconstruction arguments are unnecessary and all functionality which is needed in establishing \moccor\ is  provided in my approach by working with distinct  holomorphoids $\holt{X/L'}{\by}$. 
\erem

\subsubparat{Tensor products of Mochizuki's log-shells} Given $\bz=(\bz_0,\bz_1,\ldots,\bz_\ells)\in \yadlp'$  and if $A$ is one of the sets $\bS_{j+1}(\bz)$ for $1\leq j \leq \ell-1$, obtained using this $\ells+1$-tuple of arithmeticoids, let 
\be\label{eq:tensor1}{}^A{\sIp}(L')=\bigotimes_{a\in A}\sIp(L')=\bigotimes_{a\in A}\left(\bigoplus_{w\in\ubblv_p}  \logsh{L_w'}_{\bz_a}\right)\ee
and 
\be\label{eq:tensor2}{}^A\sIpq(L')=\bigotimes_{a\in A}\sIpq(L')=\bigotimes_{a\in A}\left(\bigoplus_{w\in\ubblv_p}L_w'\right)_{\bz_a}.\ee

Along with \eqref{eq:tensor1} and \eqref{eq:tensor2} it is also convenient to consider product versions instead of tensor products. So the product variant is this:
\be\label{eq:prod-tensor1}\widetilde{{}^A\sIp}(L')=\prod_{a\in A}\sIp(L')=\prod_{a\in A}\left(\bigoplus_{w\in\ubblv_p}  \logsh{L_w'}_{\bz_a}\right)\ee
and \be\label{eq:prod-tensor2}\widetilde{{}^A\sIpq}(L')=\prod_{a\in A}\sIpq(L')=\prod_{a\in A}\left(\bigoplus_{w\in\ubblv_p}L_{w';\bz_a}\right).\ee
These come equipped with canonical homomorphism given by the construction of tensor products:
\be\label{eq:prod-tensor1-mapping}\widetilde{{}^A\sIp}(L')\to{}^{A}\sIp(L')\ee
and \be\label{eq:prod-tensor2-mapping}\widetilde{{}^A\sIpq}(L')\to{}^{A}\sIpq(L').\ee

\bdefn\label{def:mochizuki-adelic-tensor-log-shl} 
\textit{Mochizuki's adelic tensor product log-shell} is the adelic product \cite[Theorem 3.11(i)(c)]{mochizuki-iut3}:
\be\label{eq:codmain1}\sIm(L')=\prod_{p}\left(\prod_{j=1}^\ells {}^{\bS_{j+1}}\sIp(L')\right)  \subset  \prod_{p}\left(\prod_{j=1}^\ells {}^{\bS_{j+1}}\sIpq(L')\right)=\sImq(L') .\ee
I will often suppress the explication of $L'$ in the notation and simply write  $\sIm,\sImq$.
\edefn

\brem Let me remark that one has the restriction and corestriction homomorphisms for the extension $L'/\lmod$ given using the bijection $\ubblv\isom \V_{\lmod}$: 
\be 
\sIm(\lmod)\mapright{res}\sIm(L')\mapright{cores} \sIm(\lmod).
\ee  The composite, by  \Cref{pr:galois-cohomology-section}, leads to the appearance of certain powers in each factor which will become important in the definition of weighted adelic volumes here in \cref{ss:weighted-vol} and in \cite[Section 3]{mochizuki-iut3}.
\erem

\subsubparat{}\label{ss:mochizuki-prod-var}\nwsss Along with $\sIm$ and $\sImq$ one also has their product variants in the style of \eqref{eq:prod-tensor1} and \eqref{eq:prod-tensor2}:
\be\label{eq:codmain2}\tsIm=\prod_{p}\left(\prod_{j=1}^\ells \widetilde{{}^{\bS_{j+1}}\sIp}\right)  \subset  \prod_{p}\left(\prod_{j=1}^\ells \widetilde{{}^{\bS_{j+1}}\sI_p^\Q}\right)=\tsImq\ee
One also has the homomorphisms $$\tsIm\to \sIm$$ and $$\tsImq\to \sImq$$
which are obtained using the natural  homomorphism given by the tensor product construction.

\subsubsection{}\label{ss:joshi-var}\nwsss It will be useful to work with  simpler versions of $\sIm$ and $\tsIm$ given as follows. Let $\bz_{\Theta}=( \by_1,\ldots,\by_{\ells})$ be the standard point of Mochizuki's Ansatz constructed in \ref{ss:stand-theta-link}. This point is a $\Theta_{gau}$-Link in the terminology of \iutthr.
Taking $A=\{1,\ldots,\ells\}$ in Mochizuki's construction:
\be\label{eq:codmain-joshi}\sImj=\prod_{p}\left(\bigotimes_{j=1}^\ells \sI_{p,\by_j}\right)  \subset  \prod_{p}\left(\bigotimes_{j=1}^\ells \sI_{p,\by_j}^\Q\right)=\sImjq \ee

\brem By \Cref{pr:fontaine-subspaces-identified} the log-shells are amphoric. Hence one has abstract isomorphisms
$$\sI_{p,\by_j} \isom \sI_{p,\by_\ells} \qquad \text{ for } j=1,\ldots,\ells,$$
so one can write 
\be  
\sImj=\prod_{p}\left(\bigotimes_{j=1}^\ells \sI_{p}\right)
\ee
obtained by forgetting the $\by_j$ (or any  arithmeticoids giving rise to a specific log-shell) altogether (as \iutthr\ does). But keeping the $\bz_\Theta$ components in the above notation allows us to remember that the individual factors come equipped with the valuation scaling property given by \Cref{th:adelic-theta-link}. This property is crucial for establishing \moccor\ and also for the main theorem of \cite{mochizuki-iut4}.
\erem 

\subsubsection{}\label{ss:joshi-prod-var}\nwsss One also has the product (instead of tensor product) variants of $\sImj$ and $\sImjq$ in the style of \eqref{eq:prod-tensor1} and \eqref{eq:prod-tensor2}:
\be\label{eq:codmain-joshi2}\tsImj=\prod_{p}\left(\prod_{j=1}^\ells \sI_{p,\by_j}\right)  \subset  \prod_{p}\left(\prod_{j=1}^\ells \sI_{p,\by_j}^\Q\right)=\tsImjq\ee
or in Mochizuki's style as:
\be\label{eq:codmain-joshi3}\tsImj=\prod_{p}\left(\prod_{j=1}^\ells \sIp\right)  \subset  \prod_{p}\left(\prod_{j=1}^\ells \sI_p^\Q\right)=\tsImjq\ee

One also has natural projection $$\tsIm\onto \tsImj$$
and $$\tsImq\onto \tsImjq$$
obtained by projection onto the $(j+1)^{th}$ factor of each $\bS_{j+1}$ for $j=1,\ldots,\ells$.

This discussion together with \Cref{pr:galois-cohomology-section} and \Cref{th:lmod-arith} leads to the following explicit description:
\bpro Let $\bz_{\Theta}=(\by_1,\ldots,\by_{\ells})$ be the standard point of Mochizuki's Ansatz constructed in \ref{ss:stand-theta-link}. Using \Cref{pr:galois-cohomology-section} one may identify $\tsImj$ and $\tsImjq$ as $\ells$-tuples of Galois cohomology classes of the $\ells$-tuple of arithmeticoids $\arith{L'}_{\bz_\Theta}=(\arith{L'}_{\by_{1}}, \ldots,\arith{L'}_{\by_{\ells}})$ with respect to the bijection $\ubblv\isom \V_{\lmod}$ (obtained from \inithtdata) explicitly as
\begin{align}
\tsImj & =  \prod_{j=1}^\ells \ul H_e^1(\arith{L'}_{\by_{j}}, \Z(1)) \isom H_e^1(\arith{L'}_{\by_{\ells}}, \Z(1))^\ells,\\
\tsImjq & =  \prod_{j=1}^\ells \ul H_e^1(\arith{L'}_{\by_{j}}, \Q(1))^\ells \isom H_e^1(\arith{L'}_{\by_{\ells}}, \Q(1))^\ells.
\end{align}
\epro

\brem  
In \iutthr\ Mochizuki only considers the tensor product version $\sIm$ and $\sImq$ to define his theta-values locus. But the theta-values locus may also be defined in the product versions $\tsIm$ and $\tsImq$ \eqref{eq:prod-tensor1} and \eqref{eq:prod-tensor2}.
\erem

\subsubparat{} Notably, for $j=\ells$ expanding out the tensor product \eqref{eq:tensor2}  for ${}^{\bS_{j+1}}\sIp$
one gets 
\be{}^{\bS_{\ells+1}}\sIp=\bigotimes_{a\in \bS^{\ells+1}}\sIpq=\bigotimes_{a\in \bS^{\ells+1}}\left(\bigoplus_{w|p}L_w'\right)=\bigoplus_\alpha E_\alpha''\ee
where $E_\alpha''$ runs over a multi-set of $p$-adic fields. This is similar to the multi-set of $p$-adic fields which appear in \Cref{le:direct-sum-decomp2}. If one works with ${}^{\bS_{\ells+1}^\divideontimes}\sIp$ (instead of ${}^{\bS_{\ells+1}}\sIp$) then one obtains the multi-set of $p$-adic fields which appear in \Cref{le:direct-sum-decomp2}.

\brem In particular, this remark proves my assertion that tensor packet structure in \cite[Section 3]{mochizuki-iut3} is to compensate for the missing ring structure (in \iut) needed to take products of  the theta-values at distinct primes $w|p$. From my point of view this is directly available through the product structure of the period ring $B_p$ of \cite{fargues-fontaine}.
\erem

\subsubparat{} These construction lead to the following definition
\bdefn
Mochizuki's adelic tensor product log-shells  $\sIm(L')$, given in \eqref{eq:codmain1} as a subgroup of $\sImq(L')$ (resp. the group $\tsIm(L')$, given in  \eqref{eq:codmain2} as a subgroup of $\tsImq$) equipped with the topology induced from adeles (topology on each factor is the obvious one), are \textit{Mochizuki's tensor product codomains for  the theta-values  locus} (resp. \textit{Mochizuki's  product codomains for  the theta-values  locus}) (see \moccor\ or for an explicit assertion see \cite[Definition 7.16]{fucheng} of the tensor product log-shell). One has a homomorphism 
$$\tsImq\to \sImq.$$ 
If one wants to emphasize the role of the choice of the arithmeticoids $\arith{\lmod}_{\by_{0}}$ and especially $\arith{\lmod}_{\ubyz}$, then I will write this as
\edefn
 
 \brem\ 
 \benumlab
 \item This  construction is carried out in \cite[Proposition 3.1]{mochizuki-iut3} and $\sIm$ is called the \textit{holomorphic log-shell}. 
 \item One may also carry out these constructions for arbitrary $\arith{L'}_{\by}$ and the $\arith{\lmod}_{\underline{\by}}$. In this case, in \cite[Proposition 3.2]{mochizuki-iut3}, this is called \textit{mono-analytic log-shell}.
 \item In both the instances, Mochizuki's notation is a bit more complex and I will not use it here.
 \eenum
 \erem

\subparat{Mochizuki's evaluation of theta-values in log-shells} 
The assumptions \cref{ss:elliptic-curve-assumptions} and \inithtdata\ will now be in force. Let $\holt{X/L'}{\by}$ be a holomorphoid of $X/L'$. Let $w\in\vlp$. Let $\by=(y_w)_{w\in\vlp}$ and let $\holt{X/L_w'}{y_w}$ be the local holomorphoid of $X/L$ at $w$ given by $\by$. 

\subsubpara  By properties of $L'/L$ (available by \cref{ss:elliptic-curve-assumptions} and \inithtdata) one sees that the entire $2\ell$-torsion of $C/L'$ is defined over $L'$ and hence  the $2\ell$-torsion of $C/L_w'$ is  defined over $L_w'$. So one can apply Mochizuki's Galois theoretic evaluation of theta-values (\cite[Proposition 1.4]{mochizuki-theta}) to deduce that  these classes, say one of them is $\xi_w$, lands in 
$$\xi_w \in H^1(G_{L_w';K_w},\Z_{p_w}(1)_{K_w})$$ where $\Z_{p_w}(1)_{K_w}$ is the $G_{L_w';K_w}$-module of $p_w$-power roots of unity contained in $K_w$. However it is more useful for \iutthr\ to construct cohomology classes in  $H^1_e(G_{L_w';K_w},\Z_{p_w}(1)_{K_w})$ (with $\bcris$ constructed using the completed algebraic closure of $L_w'$ provided by $K_w$). 

\subsubpara Let me say that \textcolor{red}{one cannot view the classes $\xi_w$ provided by the Tate parameter in $H^1_e$. So \cite{mochizuki-iut3} resorts to another approach to this which consists of considering multiplicative action of $\xi_w$ on the log-shell see for instance \cite[Figure 3.1]{mochizuki-iut3} and \cite[Theorem 3.11(b)]{mochizuki-iut3}. This action emerges naturally in my approach in the proof of \Cref{th:moccor-Mochizuki-form1}}. 

This is remedied (in \ssep\ref{ss:const-galois-classes} below) by explicitly constructing Frobenius trivial, crystalline cohomology classes associated to any Tate parameter and hence also theta-values on any holomorphoid of any Tate elliptic curve. These  classes do live in Mochizuki's log-shells!

\subparat{Construction of crystalline, Frobenius invariant Galois cohomology class associated to a Tate elliptic curve}\label{ss:const-galois-classes}
Here by \textit{crystalline and Frobenius invariant Galois cohomology classes}, I mean elements of $H^1_e(-,-)$. One may use this term `crystalline, Frobenius invariant' because of the definition of $H^1_e$.

\subsubsection{}\nwsss Assume $C,L,X,L,\ell$ are as in \cref{ss:elliptic-curve-assumptions}, \inithtdata. Let $\holt{X/L'}{\by}$ be a holomorphoid of $X/L'$ and write $\by=(y_w)_{v\in\V_{L'}}$.  Assume $w\in \ubblvossp\neq\emptyset$. Let $K_w$ be the residue field of $y_w$ and let $L_w'\subset \lbh_w'\subset K_w$ be the algebraic closure of $L_w'$ in $K_w$. 
I want to construct a Galois cohomology class associated to $X/L_w'$ which lives in $H_e^1(G_{L_w';K_w},\Z_p(1))$. One can do this in many ways. My approach here is adapted from \cite{colmez1992}, \cite{colmez-p-adic-integration}. 

Let $\hgm/\Z_p$ be the multiplicative formal group over $\Z_p$. Then one has
\be\hgm(\O_{L_w'})=1+\fm_{\O_{L_w'}}\ee
is the group of $1$-units of $\O_{L_w'}$. 
By assumption $w\in\ubblvossp$ means that $C/L_w'$ is a Tate elliptic curve over $L_w'$ and hence it has semi-stable reduction at $w$. Choose a $\O_{L_w'}$-minimal model for $C/L_w'$. Let  $\hat{C}/\O_{L_w'}$ be the formal group  obtained from $C/\O_{L_w'}$ (by formal completion along the identity). Then one has an isomorphism
$$\hat{C}\isom \hgm$$
of $\O_{L_w'}$-formal groups. 

For each prime $p<\infty$, define $p^*\in \Z_p$ by setting
\be 
p^*=\begin{cases*}
	p & \text{ if } $p\geq 3$, \\
	4 & \text{ if } $p=2$.
\end{cases*}
\ee
In particular, if $w|p_w$ is the rational prime lying under $w\in\vlnonp$, then $p_w^*$ is defined by the above equation.

\blem 
The quantity $p_w^*$ is $G_{L_w'}$-amphoric.
\elem
\bp 
This is proved in \cite{mochizuki-topics3} and also in \cite[Proposition 3.6]{hoshi-mono}.
\ep

\subsubsection{}\nwsss Now one uses \assumptions. By the definition of $L'$, the $2\ells$-torsion points of $C/L$ are defined over $L'$ and hence by the inclusion $L'\subset L_w'$, these $2\ell$-torsion points are also defined over $L_w'$. In particular $L_w'$ contains the $2\ell^{th}$-root of any Tate parameter of $C/L$. The choice of a $2\ell^{th}$-root of the Tate-parameter $q_{(X/L_w';K_w)}^{1/2\ell}\in \O_{L_w'}$ provides a choice of a $1$-unit 
$$1+p^*\cdot q_{(X/L_w';K_w)}^{1/2\ell}\in \hgm(\O_{L_w'})=1+\fm_{\O_{L_w'}}\subset \O_{L_w'}^*.$$

By Tate's uniformization theorem one has $$L_w^{'*}/q_{(X/L_w';K_w)}^\Z=C(L_w')$$ and hence the class  $(1+p^*\cdot q_{(X/L_w';K_w)}^{1/2\ell})\bmod{q_{(X/L_w';K_w)}^\Z}$ provides an $L_w'$-rational point on $C/L_w'$:
$$(1+p^*\cdot q_{(X/L_w';K_w)}^{1/2\ell})\bmod{q_{(X/L_w';K_w)}^\Z}\in L_w^{'*}/q_{(X/L_w';K_w)}^\Z=C(L_w').$$
Choose a compatible collection of $$\xi_w=\left\{\sqrt[p^n]{1+p^*\cdot q_{(X/L_w';K_w)}^{1/2\ell}}\right\}_{n=0,1,2,\ldots}$$ of $p^n$-th roots of $1+p^*\cdot q_{(X/L_w';K_w)}^{1/2\ell}$. 

\brem 
By taking the $1$-units $1+p^*\cdot q^{j/2\ell}$ for $j=1,\ldots,\ells$, one may similarly construct crystalline, Frobenius invariant cohomology classes arising from Mochizuki's theta-values $q^{j/2\ell}$ \cite[Proposition 1.4]{mochizuki-theta}. 
\erem

\blem
In the above notation and hypothesis. The compatible system $$\xi_w=\left\{\sqrt[p^n]{1+p^*\cdot q_{(X/L_w';K_w)}^{1/2\ell}}\right\}_{n=0,1,2,\ldots}$$ defines a crystalline, Frobenius invariant Galois cohomology class
$$\xi_w\in H^1_e(G_{L_w';K_w},\Z_p(1)).$$
\elem
\bp This can be proved in many different ways. For example, one can use $p$-adic integration methods of \cite[Theorem II.3.21]{colmez-p-adic-integration} (and \cite{colmez1992}), this choice determines a cohomology class (denoted again by $\xi$)
$$\xi_w\in H^1_e(G_{L_w';K_w},\Z_p(1)).$$
[\textcolor{red}{The important point here is that this class lands in $H^1_e$.}] Alternately one can also use \cite{perrin-riou1994}, where it is established that $1$-units provide crystalline cohomology classes (and Frobenius invariance follows from the fact that  agree in the present situation $H_e^1(G_{L_w';K_w},\Z_p(1))=H_f^1(G_{L_w';K_w},\Z_p(1))$.
\ep

\subparat{The Bloch-Kato exponential and Mochizuki's log-shell}\label{ss:bloch-kato-exp}
\subsubpara 
Let \be\label{eq:bloch-kato-exp} e_{BK}:L_w'\mapright{\isom} H^1_e(G_{L_w'},\Q_p(1))\ee
be the Bloch-Kato exponential isomorphism of topological groups  (see \cite[10.6.19]{fargues-fontaine}) in which cohomology on the right is viewed as a multiplicative group. Explicitly $e_{BK}$  is the homomorphism $$x\mapsto \frac{e^{p^m\cdot x}}{p^m} \qquad m\gg 0.$$ Let $\log_p$ denote the principal branch of the $p$-adic logarithm (i.e. $\log_p(p)=0$).  

\brem\ 
\benumlab
\item  It is important to recall a basic fact about the $p$-adic exponential function. Let me remind the reader that the power series giving the $p$-adic exponential function $e^x=\sum\limits_{n=0}^{\infty}\frac{x^n}{n!}$  (here one can take $x\in\C_p$, with the usual normalization $\abs{p}_{\C_p}=\frac{1}{p}$) has radius of convergence $\abs{x}_{\C_p}<\abs{p}^{1/(p-1)}_{\C_p}=\left(\frac{1}{p}\right)^{1/(p-1)}$ (see \cite[Chapter IV, Section 1]{koblitz-p-adic-book}).  Hence for $p=2$ the disk of convergence is $\abs{x}_{\C_2}<\abs{2}_{\C_2}=\frac{1}{2}$ while for $p$ odd, the disc of convergence strictly contains the open disk $\abs{x}_{\C_p}<\frac{1}{p}$. This should help clarify the difference in the formulae given below when $w|2$. \item Importantly, the $p$-adic exponential does not converge  on $p$-adic units or even on $1$-units but the series defining $e^{p^m\cdot u}$ converges for any $p$-adic unit $u$ for suitably large $m$. This becomes important in calculating the inverse of the Bloch-Kato exponential.
\eenum
\erem

\subsubparat{The Bloch-Kato logarithm}\label{ss:bloch-kato-logarithm} Choose for all primes $p$ the branch of the $p$-adic logarithm $\log_p$ such that $\log_p(p)=0$.
Let 
\be\label{eq:bloch-kato-log} 
\log_{BK}:H^1_e(G_{L_w'},\Q_p(1))\mapright{\isom} L_w'
\ee
be the isomorphism inverse to the Bloch-Kato exponential $e_{BK}$  \cref{ss:bloch-kato-exp}  (eq.~\eqref{eq:bloch-kato-exp}). Then the relation between Mochizuki's log-shell \eqref{eq:mono-analytic-log-shell}, the Bloch-Kato logarithm of  $1$-units given by these constructions is given by \be\label{eq:bloch-kato-log-class} \log_{BK}(\xi_{(X/L_w';K_w)})=\frac{\log_p(p_w^*\cdot(1+p^*_w\cdot q_{(X/L_w';K_w)}^{1/2\ell}))}{p_w^*}=\frac{\log_p(1+p_w^*\cdot q^{1/2\ell}_{(X/L_w';K_w)})}{p_w^*}\in \logsh{(X/L_w',K_w)}.\ee
More generally if $\xi\in\O_{L_w'}^*$ with $\xi\cong1 \bmod (p^*\cdot{\O_{L_w'}})$  then $$\xi\in H^1_e(G_{L_w'},\Z_p(1))=H^1(G_{L_w'},\Z_p(1))\cap H^1_e(G_{L_w'},\Q_p(1)) \subset H^1(G_{L_w'},\Q_p(1)),$$
and its Bloch-Kato logarithm is  $$\log_{BK}(\xi)\in \logsh{(X/L_w',K_w)}.$$ Let me record this observation in the following:
\bpro 
Let $L,C,X,L'$ be as in \cref{ss:elliptic-curve-assumptions}, \inithtdata. Let $w\in\ubblvossp\neq\emptyset$ (so $p$ is odd). Then Mochizuki's log-shell  $\logsh{\holt{X/L'_w}{y_w}}$ \eqref{eq:mono-analytic-log-shell} can be identified with the image, under the Bloch-Kato logarithm $\log_{BK}$ \eqref{eq:bloch-kato-log}, of the Bloch-Kato $\Z_p$-module $$H^1_e(G_{L_w},\Z_p(1))=H^1(G_{L_w},\Z_p(1))\cap H^1_e(G_{L_w},\Q_p(1))\subset H^1(G_{L_w},\Q_p(1))$$
and   $$\log_{BK}(\xi_{(X/L_w';K_w)})\in \logsh{(X/L_w';K_w)}.$$ 
\epro
\bp 
The only point which needs to be established (given that $\log_{BK}$ is injective) is the surjectivity, but clearly any element of the form $\log(1+p_w^*\O_{L_w'})/p_w^*$ is evidently in the image of $\log_{BK}$.
\ep

\subsubpara \label{ss:bloch-kato-isom-extended} Consider the set of prime $w\in\V_{L',p}$. One can obviously extend the above isomorphism $e_{BK}$ to the product topological groups $\prod_{w\in \vlp} L_w'$ and a similar assertion for $\log_{BK}$  and deduce the isomorphism $$e_{BK}:\prod_{w\in\V_{L',p}}L_w' \isom  \prod_{w\in\V_{L',p}} H^1_e(G_{L_w'},\Q_p(1)),$$
$$\log_{BK}: \prod_{w\in\V_{L',p}} H^1_e(G_{L_w'},\Q_p(1)) \isom \prod_{w\in\V_{L',p}}L_w',$$
which I will refer to, again, as the isomorphism given by the Bloch-Kato exponential and the Bloch-Kato logarithm respectively.

\subsubpara Recall from \constrtwoh{\ssep}{se:arith-and-galois-cohomology}, the definition of the Galois cohomology of an arithmeticoid. Then one can assemble these local classes  $\{\xi_{(X/L_w';K_w)}\}_{w\in\vlp}$ into a class  
\be 
\boldsymbol{\xi}_{\holt{X/L'}{\by}}=(\xi_{(X/L_w';K_w)})_{w\in\vlp}\in H^1_e(\arith{L'}_\by,\Z(1))
\ee in the galois cohomology  of the arithmeticoid $\arith{L'}_\by$ (see \constrtwoh{\ssep}{se:arith-and-galois-cohomology}). This cohomology group depends on $\arith{L'}_{\by}$. For notational simplicity, I will hence forth write $\boldsymbol{\xi}_{\holt{X/L'}{\by}}$ as
$$
\boldsymbol{\xi}_{\by}.
$$

Explicitly $\boldsymbol{\xi}_{\by}=\boldsymbol{\xi}_{\holt{X/L'}{\by}}=\left(\xi_w\right)_{w\in\ubblv}\in H^1_e(\arith{L'}_\by,\Z(1))$ where
\be  
\xi_w=\begin{cases}
	\left\{(1+p_w^*)^{1/p_w^n}\right\}_{n\geq0} \in H^1_e(G_{L_w'},\Z_{p_w}(1)) & \text{if } w\in\vlp, \text{ and } w\not\in \ubblvoss,\\
	\left\{(1+p_w^*\cdot q^{1/2\ell}_{w})^{1/p_w^n}\right\}_{n\geq0}  \in H^1_e(G_{L_w'},\Z_{p_w}(1)) & \text{if } w\in\ubblvoss\\
	q_w \in Ext^1_{MHS}(\Z(0),\Z(1)) & \text{if } w\in\ubblv\cap \vlp^{arc}.
\end{cases}
\ee

\subsubparat{Collation of Galois cohomology classes}\label{ss:collation-of-classes} One would like to collate and compare the cohomology classes $\boldsymbol{\xi}_{\by}$ arising from different arithmeticoids $\by$ in the cohomology of  one arithmeticoid. This is where our choice of the standard arithmeticoid \cref{ss:stand-point} comes in and the following assertion establishes this collation process.

\bpro\label{pr:general-cohomology-set-const} 
The assumptions \cref{ss:elliptic-curve-assumptions} and \inithtdata\ on $C/L,X/L,L'\supset L$ hold. Let $A$ be a non-empty, finite indexing set. For any arithmeticoid $\arith{L'}_{\by}$ let $\arith{L'}_{\uby{\by}}$ be the arithmeticoid $\arith{\lmod}$ of $\lmod$ given by \Cref{th:lmod-arith}. Suppose $\Sigma=\left\{\holt{X/L'}{\by}:\by \in\yadlp' \right\}$ is a  collection of holomorphoids of $X/L'$ and suppose for each $\by\in\Sigma$ one is given a set of 
$$\left\{ \boldsymbol{\xi}_{\uby{\by},a}:a\in A \right\} \subset H^1_e(\arith{L'}_{\uby{\by}},\Z(1))$$ cohomology classes arising from $\holt{X/L'}{\by}$ in some fashion.  Let $\arith{L'}_{\by_0}$ be a standard arithmeticoid \cref{ss:stand-point}.  Then there is a set 
$$\boldsymbol{\Psi}_{\Sigma} \subset \prod_{a\in A} H^1_e(\arith{L'}_{\ubyz},\Z(1))_a,$$
(here $H^1_e(\arith{L'}_{\ubyz},\Z(1))_a=H^1_e(\arith{L'}_{\ubyz},\Z(1))$ for each $a\in A$ is simply a reminder that the product is over the index $a\in A$), 
consisting of the images of $\boldsymbol{\xi_{\uby{\by},a}}$ (for  each $a\in A$ and for each $\by\in \Sigma$) under all   (topological) isomorphisms 
$$H^1_e(\arith{L'}_{\uby{\by}},\Z(1)) \isom H^1_e(\arith{L'}_{\ubyz},\Z(1))$$ given by  \constrtwoh{Proposition }{pr:arith-isom-gal-cohom}.
\epro
\bp 
This is proved exactly as \constrtwoh{Proposition }{pr:collation-principle}.
\ep

\subparat{The definition of $\thetami$ and $\thetaji$}\label{ss:def-thetam}
My definition of the theta-values locus $\thetaj$ is here \cref{ss:definition-of-thetaj-using-ansatz} and of $\thetam$ is here \cref{ss:def-thetam}.

Let me make the following remarks which will help clarify the construction:
\benumlab
\item One important point to be noted is that $\thetami$ and $\thetaji$ refer respectively to slightly different codomains $\sIm$ and $\sImj$. Mochizuki uses Anabelian Reconstruction Theory in \iut, and this approach requires tracking a lot of additional data and this reflected in $\sIm$. My approach is construction as opposed to reconstruction and requires less information be tracked and leads to a simpler codomain  $\sImj$ (with a homomorphism  $\sIm\to\sImj$) for defining $\thetaji$.
\item The substantively non-trivial portion of \moccor\ is the construction of the Theta-values set $\thetami$ in $\sIm$ (or $\tsIm$). 
\item The bound asserted by \moccor\ is natural  modulo the existence of this set.
\item Both \cite{joshi-teich-estimates} and \cite{mochizuki-iut3} require working with a certain number of copies  of a fixed value group (explicitly $\abs{\C_p^\flat}\subset \R$) for each prime $p$. 
\item As discussed in \cite{joshi-teich-estimates}, in my theory one  works with copies of the value group of the tilt $|\cpt|$ (one can do this by \cite{scholze12-perfectoid-ihes} or \cite{fargues-fontaine}) and all the valuation calculations can be made in this chosen, fixed value group $\abs{\cpt}$ for each prime $p$. One may further assume that $\abs{\cpt}\subset \R$ if one needs to work with realified value groups (remember that \cite{fargues-fontaine} works with valuations taking values in $\R$).
\item In \cite{mochizuki-iut3}, neither tilts nor algebraically closed complete valued fields are available, and value groups are tracked separately by means of perfect and realified Frobenioids of $\Q_p$ (for each prime $p$). My discussion of perfect and realified Frobenioids is in \cref{se:frobenioids} especially \cref{ss:perfect-realified-frobenioids}. 
\item This provides the value group equal to the value group of $\bQ_p$ (which is equal to that of  $\C_p$ as $\C_p$ is an immediate extension of $\bQ_p$ by \cite{poonen93}) and $\C_p$ has value group equal to that of $\cpt$ by \cite[Lemma 3.4]{scholze12-perfectoid-ihes}.

\item The construction of Theta-locus provides the valuation data  arising from distinct arithmetic holomorphic structures which are comparable in appropriate common value group at each prime. 
\item Without the existence of distinct arithmetic holomorphic structures \cite{joshi-teich,joshi-untilts} and the construction of \cite{joshi-teich-estimates} this idea is difficult to understand and one can arrive at incorrect conclusion regarding \moccor.
\eenum

\subsubpara \label{ss:construction-thetam-thetaj} Now one can define $\ttheta_{\scriptscriptstyle{Mochizuki}}$. The construction is similar to that of $\ttheta_{\scriptscriptstyle{Joshi}}$ (in fact the latter is modeled on the former). \Cref{th:thetam-construction} is essentially \cite[Theorem 3.11]{mochizuki-iut3}.
Write  ${\bf 1+p^*}=\left\{\{(1+p_w^*)^{1/p_w^n}\}_{n\geq0}\right\}_{w\in\ubblv}$ is a compatible  choice of $p^n_w$-th roots of $\{1+p^*_w\}_{w\in\ubblv}$.
\newcommand{\sbxi}{\boldsymbol{\xi}}
\bthmdef\label{th:thetam-construction}
Let the assumptions regarding $L,C,X,L'$ made  in \cref{ss:elliptic-curve-assumptions}, \inithtdata\ be in force.
Let $\tSigma_{L'}$ be Mochizuki's Ansatz constructed in \eqref{eq:adelic-ansatz}. Then 
\benumlab
\item For each $\bz=(\by_1,\ldots,\by_{\ells})\in\tSigma_{L'}$ one obtains
\begin{enumerate}
	\item for any choice of decoration $?$ (as specified in \Cref{tab:rosetta-stone3}),  a tuple of prime-strips 
	$$(\sF_{\scriptscriptstyle{Joshi},v}^?(\holt{X/L'}{\by_1}),\ldots,\sF_{\scriptscriptstyle{Joshi},v}^?(\holt{X/L'}{\by_\ells}))_{v\in\V_{L'}},$$
	along with a corresponding tuple of Mochizuki's prime-strips (\Cref{tab:rosetta-stone3})
	$$(\sF_{\scriptscriptstyle{Mochizuki},v}^?,\ldots,\sF_{\scriptscriptstyle{Mochizuki},v}^?)_{v\in\V_{L'}},$$
	\item and a tuple of tempered Frobenioids (at each $v\in\V_{L'}$) (see  \cref{ss:mochizuki-temp-frob})
	$$\frob_{temp}(X/L')_{\bz}:=(\frob_{temp}(\holt{X/L'}{\by_1,v}),\ldots,\frob_{temp}(\holt{X/L'}{\by_\ells,v}))_{v\in\V_{L'}},$$
	\item and  for any choice of decoration $?$ (as specified in \iutthr),  the Hodge-Theaters (given by  \cref{ss:hodge-theater}, \eqref{eq:hodge-theater-def}, and  especially \cref{th:hodge-theaters-exist}) 
	$$ \mathcal{HT}^?(X/\arith{L'}_{\bz});$$
	\item and a tuple of Galois cohomology classes
	$$\sbxi_{\bz}^{\scriptscriptstyle{Joshi}}=(\sbxi_{\uby{\by}_1},\ldots,\sbxi_{\uby{\by}_\ells}) \in \prod_{j=1}^\ells H^1_e(\arith{\lmod}_{\uby{\by}_j},\Z(1)).	$$
	\item Write $$ {}^{\bS_{j}}H^1_e(\arith{L},\Z(1))=\prod_{i=0}^{j+1} H^1_e(\arith{L}_{\by_i},\Z(1))$$
	and write $(\by_0,\by_1,\ldots,\by_{\ells})$ for the extended tuple \cref{ss:extended-theta-link} associated to $\bz=(\by_1,\ldots,\by_{\ells})\in\tSigma_{L'}$ and write
	\begin{multline*}
	\sbxi_{\bz}^{\scriptscriptstyle{Mochizuki}}=(\overbrace{({\bf 1+p^*})}^{\in \bS_1},\overbrace{({\bf 1+p^*},\sbxi_{\uby{\by}_1})}^{\in\bS_2},\ldots,\overbrace{({\bf 1+p^*},{\bf 1+p^*},\ldots,{\bf 1+p^*},\sbxi_{\uby{\by}_\ells})}^{\in\bS_{\ells+1}}) \\ \in \prod_{j=1}^\ells {}^{\bS_{j}}H^1_e(\arith{L}_{\by_j},\Z(1));	
	\end{multline*} 
	\item  Write $$\log_{BK}(\sbxi_{\bz}^{\scriptscriptstyle{Mochizuki}})\in \tsIm$$ for the image of $\sbxi_{\bz}^{\scriptscriptstyle{Mochizuki}}$ under the Bloch-Kato logarithm \cref{ss:bloch-kato-logarithm}, \cref{ss:bloch-kato-isom-extended} and again write 
	$$\log_{BK}(\sbxi_{\bz}^{\scriptscriptstyle{Mochizuki}})\in\sIm$$ for its image in $\tsIm$ under the homomorphism $$\tsIm\to\sIm.$$
\end{enumerate} 
\item For each  $\sigma\in\mathbf{G}_{L'}$ (see \constrtwoh{Theorem }{th:galois-action-on-adelic-ff}) one has 
$$\sigma(\bz)=(\sigma({\by}_1),\ldots,\sigma(\by_\ells))\in \tSigma_{L'},$$
and 
$$\sigma(\sbxi_{\bz}^{\scriptscriptstyle{Joshi}})=\sigma((\sbxi_{\uby{\by}_1},\ldots,\sbxi_{\uby{\by}_\ells})):=(\sbxi_{\sigma(\uby{\by}_1)},\ldots,\sbxi_{\sigma(\uby{\by}_\ells)}) \in \prod_{j=1}^\ells H^1_e(\arith{L}_{\sigma(\by_j)},\Z(1)),$$
and a similar action of $\sigma\in\mathbf{G}_{L'}$ on $\sbxi_{\bz}^{\scriptscriptstyle{Mochizuki}}\in\tsIm$. [Note that action of $\mathbf{G}_{L'}$ is through its action on $\yadlp$, the groups on the right do not carry any action of $\mathbf{G}_{L'}$.]
\item By construction, $\tSigma_{L'}$ is stable under global Frobenius morphism $\bvarphi$ and 
$$\bvarphi(\sbxi_{\bz}^{\scriptscriptstyle{Joshi}})=\bvarphi((\sbxi_{\uby{\by}_1},\ldots,\sbxi_{\uby{\by}_\ells}))
:=(\sbxi_{\bvarphi(\uby{\by}_1)},\ldots,\sbxi_{\bvarphi(\uby{\by}_\ells)}) \in \prod_{j=1}^\ells H^1_e(\arith{L}_{\bvarphi(\by_j)},\Z(1));$$
and a similar action of $\bvarphi$ on $\sbxi_{\bz}^{\scriptscriptstyle{Mochizuki}}$ (which is left for the reader to explicate).
\item  Taking  the sets $$\Sigma=\tSigma_{L'},\ \{\log_{BK}(\sbxi_{\bz}^{\scriptscriptstyle{Joshi}}):\bz\in\tSigma_{L'} \} \text{ and } A=\{1,2,\ldots,\ells\}$$ in \Cref{pr:general-cohomology-set-const}, one obtains from \Cref{pr:general-cohomology-set-const},  the subset $\thetajti$ given as the convex closure of $$\thetajti:=\text{Convex Closure of }\boldsymbol{\Psi}_{\tSigma_{L'}}^{\scriptscriptstyle{Joshi}}\subset \tsImj.$$
\item Taking  the sets $$\Sigma=\tSigma_{L'},\ \{\log_{BK}(\sbxi_{\bz}^{\scriptscriptstyle{Mochizuki}}):\bz\in\tSigma_{L'} \} \text{ and } A=\prod_{j=1}^\ells \bS_{j}$$
one obtains from \Cref{pr:general-cohomology-set-const} the subset $\thetamti$ given as the convex closure of:
$$\thetamti:=\text{Convex Closure of }\boldsymbol{\Psi}_{\tSigma_{L'}}^{\scriptscriptstyle{Mochizuki}}\subset \tsIm\subset \tsImq.$$
\item Let $$\thetaji$$
be the image of $\thetajti$ under the homomorphism $\tsImjq\to\sImjq$; and let $$\thetami$$
be the image of $\thetamti$ under the homomorphism $\tsImq\to\sImq$.
\item Then $$\thetajti \subset \tsImj\subset \tsImjq$$ and  
$$\thetaji \subset \sImj\subset \sImjq$$
and especially 
$$\thetamti \subset \tsIm\subset \tsImq$$ and 
$$\thetami \subset \sIm\subset \sImq$$
are Mochizuki's multi-radial representations of Theta-values i.e. the  Theta-values locii of \cite[Theorem 3.11, Corollary 3.12]{mochizuki-iut3} contained in $\tsImjq$ and $\sImjq$  (resp. $\tsImq$ and $\sImq$).
\eenum
\ethmdef
\bp
 The rest of the assertions are easily assembled from already established facts and properties.  
\ep

\brem 
Despite the notation $\sbxi_{\bz}^{\scriptscriptstyle{Mochizuki}}$ these crystalline Galois cohomology classes are not used  by Mochizuki in \iutthr. His method uses a multiplicative action of theta-values on log-shells which is described in detail in \cite[Figure 3.1]{mochizuki-iut3}. As is clarified in \Cref{re:multiplicative-action}, this multiplicative action emerges in the proof of \Cref{th:moccor-Mochizuki-form1} given below.  However, 
one can obviously follow Mochizuki's approach to constructing the theta-values locus $\ttheta_{\scriptscriptstyle{Mochizuki}}^{{\mathscrbf{I}}}$. 
\erem

\bcor\label{cor:compactness-of-locii}
The subset $\thetaji$ (resp. $\thetami$) is contained in a compact subset of $\sImjq$ (resp. $\sImq$).
\ecor
\bp 
Since $\thetaji$ (resp. $\thetami$) are the images of  $\thetajti$ (resp. $\thetamti$) under the (continuous) homomorphisms, it is enough to prove   the assertion for the latter sets. This, in turn, is immediate from the fact that each of factors of the Galois cohomology of an arithmeticoid $H^1_e(\arith{L'},\Z(1))$ is compact and  the topology on $\tsImjq$ (resp. $\tsImq$) is the product topology. 
\ep

\subsubpara For an arithmeticoid $\arith{L'}_{\by}$, let $\by=(y_w)_{w\in\vlp}$ with $K_{y_w}$ the residue field of $y_w$. Then for any element $\sbxi_{\bz}^{\scriptscriptstyle{Joshi}}\in\thetajti$ (resp. $\sbxi_{\bz}^{\scriptscriptstyle{Mochizuki}}\in\thetamti$) one can define 
\be 
\abs{\log_{BK}(\sbxi_{\by})}=\prod_{w\in\vlp} \abs{\log_{BK}(\sbxi_{\by,w})}_{K_{y_w}}=\prod_{w\in\ubblvoss} \abs{\log_{BK}(\sbxi_{\by,w})}_{K_{y_w}}.
\ee

For a tuple $\sbxi_{\bz}^{\scriptscriptstyle{Joshi}}\in\thetajti$ (resp. $\sbxi_{\bz}^{\scriptscriptstyle{Mochizuki}}\in\thetamti$) one can define 
\be 
\abs{\log_{BK}(\sbxi_{\bz}^{\scriptscriptstyle{Joshi}})}=\prod_{j=1}^\ells\abs{\log_{BK}(\sbxi_{\uby{\by}_j})}.
\ee
and similarly define $\abs{\log_{BK}(\sbxi_{\bz}^{\scriptscriptstyle{Mochizuki}})}$.

\brem 
Note that by \Cref{le:elem-lemma}, $\abs{\log_p(1+p_w^*)}_{L_w'}=\abs{p^*}_{L_w'}$ for all $w$. So $\log_{BK}({\bf 1+p^*})$ contributes an element of unit norm. Hence the product defining $\abs{\log_{BK}(\sbxi_{\bz}^{\scriptscriptstyle{Mochizuki}})}$ and is finite for each $\sbxi_{\bz}^{\scriptscriptstyle{Mochizuki}}$.
\erem

Next define
\be 
\abs{\thetajti}=\sup_{\bz\in\tSigma_{L'}}\left\{ \abs{\log_{BK}(\sbxi_{\bz}^{\scriptscriptstyle{Joshi}})} \right\}
\ee
and likewise define  
 \be 
 \abs{\thetamti}=\sup_{\bz\in\tSigma_{L'}}\left\{ \abs{\log_{BK}(\sbxi_{\bz}^{\scriptscriptstyle{Mochizuki}})} \right\}.
 \ee
  
Similarly one defines, using tensor product norms on local factors, the quantities  $\abs{\thetaji}$, $\abs{\thetami}$, as the suprema of the (tensor product) norms of the images of elements of $\thetajti$ (resp. $\thetamti$).

One sees immediately that:
\bcor\label{cor:bounded-suprema}
The four suprema $\abs{\thetajti}$, $\abs{\thetaji}$, $\abs{\thetamti}$ and $\abs{\thetami}$ are bounded.
\ecor
\bp 
This is clear from \Cref{cor:compactness-of-locii}. 
\ep 

\blem\label{le:elem-lemma}
Let $E$ be a $p$-adic field and let $\abs{-}_E$ be a $p$-adic absolute value on $E$. Then for all $x\in (p^*\cdot\O_E)-\{0\}$ one has
$$\abs{\log_E(1+x)}_E=\abs{x}_E.$$ 
\elem
\bp 
Let $\pi\in\O_E$ be a uniformiser of $\O_E$. 
By \cite{koblitz-p-adic-book}, $\log_E(1+x)$ converges for $\abs{x}_E<1$. As $E$ is a $p$-adic field and $\abs{-}_E$ is a $p$-adic absolute value on $E$ i.e. $E$ has a discrete value group, so $\abs{x}_E\leq \abs{\pi}_E$ for all $x\in \fm_E$. Thus this series converges for all $x$ such that $\abs{x}_E\leq \abs{\pi}_E$ and hence the $n^{th}$ term $\abs{\frac{x^{n}}{n}}_E\to 0$ as $n\to \infty$. So one may approximate $\log_E(1+x)$ by its partial sum $$\sum_{n=1}^N\frac{(-1)^{n-1}\cdot x^n}{n}=x + \text{ higher order terms}$$ for some sufficiently large $N$ and estimate the absolute value of this partial sum.  Hence, by the standard properties of $p$-adic absolute values, it is enough to show that 
$$\abs{\frac{x^{n}}{n}}_E<\abs{x}_E\qquad n\geq 2.$$
This is true if $\abs{x^{n-1}}< \abs{n}_{E}$ for all $n\geq 2$ or equivalently
$$\abs{x}_E<\abs{n}_E^{1/(n-1)}.$$
So it is enough to verify that $$\abs{x}_E<\inf_{n\geq1}\left\{ \abs{n}_E^{1/(n-1)}\right\}.$$
If $p\nmid n$ then $\abs{n}_E=1$, and if $p|n$ then write $n=p^k\cdot r$ with $k\geq1$, $p\nmid r$, so it is enough to check that 
$$\abs{x}_E<\inf_{k\geq 1, r\geq1}\left\{1, \abs{p^k}_E^{1/(p^k\cdot r-1)}\right\} =\inf_{k\geq 1, r\geq1}\left\{1,\abs{p}_E^{k/(p^k\cdot r-1)}\right\}.$$
For any fixed $0<t<1$ the continuous function of $y$ on the open interval $(0,1)$ given by $f(t)=t^y$ is a decreasing function of $y$. As $\frac{k}{p^k\cdot r-1}<\frac{k}{p^k-1}$ for $r\geq 2$, one may replace the rightmost infimum by
$$\abs{x}_E<\inf_{k\geq 1, r\geq1}\left\{1, \abs{p^k}_E^{1/(p^k\cdot r-1)}\right\} =\inf_{k\geq 1}\left\{1,\abs{p}_E^{k/(p^k-1)}\right\}.$$
Now for $k\geq 1$, $k\mapsto k/(p^k-1)$ is an exponentially decaying, decreasing function and hence 
$$k/(p^k-1)<\begin{cases} 1/(p-1) & \text{ if } p \neq 2,\\
1 & \text{ if } p=2,
\end{cases}
$$
equivalently (for any $p$ and each $k\geq2$)
$$k/(p^k-1)< 1/(p-1).$$
So one has $$\abs{x}_E<\inf_{k\geq 1, r\geq1}\left\{1, \abs{p^k}_E^{1/(p^k\cdot r-1)}\right\} =\inf_{k\geq 1}\left\{1,\abs{p}_E^{1/(p-1)}\right\}=\abs{p}_E^{1/(p-1)}.$$

Now by the hypothesis of the theorem $x\in(p^*\O_E)-\{0 \}$. Thus
$$\abs{x}_E\leq\abs{p^*}_E<\abs{p}_E^{1/(p-1)}.$$
Hence under the hypothesis of the lemma, the asserted estimate holds. This proves the lemma.
\ep

\subparat{Proof of the fundamental estimate for $\abs{\thetami}$ and $\abs{\thetaji}$.}
Now one is ready to prove the fundamental estimate \moccor.

\bthm\label{th:moccor-Mochizuki-form1} 
For $X,C,L,L'$ satisfying \cref{ss:elliptic-curve-assumptions}, \inithtdata, and for an odd prime $\ell\gg 0$ one has
$$  \abs{\thetami} \geq \prod_{w\in\ubblvoss} \abs{q_w^{1/2\ell}}^{\ells},$$
and also
$$  \abs{\thetaji} \geq \prod_{w\in\ubblvoss} \abs{q_w^{1/2\ell}}^{\ells}.$$
\ethm

\bp[Proof of \Cref{th:moccor-Mochizuki-form1}]
The strategy of the proof is the same as the local strategy of \cite{joshi-teich-estimates}.  One exhibits an element, namely  $\sbxi_{\bz}=\sbxi_{\bz_{\Theta}}^{\scriptscriptstyle{Joshi}}$ (resp. $\sbxi_{\bz}=\sbxi_{\bz_{\Theta}}^{\scriptscriptstyle{Mochizuki}}$), of $\thetajti$, $\thetamti$, $\thetaji$, $\thetami$ for which the absolute value can be computed and one shows that $\abs{\sbxi_{\bz}}$ exceeds the number on the right of the stated inequality. Since $\abs{\thetaji}$ (resp. $\abs{\thetami}$) on the left of the asserted inequality is a supremum, the inequality claimed in the theorem will follow if the absolute value of this chosen element  bounds above the quantity on the right of the asserted inequality in the theorem.

 Since $\thetami, \thetamti, \thetaji, \thetajti$ as  defined in \cref{ss:construction-thetam-thetaj},  \Cref{th:thetam-construction} are adelic, and $\abs{\thetami}$ and $\abs{\thetaji}$ are defined multiplicatively, one may work with components ${\ttheta^{\mathscrbf{I}}}_{\scriptscriptstyle{Joshi},w}$,  ${\ttheta^{\mathscrbf{I}}}_{\scriptscriptstyle{Mochizuki},w}$ for each $w\in\ubblv$. So it is enough to estimate $\abs{{\ttheta^{\mathscrbf{I}}}_{\scriptscriptstyle{Mochizuki},w}}$ and then take product over all $w\in\ubblv$.
 
First I will work with the product versions $\thetajti$ and $\thetamti$ and then deduce the inequality for the tensor product versions $\thetaji$ and $\thetami$ from this.

Note that these products are finite by \cref{cor:compactness-of-locii}. [As remarked earlier, by \cref{le:elem-lemma}, $\abs{\log_p(1+p_w^*)}_{L_w'}=\abs{p^*}_{L_w'}$. So $\log_{BK}({\bf 1+p^*})$ contributes an element of unit norm in the product.]

So it suffices to establish the inequality (for $w\in\ubblvoss$)
$$\abs{{\ttheta^{\mathscrbf{I}}}_{\scriptscriptstyle{Mochizuki},w}}\geq \abs{q_w^{1/2\ell}}^{\ells}.$$ 
For this purpose it will be convenient to work with the Bloch-Kato logarithms \cref{ss:bloch-kato-logarithm} of the cohomological classes $\sbxi_{\by}$. Then  \cref{eq:bloch-kato-log} asserts that the codomain of  Bloch-Kato logarithm $\log_{BK}$ is, for each $w\in\vlpp$, the field $L_w'$ considered as a $\Q_{p_w}$-vector space.  This means that in the product setting i.e. in $\tsImjq$ and $\tsImq$, each class $\sbxi_{\uby{\by}_j}$ contributes $\log_{BK}(\sbxi_{\uby{\by}_j})$ in $L'_{w;\by_j}$ to the corresponding factor of $\tsImjq$ (resp. $\tsImq$).

Now let me say how one may deduce the inequality for $\thetaji$ from $\thetajti$ (resp. $\thetami$ from $\thetamti$). One may pass from the finite direct product $\prod_{w|p} L'_w$ (equivalently from the finite direct sum) of $\Q_p$-vector spaces to the tensor product $\tensor_{w|p} L_w'$ of $\Q_p$-vector spaces as follows. Consider the mapping $$\prod_{w|p} L'_w\ni(a_w)\mapsto \otimes_{w|p} (a_w) \in \bigotimes_{w|p} L'_w \mapsto \abs{\mathop{\otimes}\limits_{w |p} a_w}_{\tensor_{w |p} L_w'} \in\R,$$ where first  mapping is the natural  homomorphism  of  $\Q_p$-vector spaces: $$\prod_{w|p} L'_w\to \otimes_{w|p} L'_w,$$  and the second mapping is the tensor norm. In other words, one considers the image of the theta-values locii in the tensor product (via the first mapping) and then takes supremum over the tensor norms of all the images under the first homomorphism. This is where  the cross-norm property \cref{ss:cross-norm-prop} of the tensor product norms comes in to play. This property says one has 
$$\abs{\mathop{\otimes}\limits_{w |p} a_w}_{\tensor_{w |p} L_w'}=\prod_{w |p}\abs{a_w}_{L_w'}.$$
Hence one may also work with suprema of absolute values of the images of the theta-locii in the tensor products i.e. one can work with $\thetajti$ instead of $\thetaji$ (resp. work with $\thetamti, \thetami$) and still arrive at the required estimates in the tensor product settings (i.e. for $\thetami$ (resp. $\thetaji$) of \cite{mochizuki-iut3}. \textit{Notably, since the tensor product $\tensor_{w |p} L_w'$ is constructed using the product $\prod_{w |p} L_w'$ (more precisely the tensor product is the quotient of a suitable free $\Q_p$-module on the set underlying this product), one can expect that the theta-values locus in the tensor product setting (i.e. $\thetaji,\thetami$) should provide tighter upper bounds than the theta-values locus in the product setting i.e. $\thetajti,\thetamti$.}

In any of these cases,  each class $\sbxi_{\uby{\by}_j}$ contributes $\log_{BK}(\sbxi_{\uby{\by}_j})$ in $L'_{w;\by_j}$ and by the definition of this class and, at each $w$ one has, by  \eqref{eq:bloch-kato-log-class} one has:
\be\label{eq:multiplicative-coh-class}\abs{\log_{BK}(\sbxi_{w,j})}_{L_w',j}= \abs{\frac{\log_{L_w'}(1+p^*\cdot q_{w;j}^{1/2\ell})}{p^*}}_{L_w',j} \qquad  \text{ for }j=1,2,\ldots,\ells\ee
 By   \Cref{le:elem-lemma} one has:
\be\abs{\log_{BK}(\sbxi_{w,j})}_{L_w',j}= \abs{\frac{\log_{L_w'}(1+p^*\cdot q_{w;j}^{1/2\ell})}{p^*}}_{L_w',j} =\abs{q_{w;j}^{1/2\ell}}_{L_w',j}\qquad  \text{ for }j=1,2,\ldots,\ells\ee
Now by the valuation scaling property \Cref{th:adelic-theta-link}{\bf(4)}  one obtains 
\be\abs{\log_{BK}(\sbxi_{w,j})}_{L_w',j}= \abs{q_{w;j}^{1/2\ell}}_{L_w',j}=\abs{q_{w;1}^{1/2\ell}}^{j^2/\ells^2}_{L_w',1} \text{ for }j=1,2,\ldots,\ells\ee
Hence (writing $q_w$ in place of $q_{w;1}$)
\be\prod_{j=1}^\ells\abs{\log_{BK}(\sbxi_{w,j})}_{L_w',j}= \prod_{j=1}^\ells\abs{q_w^{1/2\ell}}_{L_w',j}=\prod_{j=1}^\ells\abs{q_w^{1/2\ell}}^{j^2/\ells^2}_{L_w',1}=\abs{q_w^{1/2\ell}}_{L_w',1}^{\sum_{j=1}^{\ells}\frac{j^2}{\ells^2}}.\ee
Hence one may pass to product over $w\in\ubblvossp\neq\emptyset$ and over all rational primes $p$ for which $\ubblvossp\neq\emptyset$:
\be\prod_p\prod_{w \in\ubblvossp\neq\emptyset}\prod_{j=1}^\ells\abs{\log_{BK}(\sbxi_{w,j})}_{L_w',j}=\prod_p\prod_{w \in\ubblvossp\neq\emptyset}\abs{q_w^{1/2\ell}}_{L_w',1}^{\sum_{j=1}^{\ells}\frac{j^2}{\ells^2}}.\ee

Now the proof can be completed by estimating, for each $w\in\ubblvossp\neq \emptyset$, the factor $\abs{q_w^{1/2\ell}}_{L_w',1}^{\sum_{j=1}^{\ells}\frac{j^2}{\ells^2}}$ contributed by $w$ and this estimate is provided by the proof of \cite[Theorem 9.2.1]{joshi-teich-estimates}. This proves the theorem.
\ep

\brem\label{re:multiplicative-action} 
It is precisely because of \Cref{le:elem-lemma} and the property established by \Cref{eq:multiplicative-coh-class} in the proof of \Cref{th:moccor-Mochizuki-form1} that one recovers Mochizuki's multiplicative action  of theta-values on log-shells depicted in \cite[Figure 3.1]{mochizuki-iut3}, i.e. multiplication by theta-value on the log-shell, used in \cite[Theorem 3.11]{mochizuki-iut3} for defining $\thetami$ instead of constructing crystalline Galois cohomology classes as has been done here.
\erem

\newcommand{\logVol}{{\rm LogVol}}
\subparat{Relationship with Mochizuki's log-volumes and hulls}
Since $\thetaji$ (resp. $\thetami$) are contained in adelic compact sets, one can compute volumes of the smallest compact subset containing these sets. So one needs to define a volume of sets of this type which arise in realizing this idea. This theory is detailed in  \cite[Proposition 3.1, Proposition 3.9]{mochizuki-iut3} and this theory is applied in \moccor\ and in \cite{mochizuki-iut4} for estimating the volume $\Vol(\thetam^{\bsi})$ of $\thetam^{\bsi}$. More precisely, $\Vol(\thetam^{\bsi})$ is the volume of the smallest compact set (of a suitable sort chosen below) containing $\thetaji$. In \iut\ one works with  $$\logVol(\thetam^{\bsi})=\log\Vol(\thetam^{\bsi})\in \R\cup\{\pm\infty\},$$ which is the natural logarithm of $\Vol(\thetam^{\bsi})$. 

The fundamental inequality of \moccor\ is an assertion about the log-volume $\logVol(\thetam^{\bsi})$ of  $\thetam^{\bsi}$. Equivalently, it is an inequality about the volume $\Vol(\thetam^{\bsi})$. As I have remarked in \cite[\ssep 1.12]{joshi-teich-quest}, the volume computation of \cite{mochizuki-iut3} is quite similar to the sort of volume computations in classical Teichmuller Theory which, for instance, occur in \cite{mirzakhani07}.

\subsubpara{} As explained in \cref{se:intro}, my approach to establishing \moccor\ i.e. establishing the fundamental estimate for $\logVol({\thetam^{\bsi}})$ is different from that of Mochizuki. Though at this point, using the results of the preceding sections along with the \rosettastone\ and \Cref{th:thetam-construction}, it is possible to proceed exactly as Mochizuki does, to arrive at Mochizuki's proof of \moccor.

However once the principles underlying \moccor\ are understood, one can formulate approaches different from that of Mochizuki. I will take this approach to the proof rather than repeat Mochizuki's approach.

\subsubpara{} Let me begin with the theory of volumes described in \cite{mochizuki-iut3}. The basic properties are  summarized in \cite[Proposition 3.9]{mochizuki-iut3} and also in \cite[Proposition 5.7]{mochizuki-topics3}.

Let $L'$ be a number field, and for $w|p$ in $L'$, one has the  $p$-adic field $L_w'$.  Let $\O_{L_w'}\subset L_w'$ be the ring of integers. Since $L_w'$ is a locally compact abelian group,  there is a unique normalized translation invariant (real valued) Haar measure on $L_w'$ which I will denote by  $\Vol_{L_w'}$. 

If $S\subset L_w$,  is a measurable subset, for some $w$, I will write $\Vol_{L_{w}'}(S)$ for its normalized Haar measure. Here normalization means the following:
$$\Vol_{L_w'}(\O_{L_w'})=1.$$ 

Let \be\log:(0,\infty)\to \R\ee be the natural real-valued logarithm function and extend $\log$ by setting $$\log(0)=-\infty\text{ and } \log(\infty)=\infty.$$ This extended logarithm function takes values in $\R\cup\{\pm\infty \}$. By the \textit{log-volume $S\subset  L_w'$}, I will mean  
$$\mulog {L_w'} S=\log(\Vol_{L_w'}(S))\in \R\cup\{\pm\infty \}.$$ In particular, this means that one has $$\mulog{L_w'}{\O_{L_w'}}=0.$$
Let $\pi\in\O_E$ be a uniformizer. Then by translation invariance of the Haar measure $\Vol_E$ (or otherwise)  one has
$$\Vol(\pi\O_{L_w'})=\abs{\pi}_{L_w'}$$
and so 
$$\mulog{L_w'}{\pi\O_{L_w'}}=\log(\abs{\pi}_{L_w'}).$$
The following lemma is now elementary.
\blem\label{le:volume-lemma1} 
Let $E$ be a $p$-adic field with absolute value $\abs{-}_E$ and let $S= \alpha+\lambda\O_E$ be a non-empty open subset for some $\alpha\in E$ and $\lambda\in E^*$. Then 
\benumlab
\item $$\Vol_E(S) =\Vol_E(\alpha+\lambda\cdot\O_E)=\Vol(\lambda\O_E),$$ and hence
\item $$\Vol_E(\alpha+\lambda\cdot\O_E)=\Vol(\lambda\O_E)=\abs{\lambda}_E,$$ and also
\item $$\Vol(\lambda\O_E)=\sup\{\abs{s}_E:s\in \lambda\O_E\}=\abs{\lambda}_E,$$ and, in particular, 
\item for any $s\in\lambda\cdot\O_E$ (note that $s\in S$ if and only if $\alpha=0$) one has the tautology
$$\abs{\lambda}_E=\Vol_E(S)\geq \abs{s}_E.$$
\eenum
\elem
\bp 
The first equality is the translation invariance of $\Vol$ and the second equality is clear from the definition of $\Vol_E$. The remaining assertions are clear.
\ep

\subsubparat{Weighted volumes}\label{ss:weighted-vol}
These considerations can be extended as follows. Let $E_1,\ldots,E_n$ be finite extensions of $\Q_p$. The theory of volumes and log-volumes on $E_1,\ldots,E_n$ extends to a theory of weighted volumes on the tensor product 
$E=E_1\tensor_{\Q_p}E_2\tensor_{\Q_p}\cdots\tensor_{\Q_p}E_n$. Let $\Gamma=\{\gamma_1, \gamma_2, \ldots, \gamma_n\}\subset (0,1]\subset \R$ be a set of weights. Then the $\Gamma$-weighted volume $\Vol^\Gamma_{E_1\tensor_{\Q_p} \cdots \tensor_{\Q_p}E_n}$ is a function on certain measuable subsets of the $\Q_p$-vector space $E$. This is defined as follows.

Recall from \cite{schneider-book}, that $E_1\tensor_{\Q_p}E_2$ is equipped with topology which has a fundamental system of neighborhoods of zero consisting of sets of the form $V_1\tensor V_2 \tensor\cdots\tensor V_n$ where $V_i\subset E_i$ is a $\Z_p$-lattice in each $E_i$.  

So it suffices to define
\be\Vol_{E_1\tensor_{\Q_p}\cdots  \tensor_{\Q_p}E_n}(V_1\tensor_{\Z_p}\cdots \tensor_{\Z_p}V_n)=\Vol^{\gamma_1}_{E_1}(V_1)\cdots \Vol_{E_n}^{\gamma_n}(V_n).\ee

This leads to the following version of \Cref{le:volume-lemma1}:
\blem\label{le:volume-lemma2}
Let, for $i=1,\ldots,n$, $V_i=\alpha_i+\lambda_i\O_{E_i}\subset \O_{E_i}$ with $\alpha_i\in E_i$ and $\lambda_i\in E_i^*$. Then
\benumlab
\item $$\Vol^\Gamma_{E_1\tensor_{\Q_p} \cdots \tensor_{\Q_p}E_n}(V_1\tensor_{\Z_p} V_2\tensor_{\Z_p}\cdots \tensor_{\Z_p}V_n)=\prod_{i=1}^n\Vol^{\gamma_i}_{E_i}(V_i).$$
\item Hence $$\Vol^\Gamma_{E_1\tensor_{\Q_p} \cdots \tensor_{\Q_p}E_n}(V_1\tensor_{\Z_p} V_2\tensor_{\Z_p}\cdots \tensor_{\Z_p}V_n)=\prod_{i=1}^n\Vol_{E_i}^\Gamma(V_i)=\prod_{i=1}^n \abs{\lambda_i}^{\gamma_i}_{E_i},$$
\item and hence if $s_i\in \lambda_i\O_{E_i}\cap E_i^*$ then 
$$\Vol^\Gamma_{E_1\tensor_{\Q_p} \cdots \tensor_{\Q_p}E_n}(V_1\tensor_{\Z_p} V_2\tensor_{\Z_p}\cdots \tensor_{\Z_p}V_n)=\prod_{i=1}^n\Vol_{E_i}^{\gamma_i}(V_i)\geq \prod_{i=1}^n \abs{s_i}^{\gamma_i}.$$
\eenum
\elem

\subsubparat{} The sets of interest to us, arising in \cref{th:moccor-tensor-product1}, contain tensor products of sets of the form $$\prod_{w\in\ubblv}(\lambda_w\cdot \O_{L_w'})$$
for  $\lambda_w\in L_w'^*$ satisfying the condition: for all  $w\notin\cup_{p}\ubblvossp$, $\lambda_w\in \O_{L_w'}^*$. 

For tensor products sets, one may define weighted adelic volumes using the preceding formalism:
\be \Vol\left(\prod_{w\in\ubblv}(\lambda_w\cdot \O_{L_w'})\right)=\prod_{w\in\ubblv}\Vol_{L_w'}(\lambda_w\cdot \O_{L_w'})
\ee
with some choice of weights. The condition that  $\lambda_w\in \O_{L_w'}^*$ for all but finitely many $w$ ensures that this volume is finite. Hence the sets of interest to us have a finite volume. Notably, for sets of this type, the volume is related to absolute values through the above elementary lemma.

\subsubsection{}
Fix a prime $p$ such that $\ubblvoss_p\neq\emptyset$. Let $\{w_1,\ldots, w_m\}=\ubblvossp$ be all the distinct   primes of $L'$ lying over $p$ and also contained in $\ubblvoss_p$. Take 
\be\label{eq:weights}\Gamma_p=\{1/[L'_{w_i}:{\lmod}_{,v}]: w_i \in\ubblvossp \}
\ee as the set of weights. This sets of weights is determined by \Cref{pr:galois-cohomology-section} and coincides with the choice of weights in \cite[Remark 3.1.1]{mochizuki-iut3}. These are the weights and  sets of interest to us here, and in \cite{mochizuki-iut3} and \cite{mochizuki-iut4}, are the tensor product subsets which arise in the proof of \Cref{th:moccor-tensor-product1}.

This allows us to compute $\Vol_E$ for suitable subsets using absolute values and hence the two calculations will yield the same lower estimates.

\subsubparat{} 
In particular, the notion of weighted volumes and log-weighted-volumes extends verbatim to the $\Q_p$ vector spaces obtained from the $p$-component of the log-shells $\sI^{\Q_p}_{Mochizuki,p}$, and ${}^{\bS_{j+1}}\sI^{\Q_p}_{Mochizuki,p}$. In particular, for any measurable set $S\subset {}^{\bS_{j+1}}\sI^{\Q_p}_{Mochizuki,p}$, one can talk of its volume $\Vol(S)$.

\subsubparat{} Let $$S\subset E_1\tensor_{\Q_p} \cdots \tensor_{\Q_p}E_n=\bigoplus_\alpha E'_\alpha.$$ Then the \textit{hull of $S$}, denoted here by $\text{hull}(S)$,  is the smallest subset $hull(S)\subset E_1\tensor_{\Q_p} \cdots \tensor_{\Q_p}E_n$ of the tensor product such that $\text{hull}(S)\supset S$ and whose image  in the direct sum decomposition on the right is of the form
$$\bigoplus_\alpha\lambda_\alpha \O_{E_\alpha'}$$
where $0\neq\lambda_\alpha\in \O_{E_\alpha'}$.

Properties of hulls are established in \cite[Remark 3.9.5(ii)]{mochizuki-iut3}.

The following lemma is now clear:
\blem\label{le:volume-lemma3} 
Let $S'\supset S=(\tau_1\O_{E_1})\tensor_{\Z_p}(\tau_2\O_{E_2})\tensor_{\Z_p}\cdots\tensor_{\Z_p}(\tau_n\O_{E_n})$ with $\tau_i\in \O_{E_i}$ for $i=1,\ldots,n$ and suppose $\Gamma_p$ are chosen as above and suppose that one has an element $s_1\tensor s_2\tensor \cdots \tensor s_n\in S$ with $s_i\in \tau_i\O_{E_i}$ for $i=1,\ldots,n$.
Then $$\Vol^\Gamma(hull(S'))\geq\Vol^{\Gamma_p}(hull(S))\geq \Vol^{\Gamma_p}(S)\geq \prod_{i=1}^n \abs{s_i}^{\gamma_i}_{E_i}\geq \prod_{i=1}^n \abs{s_i}_{E_i}.$$
\elem
\bp
The first inequality is immediate from the properties of the hull in \cite[\ssep 3]{mochizuki-iut3}. So
the only point which needs to be checked is that $0<\gamma_i\leq 1$, so the last inequality follows from the fact that all the $\tau_i$, and hence all the $s_i$ have absolute values at most one.
\ep

\subsubsection{Mochizuki's Hulls and Convex closures}
In this subsection I want to provide the following proposition which clarifies my usage of convex closures and Mochizuki's \textit{holomorphic hulls} detailed in \cite[Remark 3.9.5]{mochizuki-iut3}. The notion of convex closure can be made in any Banach space \cite{schneider-book}. Mochizuki's notion of holomorphic hull is a bit more specific. To establish the relationship assume, for notational simplicity, that $E_1,E_2$ are two $p$-adic fields considered as finite dimensional vector spaces over $\Q_p$. 

I work with the tensor product norm $\abs{-}_{E_1\tensor_{\Q_{p}} E_2}$ and the normed $\Q_p$-vector space $$(E_1\tensor_{\Q_{p}} E_2,\abs{-}_{E_1\tensor_{\Q_{p}} E_2})$$
equipped with this tensor norm.

On the other hand Mochizuki's approach to holomorphic hulls detailed in \cite[Remark 3.9.5]{mochizuki-iut3} is based on the following decomposing the tensor product vector space $E_1\tensor_{\Q_{p}}E_2$:
$$f:E_1\tensor_{\Q_{p}}E_2\mapright{\isom}\bigoplus_{\alpha} F_\alpha$$
into a direct sum of $p$-adic fields $F_\alpha\supset \Q_p$  (obviously the set $\{\alpha\}$ is also finite) and work, according to \cite[Remark 3.9.5]{mochizuki-iut3}, relative to this decomposition. 
Since $E_1\tensor_{\Q_p}E_2$ is finite dimensional over $\Q_p$, $f$  is an  isomorphism of finite dimensional Banach spaces over $\Q_p$.

\bpro\label{pr:convexes-and-hulls}
In the notation of this paragraph, suppose $U\subset E_1\tensor_{\Q_{p}}E_2$ is a relatively compact subset. Let $H=H(U)=hull(U)$ be the hull of $U$ as defined by \cite[Remark 3.9.5, Page 543]{mochizuki-iut3}. Then
\benumlab
\item $H$ is a convex subset of $\oplus_{\alpha} F_\alpha$ containing $f(U)$ and by \cite[Properties P1, P2, P3, Page 534]{mochizuki-iut3} it is the minimal convex subset containing $f(U)$ with this property.
\item If $V\subset E_1\tensor_{\Q_{p}}E_2$ is a convex subset of this normed vector space then its image under $f$ is also convex.
\item In particular, if $P$ is a tensor product region of $P\subseteq E_1\tensor_{\Q_{p}}E_2$ in the sense of \cite[Remark 3.1.1]{mochizuki-iut3}, then $H(P)\subseteq \oplus_{\alpha} F_\alpha$ is the image of the convex closure of $P\subseteq E_1\tensor_{\Q_{p}}E_2$.
\eenum
\epro
\bp 
By the definition of the hull given in \cite[Remark 3.9.5]{mochizuki-iut3} the hull is an $\Z_p$-module. By the definition of convex subsets \cite[Chap I, \ssep 2]{schneider-book}, the hull is a convex subset as it is an $\Z_p$-module. By \cite[Remark 3.9.5, Properties P1, P2, P3]{mochizuki-iut3} it is the smallest subset with the said properties. This proves {\bf(1)}. The assertion {\bf(2)} is clear from the definitions \cite{schneider-book}. The last assertion is a consequence of the fact that holomorphic hull is the convex closure of $f(P)$.
\ep

\brem 
Because of \Cref{pr:convexes-and-hulls}, one can work with hulls and convex closures on an equal footing and this allows me to translate my results to Mochizuki's context.
\erem

\subsection{Proof of the fundamental estimate for $\Vol(\thetami)$ and $\Vol(\thetaji)$}\nwss
Now one is ready to prove the fundamental estimate \moccor.

\bthm\label{th:moccor-Mochizuki-form2} 
For $X,C,L,L'$ satisfying \cref{ss:elliptic-curve-assumptions}, \inithtdata, and for an odd prime $\ell\gg 0$ one has
$$  \Vol(\thetami) \geq\prod_{w\in\ubblvoss} \abs{q_w^{1/2\ell}}^{\ells},$$
and also
$$   \Vol(\thetaji) \geq \prod_{w\in\ubblvoss} \abs{q_w^{1/2\ell}}^{\ells}.$$
\ethm
\bp 
By \Cref{pr:convexes-and-hulls}, one can work with Mochizuki's holomorphic hulls or convex closures (in Banach spaces). It will be enough to prove this for one of the two $\thetami,\thetaji$. I will prove the assertion for $\thetami$. The proof of \Cref{th:moccor-Mochizuki-form1} will serve as a template for this assertion. By the construction of these sets,  $\thetami$ contains the classes $\bxi_{\bz}$ given by \Cref{th:thetam-construction} obtained from each choice of Tate parameters obtained using the standard point of the $\tSigma_{L'}$ and hence it contains  a subset of the type considered in \Cref{le:volume-lemma3} 
$$\prod_{w\in\ubblvoss}(\tau_1\O_{L'_{w,1}})\tensor_{\Z_p}(\tau_2\O_{L'_{w,2}})\tensor_{\Z_p}\cdots\tensor_{\Z_p}(\tau_\ells
\O_{L'_{w,\ells}}),$$
with $\tau_j=\log_{BK}(\bxi_{0,\bz_j,w})$ (for suitable $j$).
So the volume of this set can be bounded by \Cref{le:volume-lemma3}.

Since $\thetami$ defined in \cref{ss:construction-thetam-thetaj} is adelic, and weighted volumes are  defined multiplicatively, one may work with components ${\ttheta^{\mathscrbf{I}}}_{\scriptscriptstyle{Mochizuki},w}$ for each $w\in\ubblv$. So it is enough to estimate $\Vol({\ttheta^{\mathscrbf{I}}}_{\scriptscriptstyle{Mochizuki},w})$ and then take product over all $w\in\ubblv$. Note that the product is finite by \cref{cor:compactness-of-locii}. 
So it suffices to establish the inequality  $\Vol({\ttheta^{\mathscrbf{I}}}_{\scriptscriptstyle{Mochizuki},w})\geq \abs{q_w^{1/2\ell}}^\ells$. This in turn is clear from the construction of these sets: they contain the classes $\bxi_{\bz}$ given by \Cref{th:thetam-construction} obtained from each choice of Tate parameters obtained using the standard point of the $\tSigma_{L'}$ and so contain a subset of the type  \Cref{le:volume-lemma3} with $\tau_j=\log_{BK}(\bxi_{0,\bz_j,w})$ (for suitable $j$). Using this and the relationship between weighted volumes absolute values established earlier and arguing as in the proof of \Cref{th:moccor-Mochizuki-form1}, one obtains the assertion.
\ep

\subsubsection{Proof of Mochizuki's Corollary~3.12}\label{ss:moccor-proof}\nwsss 
Mochizuki's notational conventions \cite[Page 420 and Page 608]{mochizuki-iut3} sign conventions are a bit awkward to work with. The notation convention of \cite[Page 420 and Page 608]{mochizuki-iut3} is designed so that the ration of the LHS and RHS in \moccor\ is positive while both the LHS and RHS are negative. I will define 
$$
-\abs{\logVol(\thetami)}=\sum_p -\abs{\logVol(\thetamp)}=-\sum_p \abs{\logVol(\thetamp)}.
$$
This is perfectly reasonable as the $\logVol(\thetamp)\leq0$ and this notation convention is compatible with the following standard properties of real numbers and logarithms:  if $0<x\leq1$ is a real number then $$\log(x)\leq 0\text{ and so }\abs{\log(x)}=-\log(x)\geq0\text{ and hence }-\abs{\log(x)}=\log(x)\leq0.$$ Moreover if $0<x, y\leq1$ are real numbers then 
\begin{align*}
-\abs{\log(x\cdot y)}&=\log(x\cdot y)\\
&=-\abs{\log(x)}+ (-\abs{\log(y)})\\
&=-(\abs{\log(x)}+\abs{\log(y)}).
\end{align*}

To provide a translation of \Cref{th:moccor-Mochizuki-form2} into  \moccor, one must work with Mochizuki's notational convention and the above observation regarding logarithms of real numbers.  From \Cref{th:moccor-Mochizuki-form2} one obtains, using Mochizuki's notational conventions \cite[Page 420 and Page 608]{mochizuki-iut3}, the assertion known as  \moccor:
\bcor\label{cor:cor312}
With notations and assumptions of \Cref{th:moccor-Mochizuki-form2}, and in Mochizuki's notational conventions \cite[Page 420 and Page 608]{mochizuki-iut3}, one has
$$-\frac{1}{\ells}\abs{\logVol(\thetami)}\geq -\sum_{p,w\in\ubblvossp\neq\emptyset}\log\abs{q_w^{1/2\ell}}_{L_w'}.$$
\ecor

\section{Appendix: Frobenioids and Hodge-Theaters \`a la Mochizuki via \cite{joshi-teich}}\label{se:frobenioids}
\iutthr\ uses Mochizuki's notion  of Frobenioids and Hodge-Theaters. Let me explicate how these objects of \iutthr\ arise from the perspective  of \cite{joshi-teich}. I will not use the broad definition of Frobenioids given in \cite{mochizuki-frobenioid1} but use \cite[Example 6.1 and Example 6.3]{mochizuki-frobenioid1} which are more than adequate  for my discussion and for \iutthr.

\renewcommand{\deg}{{\rm deg}}
\subparat{Frobenioid of an archimedean valued field}\label{ss:frobenioid-archimedean}
The simplest example of a Frobenioid is the Frobenioid of an archimedean valued field. Let me recall it here. 

Let $K$ be a field equipped with an archimedean absolute value $\abs{-}_K$. I will assume  that $K$ contains $\sqrt{-1}$ (this is also the assumption throughout most of \iut) and hence by Ostrowski's Theorem, $K\isom \C$ as  valued fields, where $\C$ is equipped with the standard complex absolute value $\abs{-}_{\C}$.

Let \be\Phi(K)=\R_{\geq 0},\ee
where $\R_{\geq 0}$ is the additive monoid of non-negative real numbers. Here and in 
what follows, let $\Phi(K)^{gp}$ be the group associated to the monoid $\Phi(K)$. Then as $\Phi(K)=\R_{\geq 0}$, so one has $\Phi(K)^{gp}=\R$. Let 
\be  
K^*\to \Phi(K)^{gp}=\R,
\ee
be a homomorphism of groups.

Note that if $(K,\abs{-}_K)=(\C,\abs{-}_{\C_p})$ then one has the familiar homomorphism of monoids 
\be 
\C^*\mapsto \R \text{ given by } z\mapsto  \log(\abs{z}_{\C}),
\ee
where $\log:\R_{> 0}\to \R$ is the natural logarithm of a positive real number. This will be referred to as the archimedean (arithmetic) degree homomorphism.

The \emph{Frobenioid associated to the archimedean valued field $(K,\abs{-}_K)$} (see \cite[Example 6.3]{mochizuki-frobenioid1}) is the triple 
\be 
\frob(K,\abs{-}_K)=(\Phi(K)=\R_{\geq 0}, K^*\to \Phi(K)^{gp}=\R).
\ee

\subparat{Frobenioid of a non-archimedean valued field}\label{ss:frobenioid-non-archimedean}
Let me now describe Frobenioid associated to  a non-archimedean valued field as it will be used in what follows. Let
$K$ be a field equipped with a non-archimedean absolute value $\abs{-}_K$. Typically $K$ will be either perfectoid with tilt $\cpt$ or a discretely valued for this discussion.

To fix the notation and sign conventions, for $K=\C_p$, I will require that 
\be 
\abs{p}_{\C_p}=p^{-1}.
\ee

Let $\O_K\subset K$ be its ring of integers, let $\ot_K=\O_K-\{0\}$ be considered as a multiplicative monoid and let $\O_K^*$ be the group of units of $\O_K$. Let $\abs{K^*}_K=\{\abs{x}_K: x\in K^*\}$ and similarly define $\abs{\O_K}_K$. Let 
\be\Phi(K)=\ot_K/\O_K^*
\ee
be the quotient of the multiplicative monoid $\ot_K$ by its submonoid of units $\O_K^*$, this is  equipped with the canonical isomorphism of monoids \be\Phi(K)\mapright{\abs{\cdot}_K} \abs{\ot_K}_K.\ee Let $\Phi(K)^{gp}$ be the group associated with the monoid $\Phi(K)$, one has a homomorphism of monoids 
\be 
\Phi(K)\to \abs{\ot_K}_K\to \abs{K^*}_K.
\ee
which factors as 
\be\Phi(K)\to \Phi(K)^{gp}\to \abs{K^*}.\ee

Then according to \cite{mochizuki-frobenioid1}, the data
\be \frob(K,\abs{-}_K)=(\Phi(K), K^*\to \Phi(K)^{gp}),\ee
consisting of the monoid $\Phi(K)$ and the natural homomorphism of $K^*\to \Phi(K)^{gp}$ given by mapping \be x\mapsto x\bmod{\O_K^*}\ee
is the \emph{Frobenioid associated to the non-archimedean valued field  $(K,\abs{-}_K)$}.

For Diophantine applications, one needs to work with the notion of (local) arithmetic degrees. This is done by means of the natural homomorphism
\be
\Phi(K)\to \R
\ee 
given by  \be x\mapsto  \log(\abs{x}_K)\in \R,\ee which will be called \emph{the (local) arithmetic degree}. Since $\R$ is a group, by the universal property of $\Phi(K)^{gp}$, the arithmetic degree factors as a homomorphism
\be 
\Phi(K)\to \Phi(K)^{gp} \to \R,
\ee
and so composing with $K^*\to \Phi(K)^{gp}$ one gets the
arithmetic degree homomorphism
\be 
\deg_{K}:K^*\to \R
\ee
given by $x\mapsto \log(\abs{x}_K)$.

If $K=\Q_p$ with the standard valuation $v_{\Q_p}$ on $\Q_p$, then \be\Phi(\Q_p)=\Z_p^\triangleright/\Z_p^*\isom \Z_{\geq0}\text{ and }\Phi(\Q_p)^{gp}=\Z\ee and the arithmetic degree homomorphism is
a homomorphism 
\be\xymatrix{\deg_{\Q_p}:\Q_p^*\ar[rr]^{x\mapsto \log(\abs{x}_{\Q_p})}&& \R}\ee
which on $p\in\Q_p^*$ takes the value \be\deg_{\Q_p}(p)=\log(\abs{p}_{\Q_p}).\ee 

\subparat{Induced Frobenioid structures}\label{ss:induced-frobenioids}
In the context of Diophantine applications of the Arithmetic Teichmuller Theory of \joshiros\ and the theory of \iutthr, it is important to recognize that the Frobenioid of a $p$-adic field $E$
\be 
(\frob(E,\abs{-}_E),\deg_{E})
\ee 
depends on the absolute value $\abs{-}_E$. This is best understood using the notion of induced Frobenioid structures. Let me discuss this, as this is not discussed in \iutthr.

Suppose $E\to K$ is a homomorphism of valued fields, then one has  Frobenioid associated to valued field $\frob(E,\abs{-}_E)$ and as  $E\into K$ is a valued subfield, one can also  consider $E$ as a subfield of $K$ with induced valuation, and so, one has a \emph{sub-Frobenioid} $\frob(E,\abs{-}_K)$ of the Frobenioid $\frob(K,\abs{-}_K)$ of $K$ obtained by considering the valuation induced on $E$ from $K$. I will write this as \be
\frob(E,\abs{-}_K)=(\Phi(E),E^*\to \Phi(E)^{gp})_K\ee and which I will call  the \emph{Frobenioid structure on $E$ induced by $K$.}

Mochizuki does not work with induced structures in \iutthr, as Mochizuki's theory, especially \cite{mochizuki-iut3}, does not use valued fields as  such, and works with the abstract notion of Frobenioids. So my approach is slightly different from that taken in \iutthr\ and  will be used in what follows.

Typical example of this situation is the following: let $E$ be a $p$-adic field and $E\into K$ is an embedding of $E$ into an untilt $K$ of $\cpt$. Since $E$ is a discretely valued field $(E,\abs{-}_E)\isom (E,\abs{-}_K)$ is an isomorphism of discretely valued fields and hence the two Frobenioids $\frob(E,\abs{-}_E)$ and $\frob(E,\abs{-}_K)$ are abstractly isomorphic, but one wants to compare arithmetic degrees of the data $(\frob(E,\abs{-}_E),\deg_E)$ and $(\frob(E,\abs{-}_K),\deg_{K}\big|_E)$.

\bpro\label{pro:comparison-of-induced-frobenioids} 
Let $E/\Q_p$ be a finite extension of $\Q_p$ contained in $\C_p$. Let $K$ be an untilt of $\cpt$. Let $\iota:E\into K$ be an isometric embedding of $E$ into $K$ providing on $E$ a valuation given on $\iota(E)$ by the restriction of the valuation of $K$. Then one has
\benumlab 
\item The two Frobenioid structures on $E$ induced by $E\into \C_p$ and $\iota:E\into K$ are abstractly isomorphic.
\item In general, there may not exist any isomorphism of the two Frobenioid structures on $E$  provided by $\C_p$ and $K$ which is compatible with identity homomorphism on $\R$ in  the  arithmetic degree homomorphisms $(x\in E\into \C_p) \mapsto \log(\abs{x}_{\C_p})\in\R$ (resp. $(x\in E\into K) \mapsto \log(\abs{x}_{K})\in\R$).
\eenum  
\epro
\bp 
The fields $E$ and $\iota(E)$ are isomorphic as valued fields, as both are discretely valued and one has an isomorphism given by the homomorphism $\iota|_E:E\to \iota(E)$. So the Frobenioids of $E\into \C_p$ and $E\into K$ are abstractly isomorphic. This proves the first assertion.

On the other hand, under the (local) arithmetic degree homomorphism of $E\into \C_p$, $p\mapsto \log(\abs{p}_{\C_p})\in\R$ and $\iota(E)\into K$, one has $p\mapsto \log(\abs{p}_{K})\in\R$ and $v_K(p)$ which may not be equal to $v_{\C_p}(p)$ and hence there may  not exist any isomorphism between the Frobenioids of $E\into \C_p$ and $\iota(E)\into K$ which is compatible with the respective (local) arithmetic degrees and which is identity on $\R$ i.e. there may not exist any commutative diagram of groups:
\be 
\xymatrix{ E^* \ar[rr]^{x\mapsto \log(\abs{x}_{\C_p})}\ar[d]^\iota && \R\ar@{=}[d]\\
	\iota(E)^* \ar[rr]^{x\mapsto \log(\abs{x}_{K})} && \R}
\ee
if the valuation of $p$ in $\C_p$ and valuation of $p$ in $K$ do not coincide (and in general it is possible to find untilts $K$ such that this happens) i.e. $\log(\abs{x}_{\C_p})\neq \log(\abs{\iota(x)}_K)$ in $\R$ in general.
\ep

\subparat{Perfect Frobenioids and Perfectoids}\label{ss:perfect-realified-frobenioids} An important point for \iutthr\ is this: one wants to work with induced Frobenioids in \iutthr\ because one wants to allow valuations to be scaled by rationals (say $\abs{-}_{\Q_p}\mapsto \abs{-}_{\Q_p}^{1/j^2}$) (to be precise one wants a non-trivial scaling operation whose effect at the level of discretely valued $p$-adic fields is this type of scaling). But such a scaling is not possible with Frobenioids of $p$-adic fields, so Mochizuki remedies this with the introduction of \emph{perfect and realified Frobenioids} where the value group is allowed to be $\Q$ (or $\R$). This adds an additional notational and conceptual burden to his proof. Working with perfectoid fields allows one to overcome this difficulty quite naturally.

\subsubpara  My approach to Mochizuki's perfect and realified Frobenioids  is summarized by the following proposition.
\bpro
Let $\hol{X/E}{y}{F}$ be an arithmetic holomorphic structure with $y=(K,\iota_K:K^\flat\isom F)$. Then 
\benumlab
\item $\frob(K)$ is a perfect Frobenioid in the sense of \cite[Definition 1.1]{mochizuki-frobenioid1}.
\item Suppose additionally that the value group of $K$ is equal to $\R$ (equivalently $F$ has value group equal to $\R$), then $\frob(K)$ is a realified Frobenioid.
\item If $E$ is a $p$-adic field and $\iota: E\into K$ is an embedding. The inclusion of Frobenioids $\frob(E)\into \frob(K)$ provides a natural perfection, $\frob(K)$, of $\frob(E)$ and if $K$ has value group $\R$, then $\frob(K)$ is even a realification of $\frob(E)$ in the sense of \cite[Definition 1.1]{mochizuki-frobenioid1}.
\eenum 
\epro
\bp 
The proofs are clear from the definitions \cite[Definition 1.1]{mochizuki-frobenioid1}.
\ep

\brem 
Hence one could say that the above result demonstrates that there are many topologically inequivalent arithmetic models for the perfection of $\frob(E)$. This point is central even for \iutthr, as one wants a non-trivial variation of arithmetic-geometric data for the very statement of \moccor.
\erem
\subparat{Frobenioid of a number field}\label{ss:frobenioid-number-field}
Let $L$ be a number field i.e. a finite extension of $\Q$, assume that $L$ has no real embeddings for simplicity (this will always be the  case in \iutthr--for any choice of ``initial Theta data''). Then for every prime $v$ of $L$, archimedean or non-archimedean, one can associate to the completion $L_v$, a (local) Frobenioid for the local Field $L_v$ (defined as above). Following Mochizuki, \cite[Example 6.3]{mochizuki-frobenioid1},

If $v$ is a complex archimedean prime of $L$ then let $\bar{v}$ be its complex conjugate.

As discussued in \constrtwoh{\ssep}{ss:frobenioids}, the Frobenioid $\frob(L)=(L^*,\Phi(L), L^*\to \Phi(L)^{gp})$ of the chosen number field $L$ is defined as follows. It consists of the data of the group $L^*$ and the monoid
\be 
\Phi(L)=\bigoplus\limits_{v} \Phi(L_v), 
\ee
where the sum is over all primes $v$  (archimedean and non-archimdean) of $L$.
This comes equipped with the homomorphism 
\be 
L^*\to \Phi(L)^{gp}=\bigoplus\limits_{v} \Phi(L_v)^{gp}
\ee
given in the obvious way, and also a global arithmetic degree homomorphism
\be 
\Phi(L)\to \R
\ee
given 
\be(x_v)\mapsto \sum_{v}\log(\abs{x_v}_{L_v})\ee 
i.e. by the sum of local arithmetic degrees, which factors through $\Phi(L)^{gp}$ as $\R$ is a group and notably for any $x\in L^*$, one has the (logarithmic form of the) product formula \eqref{eq:prod-formula}:
\be
\sum_{v}\log(\abs{x}_{L_v})=0.
\ee

\brem As established in \constrtwoh{\ssep}{ss:frobenioids}, each arithmeticoid $\arithl_\by$ provides a natural Frobenioid $\frob(\arithl)$ and the product formula provides an hyperplane $H_{\by}$ (\cref{th:hyperplane-theta-lnk}) and in fact one has a period mapping in $\by\mapsto H_{\by}$ constructed by   \constrtwoh{Theorem }{th:hyperplane}. 
\erem
\subparat{Frobenioid associated to a variety}\label{ss:frobenioid-variety}
Let $E$ be  a field and let $X/E$ be a geometrically connected, smooth, projective variety over $E$. Let 
\be 
\Phi(X)=\divp(X)
\ee
be the monoid of effective divisors on $X$. Let $E(X)$ be the field of meromorphic functions on $X/E$. Then one has a Frobenioid associated to $X/E$ given by
the triple
\be 
(\Phi(X/E),E(X)^*\to \Phi(X/E)^{gp})
\ee
where $E(X)^*\to \Phi(X/E)^{gp}$ is the natural homomorphism  given by \be f\mapsto div(f)\in \Phi(X/E)^{gp}.\ee

\subparat{A Frobenioid associated to a rigid analytic curve}\label{ss:frobenioid-analytic-curve}
I want to now describe a Frobenioid associated to  the analytic space arising from a geometrically connected smooth quasi-projective curve over a $p$-adic field. This will lead (in \cref{ss:tempered-frobenioid}) to the construction of Mochizuki's tempered Frobenioid over an hyperbolic curve, which may be intuitively understood as the data of divisors descended (by Galois descent) from the universal cover of each tempered coverings of a fixed hyperbolic curve. So this  construction is a little more elaborate than the ones described above but the key ideas are the same as above.

\subsubpara Let $Z/E$ be a connected analytic space arising from the analytification of a geometrically connected, smooth quasi-projective variety of dimension one over a $p$-adic field $E$. I want to describe the construction of the tempered fundamental group given in \cite{andre-book} as this will be used in the construction of Mochizuki's Frobenioid associated to a hyperbolic curve over $E$. 

Let $\left\{(Z_i,z_i): i\in I\right\}$ be a countable, cofinal system of finite Galois coverings of $(Z,z_K)$ (with given  geometric base-points dominating the  geometric base-point $z_K:\sM(K)\to Z$.

\subsubpara Let $Z_i^\infty/Z_i$ be the topological universal covering of $Z_i$ \cite{berkovich-book}. Let $Z_i^{\infty,rig}$ be the associated (classical) rigid analytic space (see \cite{berkovich-book}). Then one can describe tempered fundamental group $\pit{Z/E}$ explicitly as follows:
\be 
\pit{Z/E}=\pi_1^{temp}(Z,z_K)=\invlim_{i\in I} \gal(Z_i^\infty/Z).
\ee

Then one has the following data
\be \Phi_{temp}(Z/E)=\divt(Z/E)=\invlim_{i\in I} Div_+(Z_i^{\infty,rig})^{\gal(Z_i^\infty/Z)},
\ee
where $\divt(Z_i^{\infty,rig})$ is the group of effective divisors of the (classical) rigid analytic space $Z_i^{\infty,rig}$;
and let 
\be 
\pdivt(Z/E)=\invlim_{i\in I} PDiv(Z_i^{\infty,rig})^{\gal(Z_i^{\infty,rig}/Z)},
\ee  
where $PDiv(Z_i^\infty)$ is the group of principal divisors on $Z_i^{\infty,rig}$ i.e. the group of divisors of non-zero meromorphic functions on $Z_i^{\infty,rig}$. One has a natural homomorphism 
\be \pdivt(Z/E)\to \divt(Z/E)^{gp}=\Phi_{temp}(Z/E)^{gp}.\ee

\subsubpara\label{ss:mochizuki-temp-frob} Mochizuki's Tempered Frobenioid of $Z/E$ is the data
\be 
\frob(Z/E)=\left(\Phi_{temp}(Z/E), \pdivt(Z/E)\to \Phi_{temp}(Z/E)^{gp}\right).
\ee
\subparat{The instructive example of Tate elliptic curves}
Let $K$ be an algebraically closed perfectoid field with an isometric embedding of $\Q_p$. It is instructive to understand this for the concrete case of a (projective) Tate elliptic curve $C/K$. Here by \cite{andre-book}, the (topological) universal cover of $C/E$ is $\Gm\to C$. Write $\hat{\Z}=\invlim_{n}\frac{\Z}{n\Z}$. Then one has
\be 
\pit{C/E}=\hat{\Z}\times \Z.
\ee
So for every integer $n\geq 1$ one has a Galois \'etale cover $C_n\to C$ with Galois group $\Z/n\Z$ and $C_n/K$ is again a Tate elliptic curve. So its universal covering is (again) $\Gm\to C_n$ and the tempered covering $\Gm\to C_n\to C$ has Galois group \be G_n=\Z/n\Z \times \Z\ee and 
one has 
\be \pit{C/\C_p}=\invlim_{n} G_n=\invlim_{n} ((\Z/n\Z) \times \Z)=\hat{\Z}\times \Z.
\ee

Let $\Gm^{an}$ be the analytic space over $K$ associated to the affine variety $\Gm/K$ where $\Gm/\C_p$ is the multiplicative group. By \cite{berkovich-book}, the topological universal cover of $C_n$ is \be
\Gm^{an}=C_n^\infty\to C_n
\ee given by the Berkovich's version of the Tate construction (\cite{berkovich-book}). Then the rigid analytic space associated to $\Gm^{an}$ is $\Gm^{an,rig}=\Gm^{rig}$. The rigid analytic function picture for Tate elliptic curves is discussed extensively in \cite{roquette-book}.

Let $K(C_n^{\infty,rig})$ be the function field of $C_n^{\infty,rig}$ (over $K$). Then Mochizuki's tempered Frobenioid which I describe below in \cref{ss:tempered-frobenioid} can be described in terms of the divisors on the rigid analytic spaces $\{C^{\infty,rig}_n\}$ and of functions 
\be 
\invlim_{n} K(C_n^{\infty,rig})^{G_n},
\ee
and one can think of $K(C_n^{\infty,rig})^{G_n}$ as analytic  functions on $C_n^{\infty,rig}$ with descent data with respect to the action of $G_n$.

\subparat{Mochizuki's construction of tempered Frobenioids}\label{ss:tempered-frobenioid}
Now let $X/E$ be a geometrically connected, smooth, hyperbolic curve over a $p$-adic field $E$. Let $K$ be an algebraically closed perfectoid field equipped with an isometric embedding of $E$. Let $\xan_E$ (resp. $\xan_K$) be the associated Berkovich spaces over $E$ and $K$ respectively).

\subsubparat{} Let me remark that by \cite[III, 1.4.10 and Remark after Definition 2.1.1]{andre-book} any tempered covering of $Y\to \xan_K$ is definable over a finite extension $E'$  (i.e. of the form $Y=Y_{E'}$) of $E$ contained in $K$ and by base extension any tempered covering of $Y_{E'}\to \xan_E$ provides a tempered covering of $\xan_K$ . Hence in what follows, one can work with tempered coverings of $\xan_K$.

Then for each tempered covering $Y_{E'}\to \xan_E$ one has the triple \be (\Phi_{temp}(Y/E')=\divt(Y/E'),\pdivt(Y_{E'})\longrightarrow \Phi_{temp}(Y/E')^{gp})\ee
defined in \cref{ss:frobenioid-analytic-curve}.

\subsubpara\label{ss:mochizuki-tempered-frobenioid} Mochizuki defines the tempered Frobenioid of $X/E$ as a category whose objects are Frobenioids of tempered coverings of $X/E$ defined  in \cref{ss:frobenioid-analytic-curve} and with obvious morphisms between ordinary Frobenioids arising from morphisms between the curves.  Mochizuki's tempered Frobenioid is the category whose objects are (Mochizuki considers the case $K=\ebh\isom \C_p$)
\begin{multline}
	\frob_{temp}(X/E;K)=\left\{  (\Phi_{temp}(Y/E')=\divt(Y_{E'}),\pdivt(Y_{E'})\longrightarrow \Phi_{temp}(Z/E)^{gp}): \text{ where } \right.\\
	\qquad \left. \vphantom{\pdivt(Y_{E'})} Y_{E'}\to \xan_K \text{ is a temp. covering defined over } E' \subset K
	\right\}\end{multline}
(with obvious morphisms between objects).
This is equipped with a functor \be\frob_{temp}(X/E;K)\to {\rm Temp}(X/E;K)\ee
from $\frob_{temp}(X/E;K)$ to the category of tempered coverings of $X/E$ defined over finite extensions of $E$ contained in $K$. This is given by the rule
\be 
(\divt(Y_{E}),\divt(Y_{E})^{gp}\longleftarrow\pdivt(Y_{E'}))\mapsto Y_{E'}
\ee
in which one forgets the divisors but remembers the tempered covering $Y_{E'}$ of $\xan_E$. Two remarks must be made here:

\brem
Note that for the construction of this Frobenioid one can work with the category of tempered coverings of $\xan_K$   because by  \cite[Proof of Proposition 5.1.1]{andre-book} every tempered covering of $\xan_K$ is defined over a finite extension $E'$ of $E$ contained in $K$.
\erem

\brem 
In Mochizuki's terminology , $\frob_{temp}(X/E;K)$ is  a Frobenioid over the base category of tempered coverings of $\xan_E$ and  is the Frobenioid of \cite[Example 3.2]{mochizuki-iut1} (I do not use Mochizuki's notation here). Mochizuki's definition of this Frobenioid in \cite{mochizuki-theta} does not use Berkovich's topological universal coverings of $\xan_K$ and instead works with ``combinatorial universal coverings'' defined using formal schemes associated to the tempered coverings of the analytic space $\xan_E$. But by \cite{berkovich-book}, one can pass from analytic spaces to formal schemes to go back and forth between my view point and Mochizuki's.
\erem

\subparat{Prime-strips from Frobenioids and gluing Frobenioids by prime-strips} The following theorem is now a consequence of various results proved earlier.
\bthm\label{th:gluing-frobenioids}
Let $X/L$ be as in \cref{ss:elliptic-curve-assumptions}. Let $\arith{L}_\by$ be an arithmeticoid given by $\by=(y_v)_{v\in\vl}\in\yadl$.  For each $v\in\vl$, let $\frob_{temp}(X/L_v;K_v)$ be the tempered frobenioid (see \cref{ss:mochizuki-tempered-frobenioid}) provided by the local component of $\arith{L}_\by$ at $v$. Then
$\frob_{temp}(X/L_v;K_v)$ determines each of the  prime-strips $\sF^{?}$ in the \Cref{tab:rosetta-stone3}.  Notably, the prime-strips of \iutthr\ crucial to \moccor:
\begin{align}
	\sF^{\vdash\times\mu}_{\scriptscriptstyle{Joshi}}(X/L_v;K_v) & = G_{L_v;K}\act \O_{\bL_v}^{\times \mu} \\ 
	\sF^{\vdash\blacktriangleright\times\mu}_{\scriptscriptstyle{Joshi}}(X/L_v;K_v) &=G_{L_v;K_v}\act \O_{\bL_v}^{\times \mu} \times q_{(X/L_v,X/K_v)}^\bN \\ 
	\sF^{\vdash\blacktriangleright}_{\scriptscriptstyle{Joshi}}(X/L_v;K_v)  &=   G_{L_v}\act q_{(X/L_v,X/K_v)}^\bN \\ 
	\sF^{\vdash\perp}_{\scriptscriptstyle{Joshi}}(X/L_v;K_v) &=G_{L_v}\act \O_{\bL_v}^{\mu_{2\ell}}\times q_{(X/L_v,X/K_v)}^\bN.
\end{align}
If $\frob_{temp}(Y/L_v;K_{v,1})$ and $\frob_{temp}(X/L_v;K_{v,2})$ are two tempered frobenioids arising from two holomorphoids $X/\arith{\by_1}$ and $X/\arith{\by_1}$ then one has an isomorphism between the tempered fundamental groups
$$\pit{X/L_v;K_{v,1}}\isom \pit{X/L_v;K_{v,2}},$$ and an isomorphism $$\sF_{\scriptscriptstyle{Joshi}}^?(X/L_v;K_{v,1})\isom \sF_{\scriptscriptstyle{Joshi}}^?(X/L_v;K_{v,2})$$
between each  prime-strip, given in \Cref{tab:rosetta-stone3}, corresponding to these two frobenioids.
\ethm

\bp 
The proof is immediate from \Cref{th:gluing-prime-strips}.
\ep

\begin{rmk*}
	One may view \Cref{th:gluing-frobenioids,th:gluing-prime-strips} as providing a ``gluing'' of the Frobenioids $\frob_{temp}(Y/L_v;K_{v,1})$ and $\frob_{temp}(X/L_v;K_{v,2})$  and Hodge-Theaters (see \cref{ss:hodge-theater}) by along their respective prime strips (see \cite[\ssep I2, page 21]{mochizuki-iut1}). There are other instances of similar ``gluing'' assertions in \iutthr.
\end{rmk*}

\subparat{Mochizuki's Hodge Theaters via \cite{joshi-teich}}\label{ss:hodge-theater}
Another device which is used in \iutthr\ is the notion of a \textit{Hodge-Theater} which is introduced in \cite{mochizuki-iut1}. Using \Cref{tab:rosetta-stone3}, the discussion of tempered Frobenioids in \cref{ss:tempered-frobenioid}, and the theory of arithmeticoids \cite{joshi-teich-def},  one can translate this notion into more transparent objects of \cite{joshi-teich}. Notably, this demonstrates the existence of a plurality of Hodge Theaters of all types considered in \iutthr. There are many different Hodge Theaters used in \iutthr, with increasing complexity but built roughly on the same principles. I will provide the simplest version but let me remark that such devices are unnecessary in \joshiros\ as more geometric objects are available.

\subsubpara Let $L,C/L,X/L, L',\lmod$ be as in \cref{ss:elliptic-curve-assumptions}, \inithtdata.  Let $\lmod$ be the field of moduli of $X/L$.  So $\lmod$ is  the field of moduli of $X/L$ and explicitly  $\lmod$ is the field generated over $\Q$ by the $j$-invariant of the projective elliptic curve giving rise to $X$. So one has $L\supset \lmod$.  

\subsubpara One has the adelic Arithmetic Teichmuller Space $\widetilde{\fJ(X/L)}$ constructed by \constrone{Theorem }{th:main4}. This is a product of local arithmetic Teichmuller spaces $\fjxlp$ for all $v$.  At archimedean primes one takes the usual Teichmuller spaces associated to the Riemann surface $X(L_v)$ (here $v=\infty$ and by hypothesis $L$ has no real embeddings). 

By \constrone{Theorem }{th:field-of-def}, for every prime $v$ of $L$,  every object of $\fja{X,L_v}$ is definable over $\lmod$.

\subsubpara Suppose that one is given a point $((X/L_v, X/K_{v})_{v})\in \widetilde{\fJ(X/L)}$ of the adelic arithmetic Teichmuller space $\widetilde{\fJ(X/L)}$ constructed in \cite{joshi-teich}.  Then by \constrtwoh{Proposition }{pr:arith-datum-frob} one obtains an arithmeticoid $\arith{L}_\by$ associated to $((X/L_v, X/K_{v})_{v\in\vl})\in \widetilde{\fJ(X/L)}$.

\bthm\label{th:lmod-arith}
Let $C,X,L,L'$ be as in \cref{ss:elliptic-curve-assumptions} and \inithtdata. Let $((X/L'_v, X/K_{v})_{v\in\vlp})\in \widetilde{\fJ(X/L)}\in \fJ(X/L)$ be a point of the adelic Teichmuller space of $X/L$ and $\arith{L'}_\by$ be the arithmeticoid of $L'$ provided by this point. Then 
\benumlab
\item the bijection (see \inithtdata) $$\V_{\lmod}\isom \underline{\V}\subset \V_{L'}$$ together with the arithmeticoid $\by$ provides us with an arithmeticoid $\arith{\lmod}_{\underline{\by}}$  with ${\underline{\by}}\in\yadlmod$ naturally given using $\by$ and the above bijection.
\item  In particular, every holomorphoid $\holt{X/L'}{\by}$ of $X/L'$ provides us an arithmeticoid $$\arith{\lmod}_{\underline{\by}}$$ of $\lmod$, 
\item and hence by \constrtwoh{Proposition }{pr:arith-datum-frob} a (perfect) frobenioid  $$\frob(\arith{\lmod}_{\underline{\by}}),$$  and a realified Frobenioid $$\frob(\arith{\lmod}_{\underline{\by}})^\R$$ which are isomorphic (by \constrtwoh{Proposition }{pr:arith-datum-frob}) to the perfection   $\frob(\lmod)^{pf}$ and the realification  $\frob(\lmod)^{\R}$ of the Frobenioid of the number field $\lmod$ respectively.
\eenum
\ethm
\bp To prove all the assertions, it is enough to provide a natural construction of $\underline{\by}\in\yadlmod$ from $\by\in\yadlp$ given by the data of the holomorphoid $\holt{X/L}{\by}$ and the bijection $\V_{\lmod}\isom \underline{\V}\subset \V_{L'}$ provided by \inithtdata. By definition (\constrtwoh{Definition }{def:adelic-ff-curves}) that $$\yadlp=\prod_{w\in\vlp} \abs{\syflwp}.$$ Now $\underline{\by}\in\yadlmod$ is the point obtained by projection of $\by=(y'_w)_{w\in\vlp}$ on to the coordinates corresponding to the ``index'' set identified by $\underline{\V}\isom\V_{\lmod}$ using the natural (continuous) mapping \cite[Proposition 2.3.20]{fargues-fontaine} $$f_{w|v}:\abs{\sY_{\cpt,L_w'}}\to \abs{\sY_{\cpt,{\lmod}_{,v}}}$$ 
i.e. define $$\underline{\by}=(f_{w|v}(y_w'))_{w\in\underline{\V}}\in\yadlmod'.$$
\ep

\subsubpara For each holomorphoid $\holt{X/L'}{\by}$ and for each $v\in\vlp$, let $(X/L'_v,X/K_{v})= \holt{X/L_v'}{y_v}$ be the (local) holomorphoid of $X/L_v'$. Let $\frob_{temp}((X/L'_v,X/K_{v}))_{y_v}$ be the tempered Frobenioid  of $X/L_v$ at $v$   associated to it by \cref{ss:tempered-frobenioid}, \constrtwoh{Definition }{def:frobenioid-number-field}.  

Let $\arith{\lmod}_{\underline{\by}}$  be the arithmeticoid of $\lmod$ given by \Cref{th:lmod-arith}. Let  $\frob(\arith{\lmod}_{\underline{\by}})$  (resp. $\frob(\arith{\lmod}_{\underline{\by}})^\R$)  be the  Frobenioid (resp. the realified Frobenioid) associated to  the arithmeticoid $\arith{\lmod}_{\underline{\by}}$ by \constrtwoh{Proposition }{pr:arith-datum-frob}.  Let $\frob(\lmod)^\R$ be the realification of the Frobenioid of the number field $\lmod$ (see \cref{ss:frobenioid-number-field}).

Then one can associate to $\holt{X/L'}{\by}$, a \textit{Hodge Theater}, in the sense of \cite[Definition 3.6]{mochizuki-iut1}, as follows:

\be\label{eq:hodge-theater-def} 
\mathcal{HT}((X/L'_v, X/K_{v})_{v})_{\by}=\left(\left\{\frob_{temp}((X/L'_v, X/K_{v}))\right\}_{v}, \frob(\arith{\lmod}_{\underline{\underline{\by}}})^\R\right).
\ee
I will simply write this Hodge-Theater as $$\mathcal{HT}(\holt{X/L'}{\underline{\by}})$$ or even more simply as
$$\mathcal{HT}(\holt{X/L}{})_{\underline{\by}}$$ to indicate the holomorphoid of $X/\arith{L'}_\by$ which it arises from.

By \constrtwoh{Proposition }{pr:arith-datum-frob} one may also write the above definition in the style of \cite{mochizuki-iut1} as
\be 
\mathcal{HT}((X/L'_v, X/K_{v})_{v})_{\by}=\left(\left\{\frob_{temp}((X/L'_v, X/K_{v}))\right\}_{v\in\vlp}, \frob({\lmod})^\R_{\underline{\by}}\right).
\ee
\brem 
Notably, Hodge Theaters in \iutthr\ should be labeled by the arithmetic holomorphic structures they arise from.
\erem

\subsubpara  This is an important point for \cite{mochizuki-iut3}:
\bthm\label{th:hodge-theaters-exist}
Each holomorphoid $\holt{X/L}{\by}$ provides   the Hodge-Theater $$\mathcal{HT}(\holt{X/L'}{})_{\by}$$ and all these distinct holomorphoids give rise to isomorphic Hodge Theaters. 
\ethm
\bp 
Note that by the theory of \cite{joshi-teich,joshi-untilts,joshi-teich-estimates} there are many inequivalent geometric data indexed by $\by\in\yadlp'$ and each provides the corresponding Hodge-Theater  $$\mathcal{HT}((X/L'_v, X/K_{v})_{v})=\left(\{\frob_{temp}((X/L'_v, X/K_{v}))\}_{v}, \frob^\R(\lmod)\right).$$

To prove all  these Hodge-Theaters are isomorphic, it suffices to note that each arises from the data of tempered coverings of $X/L_v'$ and the monoid of effective divisors  and principal divisors on each of these coverings. Assume that $v\in\vlp$ is non-archimedean. By \cite{andre03}, these tempered coverings of $X/L_v'$ are defined over finite extensions of $L_v'$ contained in $K_v$. For $v$ archimedean, a tempered covering is by definition a finite \'etale covering and hence  such a covering is defined over a finite extension of $L'$. So  for all $v\in\vlp$, one has a non-canonical isomorphism  between such data provided by any two holomorphoids. Hence the assertion is established.
\ep
\brem\  
\benumlab
\item Mochizuki defines other versions of Hodge Theaters $\mathcal{HT}^?$ with more refined conditions and decorations, but I will not recall these here. 
\item  Hodge Theaters of \iutthr\ consist of data coarser than the geometric data of \joshiros. All Hodge Theaters of \iutthr\ can be easily constructed by using the above methods and \present. 
\item The above proposition should be contrasted with \cite{scholze-review}. The central point which is missing in \iutthr\ (and in  \cite{scholze-review}) is the existence of geometrically inequivalent data providing the above Hodge-Theaters.  
\item Readers of \cite{scholze-review} may have drawn the (erroneous) conclusion that the equivalence of Hodge-Theaters invalidates the methodology and the claims of \iutthr. 
\item So let me say this clearly:  \moccor\ is not about Hodge-Theaters or fundamental groups  per se,  but about averaging arithmetic and geometric data arising from distinct arithmetic holomorphic structures. That distinct arithmetic holomorphic structures exist and that such an averaging is possible is asserted in \iutthr\ and unequivocally established in \present.
\eenum
\erem

\section{Appendix II: Perfect Frobenioids vs. Perfectoids}\label{se:perfect-frob-perfectoid} 
In this section, I want to provide a comparison between  the theory of perfect Frobenioids and the theory of perfectoids. Mochizuki's Frobenioid theoretic  (read multiplicative monoid theoretic) approach   requires decoupling addition and multiplication in the sense of this paper. Mochizuki asserts that this decoupling is possible by his group theoretic methods.  As noted in \constrtwoh{\ssep}{ss:formal-group-approach}, in   \cite{joshi-formal-groups}, I showed that one may also establish this decoupling algebraically keeping Mochizuki's monoid theoretic ideas as a guide;  an important observation I made early on is that the theory of perfectoid spaces of \cite{scholze12-perfectoid-ihes}  also provides such a decoupling  (as observed in this paper) via existence of general Fargues-Fontaine curves after \cite{kedlaya-liu15,kedlaya-liu19} and \cite{scholze12-perfectoid-ihes,scholze-weinstein-book}.

\subparat{A useful multiplicative sheaf}
The following will be useful in what follows:
\bpro 
Let $X/k$ be a geometrically connected, smooth, quasi-projective variety over a field $k$. For each $n\in\Z_{\geq0}$, let $\phi_n:\O_X^*\to \O_X^*$ be the homomorphism of sheaves of abelian groups given by ${x\mapsto x^n}$. Consider the (inverse) limit presheaf $$\O_X^{*pf}=\invlim_n (\O_X^*,\phi_n).$$
Then $\O_X^{*pf}$ is a sheaf of abelian groups whose sections over open subsets (and hence also stalks) are perfect monoids in the sense of \cite{mochizuki-frobenioid1}. Notably the stalks of $\O_X^{*pf}$ are perfect monoids i.e. the stalks are $\Q$-vector spaces.
\epro
\bp 
This is immediate from the well-known fact that limits of sheaves of abelian groups (on quasi-projective or more generally noetherian schemes) are sheaves. 
\ep

Along with $(X,\O_X^{*pf})$  may similarly consider the pair $(X,(\O_X,+)^{pf})$ with 
$$(\O_X,+)^{pf}=\invlim_n (\O_X,+),$$
where the inverse limit is over the homomorphism $\O_X\to \O_X$ of abelian groups given by $x\mapsto n\cdot x$ for local sections $x$ of $\O_X$. The sheaf of abelian groups $(\O_X,+)^{pf}$ is a sheaf of perfect monoids in the sense of \cite{mochizuki-frobenioid1}. 

More generally, suppose $X/k$ is an integral scheme, and let $K_X^*$ be the sheaf  of invertible, meromorphic  functions on $X$, the one may also consider $(X,K_X^{*pf})$ instead.

\brem
Notably the pairs $(X,\O_X^{*pf})$ (resp. $(X,\O_X^{*pf})$) consisting of the topological space $X$ equipped with the perfect multiplicative monoid (perfect additive monoid) (both are sheaves of abelian groups in fact) of functions on $X$   exists even if $X$ does not admit a perfectoidification in the sense of \cite{scholze12-perfectoid-ihes}.
\erem

\subparat{Multiplicatoids}
Now let me explain how one may decouple additive and multiplicative structures. Instead of working with the scheme structure $(X,\O_X)$, one can work with  the pair $(X,\O_X^{*pf})$. I call $(X,\O_X^{*pf})$ the \textit{perfect multiplicatoid of the scheme $(X,\O_X)$}.  Similarly, one calls the pair $(X,(\O_X,+)^{pf})$, the \textit{additivoid of $X$}.

One may think of the $(X,\O_X^{*pf})$ and $(X,(\O_X,+)^{pf})$ as a decoupling of the multiplication and addition structures of the scheme $(X,\O_X)$.

\brem
Simplest example of this is the case when $X=\Spec(\bL)$ where $\bL$ is the algebraic closure of a number field $L$. The case of the multiplicatoid $(\abs{\Spec(\bL)}, \bL^{*pf})$ (here $\abs{\Spec(L)}$ is the topological space underlying $X$, this of course, a one point topological space) is implicitly considered in \cite{mochizuki-frobenioid1} and importantly in \iut. Similarly, if $p$ is a prime, then the multplicatoid $(\abs{\Spec(\bQ_p)}, \bQ_p^{*pf})$ of $\Spec(\bQ_p)$ is of fundamental importance in \iut.
\erem

\subparat{A basic problem} 
An important question which comes up in the context of \iut\ is this: Given the multiplicatoid of $X$, recover or reconstruct \textit{some} scheme structure on $X$ compatible with this multiplicative structure given by the multiplicatoid on $X$. 

More precisely, consider the problem of finding schemes $(Y,\O_Y)$ with isomorphisms between their multiplicatoids
\be\label{eq:multiplicoid-eq}
(Y,\O_Y^{*pf})\isom (X,\O_X^{*pf}).
\ee 
More generally, if $M_X\subseteq K_X^{*pf}$ some suitably chosen subsheaf of monoids of the perfection of the sheaf, $K_X^*$, of invertible functions on $X$, then one may consider the same question with $(X,M_X)$ in place of $(X,\O_X^{*pf})$.

\textit{This is the precise version of the idea behind Mochizuki's decoupling of addition and multiplication and reconstruction of a new scheme structure, namely $(Y,\O_Y)$, from the multiplicatoid of $X$.} 
This point is of course not clearly enunciated in \iut\ or in \cite{mochizuki-essential-logic} and nor is the mechanism for such a claim  established.  

\subparat{How do multiplicatoids arise in the anabelian contexts?} To understand how multiplicatoids arise naturally in anabelian contexts, let me present a reformulation of the classical anabelian reconstruction results for number fields and $p$-adic fields (\cite{hoshi-number}, \cite{hoshi-mono} are great modern references for this):
\bthm 
Let $k$ be either a number field or a $p$-adic field for some prime $p$ and let $\bar{k}$ be an algebraic closure of $k$ and let $G_k$ be the absolute galois group of $k$ computed using $\bar{k}$. Then from the topological $G_k$ one can reconstruct some isomorph of the multiplicatoid $(\abs{\Spec(\bar{k})}, \bar{k}^{*pf})$ of $\bar{k}$ and this reconstruction is functorial in $k$. In other words, $(\abs{\Spec(\bar{k})}, \bar{k}^{*pf})$ is amphoric.
\ethm

\subpara In \cite{joshi-formal-groups}, I showed that one may in fact algebraize the idea of decoupling of the additive and multiplicative structures by means of monoid formal group laws and universal monoid formal group law established in \cite{joshi-formal-groups} (a precise version of this relevant to the main results of this paper  is recorded in \constrtwoh{Theorem }{th:formal-group-law}). 

\subpara My first important observation is that the existence of distinct schemes $(Y,\O_Y)$ in \eqref{eq:multiplicoid-eq}  is difficult to establish convincingly if one works with the category of schemes (as Mochizuki does). This difficulty already appears in the special case of of a geometric schematic point $X=\Spec(\bQ_p)$ which is relevant to \iut. For \iut, one needs a precise quantification of what ``some'' means in the above theorem.

\subparat{Multiplicatoids from perfectoid spaces} My second important observation is that Scholze's work \cite{scholze12-perfectoid-ihes} provides solution to a similar problem in a manner compatible with Mochizuki's anabelian view point. Here is precise formulation:

\bthm
Let $(X,\O_X)$ is a perfectoid space (for simplicity, one may assume that $X$ arises from a perfectoid algebra) over an algebraically closed perfectoid field. Suppose that  $(X^\flat,\O_{X^\flat})$ is the tilt of $(X,\O_X)$. Let $(Y,\O_Y)$ be another untilt of $(X^\flat,\O_{X^\flat})$   then 
\benumlab
\item one has a homeomorphism of pairs consisting of a topological space and a sheaf of $p$-perfect multiplicative monoids 
\be\label{eq:multiplicoid-eq2}
(\abs{Y},\invlim_{x\mapsto x^p} \O_Y)\isom (\abs{X},\invlim_{x\mapsto x^p} \O_X).
\ee 
\item and there exists pairs $(X,\O_X)$ of perfectoid spaces which are not isomorphic, but for which the above isomorphism holds.
\item Moreover, for a suitable choice of geometric basepoints, one has an isomorphism of \'etale fundamental groups
\be 
\pi_1^{et}(X)\isom \pi_1^{et}(Y).
\ee
\eenum
\ethm
\bp By \cite[Theorem 1.8]{scholze12-perfectoid-ihes} one has an isomorphism of sheaves of topological rings 
(and hence of their subsheaves of multiplicative monoids)
\be\O_{X^\flat}=\invlim_{x\mapsto x^p} \O_X,\ee
and by \cite[Theorem 6.2]{scholze12-perfectoid-ihes} one has a homeomorphism of topological spaces 
\be 
\abs{X}\isom \abs{X^\flat}.
\ee 
Thus one obtains an isomorphism of  pairs for $X,X^\flat$
\be(\abs{X^\flat},\O_{X^\flat})\isom(\abs{X},\invlim_{x\mapsto x^p} \O_X).\ee
In particular, the isomorphism class of the pair for $X$ depends only on $X^\flat$. Hence if $(X,\O_{X}),(Y,\O_{Y})$ are two untilts of $(X^\flat,\O_{X^\flat})$ then one has an isomorphism of the  stated pairs. This proves {\bf(1)}.

In general the two perfectoid spaces with isomorphic tilts  need not be isomorphic.  An explicit example of this is provided, in the case $X$ arises from a perfectoid field, by the main theorem of \cite{kedlaya18}.  This proves {\bf(2)}.

Now suppose one chooses geometric basepoints for $X,X^\flat$ then \cite[Theorem 7.12]{scholze12-perfectoid-ihes} shows that one has an isomorphism of \'etale fundamental groups of $X$ and $X^{\flat}$, and hence the isomorphism class of this group depends only on that of $X^\flat$. This proves {\bf(3)}.
\ep 

\subparat{Comparison} 
\brem\ 
\benumlab
\item So one has a non-trivial solutions to the above problem  in the category of perfectoid spaces.
\item Hence for a perfectoid space $X$ over an algebraically closed perfectoid field, the sheaf $\O_{X^\flat}=\invlim_{x\mapsto x^p} \O_X$
plays a role similar to the sheaf encapsulating the multiplicatoid of a scheme considered above.
\eenum
\erem

\brem 
My approach to the aforementioned problem, developed in the present series of papers, is to work with Berkovich analytic spaces. This has an added advantage which is not available in the category of schemes: in category of Berkovich analytic spaces, generally speaking, logarithms of invertible functions are available and  one can recover the additive structure  of functions (i.e. recover the sheaf of abelian groups $(\O_X,+)$--or some subsheaf thereof) from the multiplicative datum $\O_X^{*pf}$ by passage to logarithms of suitably invertible functions i.e. one may treat the additive and multiplicative structure as being interchangeable.
\erem

\brem Now to compare the two notions:  the theory of perfect Frobenioids or more precisely perfect multiplicatoids is available for any geometrically connected smooth hyperbolic curve over a $p$-adic field (and more generally for any smooth variety over a $p$-adic field). More precisely the perfect multiplicatoid \eqref{eq:multiplicoid-eq} exists for any scheme $X$. But as has been remarked in \cite{scholze14-icm2}, the `perfectoidification' of a hyperbolic curve does not exist (even for hyperbolic curves over algebraically closed perfectoid fields).

Hence, the theory of pairs \eqref{eq:multiplicoid-eq} implicit in \cite{mochizuki-frobenioid1} is at once more general but at the same time weaker than the theory of perfectoid spaces of \cite{scholze12-perfectoid-ihes}. 
\erem

\newcommand{\glqp}{{\rm GL}_2^+(\Q)}
\newcommand{\tslq}{\widetilde{{\rm SL}_2(\Q)}}
\newcommand{\tslql}{{\widetilde{{\rm SL}_2(\Q)}}^\ells}
\newcommand{\slq}{{\rm SL_2}(\Q)}
\newcommand{\tglqp}{\widetilde{{\rm GL}_2^+(\Q)}}
\newcommand{\tglqpl}{{\widetilde{{\rm GL}_2^+(\Q)}}^\ells}
\newcommand{\slr}{{\rm SL}_2(\R)}
\newcommand{\tslr}{\widetilde{{\rm SL}_2(\R)}}
\newcommand{\slz}{{\rm SL}_2(\Z)}
\newcommand{\tslz}{\widetilde{{\rm SL}_2(\Z)}}
\newcommand{\fH}{\mathfrak{H}}
\newcommand{\syi}{\sY_{\infty}}
\newcommand{\syis}{\sY_{\infty,s}}
\newcommand{\syils}{\sY_{\infty,s}^\ells}
\newcommand{\sxi}{\sX_\infty}

\section{Appendix III: The proof of geometric Szpiro inequality revisited}\label{se:appendix-geom-case}
\subparat{The proofs of geometric Szpiro inequality} In \cite{bogomolov00}, the authors gave a proof of geometric Szpiro inequality for elliptic fibrations over a genus zero curve and in \cite{zhang01} this was generalized to elliptic fibrations over base curves of genus $g\geq 0$. In Spring 2018, I was on a sabbatical to RIMS, Kyoto with Shinichi Mochizuki and he gave me a few private lectures on his approach ( \cite{mochizuki-bogomolov}) to these works. As Mochizuki informs us in \cite{mochizuki-bogomolov}, Ivan Fesenko had pointed out to him the existence of \cite{bogomolov00}, \cite{zhang01}.

Note that \cite{bogomolov00}, \cite{zhang01} and \cite{mochizuki-bogomolov} deal with  the same theorem (namely the geometric Szpiro inequality) from three different perspectives. 

Here I will record a fourth distinct approach to understanding these three approaches to geometric Szpiro inequality which I discovered in the course of my work on arithmetic Teichmuller spaces. Notably the perspective presented below exhibits all the key features of \iut, and also the perspective of this paper and a reading of this is essential for anyone who wishes to understand \iut, my work and the relationship between these two theories. In particular,  I provide a formulation (and a proof) of \cite[Corollary 3.12]{mochizuki-iut3} in the context of the geometric Szpiro inequality (this is not formulated in \cite{mochizuki-bogomolov}).

\textcolor{red}{This section should be read in tandem with \cite{zhang01}, \cite{bogomolov00}.} (\cite{mochizuki-bogomolov} largely follows these references with added commentary).

\subparat{The setup} Let $C$ be a connected, smooth, projective curve over $\C$ of genus $g\geq 0$. Let $f:X\to C$ 
be a proper, generically smooth morphism with geometric generic fiber of genus one.  This fibration will be fixed for the remainder of this paper. Additionally I will assume that this is not an isotrivial fibration.

\subparat{The universal cover of $\slr$ and its properties} Let $\slr$ be the topological group of $2\times2$ matrices with real entries and of determinant one. Let $\tslr$ be the universal cover of $\slr$. Then $\tslr$ is described as a central extension (\cite[Lemma 3.5]{zhang01})
\be\label{eq:slr-cover}
1\to \langle z^2\rangle \to\tslr \to \slr \to 1
\ee
where $z\in\tslr$ is a certain element, which generates the center of $\tslr$ and which is described explicitly in \cite{zhang01} and so I will not recall that description here.

Let $$\tslz\subset \tslr$$ be the inverse image of $\slz\subset\slr$. Let $$\tslq\subset \tslr$$ be the inverse image of $\slq\subset \slr$.

\brem 
Note that by classical results or the explicit description in \cite{bogomolov00} or \cite{zhang01}, $\tslz$ is isomorphic to the braid group on three letters (see \cite[Theorem 16]{rawnsley12}).
\erem

\subparat{Classical Teichmuller space in genus one} Let $\fH$ be the complex upper half plane. Then one has a natural mapping 
$$\slr\to \fH$$
given by $$\begin{bmatrix} a & b \\ c & d \end{bmatrix}\mapsto \frac{a\cdot i+b}{c\cdot i +d}.
$$

If $\tau\in \fH$, I will write $E_\tau$ for the elliptic curve given by the lattice $[1,\tau]\subset \C$.

Recall from \cite{imayoshi-book} that $\fH$ is the Teichmuller space of genus one Riemann surfaces and the usual modular group $\slz\subset \slr$ can be identified with the Teichmuller modular group (in genus one).

\subparat{Global (geometric) Fargues-Fontaine curves $\syi$, $\sxi$, and the global Frobenius morphism $\vphi_\infty$}\label{app:global-frobenius} Let me now illustrate how the geometric analogues of $\yadl$ i.e. global Fargues-Fontaine curve type objects with a similar geometry appears in the context of \cite{zhang01}, \cite{bogomolov00}. This is an important point in understanding my approach and Mochizuki's approach to global Diophantine problems via our respective theories.
\bdefn\label{def:global-frob-classical-case} 
Let $\syi=\tslr$ and $\sxi=\slr$. I will refer to $\syi$ as the \textit{incomplete global Fargues-Fontaine curve} and $\sxi=\slr$ as the \textit{complete global Fargues-Fontaine curve.} I will write  
\be\label{eq:global-frob-classical-case}\vphi_\infty=z^2\in\tslr\ee for the generator of the central extension \eqref{eq:slr-cover} and refer to $\vphi_\infty$ as the \textit{global Frobenius morphism of $\syi$}. In particular,  with these definitions, and because of \eqref{eq:slr-cover}, one has
$$\syi/\vphi_\infty^\Z\mapright{\isom} \sxi$$
is a quotient by the global Frobenius morphism.
\edefn

\brem\
\benumlab
\item The relationship with the global arithmetic constructions of \cite{joshi-teich-def} is this: $\syi$ (resp. $\sxi$) is the analog of global curve $\yadl$ (resp. $\sxi=\yadl/\vphi^\Z$).
\item  Let me spell out the analogy  to the situation of usual (local) Fargues-Fontaine curves where one has a quotient of adic spaces
$$\syfe/\vphi^\Z\mapright{\isom}\sxfe$$
with $\vphi$ the Frobenius morphism of $\syfe$. 
\item Note however that the usual Fargues-Fontaine curve $\syfe,\sxfe$ are local objects from the the point of view of theory of the present series of papers, while $\syi,\sxi$ are global objects in the context of \cite{zhang01}, \cite{bogomolov00}.
\item There is also an alternate global analog in the present theory: $\yadl$ and the quotient (topological stack) $[\yadl/L^*]$ are also global analogs of $\syi,\sxi$ in the context of the present series of papers (and hence in the context of \iut).
\item Notably $\yadl$ equipped with the action of $L^*$  or its global Frobenius  morphism $\boldsymbol{\varphi}$ given by \constrtwoh{Theorem }{th:galois-action-on-adelic-ff} and \constrtwoh{Corollary }{cor:global-frobenius} are the global arithmetic analogs of $\syi$ and its Frobenius morphism $\varphi_\infty$.
\eenum
\erem

\subparat{Existence of $\flog$-links} In \iut\ Mochizuki works with $\flog$-links. Let me discuss how this appears in the context of \cite{zhang01}, \cite{bogomolov00}. Let $g_\tau\in\slr$ be a lift of $\tau\in \fH$ to $\slr$. Let $\tilde{g}_\tau\in\tslr$ be the lift of $g_\tau$ to $\tslr$.
\bdefn\label{de:geom-log-links}
Let $\tilde{g}_{\tau,0}\in\tslr$ be a lift of $\tau\in\fH$ to $\tslr$. If $\tilde{g}'_{\tau,0}$ is another lift of $\tau$ then I will say that $\tilde{g}_{\tau,0}$ and $\tilde{g}_{\tau,0}'$ are $\flog$-linked if $$\tilde{g}'_{\tau,0}=\tilde{g}_{\tau,0}\cdot \vphi_\infty.$$
I will call $$\left\{\tilde g_{\tau,n}=\tilde{g}_{\tau,0}\cdot \vphi_\infty^n: n\in \Z \right\}$$
a chain of $\flog$-links.
\edefn

\brem
Note that a chain of $\flog$-links is the fiber over $\syi\to\sxi$ over the common image of $\tilde{g}_{\tau,n}$ in $\sxi=\slr$.
\erem

\subparat{The Schottky parameter}\label{ss:schottky-parameterization} Let me begin with the following 
\blem\label{le:link-lemma}
Let $\tau\in\fH$. For any  $\alpha\in\Q^{>0}$,  the elliptic curves $E_1=\C/\Lambda_1$ and $E_{\tau,\alpha}=\C/\Lambda_{\tau,\alpha}$ given by the lattices $\Lambda_{\tau,1}=[1,\tau]\subset \C$ and $\Lambda_{\tau,\alpha}=[1,\alpha\cdot\tau]$ are $\C$-isogenous.
\elem
\bp 
It is sufficient to note that 
$\begin{pmatrix}
1 & 0\\
0 & \alpha
\end{pmatrix}\in\glqp$ and
$$\begin{pmatrix}
1 & 0\\
0 & \alpha
\end{pmatrix} \cdot \begin{pmatrix}
1 \\ \tau
\end{pmatrix} =\begin{pmatrix}
1 \\ \alpha\cdot\tau
\end{pmatrix}
$$ takes $\Lambda_1$ to $\Lambda_\alpha$. So the assertion is clear.
\ep

\bdefn 
For an elliptic curve $E$ given by a lattice $[1,\tau]\subset\C$ let $q_E=e^{2\pi \cdot \sqrt{-1}\tau}$ be its Schottky parameter. Then one has the Schottky uniformization
$$E\isom \C^*/q_E^\Z.$$
\edefn
\bp 
This is standard.
\ep 

\blem\label{le:link-schottky-scaling} 
For the elliptic curves $E_1$ and $E_{\tau,\alpha}$ one has the following relationship between their Schottky parameters:
$$q_{E_{\tau,\alpha}}=q_{E_1}^\alpha.$$
\elem
\bp 
This is an elementary calculation:
$$q_{E_{\tau,\alpha}}=e^{2\pi \cdot \sqrt{-1}(\alpha\cdot\tau)}=e^{2\pi \cdot \sqrt{-1}\cdot \alpha\cdot\tau}=q_{E_1}^\alpha.$$
\ep

\subparat{Existence of $\Theta_{gau}$-links}\label{pa:schottky-param}
Now I want to discuss how Mochizuki's $\Theta_{gau}$-link appears in the geometric case of \cite{zhang01}, \cite{bogomolov00}. Let $\ell\geq 5$ be a prime number, let $\ells=\frac{\ell-1}{2}$. Let $\tau\in \fH$ and  write $q=e^{2\pi\sqrt{-1}\tau}$. 
From now on, for $j=1,2,\ldots,\ells$, for simplicity of notation, I will write $$E_{\tau,j}=E_{\frac{\tau}{2\ell},j^2}.$$

Note that the lattice $[1,\frac{\tau\cdot j^2}{2\ell}]$ for $E_{\tau,j}$ can be obtained from the lattice $[1,\frac{\tau}{2\ell}]$ of $E_{\tau,1}$ by using the matrix
$$\begin{pmatrix}
j & 0\\
0 & \frac{1}{j}
\end{pmatrix}\in \slq$$
which takes operates on $\fH$ by $$\tau\mapsto \begin{pmatrix}
j & 0\\
0 & \frac{1}{j}
\end{pmatrix}\cdot \tau=\frac{j\cdot \tau}{j^{-1}}={j^2\cdot \tau},$$
and hence $\frac{\tau}{2\ell}\mapsto \frac{j^2\cdot \tau}{2\ell}$ under its action.

\brem 
The reason I want to consider $\slq$ in place of $\glqp$ is that $\slq\subset \slr$. Working with $\glqp$ would require replacing $\slr$ by ${\rm GL}_2^+(\R)$ in \cite{zhang01}, \cite{bogomolov00}--this can be done, but presently I have not written down all the changes this necessitates in \cite{zhang01}, \cite{bogomolov00}.
\erem

\blem 
Let $$\alpha_j=\begin{pmatrix}
j & 0\\
0 & \frac{1}{j}
\end{pmatrix}\text{ for }j=1,2,\ldots,\ells.$$
Then the $\ells$-tuple $(\alpha_1,\ldots,\alpha_\ells)\in\slq^\ells$ operates on $\fH^\ells$, by the usual Mobius action on each of the factors, and under this action the diagonal tuple $$(\tau,\tau,\ldots,\tau)\in\fH^\ells$$ is mapped to 
$$(1^2\cdot\tau,2^2\cdot\tau,\ldots,\ells^2\cdot\tau).$$
\elem
\bp 
This is completely clear.
\ep

\bdefn 
Let $\ell\geq5$ be a prime number.
A $\Theta_{gau}$-link is an $\ells$-tuple $$(\tilde{g}_1,\ldots,\tilde{g}_{\ells})\in\tslql\subset\syi^\ells$$ lying over a tuple of elliptic curves $(E_{\tau,1},E_{\tau,2},\ldots,E_{\tau,\ells})$ for some $\tau\in\fH$ (with $E_{\tau,j}$ given by the above convention  and $j=1,2,\ldots,\ells$). The ``link'' here is the rule
$$E_\tau\mapsto (E_{\tau,1},E_{\tau,2},\ldots,E_{\tau,\ells}).$$
\edefn
\brem 
One may think of $E_\tau\mapsto (E_{\tau,1},E_{\tau,2},\ldots,E_{\tau,\ells})$ as a correspondence i.e. a one-to-many function on the moduli of elliptic curves arising from the correspondence on $\fH$  given by the correspondence:
$$\tau\mapsto (1^2\cdot \tau,2^2\cdot \tau, \ldots,\ells^2\cdot \tau)\in\fH^\ells.$$
\erem

\bthm\label{th:theta-tuples}
Let $(\tilde{g}_1,\ldots,\tilde{g}_\ells)\in\syi^\ells$ be a $\Theta_{gau}$-link. Let $E_{\tau,1}, E_{\tau,2}, \ldots,E_{\tau,\ells}$ be the corresponding elliptic curves. Then 
\benumlab
\item The elliptic curves $E_{\tau,1}, E_{\tau,2} ,\ldots , E_{\tau,\ells}$ correspond to distinct points of the genus Teichmuller space $T_1=\fH$.
\item The elliptic curves $$E_{\tau,1}, E_{\tau,2} ,\cdots , E_{\tau,\ells}$$ are all isogenous to $E_\tau$.
\item One has the  tuple of the corresponding Schottky parameters
$$(q_{\tau,1}, q_{\tau,2},\ldots,q_{\tau,\ells}) \in(\C^*)^\ells$$ corresponding to $E_{\tau,1}, E_{\tau,2} ,\ldots , E_{\tau,\ells}$.
\item Let   $P_1$ correspond to the $\ell$-torsion point on the double cover $E_{\tau/2}$ of $E_\tau$ given by $\tau'=\frac{\tau}{\ell}$ and choose a $\Theta$-function as in \cite{mochizuki-theta}. Then
the above tuple of Schottky parameters
$$(q_{\tau}^{1^2/2\ell},q_\tau^{2^2/2\ell},\ldots,q_\tau^{\ells^2/2\ell})\in (\C^{*})^\ells$$
corresponds to the set of values of the chosen $\Theta$-function at the $\ell$-torsion points $$\left\{P_j=j\cdot P_1\in E_{\tau/2}[\ell]: j=1,2,\ldots,\ells \right\}.$$
\eenum
\ethm
\bp 
The proof is clear from the previous lemmas.
\ep

\brem 
Since $\begin{pmatrix}
j & 0\\
0 & \frac{1}{j}
\end{pmatrix}\in\slq\subset \glqp$, this theorem makes it clear that one can think of a $\Theta_{gau}$-Link as arising  from natural correspondences on the upper half plane in a manner similar to Hecke correspondences on classical modular curves. By \Cref{le:link-schottky-scaling}, these correspondences provide a (non-trivial) scaling $$q_\tau\mapsto q_\tau^{j^2}$$ of the Schottky parameters. In the $p$-adic case, \cite{joshi-teich-estimates} takes a similar approach via correspondences on suitable Fargues-Fontaine curves for the construction of $\Theta_{gau}$-links in the context of \iut.
\erem
\newcommand{\bTheta}{\boldsymbol{\Theta}}

\subparat{The $\flog$-$\Theta$-lattice}
Let me note that $\slq,\tslq$ are both countable sets and hence so are $\slq^\ells$ and $\tslql$. Choose an enumerating function to enumerate the elements of $$\slq^\ells=\{\theta_{n}:n\in\Z\}.$$
Then one may  enumerate the set $\tslql$ as 
$$\tslql=\{\theta_{n,m}:(n,m)\in\Z\times \Z\},$$
with the canonical surjection $$\tslql\to \slq^\ells$$ being given by 
$$\theta_{n,m}\mapsto \theta_n.$$
Following \cite{mochizuki-iut1,mochizuki-iut3}, one could represent the  data $\tslql=\{\theta_{n,m}\}$   as a ``lattice'' in the plane with the vertical direction corresponding to the fiber $\{\theta_{n,m}:m\in\Z \}$ of $\tslql\to\slq^\ells$ over $\theta_n\in\slq^\ells$ i.e. over a $\Theta_{gau}$-link. 
In this picture the fiber over $\theta_n$ consists of $\flog$-links i.e. 
while the horizontal direction being represented as $\Theta_{gau}$-links. In \cite{mochizuki-iut3} this sort of presentation of $\Theta_{gau}$-links and $\flog$-links is called the \textit{$\log$-$\Theta$-lattice}:
$$
\xymatrix{
	&  &  & \\ 	
	\cdots \ar@{-}[r] &\theta_{n-1,m+1}\ar@{-}[r]\ar@{-}[d]_{\flog=\vphi_\infty}\ar@{-}[u]^{\flog=\vphi_\infty} & \theta_{n,m+1}\ar@{-}[d] \ar@{-}[u]_{\flog=\vphi_\infty}\ar@{-}[r] & \cdots\\ 
	\cdots \ar@{-}[r] & \theta_{n-1,m} \ar@{-}[r]\ar@{-}[d]_{\flog=\vphi_\infty} &  \theta_{n,m}\ar@{-}[u]_{\flog=\vphi_\infty} \ar@{-}[r]\ar@{-}[d]^{\flog=\vphi_\infty} & \cdots,\\
	&  &  & 
}$$ 
\brem
This should be compared with Mochizuki's discussion of the $\log$-$\Theta$-lattice in \cite[\ssep I, Page 405]{mochizuki-iut3} and its construction  in \cite[\ssep 1, Page 427]{mochizuki-iut3}.
\erem
\subparat{Hodge structures from Schottky parametrization}\label{ss:hodge-str-schottky}
\textcolor{red}{This is an elaboration of \cite[\ssep 22]{joshi-anabelomorphy}.}
Let $E_\tau$ be an elliptic curve with period lattice $[1,\tau]\in\fH$. Then the Schottky parameterization \ssep\ref{pa:schottky-param} provides a mixed Tate Hodge structure 
$$H_{Schottky}(E_\tau)\in \Ext^1_{\Z-MHS}(\Z(0),\Z(1))\isom \C^*$$ given via \cite{deligne-local} by $$H_{Schottky}(E_\tau)=q_\tau=e^{2\pi \sqrt{-1}\tau}\in\C^*.$$  This will be useful in understanding how one may carry out constructions in cohomology theory. In \iut\ cohomology theory means Galois cohomology theory $H^1(G_{L_v},\Q_v(1))$ and one has $$H^1(G_{L_v},\Q_v(1))=\Ext^1_{G_{L_v}}(\Q_v(0),\Q_v(1))$$ for each $v\in\V_L$, so the analogy I make here should not be lost on the readers.

\brem Note that one way to understand the Schottky parametrization  is to recognize that it asserts that the elliptic curve $E_\tau$ is  ordinary at an archimedean prime (and this ordinarity is witnessed by the Schottky Hodge structure $H_{Schottky}(E_\tau)$). \textit{However note that the usual Hodge structure provided by $H^1(E_\tau,\Z)$ is pure of weight one.} Tate's parameterization in the $p$-adic context is the $p$-adic version of Schottky's paramterization.
\erem

Thus one has 
\blem 
Any tuple  $(q_{\tau,1}, q_{\tau,2},\ldots,q_{\tau,\ells}) \in(\C^*)^\ells$ in \Cref{th:theta-tuples} provides a tuple of $\Z$-mixed Hodge structures  $$\left(H_{Schottky}(E_{\tau,1}), H_{Schottky}(E_{\tau,2}) ,\ldots , H_{Schottky}(E_{\tau,\ells})\right) \in \Ext^1_{\Z-MHS}(\Z(0),\Z(1))^\ells.$$
\elem

\subparat{A prelude to the definition of the $\Theta$-values locus} 
Let $$\bTheta_{classical,\tau}=\left\{ (\tilde{g}_1,\ldots,\tilde{g}_\ells)\in\tslql\subset  \syi^\ells: (\tilde{g}_1,\ldots,\tilde{g}_\ells) \text{ provides } E_{\tau} \right\}.$$

Let $S\subset C(\C)$ be the finite set of points over which the fibers of $X\to C$ are singular. 

Let $s\in S$ be a singular point of $f:X\to C$. One would like to understand $$\bTheta_{classical,\tau_s}$$ for $\tau_s\in \fH$ in a sufficiently close i.e. in a small neighborhood of the singular point $s\in S$ of $f:X\to C$ so that the local description of \cite{zhang01} remains valid (for each singular point $s\in S$). 

\subparat{Mochizuki style Theta-values locus for \cite{zhang01},
	\cite{bogomolov00}} So to construct a $\Theta$-values locus in the style of \cite{mochizuki-iut3} one does the following. Let $\syis$ be a copy of $\syi$ indexed by an $s\in S$.

Let \be\label{eq:theta-value-locus}\bTheta_{classical}=\left\{ (\tilde{g}_{s,1},\ldots,\tilde{g}_{s,\ells})_{s\in S}\in   \prod_{s\in S} {\syils} \text{ and } (\tilde{g}_{s,1},\ldots,\tilde{g}_{s,\ells})_{s} \in \bTheta_{classical,\tau_s} \text{ for all } s\in S  \right\}\ee
Each $(\tilde{g}_{s,1},\ldots,\tilde{g}_{s,\ells})_{s\in S}\in\bTheta_{classical}$ is 
the analog of the  \textit{$\Theta$-pilot object}, in the sense of \cite{mochizuki-iut3}, in the context of \cite{zhang01,bogomolov00}. 

One may also define the theta-values set $$\bTheta_{classical}^{val}\subset \Ext^1_{\Z-MHS}(\Z(0),\Z(1))^\ells,$$ in the style of  \iut, for the context of \cite{zhang01, bogomolov00} by taking the tuples of Schottky parameters provided by each tuple $(\tilde{g}_{s,1},\ldots,\tilde{g}_{s,\ells})_{s\in S}\in\bTheta_{classical}$.

\subparat{Stablized height functions $h$ and $H$}\label{ss:height-stabilzation}\
Let $$h:\tslr \to \R$$
be the function defined in \cite[Proof of Theorem 3.3]{zhang01} (also defined in \cite{bogomolov00}) and denoted in \cite{zhang01} by $$\ell:\tslr \to \R.$$ For $g\in\tslr$, write
$$H=e^{h}:\tslr\to \R.$$
One should think of $h$ as the \textit{logarithmic Frobenius-stabilized height function} and $H$ as the \textit{Frobenius-stabilized height function} or more simply as the \textit{stablized logarithmic height function} and \textit{stablized height function} respectively (also see \ssep\ref{ss:height-stabilzation}).

\subparat{Properties of the height function} The central point in the proofs of \cite{bogomolov00,zhang01} is the following property of the height function introduced in \ssep\ref{ss:height-stabilzation}. By \cite{zhang01} one has  for any $\tilde{g}_{1},\tilde{g}_{2}\in\tslr$:
\be\label{eq:bogomolov-zhang} h(\tilde{g}_{1}\cdot\tilde{g}_2)\leq h(\tilde{g}_{1})+h(\tilde{g}_{2}).\ee 
Given a $\Theta$-link $(\tilde{g}_{s,1},\ldots,\tilde{g}_{s,\ells})_{s\in S}\in \bTheta_{classical}$, the quantity  of interest to us is:
$$\sum_{s\in S}h(\tilde{g}_{s,1}\cdot\tilde{g}_{s,2} \cdots \tilde{g}_{s,\ells} ).$$

Hence by \eqref{eq:bogomolov-zhang} one has $$\sum_{s\in S} h(\tilde{g}_{s,1}\cdot\tilde{g}_{s,2} \cdots \tilde{g}_{s,\ells} )\leq \sum_{s\in S}\sum_{j=1}^\ells h(\tilde{g}_j).$$

In particular one can also consider the average
$$\frac{1}{\ells}\sum_{s\in S} h(\tilde{g}_{s,1}\cdot\tilde{g}_{s,2} \cdots \tilde{g}_{s,\ells} )\leq \frac{1}{\ells}\sum_{s\in S}\sum_{j=1}^\ells h(\tilde{g}_j).$$

\brem 
The proof of \cite{zhang01}, \cite{bogomolov00} (and \cite{mochizuki-bogomolov}) can be understood as a sort of averaging over $\tslz$. In what I do here,  this averaging is being replaced here by averaging over $\tslql$. The innovation presented here is the idea of averaging over $\tslql$ which is not considered in \cite{mochizuki-bogomolov}--but corresponds to working simultaneously with Mochizuki's $\Theta_{gau}$-Links,  and $\flog$-links (\cite{mochizuki-iut3}).
\erem

\subparat{The height function $h$ and the global Frobenius morphism (aka $\flog$-Link)}\label{ss:geom-case-height-frobenius}
Let me explain the reason why $h$ is referred to as a stabilized height function in \cref{ss:height-stabilzation}. Let $\tilde{g}'_{\tau,0},\tilde{g}_{\tau,0}\in\tslr$ and assume that $$\tilde{g}'_{\tau,0}=\tilde{g}_{\tau,0}\cdot \vphi_\infty^m$$
for $m\in\Z$. By \Cref{de:geom-log-links}, this equation means that $\tilde{g}'_{\tau,0},\tilde{g}_{\tau,0}$ differ by an iterate of the $\flog$-Link or equivalently differ by the  $m^{th}$ iterate of the global Frobenius $\vphi_\infty$ for some $m\in\Z$. 
Then one has the following key relationship established in \cite[Proof of the Geometric Szpiro Inequality]{zhang01} between the heights  $h(\tilde{g}'_{\tau,0})$ and $h(\tilde{g}_{\tau,0})$ i.e. between the height of $\tilde{g}_{\tau,0}$ and the height of a log-linked element $\tilde{g}'_{\tau,0}$:
\be h(\tilde{g}'_{\tau,0}) = h(\tilde{g}_{\tau,0}\cdot\vphi_\infty^m) \leq h(\tilde{g}_{\tau,0}) + \pi\cdot m\in\R.\ee
(where $\pi\in\R$ is the familiar real number). 

This property plays a central role in the proof of \cite[Proof of the Geometric Szpiro Inequality]{zhang01} and Mochizuki's discussion of this point and its relationship to his proof can be found in \cite{mochizuki-bogomolov}.

\subparat{The definition of $\abs{\bTheta_{classical}(D)}$} More generally, keeping in mind \cite{mochizuki-iut4}, let $S\subset D\subset C(\C)$ be a subset containing $S$ such that the local theory of \cite{zhang01} can be applied (for example let $D$ be the disjoint union of small disk $D_s\subset C(\C)$ around each $s\in S$). Then it is possible to define 
$\bTheta_{classical}(D)\supset \bTheta_{classical}$ as follows:
\bdefn\label{def:theta-abs-val-def} 
For a compact domain $S\subset D\subset C(\C)$ chosen as above, let
$$\abs{\bTheta_{classical}(D)}=\sup\left\{\sum_{s\in S} \sum_{j=1}^\ells h(\tilde{g}_{s,j}): (\tilde{g}_{s,1},\tilde{g}_2 \cdots, \tilde{g}_{s,\ells})_{s\in S}\in \bTheta_{classical}(D) \right\}.$$
\edefn

The following lemma will be useful:
\blem\label{le:sum-product-h-prop} 
$$\abs{\bTheta_{classical}(D)}\geq \sup\left\{\sum_{s\in S}h(\tilde{g}_{s,1}\cdot\tilde{g}_2 \cdots \tilde{g}_{s,\ells}): (\tilde{g}_{s,1},\tilde{g}_2 \cdots, \tilde{g}_{s,\ells})_{s\in S}\in \bTheta_{classical}(D) \right\}$$
\elem
\bp 
This is immediate from the above discussed property \eqref{eq:bogomolov-zhang} of $h$ and the definition of $\abs{\bTheta_{classical}(D)}$.
\ep
\brem 
Working with a compact domain $D$ is analogous to working with compactly bounded domains (in the sense of \cite{mochizuki-general-pos}) in the proof of the main theorem of \cite{mochizuki-iut4}.
\erem
\subparat{The monodromy representation and its reduction modulo $\ell$} Let $\pi_1^{top}(C(\C)-S,*)$ be the topological fundamental group of $C(\C)-S$ for some choice of basepoint $*\in C(\C)-S$. This is a group with $2g+\abs{S}$ generators $$a_1,\ldots,b_1,\ldots,a_g,b_g,\{\gamma_s\}_{s\in S}$$ and one relation
$$\prod_{j=1}^g[a_j,b_j]\prod_{s\in S}\gamma_s=1$$
between these generators.

Let $t\in C(\C)-S$ be a general point and $X_t$ be the fiber of $X\to C(\C)$ over this point and consider the monodromy representation (on the singular homology)
$$\rho:\pi_1^{top}(C(\C)-S,*)\to {Aut}(H_1(X_t,\Z)),$$
and also its ``reduction modulo $\ell$'':
$$\rho_\ell:\pi_1^{top}(C(\C)-S,*)\to {Aut}(H_1(X_t,\Z/\ell)).$$

One can assume without any loss of generality that $\rho_\ell$ is irreducible otherwise by a classical argument due to Faltings, Szpiro's inequality is trivial.

At any point $s\in S$ one has the monodromy around $s$ which is the image of $\rho(\gamma_s)$ and also $\rho_\ell(\gamma_s)$. Let $\tilde{\gamma}^{(j)}_s\in\tslz$ be an arbitrary lift of $\rho_\ell(\gamma_s^{j^2})$ for $j=1,\ldots,\ells$.

\subparat{Mochizuki's Corollary 3.12 in the geometric case} \cite[Corollary 3.12]{mochizuki-iut3} now has a simple formulation with a tautological proof:
\bthm\label{th:geometric-case-of-moccor} 
For a compact domain $D$ as above, and any choice of $(\tilde{g}_1,\ldots,\tilde{g}_\ells)\in\bTheta_{classical}(D)$ , one has
$$\abs{\bTheta_{classical}(D)}\geq \sum_{s\in S} \sum_{j=1}^\ells h(\tilde{\gamma}_{s}^{(j)}) \geq h\left(\prod_{s\in S}\tilde{\gamma}^{(1)}_s\cdot\tilde{\gamma}^{(2)}_s\cdots\tilde{\gamma}^{(\ells)}_s\right).$$
\ethm
\bp 
The lower bound is immediate from \Cref{le:sum-product-h-prop} and the fact that, by \Cref{def:theta-abs-val-def}, $\abs{\bTheta_{classical}(D)}$ is the supremum of all such sums.
\ep
\brem 
One can think of the product $\prod_{s\in S}\left(\tilde{\gamma}^{(1)}_s\cdot\tilde{\gamma}^{(2)}_s\cdots\tilde{\gamma}^{(\ells)}_s\right)$ and more generally the product $\prod_{s\in S}\left(\tilde{\gamma}^{(1)}_s\cdot\tilde{\gamma}^{(2)}_s\cdots\tilde{\gamma}^{(\ell)}_s\right)$  as a ``simulation'' of a ``multiplicative'' one dimensional subspace of $\rho_\ell$. Since $\rho_\ell$ is irreducible there does not exist any natural one dimensional subspace  (globally) in $\rho_\ell$. This multiplicative subspace point of view is due to \cite{mochizuki-iut3}.
\erem 
\subparat{Additional Remarks} At this point I want to make the following important remarks. 
\benumlab
\item Working with $\slq\subset \glqp$ is essentially equivalent to  Mochizuki's idea of working with $\Theta_{gau}$-links and working with $\vphi_\infty$ is equivalent  $\flog$-links (in \cite{joshi-teich-estimates}, I prove both these assertions in the $p$-adic context). 
\item Especially working with $\tslq$ allows us to work with $\Theta_{gau}$-links and $\flog$-links simultaneously. In \iutthr, Mochizuki works with $\Theta_{gau}$-Links and $\flog$-Links separately. My construction of Mochizuki's Ansatz in \cite{joshi-teich-estimates} allows me to work with both these devices simultaneously.
\item The global Frobenius morphism \eqref{def:global-frob-classical-case} plays a central role in \cite{zhang01}, \cite{bogomolov00}.
\item The global relation between generators of the topological fundamental group $$\pi_1^{top}(C(\C)-S,*),$$ namely
$$\prod_{j=1}^g[a_i,b_i]\prod_{s\in S}\gamma_s=1,$$
plays the role of the product formula \eqref{eq:prod-formula} in the theory of arithmeticoids and their height functions.
\item The function $h$, which is denoted by $\ell$ in \cite{zhang01}, plays  the role of stabilized height function in the  proof of the geometric Szpiro Inequality in \cite{zhang01}.
\eenum

\phantomsection
\addcontentsline{toc}{section}{Bibliography}
\bibliography{../../master/masterofallbibs.bib}
\end{document}